# INTERVAL ALGEBRAIC BISTRUCTURES

W. B. Vasantha Kandasamy
Florentin Smarandache

**2011**

# INTERVAL ALGEBRAIC BISTRUCTURES

W. B. Vasantha Kandasamy
Florentin Smarandache

**2011**

# CONTENTS

# PREFACE

Authors in this book construct interval bistructures using only interval groups, interval loops, interval semigroups and interval groupoids.

Several results enjoyed by these interval bistructures are described. By this method, we obtain interval bistructures which are associative or non associative or quasi associative. The term quasi is used mainly in the interval bistructure $B = B_1 \cup B_2$ (or in n-interval structure) if one of $B_1$ (or $B_2$) enjoys an algebraic property and the other does not enjoy that property (one section of interval structure satisfies an algebraic property and the remaining section does not satisfy that particular property). The term quasi and semi are used in a synonymous way.

This book has four chapters. In the first chapter interval bistructures (biinterval structures) such as interval bisemigroup, interval bigroupoid, interval bigroup and interval biloops are introduced. Throughout this book we work only with the intervals of the form [0, a] where $a \in Z_n$ or $Z^+ \cup \{0\}$ or $R^+ \cup \{0\}$ or $Q^+ \cup \{0\}$ unless otherwise specified. Also interval bistructures of the form interval loop-group, interval group-groupoid so on are introduced and studied.

In chapter two n-interval structures are introduced. n-interval groupoids, n-interval semigroups, n-interval loops and so on are introduced and analysed. Using these notions n-interval mixed algebraic structure are defined and described.

Some probable applications are discussed. Only in due course of time several applications would be evolved by researchers as per their need.

The final chapter suggests around 295 problems of which some are simple exercises, some are difficult and some of them are research problems.

This book gives around 388 examples and 124 theorems which is an attractive feature of this book.

We thank Dr. K.Kandasamy for proof reading.

W.B.VASANTHA KANDASAMY
FLORENTIN SMARANDACHE

**Chapter One**

# INTERVAL BISTRUCTURES

This chapter has four sections. Section one is introductory. We introduce biintervals and biinterval matrices. These concepts will be used to construct bi-interval algebraic structures.

Let $I(Z_n) = \{[0, a] / a \in Z_n\}$ be the set of modulo integer intervals.
$I(Z^+ \cup \{0\}) = \{[0, a] / a \in Z^+ \cup \{0\}\}$ is the integer intervals.
$I(Q^+ \cup \{0\}) = \{[0, a] / a \in Q^+ \cup \{0\}\}$ is the rational intervals.
$I(R^+ \cup \{0\}) = \{[0, a] / a \in R^+ \cup \{0\}\}$ is the real intervals.

Section two introduces interval bigroupoids. Interval bigroups and their generalizations are given in section three. In section four interval biloops are introduced and studied.

## 1.1 Interval Bisemigroups

Now we define a biinterval. A biinterval is the union of two distinct intervals $I = I_1 \cup I_2 = [0, a] \cup [0, b]$ where $a \neq b$, $I_1$ and $I_2$ are intervals.

If a = b the biinterval will be known as the pseudo biinterval. We do not demand $[0, a] \not\subseteq [0, b]$ or $[0, b] \not\subseteq [0, a]$ or $[0, b] \subseteq [0, a]$.

We will illustrate this situation by some examples.

***Example 1.1.1:*** Let $I = \{[0, 7]\} \cup \{[0, \sqrt{2}]\}$ be the biinterval where $7, \sqrt{2}, \in R^+ \cup \{0\}$.

***Example 1.1.2:*** Let $I = [0, 9] \cup [0, \overline{2}]$ be the biinterval where $9 \in Z^+ \cup \{0\}$ and $\overline{2} \in Z_4$.

***Example 1.1.3:*** Let $I = I_1 \cup I_2 = [0, \overline{3}] \cup [0, \overline{7}]$ where $\overline{3} \in Z_6$ and $\overline{7} \in Z_{10}$ be the biinterval.

***Example 1.1.4:*** Let $I = I_1 \cup I_2 = [0, 20] \cup [0, \sqrt{17}]$ where $20 \in Z^+ \cup \{0\}$ and $\sqrt{17} \in R^+ \cup \{0\}$ be the biinterval.

Now we have seen examples of biintervals. We will define some operations on these biintervals so that they get some algebraic structures.

Further it is easy to define concepts like n-intervals; when n = 2 we get the biinterval and when n = 3 we get the triinterval.

We will define them in the following.

Let $I = I_1 \cup I_2 \cup I_3$ be such that each interval $I_j$ is distinct $1 \leq j \leq 3$, then we define I to be a triinterval and $I_j$ can take their values from $Z^+ \cup \{0\}$ or $Z_n$ or $R^+ \cup \{0\}$ or $Q^+ \cup \{0\}$ or not used in the mutually exclusive sence, $1 \leq j \leq 3$.

Suppose $I = I_1 \cup I_2 \cup \ldots \cup I_n$, $(n \geq 2)$ where $I_j$ is an interval from $Q^+ \cup \{0\}$ or $R^+ \cup \{0\}$ or $Z^+ \cup \{0\}$ or $Z_m$; $1 \leq j \leq n$, $I_j \neq I_k$ if $k \neq j$; $1 \leq j, k \leq n$ then we define I to be an n-interval.

We will give some examples before we proceed onto define further structures.

*Example 1.1.5:* Let $I = I_1 \cup I_2 \cup I_3 = [0, 3] \cup [0, \sqrt{7}] \cup [0, \overline{5}]$ be a triinterval where $3 \in Z^+ \cup \{0\}$ $\sqrt{7} \in R^+ \cup \{0\}$ and $\overline{5} \in Z_{11}$.

*Example 1.1.6:* Let $I = I_1 \cup I_2 \cup I_3 \cup I_4 \cup I_5 \cup I_6 = \{[0, 7] \cup [0, 5] \cup [0, 2] \cup [0, 7/3] \cup [0, \sqrt{71}] \cup [0, \frac{27}{31}]\}$ be a 6-interval where the entries in each interval is from $R^+ \cup \{0\}$.

*Example 1.1.7:* Let $I = I_1 \cup I_2 \cup \ldots \cup I_9 = \{[0, x_1] \cup \ldots \cup [0, x_9]\}$ where $x_i \in Z^+ \cup \{0\}$, $1 \leq i \leq 9$, be the 9-interval.

*Example 1.1.8:* Let $I = I_1 \cup I_2 \cup I_3 \cup I_4 = [0, x_1] \cup [0, x_2] \cup [0, x_3] \cup [0, x_4]$ where $x_1 \in Z_7$, $x_2 \in Z_{12}$, $x_3 \in Z_{17}$ and $x_4 \in Z_{21}$ be the 4-interval.

**DEFINITION 1.1.1:** *Let $S = S_1 \cup S_2$ where $S_1$ and $S_2$ are distinct interval semigroups under the operations '*' and 'o' respectively. ($S_1 \neq S_2$, $S_1 \not\subseteq S_2$ and $S_2 \not\subseteq S_1$) then $(S, .)$ is defined as a interval bisemigroup or biinterval semigroup denoted by $S = S_1 \cup S_2 = \{[0, a], *\} \cup \{[0, b], o\} = \{[0, a] \cup [0, b] / [0, a] \in S_1$ and $[0, b] \in S_2\}$. $I.J = ([0, a] \cup [0, b]) . ([0, x] \cup [0, y])$ where '.' is defined as $I.J = \{[0, a] * [0, x] \cup [0, b] \, o \, [0, y]\}$.*

We will illustrate this situation by some examples.

*Example 1.1.9:* Let $S = S_1 \cup S_2 = \{[0, a] / a \in Z^+ \cup \{0\}$ be an interval semigroup under addition$\} \cup \{[0, b] / b \in Z_{20}$ under multiplication modulo 20 is an interval semigroup$\}$ be the interval bisemigroup or biinterval semigroup.

*Example 1.1.10:* Let $S = S_1 \cup S_2$ where $S_1 = \{[0, a] / a \in Q^+ \cup \{0\}$ under multiplication$\}$ is an interval semigroup and $S_2 = \{[0,$ a] / a ∈ $R^+ \cup \{0\}$ under addition} be the interval semigroup. S is an interval bisemigroup.

Now we see the examples given in 1.1.9 and 1.1.10 are biinterval semigroups of infinite order.

We see how in the example 1.1.9 the operation is carried out on the interval bisemigroup. Let $I = I_1 \cup I_2 = [0, 10] \cup [0, 4]$ and $J = [0, 27] \cup [0, 12]$ be in S.

$$
\begin{aligned}
I \cdot J &= \{[0, 10] \cup [0, 4]\} \cdot \{[0, 27] \cup [0, 12]\} \\
&= \{[0, 10] + [0, 2] \cup [0, 4] \times [0, 12]\} \\
&= \{[0, 10+2] \cup [0, 4 \times 12 \ (\text{mod } 20)]\} \\
&= [0, 12] \cup [0, 8].
\end{aligned}
$$

Thus (S, '·') is also known as biinterval semigroup as its elements are biintervals.

***Example 1.1.11:*** Let $S = S_1 \cup S_2 = \{[0, a] / a \in Z_{13}$ under addition modulo 13$\} \cup \{[0, b] / b \in Z_{16}$ under multiplication modulo 16$\}$ be a biinterval semigroup.

Clearly if $x = [0, 10] \cup [0, 14]$ and $y = [0, 7] \cup [0, 10]$ are in S, then $x \cdot y = ([0, 10] \cup [0, 14]) \cdot ([0, 7] \cup [0, 10]) = ([0, 10] + [0, 7]) \cup ([0, 14] \times [0, 10]) = [0, 17 \ (\text{mod } 13)] \cup [0, 140 \ (\text{mod } 16)] = [0, 4] \cup [0, 12]$.

It is interesting to note that the biinterval semigroup given in example 1.1.11 is of finite order.

We say biinterval semigroup is commutative if both $S_1$ and $S_2$ are commutative interval semigroups. If only one of $S_1$ or $S_2$ is a commutative interval semigroup and the other is not commutative then we say S is a semi commutative biinterval semigroup.

All the examples of biinterval semigroups given are commutative to define biinterval semigroups which are semi commutative we have to use either matrix interval semigroups or symmetric interval semigroups. For more about these please refer [10, 14].

We will give examples of non commutative biinterval semigroups and semi commutative biinterval semigroups.

Clearly elements of these biinterval semigroups will not be biintervals so with some flaw we will call them as interval bisemigroups.

*Example 1.1.12:* Let $S = S_1 \cup S_2$ where

$$S_1 = \left\{ \begin{bmatrix} [0,a] & [0,b] \\ [0,c] & [0,d] \end{bmatrix} \middle| a,b,c,d \in Z^+ \cup \{0\} \right\} \text{ and}$$

$$S_2 = \left\{ \begin{bmatrix} [0,a] & [0,b] & [0,c] \\ [0,d] & [0,f] & [0,h] \\ [0,m] & [0,n] & [0,p] \end{bmatrix} \middle| a,b,c,d,f,h,m,n,p \in Z_7 \right\} \text{ are}$$

interval matrix semigroups under multiplication. $S = S_1 \cup S_2$ is an interval bisemigroup which is non commutative.
For take

$$x = \begin{bmatrix} [0,9] & [0,2] \\ [0,3] & [0,1] \end{bmatrix} \cup \begin{bmatrix} [0,1] & 0 & 0 \\ [0,3] & [0,4] & 0 \\ [0,2] & 0 & [0,6] \end{bmatrix}$$

and

$$y = \begin{bmatrix} [0,1] & [0,10] \\ [0,12] & [0,4] \end{bmatrix} \cup \begin{bmatrix} [0,2] & [0,1] & [0,3] \\ [0,1] & [0,4] & [0,1] \\ 0 & [0,1] & 0 \end{bmatrix}$$

in S.

$$x.y = \begin{bmatrix} [0,9] & [0,2] \\ [0,3] & [0,1] \end{bmatrix} \times \begin{bmatrix} [0,1] & [0,10] \\ [0,12] & [0,4] \end{bmatrix} \cup$$

$$\begin{bmatrix} [0,1] & 0 & 0 \\ [0,3] & [0,4] & 0 \\ [0,2] & 0 & [0,6] \end{bmatrix} \times \begin{bmatrix} [0,2] & [0,1] & [0,3] \\ [0,1] & [0,4] & [0,1] \\ 0 & [0,1] & 0 \end{bmatrix}$$

$$= \begin{bmatrix} [0,9][0,1]+[0,2][0,12] & [0,9][0,10]+[0,2][0,4] \\ [0,3][0,1]+[0,1][0,12] & [0,3][0,10]+[0,1][0,4] \end{bmatrix} \cup$$

$$\begin{bmatrix} [0,1][0,2]+0+0 \\ [0,3][0,2]+[0,4][0,1]+0 \\ [0,2][0,2]+0+0 \end{bmatrix}$$

$$\begin{bmatrix} [0,1][0,1]+0+0 & [0,1][0,3]+0+0 \\ [0,3][0,1]+[0,4][0,4]+0 & [0,3][0,3]+[0,4][0,1]+0 \\ [0,2][0,1]+0+[0,6][0,1] & [0,2][0,3]+0+0 \end{bmatrix}$$

$$= \begin{bmatrix} [0,9]+[0,24] & [0,90]+[0,8] \\ [0,3]+[0,12] & [0,30]+[0,4] \end{bmatrix} \cup$$

$$\begin{bmatrix} [0,2]+0+0 & [0,1]+0+0 & [0,3]+0+0 \\ [0,6]+[0,4]+0 & [0,3]+[0,2]+0 & [0,2]+[0,4]+0 \\ [0,4]+0+0 & [0,2]+0+[0,6] & [0,6]+0+0 \end{bmatrix}$$

$$= \begin{bmatrix} [0,33] & [0,98] \\ [0,15] & [0,34] \end{bmatrix} \cup \begin{bmatrix} [0,2] & [0,1] & [0,3] \\ [0,3] & [0,5] & [0,6] \\ [0,4] & [0,1] & [0,6] \end{bmatrix}.$$

***Example 1.1.13:*** Let $S = S_1 \cup S_2$ where $S_1 = \{S(X)$ where $X = ([0, 1], [0, 2], [0, 3])\}$ is the interval symmetric semigroup and $S_2 = \{$All $5 \times 5$ interval matrices with intervals of the form $[0, a]$ where $a \in Z_{12}\}$ be the interval matrix semigroup under multiplication. S is a interval bisemigroup, non commutative and is of finite order. The product on S is defined as follows. Let

$$x = \begin{pmatrix} [0,1] & [0,2] & [0,3] \\ [0,3] & [0,2] & [0,1] \end{pmatrix} \cup$$

$$\begin{pmatrix} 0 & [0,1] & 0 & [0,2] & 0 \\ [0,3] & 0 & [0,4] & 0 & [0,1] \\ 0 & [0,2] & 0 & [0,3] & 0 \\ [0,5] & 0 & [0,1] & 0 & [0,2] \\ 0 & [0,3] & 0 & [0,7] & 0 \end{pmatrix}$$

and
$$y = \begin{pmatrix} [0,1] & [0,2] & [0,3] \\ [0,2] & [0,1] & [0,3] \end{pmatrix} \cup$$
$$\begin{pmatrix} [0,1] & 0 & [0,2] & 0 & [0,3] \\ 0 & [0,4] & 0 & [0,5] & 0 \\ [0,6] & 0 & [0,7] & 0 & [0,8] \\ 0 & [0,9] & 0 & [0,10] & 0 \\ [0,11] & 0 & [0,1] & 0 & [0,2] \end{pmatrix}$$

belong to S.

Now
$$x.y = \begin{pmatrix} [0,1] & [0,2] & [0,3] \\ [0,3] & [0,2] & [0,1] \end{pmatrix} \cdot \begin{pmatrix} [0,1] & [0,2] & [0,3] \\ [0,2] & [0,1] & [0,3] \end{pmatrix} \cup$$

$$\begin{pmatrix} 0 & [0,1] & 0 & [0,2] & 0 \\ [0,3] & 0 & [0,4] & 0 & [0,1] \\ 0 & [0,2] & 0 & [0,3] & 0 \\ [0,5] & 0 & [0,1] & 0 & [0,2] \\ 0 & [0,3] & 0 & [0,7] & 0 \end{pmatrix} \times$$
$$\begin{pmatrix} [0,1] & 0 & [0,2] & 0 & [0,3] \\ 0 & [0,4] & 0 & [0,5] & 0 \\ [0,6] & 0 & [0,7] & 0 & [0,8] \\ 0 & [0,9] & 0 & [0,10] & 0 \\ [0,11] & 0 & [0,1] & 0 & [0,2] \end{pmatrix}$$

$$= \begin{pmatrix} [0,1] & [0,2] & [0,3] \\ [0,3] & [0,1] & [0,2] \end{pmatrix} \cup \begin{pmatrix} 0 & [0,10] & 0 & [0,1] & 0 \\ [0,2] & 0 & [0,5] & 0 & [0,9] \\ 0 & [0,5] & 0 & [0,4] & 0 \\ [0,3] & 0 & [0,1] & 0 & [0,3] \\ 0 & [0,3] & 0 & [0,1] & 0 \end{pmatrix}$$

is in S.

We can have any number of finite or infinite, commutative or non commutative interval bisemigroups. We have seen examples of them.

We can define substructures for them.

**DEFINITION 1.1.2:** *Let $S = S_1 \cup S_2$ be an interval bisemigroup. $P = P_1 \cup P_2 \subseteq S_1 \cup S_2 = S$ be a proper non empty bisubset of S. If each $P_i$ is an interval subsemigroup of $S_i$, i=1,2 then we call P to be an interval subbisemigroup or interval bisubsemigroup of S.*

We will illustrate this situation by some examples.

*Example 1.1.14:* Let $S = S_1 \cup S_2$ where $S_1 = \{[0, a] / a \in Z^+ \cup \{0\}\}$ is the interval semigroup under multiplication and $S_2 = \{[0, a] / a \in Z_{24}\}$ is the interval semigroup under addition be a biinterval semigroup. Consider $P = P_1 \cup P_2 = \{[0, a] / a \in 3Z^+ \cup \{0\}\} \cup \{[0, a] / a \in \{0, 2, 4, 6, 8, 10, \ldots, 20, 22\} \subseteq Z_{24}\} \subseteq S_1 \cup S_2 = S$. It is easily verified P is a biinterval subsemigroup of S.

*Example 1.1.15:* Let $S = S_1 \cup S_2$ where $S_1 = \{[0, a] / a \in Z_{11}\}$ is the interval semigroup under addition modulo 11 and $S_2 = \{[0, b] / b \in Z_{17}\}$ be the interval semigroup under addition modulo 17. S is a biinterval semigroup. It is easily verified S has no biinterval subsemigroups.

If an interval bisemigroup has no interval subsemigroups then we call S to be a simple biinterval semigroup.

We have a class of biinterval semigroups which are simple.

**THEOREM 1.1.1:** *Let $S = S_1 \cup S_2 = \{[0, a] / a \in Z_p$, p a prime, $+\} \cup \{[0, b] / b \in Z_q$, q a prime, $+\}$ $(p \neq q)$ be a biinterval semigroup. S is a simple biinterval semigroup.*

The proof is direct and left as an exercise for the reader.

*Example 1.1.16:* Let $S = S_1 \cup S_2 = \{[0, a] / a \in Z^+ \cup \{0\}\} \cup \{[0, b] / b \in Z_{43}\}$ be an interval bisemigroup with operation addition defined on both $S_1$ and $S_2$. Clearly $S_1$ has interval subsemigroups where as $S_2$ has no interval subsemigroups.

We do not call S as simple but we define such substructure as quasi interval bisubsemigroup and these interval semigroups are also known as quasi simple interval bisemigroups.

We will illustrate this situation by some examples.

*Example 1.1.17:* Let $S = S_1 \cup S_2 = \{[0, a] / a \in R^+ \cup \{0\}\} \cup \{[0, a] / a \in Z_{23}\}$ be an interval bisemigroup under interval addition. Take $P = P_1 \cup P_2 = \{[0, a] / a \in Q^+ \cup \{0\}\} \cup P_2 = (=S_2) \subseteq S_1 \cup S_2 = S$; P is a quasi interval bisubsemigroup of S.

Infact S has infinitely many quasi interval bisubsemigroups.

*Example 1.1.18:* Let $S = S_1 \cup S_2 = \{[0, a] / a \in Z_{12}, \times,$ multiplication modulo $12\} \cup \{[0, a] / a \in Z_{43},$ addition modulo $43\}$ be an interval bisemigroup.

Consider $P = P_1 \cup P_2 = \{[0, a] / a \in \{0, 3, 6, 9\} \subseteq Z_{12}\} \cup P_2 (= S_2) \subseteq S_1 \cup S_2 = S$; P is a quasi interval bisubsemigroup. S is a quasi simple interval bisemigroup. We see S has only 4 quasi interval bisubsemigroups and S is of finite order.

Now we proceed onto give examples of interval biideals or biinterval ideals of an interval bisemigroup S.

*Example 1.1.19:* Let $S = S_1 \cup S_2 = \{[0, a] / +, a \in Z^+ \cup \{0\}\} \cup \{[0, b] / \times, b \in Z_{24}\}$ be an interval bisemigroup. Consider $J = J_1 \cup J_2 = \{[0, a] / a \in 5Z^+ \cup \{0\}\} \cup \{[0, b] / b \in \{0, 3, 6, 9, 12, 15, 18, 21\} \subseteq Z_{24}\} \subseteq S_1 \cup S_2 = S$; it is easily verified J is a biideal of S called the interval biideal or biinterval ideal of S.

*Example 1.1.20:* Let $S = S_1 \cup S_2 = \{[0, a] / a \in Z_6\} \cup \{[0, b] / b \in Z_{30}\}$ be a biinterval semigroup under multiplication. Consider $P = P_1 \cup P_2 = \{[0, a] / a \in \{0, 2, 4\} \subseteq Z_6\} \cup \{[0, b] / b \in \{10, 20, 0\} \subseteq Z_{30}\} \subseteq S_1 \cup S_2$; P is a biinterval ideal of S.

***Example 1.1.21:*** Let $S = S_1 \cup S_2 = \{[0, a] / a \in Z_{11}\} \cup \{[0, b] / b \in Z_{17}\}$ be an interval bisemigroup. S has no biideals.

We call those interval bisemigroups to be ideally simple interval bisemigroups if S has no biideals. Example 1.1.21 is an ideally simple interval bisemigroup. We have an infinite class of ideally simple interval bisemigroups.

**THEOREM 1.1.2:** *Let $S = S_1 \cup S_2 = \{[0, a] / a \in Z_p, p$ a prime$\} \cup \{[0, b] / b \in Z_q, q$ a prime$\}$ $(p \neq q)$ be an interval bisemigroup. S is an ideally simple interval bisemigroup.*

The reader is left with the task of proving this result.

***Example 1.1.22:*** Let $S = S_1 \cup S_2 = \{[0, a] / a \in Z_{11}\} \cup \{[0, b] / b \in Z^+ \cup \{0\}\}$ be an interval bisemigroup under multiplication. We see $S_1$ has no ideals, but $S_2$ has infinitely many ideals. Thus S has quasi interval ideals.

If S has only quasi interval ideals then we define S to be a quasi ideally simple interval bisemigroup.

We give examples of this structure.

***Example 1.1.23:*** Let $S = S_1 \cup S_2 = \{[0, a] / a \in R^+ \cup \{0\}\} \cup \{[0, b]/b \in Z_{30}\}$ be a biinterval semigroup under multiplication. $P = P_1 (=S_1) \subseteq S_1 \cup S_2$; P is a quasi ideal hence S is a quasi ideally simple biinterval semigroup.

***Example 1.1.24:*** Let $S = S_1 \cup S_2 = \{[0, a] / a \in Z_7\} \cup \{[0, a] / a \in Z_{50}\}$ be an interval bisemigroup. Consider $H = H_1 (=S_1) \cup H_2 = S_1 \cup \{[0, a] / a \in \{0, 5, 10, 15, \ldots, 45\} \subseteq Z_{50}\} \subseteq S_1 \cup S_2$; H is a biideal and S is a quasi ideally simple biinterval semigroup.

***Example 1.1.25:*** Let $S = S_1 \cup S_2 = \{[0, a] / a \in R^+ \cup \{0\}\} \cup \{[0, b] / b \in Q^+ \cup \{0\}\}$ be a biinterval semigroup. S is an ideally simple biinterval semigroup.

We can define bizero divisors, quasi bizero divisors, biunits, quasi biunits and biidempotents and quasi biidempotents in these biinterval semigroups.

Let $S = S_1 \cup S_2$ be an biinterval semigroup. Let $\alpha = \alpha_1 \cup \alpha_2$ be a biinterval in S if there exists a $\beta = \beta_1 \cup \beta_2$ in S which is a biinterval such that $\alpha \cdot \beta = 0 \cup 0$ then we say $\alpha$ is a bizero interval zero divisor in S.

We will first illustrate this situation by some examples.

*Example 1.1.26:* Let $S = S_1 \cup S_2 = \{[0, a] / a \in Z_{12}\} \cup \{[0, b] / b \in Z_{420}\}$ be a biinterval semigroup.
Choose $\alpha = \alpha_1 \cup \alpha_2 = [0, 6] \cup [0, 60]$ in $S_1 \cup S_2$. $\beta = \beta_1 \cup \beta_2 = [0, 4] \cup [0, 7]$ in $S_1 \cup S_2$ is such that $\alpha\beta = ([0, 6] \cup [0, 60]) [[0, 4] \cup [0, 7]) = [0, 0] \cup [0, 0]$ thus $\alpha$ is a interval bizero divisor in S.

Let $x = x_1 \cup x_2 = [0, 11] \cup [0, 419]$ be in $S = S_1 \cup S_2$.

We see $x^2 = [0, 1] \cup [0, 1]$ is a biunit in S. Consider $y = y_1 \cup y_2 = [0, 4] \cup [0, 36] \in S_1 \cup S_2$ is such that $([0, 4] \cup [0, 36])^2 = [0, 4] \cup [0, 36] = y_1 \cup y_2$ is an interval biidempotent of $S = S_1 \cup S_2$.

*Example 1.1.27:* Let $S = S_1 \cup S_2 = \{[0, a] / a \in Z^+ \cup \{0\}\} \cup \{[0,b] / b \in Q^+ \cup \{0\}\}$ be an interval bisemigroup; S has no bizero divisors or biunits or biidempotents.

*Example 1.1.28:* Let $S = S_1 \cup S_2 = \{[0, a] / a \in Z_7\} \cup \{[0, b] / b \in Z_{19}\}$ be an interval bisemigroup of finite order. S has no non trivial biidempotents or bizero divisors but has a biunit given by $x = [0, 6] \cup [0, 18] \in S_1 \cup S_2$ is such that $x^2 = [0, 1] \cup [0, 1]$.

*Example 1.1.29:* Suppose $S = S_1 \cup S_2 = \{[0, a] / a \in Z_p$, p a prime$\} \cup \{[0, b] / b \in Z_q$, q a prime$\}$; $p \neq q$ be a biinterval semigroup. Then S has only one non trivial biunit given by $\alpha = [0, p-1] \cup [0, q-1] \in S$, is such that $\alpha^2 = [0, 1] \cup [0, 1]$.

Now we have seen the special elements in the interval bisemigroup. It may so happen one of $S_i$ may have units, zero divisors and idempotents and $S_j$ ($i \neq j$) may not have any of these in these situations we do define the following.

$\alpha = [0, x] + 0$ in S is a quasi biunit if there exists a $\beta = [0, y] + 0$ in S such that $\alpha\beta = [0, 1] + 0$.

Similarly we define quasi biidempotent and quasi bizero divisor in $S = S_1 \cup S_2$.

We will only illustrate these situations by some examples.

***Example 1.1.30:*** Let $S = S_1 \cup S_2 = \{[0, a] \,/\, a \in Z^+ \cup \{0\}\} \cup \{[0, b] \,/\, b \in Z_{30}\}$ be a biinterval semigroup. We see $S_1$ has no units or zero divisors or idempotents but $S_2$ has zero divisors, units and idempotents, hence S can have only quasi biunits, quasi biidempotents and quasi bizero divisors.

Consider $\alpha = 0 \cup [0, 29] \in S = S_1 \cup S_2$. Clearly $\alpha^2 = 0 \cup [0, 1]$ is a quasi biunit.

Take $\beta = 0 \cup [0, 15] \in S = S_1 \cup S_2$ we have $\gamma = 0 \cup [0, 4] \in S$ such that $\beta\gamma = 0 \cup 0$. Hence $\beta$ is a quasi bizero divisor in S.

Take $x = 0 \cup [0, 15] \in S$, we see $x^2 = 0 \cup [0, 225] = [0, 15] \in S$ is a quasi biidempotent of S.

Thus S has only quasi biunits, quasi bizero divisors and quasi biidempotents in it.

***Example 1.1.31:*** Let $S = S_1 \cup S_2 = \{[0, a] \,/\, a \in R^+ \cup \{0\}\} \cup \{[0, a] \,/\, a \in Z_{31}\}$ be a biinterval semigroup. We see S has only one quasi biunit and has no quasi bizero divisors or quasi biidempotents.

This is also a different situation from example 1.1.30.

***Example 1.1.32:*** Let $S = S_1 \cup S_2 = \{[0, a] \,/\, a \in Z_{29}\} \cup \{[0, b] \,/\, b \in Z_{35}\}$ be a biinterval semigroup. We see S has a biunit given by $x = [0, 28] \cup [0, 34] \in S = S_1 \cup S_2$. Clearly $x^2 = [0, 1] \cup$ [0, 1]. But S has no biidempotents or bizero divisors. But S has both quasi bizero divisors and quasi biidempotents.

Take $\alpha = 0 \cup [0, 5] \in S$, we see $\beta = 0 \cup [0, 7]$ is such that $\alpha\beta = 0 \cup 0$. Also $y = 0 \cup [0, 14]$ is such that $\alpha y = 0 \cup 0$. Thus S has quasi bizero divisors. Consider $x = 0 \cup [0, 15] \in S$; we see $x^2 = 0 \cup [0, 15] = x$. Thus x is a quasi biidempotent in S.

Now having seen special elements in biinterval semigroups. We see some examples of special type of interval bisemigroups.

*Example 1.1.33:* Let $S = S_1 \cup S_2$

$$= \left\{ \sum_{i=0}^{\infty} [0,a] x^i \,\bigg|\, a \in Z_{12} \right\} \cup \left\{ \begin{bmatrix} [0,a] & [0,b] \\ [0,c] & [0,d] \end{bmatrix} \,\bigg|\, a,b,c,d \in Z^+ \cup \{0\} \right\}$$

be an interval bisemigroup. S has bizero divisors.
Take

$$x = \{[0, 4] x^7\} \cup \left\{ \begin{bmatrix} [0,a] & 0 \\ 0 & 0 \end{bmatrix} \right\} \text{ in } S = S_1 \cup S_2.$$

We see

$$y = \{[0, 3] x^9\} \cup \left\{ \begin{bmatrix} 0 & 0 \\ 0 & [0,b] \end{bmatrix} \right\}$$

in S is such that

$$xy = 0 \cup \left\{ \begin{bmatrix} 0 & 0 \\ 0 & 0 \end{bmatrix} \right\}.$$

Thus S has bizero divisors. It is easily verified S has both biunits and biidempotents.

*Example 1.1.34:* Let $S = S_1 \cup S_2 = \{$All $3 \times 3$ interval matrices with intervals of the form [0, a] where $a \in Z^+ \cup \{0\}\} \cup \{$All $2 \times 2$ interval matrices with intervals of the form [0, b], $b \in Z_{47}\}$ be an interval bisemigroup. S has bizero divisors, biunits and biidemponents.

The following theorem is direct and the proof is left as an exercise for the reader.

**THEOREM 1.1.3:** *Let $S = S_1 \cup S_2$ be an interval bisemigroup (or biinterval semigroup). If S has biunits or biidempotents or bizero divisors then S has quasi biunits or quasi biidempotents or quasi bizero divisors. But if S has quasi biunits or quasi biidemponents or quasi bizero divisors then S in general need not contain biunits or bizero divisors or biidempotents.*

Now we proceed onto define Smarandache interval bisemigroup and quasi Smarandache interval bisemigroups and illustrate them by examples.

**DEFINITION 1.1.3:** *Let $S = S_1 \cup S_2$ be an interval bisemigroup. Suppose $A = A_1 \cup A_2 \subseteq S_1 \cup S_2$ is such that each $A_i$ is an interval group under the operations of $S_i$ (i = 1, 2) then we define $S = S_1 \cup S_2$ to be a Smarandache interval bisemigroup. (S-interval bisemigroup). If only one of the $A_i$'s is an interval group and other $A_j$ is not an interval group for any subset $A_j$ of $S_j$ ($i \neq j$) then we call $S = S_1 \cup S_2$ to be a quasi Smarandache interval bisemigroup (quasi S-interval bisemigroup).*

We will provide examples of these definitions.

***Example 1.1.35:*** Let $S = S_1 \cup S_2 = \{[0, a] / a \in Z_{40}\} \cup \{[0, a] / a \in Z_{24}\}$ be an interval bisemigroup. Consider $A = A_1 \cup A_2 = \{[0, 1], [0, 39]\} \cup \{[0, 1], [0, 23]\} \subseteq S_1 \cup S_2 = S$ is an interval bigroup (as $A_i$ is an interval group for i = 1,2). Hence S is a S-interval bisemigroup of finite order.

***Example 1.1.36:*** Let $S = \{S(X)$ where $X = ([0, 1] [0, 2] [0, 3] [0, 4])$ be an interval symmetric semigroup$\} \cup \{[0, a] / a \in Z_{19}\}$ be an interval bisemigroup of finite order. Consider $A = S_1 \cup S_2 = \{S_x$, the interval symmetric group of $S(X)\} \cup \{[0, a] / a \in Z_{19} \setminus \{0\}\} \subseteq S_1 \cup S_2$ is an interval bigroup. Thus S is a S-interval bisemigroup.

*Example 1.1.37:* Let $S = \{[0, a] / a \in Z_{23}\} \cup \{[0, b] / b \in Z_{41}\}$ be an interval bisemigroup. $A = A_1 \cup A_2 = \{[0, a] / Z_{23} \setminus \{0\}\} \cup \{[0, b] / Z_{41} \setminus \{0\}\} \subseteq S_1 \cup S_2 = S$ is a interval bigroup. Thus S is a S-interval bisemigroup.

We see from examples 1.1.36 and 1.1.37 that the interval bigroup which we have considered in these interval bisemigroups are the largest ones. We call such large A as Smarandache hyper bigroup of S.

**THEOREM 1.1.4:** *Let $S = S_1 \cup S_2 = \{[0, a] / a \in Z_p, p \text{ a prime}\} \cup \{[0, b] / b \in Z_q, q \text{ a prime}\}$ ($p \neq q$). S has only one large interval bigroup and it is the S-hyper bisubgroup of S.*

The proof is left as an exercise for the reader.

**THEOREM 1.1.5:** *Let $S = S_1 \cup S_2 = \{S(x) / x = ([0, a_1], ..., [0, a_n])\} \cup \{S(y) / y = ([0, b_1], ...., [0, b_m])\}$ $m \neq n$ be an interval bisemigroup. S has a S-hyper subbigroup. However S has several interval bigroups.*

This proof is also left as an exercise to the reader. For more information refer [10, 13-4]. Now having seen examples of S-interval bisemigroups we proceed onto give examples of quasi S-interval bisemigroups.

*Example 1.1.38:* Let $S = S_1 \cup S_2 = \{[0, a] / a \in Z^+ \cup \{0\}\} \cup \{[0, b] / b \in Z_6\}$ be an interval bisemigroup. We see $S_1$ has no proper interval subgroup but $S_2$ has proper interval subgroups. Hence S is only a quasi Smarandache interval bisemigroup.

*Example 1.1.39:* Let $S = S_1 \cup S_2 = \{[0, a] / a \in 3Z^+ \cup \{0\}\} \cup \{[0, b] / b \in Z_{52}\}$ be an interval bisemigroup. It is easily verified S is only a quasi S-interval bisemigroup.

Now we proceed onto define semi interval bisemigroup.

**DEFINITION 1.1.4:** *Let $S = S_1 \cup S_2$ where only one of $S_1$ or $S_2$ is an interval semigroup and the other is a semigroup. We call S a semi interval bisemigroup.*

*Example 1.1.40:* Let $S = S_1 \cup S_2 = \{[0, a] \,/\, a \in Z^+ \cup \{0\}\} \cup \{Z_{12}, \times\}$ be a semi interval bisemigroup of infinite order.

*Example 1.1.41:* Let $S = S_1 \cup S_2 = \{[0, a] \,/\, a \in Z_{49}\} \cup \{Z_{15}, \times\}$ be a semi interval bisemigroup of finite order.

*Example 1.1.42:* Let $S = \{S(X) \mid X = \{([0, 1], [0, 2], [0, 3] [0, 4])\}$ be the interval symmetric semigroup$\} \cup \{Z_{29}, \times\} = S_1 \cup S_2$ is the semi interval bisemigroup of finite order which is quasi commutative.

*Example 1.1.43:* Let $S = S_1 \cup S_2 =$

$$\left\{ \begin{bmatrix} a_{11} & a_{12} & a_{13} \\ a_{21} & a_{22} & a_{23} \\ a_{31} & a_{32} & a_{33} \end{bmatrix} \middle| a_{ij} \in Z^+ \cup \{0\} \right\} \cup \{[0, a] \,/\, a \in Q^+ \cup \{0\}\}$$

be a semi interval bisemigroup of infinite order which is quasi commutative.

Now as in case of usual biinterval semigroup we can define substructures and special elements, this task is assigned to the reader.

We will however illustrate these by some examples.

*Example 1.1.44:* Let $S = S_1 \cup S_2 = (Z^+ \cup \{0\}) \cup \{[0, a] \mid a \in Q^+ \cup \{0\}\}$ be a semi interval bisemigroup of infinite order. $T = T_1 \cup T_2 = \{5Z^+ \cup \{0\}\} \cup \{[0, a] \mid a \in Z^+ \cup \{0\}\} \subseteq S_1 \cup S_2$ is a semi interval bisubsemigroup of S. Clearly T is not an ideal of S.

*Example 1.1.45:* Let $S = S_1 \cup S_2 = \{([0, a_1], [0, a_2], [0, a_3], [0, a_4], [0, a_5]) \mid a_i \in Z_{20}$ under multiplication modulo 20$\} \cup$

S(6) be a semi interval bisemigroup of finite order. $T = T_1 \cup T_2$
= {([0, $a_1$], [0, $a_2$], 0, 0, [0, $a_3$], [0, $a_4$]) / $a_i \in Z_{20}$, $1 \leq i \leq 4$} $\cup$

$$\left\{ \begin{pmatrix} 1 & 2 & 3 & 4 & 5 & 6 \\ 2 & 1 & 3 & 4 & 5 & 6 \end{pmatrix}, \begin{pmatrix} 1 & 2 & 3 & 4 & 5 & 6 \\ 1 & 2 & 3 & 4 & 5 & 6 \end{pmatrix}, \begin{pmatrix} 1 & 2 & 3 & 4 & 5 & 6 \\ 1 & 1 & 3 & 4 & 5 & 6 \end{pmatrix}, \begin{pmatrix} 1 & 2 & 3 & 4 & 5 & 6 \\ 2 & 2 & 3 & 4 & 5 & 6 \end{pmatrix} \right\}$$

$\subseteq S_1 \cup S_2$ be a semi interval bisubsemigroup of finite order. Clearly T is not an ideal of S.

*Example 1.1.46:* Let $S = S_1 \cup S_2 =$

$$\left\{ \begin{bmatrix} [0, a_1] \\ [0, a_2] \\ \vdots \\ [0, a_{10}] \end{bmatrix} \middle| a_i \in Q^+ \cup \{0\}; 1 \leq i \leq 10 \right\}$$

$\cup$ {All 10 × 10 matrices with entries from $Z_6$}, ($S_1$ under interval matrix multiplication) is a semi interval bisemigroup of infinite order. Take $I = I_1 \cup I_2$

$$= \left\{ \begin{bmatrix} [0, a] \\ [0, a] \\ \vdots \\ [0, a] \end{bmatrix} \middle| a \in Z^+ \cup \{0\} \right\}$$

$\cup$ {All 10 × 10 upper triangular matrices with entries from $Z_6$} $\subseteq S_1 \cup S_2$ is a semi interval subbisemigroup of S and I is not an ideal of S.

Now we will give examples of ideals in a semi interval bisemigroups.

*Example 1.1.47:* Let $S = S_1 \cup S_2 =$ {[0, a] / $a \in Z^+ \cup \{0\}$} $\cup$ {All 5 × 5 matrices with entries from $Z^+ \cup \{0\}$} be a semi interval bisemigroup of infinite order. $T = T_1 \cup T_2 = $ {[0, a] / a $\in 3Z^+ \cup \{0\}\} \cup \{$All $5 \times 5$ matrices with entries from $5Z^+ \cup \{0\}\} \subseteq S_1 \cup S_2$, is an ideal of S.

***Example 1.1.48:*** Let $S = S_1 \cup S_2 = \{([0, a_1], [0, a_2], \ldots, [0, a_7]) / a_i \in Z_{24}; 1 \leq i \leq 7\} \cup \{Z_{30}, \times\}$ be a semi interval bisemigroup of finite order. Consider $I = I_1 \cup I_2 = \{([0, a_1], [0, a_2], \ldots, [0, a_7]) / a_i \in \{0, 2, 4, 6, \ldots, 22\} \subseteq Z_{24}; 1 \leq i \leq 7\} \cup \{\{0, 10, 20\}, \times\} \subseteq S_1 \cup S_2$; I is an ideal of S.

***Example 1.1.49:*** Let $S = S_1 \cup S_2 = \{$All $6 \times 6$ matrices with entries from $Z^+ \cup \{0\}\} \cup \{[0, a] / a \in Q^+ \cup \{0\}\}$ be a semi interval bisemigroup. We see S has semi interval sub bisemigroups but no ideals. However $S_1$ contains ideals and $S_2$ has no ideals. We call $I = I_1 \cup \{0\}$ where $I_1$ is an ideal of $S_1$ as quasi semi interval biideal of S.
   We will first illustrate this situation by some examples.

***Example 1.1.50:*** Let $S = S_1 \cup S_2 = (Q^+ \cup \{0\}) \cup \{[0, a] / a \in Z^+ \cup \{0\}\}$ be a semi interval bisemigroup. Consider $J = J_1 \cup J_2 = \{0\} \cup \{[0, a] / a \in 7Z^+ \cup \{0\}\} \subseteq S_1 \cup S_2$; J is a quasi semi interval biideal of S.

***Example 1.1.51:*** Let $S = S_1 \cup S_2 = \{[0, a] / a \in Z_7\} \cup \{$all $12 \times 12$ matrices with entries from $Z_{12}\}$ be a semi interval bisemigroup. $T = T_1 \cup T_2 = \{0\} \cup \{12 \times 12$ matrices with entries from $\{0, 2, 4, 6, 8, 10\} \subseteq Z_{12}\} \subseteq S_1 \cup S_2$; T is a quasi semi interval biideal of S.

***Example 1.1.52:*** Let $S = S_1 \cup S_2 = \{Z_{11}\} \cup \{[0, a] / a \in Z^+ \cup \{0\}\}$ be a semi interval bisemigroup. $P = P_1 \cup P_2 = \{0\} \cup \{[0, a] \mid a \in 3Z^+ \cup \{0\}\} \subseteq S_1 \cup S_2$ is a quasi semi interval biideal of S.

***Example 1.1.53:*** Let $S = S_1 \cup S_2 = Q^+ \cup \{0\} \cup \{[0, a] / a \in Z_{13}\}$ be a semi interval bisemigroup, S has no ideals.

In view of this we can as in case of interval bisemigroups define simple semi interval bisemigroups and ideally simple semi interval bisemigroup; we only give examples of this.

***Example 1.1.54:*** Let $S = S_1 \cup S_2 = (Z_7, \times) \cup \{[0, a] / a \in Z_{13}\}$ be a semi interval bisemigroup. Clearly S is an ideally simple interval bisemigroup. However $P = P_1 \cup P_2 = \{0, 1, 6\} \cup \{[0, 1], 0, [0, 12]\} \subseteq S_1 \cup S_2$ is a semi interval bisubsemigroup of S which is not an ideal of S.

***Example 1.1.55:*** Let $S = S_1 \cup S_2 = \left\{ \begin{bmatrix} [0,a] \\ [0,a] \\ \vdots \\ [0,a] \end{bmatrix} \middle| a \in Z_{13} \right\} \cup$

$\{([0, a] [0, a] [0, a]) / a \in Z_{43}\}$ be a semi interval bisemigroup. Clearly S is a simple as well as ideally simple semi interval bisemigroup.

We can define bizero divisors, biidempotents and biunits in semi interval bisemigroups. We also can define quasi bizero divisors, quasi biunits and quasi biidempotents for semi interval bisemigroups as in case of interval bisemigroups. This task is a matter of routine and the reader is expected to do this job.

We will only give examples of these.

***Example 1.1.56:*** Let $S = S_1 \cup S_2 = \{Z_{12}, \times\} \cup \{[0, a] / a \in Z_{24}\}$ be a semi interval bisemigroup. $X = \{6\} \cup \{[0, 12]\}$ is a zero divisor for $Y = \{4\} \cup \{[0, 2]\}$ in S is such that $XY = 0 \cup 0$. Consider $x = \{4\} \cup \{0, 16]\} \in S$ is such that $x^2 = x = \{4\} \cup \{[0, 16]\}$, thus x is a biidempotent.

Consider $m = \{11\} \cup \{[0, 23]\} \in S$ is such that $m^2 = \{1\} \cup [0, 1]$, hence m is a biunit in S.

***Example 1.1.57:*** Let $S = S_1 \cup S_2 = \{[0, a] / a \in Z_{18}\} \cup \{Z^+ \cup \{0\}\}$ be a semi interval bisemigroup. Consider $x = [0, 17] \cup 0$ is a quasi biunit. Y = [0, 9] ∪ {0} is a quasi biidempotent in S. Consider α = [0, 6] ∪ {0} ∈ S; we have β = [0, 3] ∪ {0} ∈ S is such that α.β = 0 ∪ 0.

We call a semi interval bisemigroup $S = S_1 \cup S_2$ to be S-semi interval bisemigroup (Smarandache semi interval bisemigroup) if both $S_1$ and $S_2$ are Smarandache semigroups.

We will give examples of them.

*Example 1.1.58:* Let $S = S_1 \cup S_2 = \{[0, a] / a \in Z_{17}\} \cup \{2 \times 2$ matrices with entries from $Q^+ \cup \{0\}\}$ be a semi interval bisemigroup. Consider

$$A = \{[0, a] / a \in Z_{17} \setminus \{0\}\} \cup \left\{ \begin{bmatrix} a & 0 \\ 0 & b \end{bmatrix} \middle| b, a \in Q^+ \right\} \subseteq S_1 \cup S_2,$$

is a semi interval bigroup.

Hence S is a S-semi interval bisemigroup.

*Example 1.1.59:* Let $S = S_1 \cup S_2 = Z^+ \cup \{0\} \cup \{[0, a] / a \in Q^+ \cup \{0\}\}$ be a semi interval bisemigroup. S is not a S-semi interval bisemigroup as $S_1$ is not a S-semigroup.

In view of this we have the following definition. If a semi interval bisemigroup S has only one of $S_1$ or $S_2$ to be a S-semigroup then we define S to be a quasi Smarandache semi interval bisemigroup.

*Example 1.1.60:* Let $S = S_1 \cup S_2 = \{3Z^+ \cup \{0\}\} \cup \{[0, a] / a \in Z_{48}\}$ be a semi interval bisemigroup. Clearly S is only a quasi S-semi interval bisemigroup as $S_1$ is not a S-semigroup. $S_2$ contains $A_2 = \{[0, 1], [0, 47]\} \subseteq S_2$ is an interval group.

*Example 1.1.61:* Let $S = S_1 \cup S_2 = \{[0, a] / a \in 5Z^+ \cup \{0\}\} \cup \{Z_{30}\}$ be a semi interval bisemigroup. S is only a quasi S-semi interval bisemigroup as $S_1$ has no proper interval subgroup where as in $S_2$, $A_2 = \{1, 29\}$ to be a subset which is a group. That is $A_2 = \{1, 29\} \subseteq S_2$. Hence S is only a quasi S-semi interval bisemigroup. Now we can define ideals and subsemigroups for semi interval bisemigroups also. This task is left as an exercise for the reader.

In the next section we proceed onto define the notion of interval bigroupoids and discuss the properties associated with them.

## 1.2 Interval Bigroupoids

In this section we define interval bigroupoids and quasi interval bigroupoids. We also discuss the properties associated with them and analyse their substructures.

**DEFINITION 1.2.1:** *Let $G = G_1 \cup G_2$ where both $G_1$ and $G_2$ are interval groupoids; then we define G to be a interval bigroupoid or biinterval groupoid under the operation '·' inherited from $G_1$ and $G_2$.*

For more about interval groupoids please refer [7, 11, 14]. We will illustrate this situation by some examples.

*Example 1.2.1*: Let $G = G_1 \cup G_2 = \{[0, a], *, (2, 3), a \in Z_7\} \cup \{[0, b] / b \in Z^+ \cup \{0\}, *, (7, 10)\}$ be an interval bigroupoid or biinterval groupoid under the operation '.'; we see G is of infinite biorder. Further if $x = [0, 3] \cup [0, 20]$ and $y = [0, 2] \cup [0, 1] \in G = G_1 \cup G_2$;

$$\begin{aligned}
x.y &= ([0, 3] \cup [0, 20]) . ([0, 2] \cup [0, 1]) \\
&= [0, 3] * [0, 2] \cup [0, 20] * [0, 1] \\
&= [0, 6 + 6 \pmod 7] \cup [0, 20 \times 7 + 1 \times 10] \\
&= [0, 5] \cup [0, 150] \in G.
\end{aligned}$$

*Example 1.2.2:* Let $G = G_1 \cup G_2 = \{[0, a] / *, (3, 9), a \in Z_{20}\} \cup [0, b] / *, (2, 11), b \in Z_{14}\}$ be an interval bigroupoid with an operation '.' on G. Let x = [0, 3] ∪ [0, 7] and y = [0, 1] ∪ [0, 13] be in G.

x.y  =  ([0, 3] ∪ [0, 7]) · ([0, 1] ∪ [0, 13])
     =  [0, 3] * [0, 1] ∪ [0, 7] * [0, 13]
     =  [0, 9 + 9 (mod 20)] ∪ [0, 14+13 × 11 (mod 14)]
     =  [0, 18] ∪ [0, 3].

(G, .) is an biinterval groupoid of finite order. Clearly G is non associative and non commutative.

*Example 1.2.3:* Let G = $G_1 \cup G_2$ = {[0, a]| a ∈ $Z^+ \cup \{0\}$, *, (7, 2)} ∪ {[0, b]| b ∈ $Q^+ \cup \{0\}$, *, (17, 5)} be a biinterval groupoid of infinite order. G is both non associative and non commutative.

*Example 1.2.4:* Let G = $G_1 \cup G_2$ = {[0, a] / a ∈ $Z^+ \cup \{0\}$, *, (4, 7)} ∪ {[0, b] / b ∈ $Z^+ \cup \{0\}$, *, (21, 17)} be a interval bigroupoid of infinite order.

*Example 1.2.5:* Let G = $G_1 \cup G_2$ = {[0, a] / a ∈ $Z_{14}$, *, (10, 0)} ∪ {[0, b] / b ∈ $Z_{14}$, *, (2, 5)} be an interval bigroupoid of finite order. G is non commutative.

*Example 1.2.6:* Let G = $G_1 \cup G_2$ = {[0, a] / a ∈ $Z_{11}$, *, (3, 7)} ∪ {[0, b] / b ∈ $Z_{13}$, *, (3, 7)} be an interval bigroupoid of finite order.

Now we have seen examples of interval bigroupoids.

We proceed onto give examples of substructures in them. It is important and interesting to note that it is not very easy to find substructures for the operation '*' defined on them is in a complicated way to make the binary operation non associative.

*Example 1.2.7:* Let G = $G_1 \cup G_2$ = {[0, a] / a ∈ $Z^+ \cup \{0\}$, *, (8, 12)} ∪ {[0, a] / a ∈ $Z^+ \cup \{0\}$, *, (2, 6)} be an interval bigroupoid. Consider H = $H_1 \cup H_2$ = {[0, a] / a ∈ $2Z^+ \cup \{0\}$} ∪ $\{[0, b] / b \in 2Z^+ \cup \{0\}, *, (2, 6)\} \subseteq G_1 \cup G_2$; H is an interval subgroupoid of G.

Suppose we have $G_1 = \{[0, a] / a \in 3Z^+ \cup \{0\}, *, (3/2, 9)\}$ is not a groupoid for if we consider [0, 3], [0, 4] in $G_1$

$$[0, 3] * [0, 4] = [0, 3.3/2 + 4.3]$$
$$= [0, 9/2 + 12]$$
$$= [0, 33/2] \notin G_1.$$

So while choosing the pair which operates on $G_1$ we need them to be present in $G_1$. For if this point is not taken into account we will have problem about the closure of the operation defined on $G_1$.

***Example 1.2.8:*** Let $G = G_1 \cup G_2 = \{[0, a] / a \in Z_6, *, (2, 4)\} \cup \{[0, b] / b \in Z_8, *, (2, 6)\}$ be an interval bigroupoid. Take $P = P_1 \cup P_2 = \{[0, a] / a \in \{0, 2, 4\} \subseteq Z_6, *, (2, 4)\} \cup \{[0, b] | b \in \{0, 2, 4, 6\} \subseteq Z_8, *, (2, 6)\} \subseteq G_1 \cup G_2$, P is an interval bigroupoid of G.

***Example 1.2.9:*** Let $G = G_1 \cup G_2 = \{[0, a] / a \in Z_{12}, *, (2, 10)\} \cup \{[0, b] / b \in Z_{12}, *, (10, 8)\}$ be an interval bigroupoid. $P = P_1 \cup P_2 = \{[0, a] / a \in \{0, 2, 4, 6, 8, 10\} \subseteq Z_{12}, * (2, 10)\} \cup \{[0, b] / b \in \{0, 4, 8\} \subseteq Z_{12}, (10, 8)\} \subseteq G_1 \cup G_2$ is such that P is an interval subbigroupoid.

Now we can define on an interval bigroupoid $G = G_1 \cup G_2$ the notion of P-interval bigroupoid, alternative interval bigroupoid and so on.

We say an interval bigroupoid $G = G_1 \cup G_2$ to be a P-interval bigroupoid if both $G_1$ and $G_2$ are P-interval groupoids.

If only one of $G_1$ or $G_2$ is a P-interval groupoid and the other is not a P-interval groupoid then we define $G = G_1 \cup G_2$ to be a quasi P-interval bigroupoid.

We will give examples of both the situations.

***Example 1.2.10:*** Let $G = G_1 \cup G_2 = \{[0, a] / a \in Z_{45}, *, (7, 7)\} \cup \{[0, b] / a \in Z_{15}, *, (4, 4)\}$ is an interval bigroupoid which is a P-interval bigroupoid.

***Example 1.2.11:*** Let $G = G_1 \cup G_2 = \{[0, a] / a \in Z_{19}, *, (3, 3)\} \cup \{[0, b] / b \in Z_{21}, *, (5, 5)\}$ be an interval bigroupoid. G is an interval P-bigroupoid.

We have the following theorem the proof of which is left as an exercise to the reader.

**THEOREM 1.2.1:** *Let $G = G_1 \cup G_2 = \{[0, a] / a \in Z_n, *, (t, t)\} \cup \{[0, b] / b \in Z_m, *, (u, u)\}$ $(m \neq n)$ be an interval bigroupoid. G is a P-interval bigroupoid.*

***Example 1.2.12:*** Let $G = G_1 \cup G_2 = \{[0, a] / a \in Z_{25}, *, (2, 2)\} \cup \{[0, a] / a \in Z_{19}, *, (3, 0)\}$ be an interval bigroupoid. Clearly G is not a P-interval bigroupoid. It is only a quasi interval P-bigroupoid.

***Example 1.2.13:*** Let $G = G_1 \cup G_2 = \{[0, a] / a \in Z_{23}, *, (4, 4)\} \cup \{[0, b] / b \in Z_{17}, *, (8, 0)\}$ be an interval bigroupoid. Clearly G is a quasi P-interval bigroupoid.

Now we have some results the proofs which are direct.

**THEOREM 1.2.2:** *Let $G = G_1 \cup G_2 = \{[0, a] / a \in Z_n, *, (0, t)\} \cup \{[0, b] / b \in Z_m, *, (u, 0)\}$ be an interval bigroupoid. G is a P-interval groupoid if and only if $t^2 = t \mod n$ and $u^2 = u \mod m$.*

**THEOREM 1.2.3:** *Let $G = G_1 \cup G_2 = \{[0, a] / a \in Z_p, *, (t, 0)\};$ p a prime$\} \cup \{[0, b] / b \in Z_n, *, (r, 0)\}$ be an interval bigroupoid. G is a quasi interval P-bigroupoid if and only if $r^2 \equiv r \mod n$.*

**THEOREM 1.2.4:** *Let $G = G_1 \cup G_2 = \{[0, a] / a \in Z_n, (t, t), *\} \cup \{[0, b] / b \in Z_p, *, (u, 0), p \text{ a prime}\}$ be an interval bigroupoid. G is a quasi interval P-bigroupoid.*

*Example 1.2.14:* Let $G = G_1 \cup G_2 = \{[0, a] / a \in Z_{11}, *, (8, 0)\} \cup \{[0, b] / b \in Z_{43}, *, (10, 0)\}$ be an interval bigroupoid. G is not an interval P-bigroupoid or quasi interval P-bigroupoid.

We define an interval bigroupoid $G = G_1 \cup G_2$ to be an alternative interval bigroupoid if both $G_1$ and $G_2$ are interval alternative groupoid. If only one of $G_1$ or $G_2$ is an interval alternative groupoid then we define G to be a quasi alternative interval bigroupoid.

We will give examples of this situation.

*Example 1.2.15:* Let $G = G_1 \cup G_2 = \{[0, a] / a \in Z_{24}, *, (9, 9)\} \cup \{[0, a] / a \in Z_{12}, *, (4, 4)\}$ be an alternative interval bigroupoid.

*Example 1.2.16:* Let $G = G_1 \cup G_2 = \{[0, a] / a \in Z_{24}, *, (9, 0)\} \cup \{[0, b] / b \in Z_{36}, *, (9, 9)\}$ be an alternative interval bigroupoid.

We have the following results which are direct.

**THEOREM 1.2.5:** *Let $G = G_1 \cup G_2 = \{[0, a] / a \in Z_n, *, (t, t)\} \cup \{[0, b] / b \in Z_m, *, (u, u)\}$ be an interval bigroupoid. G is an alternative interval bigroupoid if and only if $t^2 = t \mod n$ and $u^2 = u \mod m$.*

**THEOREM 1.2.6:** *Let $G = G_1 \cup G_2 = \{[0, a] / a \in Z_n, *, (t, 0)\} \cup \{[0, b] / b \in Z_m, *, (u, 0)\}$ be an interval bigroupoid G is an alternative interval bigroupoid if and only if $t^2 = t \mod n$ and $u^2 = u \mod m$.*

**THEOREM 1.2.7:** *Let $G = G_1 \cup G_2 = \{[0, a] / a \in Z_n, *, (0, t)\}$ $\cup \{[0, b] / b \in Z_m, *, (u, u)\}$ be an interval bigroupoid. G is an alternative interval bigroupoid if and only if $t^2 = t$ mod n and $u^2 = u$ mod m.*

We also have a class of interval bigroupoids which are not alternative interval bigroupoids, which is given by the following theorem.

**THEOREM 1.2.8:** *Let $G = G_1 \cup G_2 = \{[0, a] / a \in Z_p, *, (t, t); p$ a prime$\} \cup \{[0, b] / b \in Z_p, *, (u, u); q$ a prime$\}$ be an interval bigroupoid. G is not an alternative interval bigroupoid.*

Now we will proceed onto give examples and results related with quasi interval alternative bigroupoids.

*Example 1.2.17:* Let $G = G_1 \cup G_2 = \{[0, a] / a \in Z_{19}, *, (9, 0)\}$ $\cup \{[0, a] / a \in Z_{24}, *, (9, 9)\}$ be an interval bigroupoid. Clearly G is not an alternative interval bigroupoid but G is a quasi alternative interval bigroupoid.

*Example 1.2.18:* Let $G = G_1 \cup G_2 = \{[0, a] / a \in Z_{47}, *, (8, 8)\}$ $\cup \{[0, a] / a \in Z_{12}, *, (4, 4)\}$ be an interval bigroupoid. Clearly G is a quasi alternative interval bigroupoid.

**THEOREM 1.2.9:** *$G = G_1 \cup G_2 = \{[0, a] / a \in Z_p, *, (t, t), p$ a prime$\} \cup \{[0, b] / a \in Z_n, *, (u, u)\}$ is a quasi interval alternative bigroupoid if and only if $u^2 = u$ mod n.*

**THEOREM 1.2.10:** *$G = G_1 \cup G_2 = \{[0, a] / a \in Z_p, *, (t, 0), p$ a prime$\} \cup \{[0, b] / b \in Z_n, (u, 0)\}$ is a quasi alternative interval bigroupoid if and only if $u^2 = u \pmod{n}$.*

We call an interval bigroupoid $G = G_1 \cup G_2$ to be an idempotent interval bigroupoid if both $G_1$ and $G_2$ are idempotent interval groupoids. If only one of $G_1$ or $G_2$ is an idempotent interval groupoid then we define G to be a quasi idempotent interval bigroupoid.

We will first give examples of this situation.

***Example 1.2.19:*** Let $G = G_1 \cup G_2 = \{[0, a] / a \in Z_{12}, *, (7, 6)\} \cup \{[0, b] / b \in Z_{15}, *, (9, 6)\}$ be an interval bigroupoid. Clearly G is an idempotent interval bigroupoid.

***Example 1.2.20:*** $G = G_1 \cup G_2 = \{[0, a] / a \in Z_{19}, *, (16, 4)\} \cup \{[0, a] / a \in Z_{13}, (8, 6)\}$ is an idempotent interval bigroupoid.

***Example 1.2.21:*** Let $G = G_1 \cup G_2 = \{[0, a] / a \in Z_{15}, *, (2, 3)\} \cup \{[0, b] / b \in Z_{11}, *, (3, 5)\}$ be an quasi idempotent interval bigroupoid.

We state a few theorems the proofs of which are direct.

**THEOREM 1.2.11:** *$G = G_1 \cup G_2 = \{[0, a] / a \in Z_n, *, (t, u)\} \cup \{[0, b] / b \in Z_m, *, (r, s)\}$ is an idempotent interval bigroupoid if and only if $t + u = 1 \mod n$ and $r + s = 1 \mod m$.*

**THEOREM 1.2.12:** *$G = G_1 \cup G_2 = \{[0, a] / a \in Z_n, *, (t, u); t + u \not\equiv 1 \mod n\} \cup \{[0, b] / b \in Z_m, *, (r, s)\}$ is a quasi interval idempotent bigroupoid if and only if $r + s = 1 \mod m$.*

We say an interval bigroupoid is bisimple if and only if both $G_1$ and $G_2$ are simple.

We say quasi bisimple only one of $G_1$ or $G_2$ is simple.

***Example 1.2.22:*** Let $G = G_1 \cup G_2 = \{[0, a] / a \in Z_{12}, *, (5, 7)\} \cup \{[0, b] / b \in Z_{20}, *, (7, 13)\}$ be an interval bigroupoid. G is a bisimple interval bigroupoid.

***Example 1.2.23:*** Let $G = G_1 \cup G_2 = \{[0, a] / a \in Z_{19}, (2, 17)\} \cup \{[0, b] / b \in Z_{12}, *, (11, 2)\}$ be an interval bigroupoid. Clearly G is an interval bisimple bigroupoid.

*Example 1.2.24:* Let $G = G_1 \cup G_2 = \{[0, a] / a \in Z_{14}, *, (2, 3)\} \cup \{[0, a] / a \in Z_{20}, *, (17, 3)\}$ be a quasi bisimple interval bigroupoid.

*Example 1.2.25:* Let $G = G_1 \cup G_2 = \{[0, a] / a \in Z_{13}, *, (9, 4)\} \cup \{[0, b] / a \in Z_{36}, *, (31, 5)\}$ be a quasi bisimple interval bigroupoid.

We give a few results and expect the reader to supply the proof.

**THEOREM 1.2.13:** *Let $G = G_1 \cup G_2 = \{[0, a] / a \in Z_n, *, (t, u)\} \cup \{[0, b] / b \in Z_m, *, (r, s)\}$ be an interval bigroupoid. If $n = t + u$ and $m = r + s$ and $t, u, r$ and $s$ are primes then $G$ is a bisimple interval bigroupoid.*

**THEOREM 1.2.14:** *Let $G = G_1 \cup G_2 = \{[0, a] / a \in Z_n, *, (t, u); (t, u) = d; d \neq 1\} \cup \{[0, a] / a \in Z_m, *, (r, s)\}$ be an interval bigroupoid. If $r + s = m$ and $r$ and $s$ are primes then $G$ is a quasi bisimple interval bigroupoid.*

**THEOREM 1.2.15:** *Let $G = G_1 \cup G_2 = \{[0, a] / a \in Z_n, *, (t, u)\} \cup \{[0, b] / b \in Z_m, *, (r, s)\}$ be an interval bigroupoid. $G$ is an interval bisemigroupoid if and only if $(t, u) = 1$ and $t^2 \equiv t \bmod n$ and $u^2 \equiv u \bmod n$ and $(r, s) = 1$ with $r^2 = r \bmod m$ and $s^2 = s \bmod m$ ($t \neq 0, u \neq 0, r \neq 0$ and $s \neq 0$).*

Now we proceed onto define Smarandache interval bigroupoid $G = G_1 \cup G_2$. We say $G$ is a Smarandache interval bigroupoid (S-interval bigroupoid) if each $G_i$ is a S-interval groupoid; i=1,2.

We say quasi Smarandache interval bigroupoid if only one if one of $G_1$ or $G_2$ is a S-interval groupoid.

*Example 1.2.26:* Let $G = G_1 \cup G_2 = \{[0, a] / a \in Z_8, *, (2, 6)\} \cup \{[0, b] / b \in Z_5, *, (3, 3)\}$ be a S-interval bigroupoid. For A =

$\{[0, a] / a \in \{0, 4\} \subseteq Z_8, *, (2, 6)\} \cup \{[0, b] / b \in \{4\} \subseteq Z_5\} = A_1 \cup A_2 \subseteq G_1 \cup G_2$ is a bisemigroup in G.

***Example 1.2.27:*** Let $G = G_1 \cup G_2 = \{[0, a] / a \in Z_5, *, (1, 3)\} \cup \{[0, a] / a \in Z_5, *, (2, 1)\}$ be an interval bigroupoid. G is not a S-interval bigroupoid as well as not a quasi S-interval bigroupoid.

***Example 1.2.28:*** Let $G = G_1 \cup G_2 = \{[0, a] / a \in Z_{10}, *, (1, 2)\} \cup \{[0, b] / b \in Z_5, *, (2, 1)\}$ be a quasi Smarandache interval bigroupoid.

In view of this we give the following results.

**THEOREM 1.2.16:** *Let $G = G_1 \cup G_2 = \{[0, a] / a \in Z_n, *, (t, u), (t, u) = 1 ; t \neq u; t + u \equiv 1 \pmod{n}\} \cup \{[0, b] / b \in Z_m, *, (r, s), (r, s) = 1, r \neq s, r + s \equiv 1 \pmod{m}\}$ ($m > 5$ and $n > 5$) be an interval bigroupoid. G is a S-biinterval groupoid.*

Now as in case in usual bigroupoids we can define Smarandache interval bigroupoid, this task exercise for the reader.

***Example 1.2.29:*** Let $G = G_1 \cup G_2 = \{[0, a] / a \in Z_8, *, (2, 6)\} \cup \{[0, b] / b \in Z_{12}, *, (3, 9)\}$ be an interval bigroupoid $P = P_1 \cup P_2 = \{[0, a] / a \in \{0, 2, 4, 6\} \subseteq Z_8, *, (2, 6)\} \cup \{[0, b] / b \in \{0, 3, 6, 9\} \subseteq Z_{12}, *, (3, 9)\} \subseteq G_1 \cup G_2$ is a Smarandache interval subbigroupoid of G.

The notion of Smarandache interval P-subbigroupoid, Smarandache strong P-interval groupoid, Smarandache Bol interval bigroupoid, Smarandache strong Bol interval bigroupoid, Smarandache right alternative interval bigroupoid so on can be defined as in case of usual groupoids [7, 11, 14].

We will give examples and related results for the reader. The definition is a matter of routine.

***Example 1.2.30:*** Let $G = G_1 \cup G_2 = \{[0, a] / a \in Z_{12}, *, (3, 4)\} \cup \{[0, b] / b \in Z_4, *, (2, 3)\}$ be an interval bigroupoid. It is easily verified G is a Smarandache interval Bol bigroupoid.

***Example 1.2.31:*** Let $G = G_1 \cup G_2 = \{[0, a] / a \in Z_6, *, (4, 3)\} \cup \{[0, b] / b \in Z_6, *, (3, 5)\}$ be an interval groupoid. G is a Smarandache interval P-bigroupoid.

***Example 1.2.32:*** Let $G = G_1 \cup G_2 = \{[0, a] / a \in Z_{10}, *, (5, 6)\} \cup \{[0, b] / b \in Z_{12}, *, (3, 9)\}$ be an interval bigroupoid. G is a Smarandache Moufang interval bigroupoid.

***Example 1.2.33:*** Let $G = G_1 \cup G_2 = \{[0, a] / a \in Z_{11}, *, (6, 6)\} \cup \{[0, b] / b \in Z_{19}, *, (10, 10)\}$ be an interval bigroupoid. It is easily verified G is a Smarandache idempotent bigroupoid.

We will now state a theorem the proof of which is left as an exercise for the reader.

**THEOREM 1.2.17:** *Let $G = G_1 \cup G_2 = \{[0, a] / *, a \in Z_p, \left(\frac{p+1}{2}, \frac{p+1}{2}\right)\} \cup \{[0, b] / *, b \in Z_q, \left(\frac{q+1}{2}, \frac{q+1}{2}\right)\}$, where p and q are two distinct primes, then G is a Smarandache idempotent interval bigroupoid.*

**THEOREM 1.2.18:** *Let $G = G_1 \cup G_2 = \{[0, a] / a \in Z_n, *, (m, m)\} \cup \{[0, a] / a \in Z_s, *, (t, t)\}$ where $t \neq n$ and $m + m \equiv 1 \pmod n$, $t + t \equiv 1 \pmod t$ and $m^2 = m \pmod n$ and $t^2 = t \pmod s$. G is a*

1. *Smarandache idempotent interval bigroupoid.*
2. *S-strong P-interval bigroupoid .*
3. *Smarandache strong interval Bol bigroupoid.*
4. *Smarandache strong Moufang interval bigroupoid.*
5. *Smarandache strong interval alternative bigroupoid.*

**THEOREM 1.2.19:** *Let $G = G_1 \cup G_2$ be a S-strong alternative interval bigroupoid then G is a S-alternative interval bigroupoid.*

In this theorem 1.2.19 if we replace the S-strong alternative interval groupoid. G by S-strong Moufang or S-strong P-bigroupoid or S-strong idempotent bigroupoid or S-strong Bol bigroupoid then G will be a S-Moufang or S-P-groupoid or S-idempotent bigroupoid or S-Bol interval bigroupoid.

**THEOREM 1.2.20:** *$G = G_1 \cup G_2 = \{[0, a] / a \in Z_n, *, (t, u), (t + u) \equiv 1 \bmod n\} \cup \{[0, b] / b \in Z_m, *, (r, s), (r + s) \equiv 1 \bmod m\}$, $m \neq n$ is a S-strong Moufang interval bigroupoid, S-interval idempotent bigroupoid, S-strong Bol interval bigroupoid and S-strong alternative interval bigroupoid if and only if $t^2 \equiv t \pmod{n}$, $u^2 \equiv u \pmod{n}$, $r^2 \equiv r \pmod{m}$ and $s^2 \equiv s \pmod{m}$.*

**THEOREM 1.2.21:** *Let $G = G_1 \cup G_2 = \{[0, a] / a \in Z_n, *, (m, 0)\} \cup \{[0, b] / b \in Z_t, *, (r, 0)\}$, $t \neq n$ with $m^2 \equiv m \pmod{n}$ and $r^2 = r \pmod{t}$. Then G is a S-strong Bol interval bigroupoid, S-strong Moufang interval bigroupoid, S-strong P-interval bigroupoid and S-strong alternative interval bigroupoid.*

Several interesting results in this direction can be obtained by any interested reader.

Now we can define quasi interval bigroupoid, as a bigroupoid $G = G_1 \cup G_2$ where only one of $G_1$ or $G_2$ is an interval groupoid and the other is just a groupoid.

We will illustrate this situation by some examples.

*Example 1.2.34:* Let $G = G_1 \cup G_2 = \{[0, a] / a \in Z_7, *, (3, 2)\} \cup Z_{11}(5, 4)$ be a quasi interval bigroupoid. Clearly G is of finite order.

***Example 1.2.35:*** Let $G = G_1 \cup G_2 = \{[0, a] / a \in Z^+ \cup \{0\}, *, (7, 8)\} \cup Z_{17} (3, 8)$ be a quasi interval bigroupoid of infinite order.

***Example 1.2.36:*** Let $G = G_1 \cup G_2 = \{[0, a] / a \in Z_{27}, *, (3, 8)\} \cup \{Z_{27} (19, 4)\}$ be a quasi interval bigroupoid of finite order. Now we can find substructure in them.

***Example 1.2.37:*** Let $G = G_1 \cup G_2 = \{[0, a] / a \in Z_6, *, (2, 2)\} \cup \{Z_6 (0, 2)\}$ be a quasi interval bigroupoid. $P = P_1 \cup P_2 = \{[0, a] / a \in \{0, 2, 4\} *, (2, 2)\} \cup \{\{0, 2, 4\} \subseteq Z_6, *, (0, 2)\}$ is a quasi interval subbigroupoid of G.

***Example 1.2.38:*** Let $G = G_1 \cup G_2 = \{[0, a] / a \in Z_{12}, *, (3, 9)\} \cup \{Z_6 (5, 5)\}$ be a quasi interval bigroupoid. $S = S_1 \cup S_2 = \{[0, a] / a \in \{0, 3, 6, 9\}, *, (3, 9)\} \cup \{\{4\} \subseteq Z_6, *, (5, 5)\} \subseteq G_1 \cup G_2$ is a quasi interval subbigroupoid of G.

***Example 1.2.39:*** Let $G = G_1 \cup G_2 = \{[0, a] / a \in Z_{12}, *, (2, 10)\} \cup \{Z_{12} (10,8)\}$ be a quasi interval bigroupoid. Consider $P = P_1 \cup P_2 = \{[0, a] / a \in \{0, 2, 4, 6, 8, 10\} \subseteq Z_{12}, *, (2, 10)\} \cup \{\{2, 6, 10\} \subseteq Z_{12}, *, (10, 8)\} \subseteq G_1 \cup G_2$, P is a quasi interval subbigroupoid of G.

Now we say a quasi interval bigroupoid is a Smarandache quasi interval bigroupoid if it has quasi interval bisemigroup. Likewise one can derive results for other S-structures associated with identities. We will illustrate these situations by some examples.

***Example 1.2.40:*** Let $G = G_1 \cup G_2 = \{[0, a] / a \in Z_{18}, *, (0, 3)\} \cup \{Z_{18} (0, 11)$ be a quasi interval bigroupoid. Clearly G is not a quasi P-interval bigroupoid or a quasi alternative interval bigroupoid.

***Example 1.2.41:*** Let $G = G_1 \cup G_2 = \{[0, a] / a \in Z_{12}, *, (0, 4)\}$ $\cup \{Z_6 (0, 3)\}$ be a quasi P interval bigroupoid as well as quasi alternative interval bigroupoid.

***Example 1.2.42:*** Let $G = G_1 \cup G_2 = \{[0, a] / a \in Z_{12}, *, (9, 0)\}$ $\cup \{Z_9 (0, 10)\}$ be a quasi P-interval bigroupoid as well as quasi alternative interval bigroupoid.

We will give a theorem which gurantees the existence of quasi P-interval bigroupoid and quasi alternative interval bigroupoid.

**THEOREM 1.2.22:** *Let $G = G_1 \cup G_2 = \{[0, a] / a \in Z_n, *, (0, t)\} \cup \{Z_m (0, r)\}$ be a quasi interval bigroupoid. G is a quasi P-interval bigroupoid and quasi interval alternative bigroupoid if and only if $t^2 \equiv t \pmod{n}$ and $r^2 \equiv r \pmod{m}$, $m \neq n$.*

**COROLLARY 1.2.1:** *In the above theorem if n and m are primes then the statement of the theorem is not true.*

We also have a class of quasi normal interval bigroupoid.

***Example 1.2.43:*** Let $G = G_1 \cup G_2 = \{[0, a] / a \in Z_{13}, *, (7, 7)\}$ $\cup \{Z_{19} (2, 2)\}$ be a quasi interval bigroupoid. Clearly G is a quasi interval normal bigroupoid.

***Example 1.2.44:*** Let $G = G_1 \cup G_2 = \{[0, a] / a \in Z_{43}, (17, 17)\}$ $\cup Z_{23} (4, 4)$ be a quasi normal interval bigroupoid.
   In view of this we have the following theorem.

**THEOREM 1.2.23:** *Let $G = G_1 \cup G_2 = Z_p (t, t) \cup \{[0, a] / a \in Z_q, *, (u, u)\}$ p and q are primes with $t < p$ and $u < q$, then $G = G_1 \cup G_2$ is a quasi normal interval bigroupoid.*

The proof is easy and hence is left as an exercise for the reader.

**THEOREM 1.2.24:** *Let $G = G_1 \cup G_2 = Z_n (t, t) \cup \{[0, b] / b \in Z_m, *, (u, u)\}$ be a quasi interval bigroupoid. G is a quasi interval P-bigroupoid.*

**THEOREM 1.2.25:** *Let $G = G_1 \cup G_2 = Z_n (t, t) \cup \{[0, a] / a \in Z_m, *, (u, u)\}$, $1 < t < n$ and $1 < u < m$ be a quasi interval bigroupoid, G is a not a quasi alternative interval bigroupoid if m and n are primes.*

In view of this we have the following corollary.

**COROLLARY 1.2.2:** *Let $G = G_1 \cup G_2 = Z_n (t, t) \cup \{[0, b] / b \in Z_m, *, (u, u)\}$, n and m are not primes be a quasi interval bigroupoid. G is a quasi alternative interval bigroupoid if and only if $t^2 = t \pmod{n}$ and $u^2 = u \pmod{m}$.*

*Proof:* Follows from the fact
$$t^2x + t^2y + tx = tx + t^2y + t^2x \bmod n$$
and
$$[0, u^2b] + [0, u^2d] + [0, ub] = [0, ub] + [0, u^2d] + [0, u^2b],$$
$$\text{i.e., } [0, u^2b + u^2d + ub] = [0, ub + u^2d + u^2b]$$
where
$$([0, b] * [0, d]) * ([0, b]) = ([0, b] * ([0, d] * [0, b]).$$
For more about these properties refer [7, 11, 14].

*Example 1.2.45:* Let $G = G_1 \cup G_2 = \{[0, a] / a \in Z_{12}, *, (7, 5)\} \cup \{Z_{24} (13, 11)\}$ be a quasi interval bigroupoid. Clearly G is a simple quasi interval bigroupoid.

*Example 1.2.46:* Let $G = G_1 \cup G_2 = \{[0, a] / a \in Z_7 (2, 5)\} \cup \{Z_{20} (7, 13)\}$ be a simple quasi interval bigroupoid.

**THEOREM 1.2.26:** *Let $G = G_1 \cup G_2 = Z_n (t, p) \cup \{[0, a] / a \in Z_m; *, (q, r)\}$ be a quasi interval bigroupoid. If $t + p = n$ and $q + r = m$ where t, p, q and r are primes then G is a simple quasi interval bigroupoid.*

**THEOREM 1.2.27:** *Let $G = G_1 \cup G_2 = \{[0, a] / a \in Z_n, *, (t, u) = 1; t, u \in Z_n \setminus \{0\}\} \cup \{Z_m (r, s) / (r, s) = 1 \; r, s, \in Z_m \setminus \{0\}\}$ be a quasi interval bigroupoid. $\{0\} \cup \{0\}$ is not a biideal of G.*

**THEOREM 1.2.28:** *Let $G = G_1 \cup G_2 = \{[0, a] / a \in Z_n, *, (t, u); t, u \in Z_n \setminus \{0\} \text{ and } (t, u) = 1\} \cup \{Z_m (r, s) / r, s \in Z_m \setminus \{0\} \text{ and } (r, s) = 1\}$ be a quasi interval bigroupoid. G is a quasi interval bisemigroup if and only if $t^2 \equiv t \pmod{n}$, $u^2 \equiv u \pmod{n}$, $r^2 \equiv r \pmod{n}$ and $s^2 \equiv s \pmod{n}$.*

Now we give examples of these situations.

*Example 1.2.47:* Let $G = G_1 \cup G_2 = \{[0, a] / a \in Z_{20}, *, (8, 12)\} \cup \{Z_{15}, (7, 8)\}$ be a quasi interval bigroupoid. G is a quasi idempotent interval bigroupoid.

*Example 1.2.48:* Let $G = G_1 \cup G_2 = \{[0, a] / a \in Z_{10}, *, (5, 6)\} \cup \{Z_{12} (4, 9)\}$ be a quasi interval bigroupoid. Clearly G is a quasi interval bisemigroup.

*Example 1.2.49:* Let $G = G_1 \cup G_2 = \{[0, a] / a \in Z_4, *, (2, 3)\} \cup \{Z_4 (3, 2)\}$ be the quasi interval bigroupoid. G has both left ideals as well as right ideals.

In view of this we have the following theorem.

**THEOREM 1.2.29:** *Let $G = G_1 \cup G_2 = \{[0, a] / a \in Z_n, *, (t, u)\} \cup \{Z_m (p, q)\}$ be a quasi interval bigroupoid. $P = P_1 \cup P_2$ is a right biideal of G if and only if $P = P_1 \cup P_2$ is a left biideal of $G' = G'_1 \cup G'_2 = \{[0, a] / a \in Z_n, *, (u, t)\} \cup \{Z_n (q, p)\}$.*

We as in case of interval bigroupoid can define S-interval quasi bigroupoid or S-quasi interval bigroupoid. The definition is a matter of routine, so we will give only examples of them.

*Example 1.2.50:* Let $G = G_1 \cup G_2 = \{[0, a] / a \in Z_{12}, *, (1, 3)\} \cup \{Z_5 (3, 3)$ be a quasi interval bigroupoid. Clearly G is a Smarandache quasi interval bigroupoid as $P = P_1 \cup P_2 = \{[0, a] / a \in \{0, 6\}, *, (1, 3)\} \cup \{a = \{4\}, *, (3, 3), 4 \in Z_5\} \subseteq G_1 \cup G_2$ is a quasi interval bisemigroup. Hence the claim.

*Example 1.2.51:* Let $G = G_1 \cup G_2 = \{[0, a] / a \in Z_6, *, (4, 5)\} \cup \{Z_6, (2, 4)\}$ be a quasi interval bigroupoid. Consider $P = P_1 \cup P_2 = \{[0, a] / a \in \{1, 3, 5\} \subseteq Z_6, *, (4, 5)\} \cup \{(0, 2, 4) \subseteq Z_6, *, (2, 4)\} \subseteq G_1 \cup G_2$, P is an interval biideal of G.

However P is not a S-biideal of G.

Several other concepts like S-semi normal S-semiconjugate and S-normal can be defined for quasi interval bigroupoids. We can define identities in these special type of bigroupoids also. This task is also a matter of routine hence left to the reader. However we will give examples of them.

*Example 1.2.52:* Let $G = G_1 \cup G_2 = \{[0, a] / a \in Z_{10}, *, (5, 6)\} \cup \{Z_{12}, (3, 9)\}$ be a quasi interval Smarandache Moufang bigroupoid.

*Example 1.2.53:* Let $G = G_1 \cup G_2 = \{[0, a] / a \in Z_{12}, *, (3, 4)\} \cup \{Z_4 (2,3)\}$ be a quasi Smarandache Bol interval bigroupoid.

*Example 1.2.54:* Let $G = G_1 \cup G_2 = \{[0, a] / a \in Z_6, *, (4, 3)\} \cup \{Z_6 (3, 5)\}$ be a quasi Smarandache P-interval bigroupoid of finite order.

*Example 1.2.55:* Let $G = G_1 \cup G_2 = \{[0, a] / a \in Z_{14}, *, (7, 8)\} \cup \{Z_{12} (1, 6)\}$ be a Smarandache quasi alternative interval bigroupoid.

Now results related with interval bigroupoids can be got in an analogous way for quasi interval bigroupoids also. This task is left as an exercise to the reader. Now we define interval semigroup interval groupoid which we term as biinterval groupoid-semigroup or biinterval semigroup-groupoid.

**DEFINITION 1.2.2:** *Let $G = S \cup G_1$ where $(S, .)$ is an interval semigroup and $(G_1, *)$ is an interval groupoid. $(G, o)$ which inherits the operations of $S$ and $G_1$ is defined to be a biinterval semigroup groupoid. The operation 'o' is defined as follows.*

*First of all $G = S \cup G_1 = \{[0, s] \, / \, [0,s] \in S\} \cup \{[0, g_1] \, / \, [0,g_1] \in G_1\} = \{[0, s] \cup [0, g_1] \, / \, [0, s] \in S$ and $[0, g_1] \in G_1\}$.*

*Let $x = [0, s] \cup [0, g]$ and $y = [0, t] \cup [0, h]$ be in G.*

$$x \, o \, y = ([0, s] \cup [0, g]) \, o \, ([0, t] \cup [0, h])$$

*$= \{([0, s] . [0, t] \cup [0, g] * [0, h])\} = [0, st] \cup [0, g*h] \in G$ as $[0, st] \in S$ and $[0, g*h] \in G_1$.*

Elements of G are biintervals so we define G as biinterval semigroup- groupoid. We will first illustrate this situation by some examples.

***Example 1.2.56:*** Let $G = G_1 \cup S = \{[0, a] \, / \, a \in Z_{14}, *, (3, 7)\} \cup \{[0, t] \, / \, t \in Z^+ \cup \{0\}, \times \}$ be a biinterval groupoid-semigroup.

***Example 1.2.57:*** Let $G = S \cup G_2 = \{[0, a] \, / \, a \in Z_{15}$, under multiplication modulo 15$\} \cup \{[0, b] \, / \, b \in Z_{18}, *, (5, 7)\}$ be a biinterval semigroup-groupoid. We see G is of finite order and is non commutative.

***Example 1.2.58:*** Let $G = G_1 \cup G_2 = \{[0, a] \, / \, a \in Z_{20}$, multiplication under modulo 20$\} \cup \{[0, b] \, / \, b \in Z_{20}, *, (3, 7)\}$ be a biinterval semigroup - groupoid of finite order; order of G is just (20, 20) so $20^2$ elements are in G.

***Example 1.2.59:*** Let $G = S \cup G_1 = \{$interval symmetric semigroup S(X) where $X = \{([0, a_1], [0, a_2], [0, a_3]) \} \cup \{[0, a] \, / \, a \in Z_{49}, *, (0, 9)\}$ be a biinterval semigroup groupoid of finite order. Clearly G is non commutative.

Now we proceed onto define substructures in them.

We call a proper bisubset $P = P_1 \cup P_2 \subseteq S \cup G_1 = G$ a biinterval semigroup-groupoid to be a biinterval subsemigroup-subgroupoid if $P_1$ is an interval subsemigroup of S and $P_2$ is an interval subgroupoid of $G_1$.

We will illustrate this situation by some examples.

*Example 1.2.60:* Let $G = S \cup G_1 = \{[0, a] / a \in Z_{40}$ under multiplication modulo 40$\} \cup \{[0, b] / b \in Z_{14}, *, (8, 7)\}$ be a biinterval semigroup groupoid. Consider $P = P_1 \cup P_2 = \{[0, a] / a \in \{0, 2, 4, 6, …, 38\} \subseteq Z_{40}\} \cup \{[0, a] / a \in \{0, 4\} \subseteq Z_{14}\} \subseteq G = S \cup G_1$ is a biinterval subsemigroup - subgroupoid of G.

*Example 1.2.61:* Let $G = G_1 \cup S = \{[0, a] / a \in Z_{12}, *, (1, 3)\} \cup \{[0, b] / b \in Z_{16}\}$ be a biinterval groupoid - semigroup. $T = T_1 \cup T_2 = \{[0, a] / a \in \{0, 3, 6, 9\} \subseteq Z_{12}\} \cup \{[0, b] / b \in \{0, 2, 4, 6, …, 14\} \subseteq Z_{16}\} \subseteq G_1 \cup S$ is a biinterval subgroupoid - subsemigroup of G.

In a similar way we can define ideals, S-subgroupoid S-subsemigroups and so on.

We will call a biset $P = P_1 \cup P_2 \subseteq G_1 \cup S = G$ to be a Smarandache sub biinterval groupoid - semigroup if $P_1$ is an interval S-subgroupoid and $P_2$ is an interval S-semigroup.

We will illustrate this situation by some examples.

*Example 1.2.62:* Let $G = G_1 \cup S_2$ be a biinterval groupoid-semigroup where $G_1 = \{[0, a] / a \in Z_{12}, *, (1, 3)\}$ the interval groupoid and $S_2 = \{[0, a] / a \in Z_{20}\}$ an interval semigroup under multiplication modulo 20. Choose $P = P_1 \cup P_2 = \{[0, a] / a \in \{0, 3, 6, 9\} \subseteq Z_{12}, *, (1, 3)\} \cup \{[0, a] / a \in \{0, 5, 10, 15\} \subseteq Z_{20}\} \subseteq G_1 \cup S_2 = G$. P is a Smarandache subinterval groupoid - semigroup of S.

It is pertinent to mention here that every biinterval subgroupoid-subsemigroup need not in general be a Smarandache biinterval subgroupoid - subsemigroup.

Infact every S-sub biinterval groupoid-semigroup (Smarandache biinterval subgroupoid - subsemigroup) is also a biinterval subgroupoid - subsemigroup and not conversely. We also have the following theorem the proof of which is left to the reader.

**THEOREM 1.2.30:** *Let $S_1 \cup G_1 = G$ be a biinterval semigroup - groupoid. If G has a S-subbiinterval semigroup-groupoid then G is a S-biinterval semigroup - groupoid.*

S-interval biideals can be defined in an analogous way. We now define the notion of quasi biinterval semigroup - groupoid (groupoid - semigroup). Let $G = S_1 \cup G_1$, if $S_1$ is a semigroup and not an interval semigroup and $G_1$ and interval groupoid or $S_1$ an interval semigroup and $G_1$ an ordinary groupoid not an interval groupoid then we define G to be a quasi biinterval semigroup - groupoid (quasi biinterval groupoid - semigroup).

We will illustrate this situation by some examples.

***Example 1.2.63:*** Let $G = S \cup G_1 = \{[0, a] / a \in Z_{45}$ under multiplication modulo 45$\} \cup \{Z_{12} (3, 4)\}$ be a quasi biinterval semigroup - groupoid.

***Example 1.2.64:*** Let $G = S \cup G_1 = \{S (7)$; the symmetric semigroup on $(1, 2, …, 7)\} \cup \{[0, a] / a \in Z_{43}, *, (10, 11)\}$ be a quasi biinterval semigroup - groupoid of finite order which is non commutative.

***Example 1.2.65:*** Let $G = S \cup G_1 = \{$All $3 \times 3$ matrices with entries from $Z_{40}\} \cup \{[0, a] / a \in Z_{19}, *, (2, 3)\}$ be a quasi biinterval semigroup - groupoid.

***Example 1.2.66:*** Let $G = S \cup G_1 = \{[0, a] / a \in Z^+ \cup \{0\}$ under multiplication$\} \cup \{Z_{45} (3, 8)\}$ be a quasi biinterval semigroup-groupoid.

We can define substructures and S-structures and S-substructure in an analogous way. This is direct and the interested reader can do the job.

*Example 1.2.67:* Let $G = G_1 \cup S = \{Z_{14} (7, 8)\} \cup \{[0, a] / a \in Z^+ \cup \{0\}\}$ be a quasi biinterval groupoid - semigroup. $P = P_1 \cup P_2 = \{0, 2\} \subseteq Z_{14}, *, (7, 8)\} \cup \{[0, a] / a \in 3Z^+ \cup \{0\}\} \subseteq G_1 \cup S$ is a quasi biinterval subgroupoid-subsemigroup (quasi sub biinterval groupoid - semigroup).

*Example 1.2.68:* Let $G = G_1 \cup S = \{Z_7 \setminus \{0\}, \times\} \cup \{[0, a] / a \in Z_7, *, (3, 2)\}$ be a biinterval quasi semigroup - groupoid. G is a simple biinterval quasi semigroup - groupoid as G has no subbiinterval quasi semigroup - groupoid.

Now we can also define bizerodivisors, biunits, biidempotents and their Smarandache analogue for these structures.

This task is also left for the reader.

We can define identities as semigroups can satisfy these identities as it holds an associative operation on it.

1.3 Interval Bigroups and their Generalization

In this section we introduce interval bigroups and generalize them to biinterval group-semigroup and biinterval group - groupoid. We will illustrate this by examples.

**DEFINITION 1.3.1:** *Let $G = G_1 \cup G_2$ where $(G_1, o)$ and $(G_2, *)$ be two distinct interval groups; G with the inherited operations of $G_1$ and $G_2$ is a bigroup called the interval bigroup. $G = G_1 \cup G_2 = \{[0, a] \cup [0, b] / [0,a] \in G_1$ and $[0,b] \in G_2\}$ with operation '·' defined by, for $g = [0, a] \cup [0, b]$ and $h = [0, c] \cup [0, d]$ we have $g.h = ([0, a] \cup [0, b]) . ([0, c] \cup [0, d]) = [0, a]$ o $[0, c] \cup [0, b] * [0, d] = [0, a \text{ o } c] \cup [0, b * d]$ is in G.*

It is a matter of routine to prove (G, .) is a bigroup. This task is left as an exercise for the reader.

***Example 1.3.1:*** Let $G = G_1 \cup G_2 = \{[0, a] / a \in Z_{10}$ under addition modulo 10$\} \cup \{[0, b] / b \in Z_{11} \setminus \{0\}$ under multiplication modulo 11$\}$ is an interval bigroup. $[0, 0] \cup [0, 1]$ is the bidentity of G. We will describe the operations on the elements of G. Let $x = [0, 6] \cup [0, 3] \in G$ the inverse of x is $[0, 4] \cup [0, 4] = x^{-1}$. For we see

$$\begin{aligned}
xx^{-1} &= ([0, 6] \cup [0, 3]) \cdot ([0, 4] \cup [0, 4]) \\
&= [0, 10 \ (\text{mod } 10)] \cup [0, 12 \ (\text{mod } 11)] \\
&= [0, 0] \cup [0, 1] = \text{the biidentity in G.}
\end{aligned}$$

Now consider $x = [0, 9] \cup [0, 2]$ and $y = [0, 4] \cup [0, 7]$ in G.

$$\begin{aligned}
x \cdot y &= ([0, 9] \cup [0, 2]) \cdot ([0, 4] \cup [0, 7]) \\
&= ([0, 9] + [0, 4]) \cup [0, 2] \times [0, 7]) \\
&= [0, 13 \ (\text{mod } 10)] \cup [0, 14 \ (\text{mod } 11)] \\
&= [0, 3] \cup [0, 3] \in G.
\end{aligned}$$

We will give more examples of these interval bigroups.

***Example 1.3.2:*** Let $G = G_1 \cup G_2 = \{[0, a] / a \in Z_{23} \setminus \{0\}$, multiplication modulo 23$\} \cup \{[0, b] / b \in Z_{19} \setminus \{0\}$, multiplication modulo 19$\}$ be a interval bigroup. Clearly G is commutative and is of finite order.

***Example 1.3.3:*** Let $G = G_1 \cup G_2 = \left\{\sum_{i=0}^{\infty}[0,a]x^i \bigg| a \in Z_{45}\right\} \cup \{[0, a] / a \in Z_{29} \setminus \{0\}\}$ be an interval bigroup of infinite order.

***Example 1.3.4:*** Let $G = G_1 \cup G_2$

$$= \left\{\sum_{i=0}^{7}[0,a]x^i \bigg| a \in Z_{25}\right\} \cup \left\{\sum_{i=0}^{9}[0,a]x^i \bigg| a \in Z_{47}\right\}$$

be a finite interval bigroup, both $G_1$ and $G_2$ are groups under polynomial modulo addition.

***Example 1.3.5:*** Let $G = G_1 \cup G_2 = \{G_X$ where $X = \{([0, a_1], [0, a_2] \ldots [0, a_7])$ under composition of mapping$\} \cup \{S_Y$ where $Y =$ $\{([0, b_1], [0, b_2], [0, b_3])\}$ under composition of mapping} be an interval bigroup. Clearly G is a non commutative interval bigroup of finite order.

Now we proceed onto give examples of substructures in them.

***Example 1.3.6:*** Let $G = G_1 \cup G_2 = \{[0, a] / a \in Z_{28}$, under addition modulo 28$\} \cup \{[0, b] / b \in Z_{19} \setminus \{0\}$, under multiplication modulo 19$\}$ be an interval bigroup. Choose $P = P_1 \cup P_2 = \{[0, a] / a \in \{0, 7, 14, 21\} \subseteq Z_{28}\} \cup \{[0, b] / b \in \{1, 18\} \subseteq Z_{19} \setminus \{0\} \subseteq G_1 \cup G_2$ is an interval subbigroup of G.

***Example 1.3.7:*** Let $G = G_1 \cup G_2 = \{[0, a] / a \in Z_{200}\} \cup \{[0, a] / a \in Z_{144}\}$ be an interval bigroup. $P = P_1 \cup P_2 = \{[0, a] / a \in \{0, 10, 20, \ldots, 190\} \subseteq Z_{200}\} \cup \{[0, a] / a \in \{0, 2, \ldots, 142\} \subseteq G_1 \cup G_2$ be an interval subbigroup of G.

All these interval bisubgroups are also normal as G is a commutative interval bigroup.

***Example 1.3.8:*** Let G be an interval bigroup where $G_1$ is an interval symmetric group on $X = \{([0, a_1], \ldots, [0, a_{11}])\}$ and $G_2$ is also an interval symmetric group on $Y = \{([0, b_1], \ldots, [0, b_7])\}$ i.e., $G = G_1 \cup G_2$.

Let $A = A_x \cup A_y$, the alternative interval bisubgroup of G.

Clearly A is an interval binormal subgroup of G. Infact all interval bisubgroups of G are not binormal in G.

***Example 1.3.9:*** Let $G = G_1 \cup G_2 = \{[0, a] / a \in Z_{10}\} \cup \{[0, b] / b \in Z_{42}\}$ be an interval bigroup. All interval bisubgroups of G are normal as G is a commutative interval bigroup.

We define the biorder of G to be $|G_1| \cdot |G_2|$. We see all the properties enjoyed by usual groups are true in case of interval bigroups provided this is taken as the order, for $G = G_1 \cup G_2 = \{[0, a] \cup [0, b] / a \in G_1$ and $b \in G_2\}$; if both $G_1$ and $G_2$ are assumed to be of finite order all classical theorems for finite groups are true in case of interval bigroups. This is only a matter of routine and can be easily derived as a regular simple exercise.

We will illustrate this by some simple examples.

***Example 1.3.10:*** Let $G = G_1 \cup G_2 = \{[0, a] / a \in Z_{12}\} \cup \{[0, a] / a \in Z_5 \setminus \{0\}\}$ be an interval bigroup. Clearly G is of finite biorder. $|G| = |G_1| \cdot |G_2| = 12.4 = 48$.

Now consider $H = H_1 \cup H_2 = \{[0, a] / a \in \{0, 2, 4, 6, 8, 10\} \subseteq Z_{12}\} \cup \{[0, a] / a \in \{1, 3\} \subseteq G_2 = Z_5 \setminus \{0\}\} \subseteq G_1 \cup G_2$, the interval subbigroup of G. Order of H is $|H_1| \cdot |H_2| = 6.2 = 12$. We see 12/48. Now we define for any $x = [0, a] \cup [0, b] \in G$ we define $x^n = [0, 1] \cup [0, 1]$ if $[0, a]^t = 1$ and $[0, b]^s = 1$ then n = st. We will illustrate this by some examples.

Let $x = [0, 6] \cup [0, 4] \in G$ then $x^2 = [0, 0] \cup [0, 1]$ biidentity element of G. If $x = [0, 9] \cup [0, 3] \in G$ then $x^{16} = [0, 9]^4 \cup [0, 3]^4 = [0, 0] \cup [0, 1]$.

Likewise Cauchy theorem will be true in case of finite interval bigroups. Sylow theorems can be easily proved for finite interval bigroups. Cayley theorem can be extended by using two suitable symmetric interval bigroup in which all interval bigroups can be embedded.

Now we proceed onto define quasi interval bigroups.

**DEFINITION 1.3.2:** *Let $G = G_1 \cup G_2$, where only one of $G_1$ or $G_2$ is an interval group and the other is just a group then we call G to be a quasi interval bigroup. The operations are defined on G as in case bigroups and interval bigroups.*

We will give examples of these structures.

***Example 1.3.11:*** Let $G = G_1 \cup G_2 = \langle g / g^{12} = 1 \rangle \cup \{[0, a] / a \in Z_{32}\}$ be a quasi interval bigroup of finite order. $|G| = 12 \times 32$.

***Example 1.3.12:*** Let $G = G_1 \cup G_2 = \{[0, a] / a \in Z_{23} \setminus \{0\}\} \cup$ {all $5 \times 5$ matrices with entries from R such that their determinant is non zero} be a quasi interval bigroup. Clearly G is of infinite order and G is non commutative.

***Example 1.3.13:*** Let $G = G_1 \cup G_2 = \{S_3\} \cup \{[0, a] / a \in Z_{45}\}$ be a quasi interval bigroup of order $6 \times 45$.

***Example 1.3.14:*** Let $G = G_1 \cup G_2 = S_5 \cup \{([0, a_1] \ldots [0, a_6])$ the symmetric interval group on 6 intervals} be a quasi interval bigroup of order $5! \times 6!$.

Clearly G is a non commutative quasi interval bigroup.

We can define substructures in them in an analogous way which is direct.

***Example 1.3.15:*** Let $G = G_1 \cup G_2 = S_3 \cup \{[0, a] / a \in Z_{40}\}$ be a quasi interval bigroup.
Consider $P = P_1 \cup P_2 =$
$$\left\{\begin{pmatrix} 1 & 2 & 3 \\ 1 & 2 & 3 \end{pmatrix}, \begin{pmatrix} 1 & 2 & 3 \\ 2 & 3 & 1 \end{pmatrix}, \begin{pmatrix} 1 & 2 & 3 \\ 3 & 1 & 2 \end{pmatrix}\right\}$$
$\cup \{[0, a] / a \in \{0, 2, 4, \ldots, 38\} \subseteq Z_{40}\} \subseteq G_1 \cup G_2$. P is a quasi interval bisubgroup of G and $|P| = 3.20 = 60$. We see P is also the quasi interval normal subbigroup of G. We see G has quasi interval subbigroups which are not normal in G.
Choose
$$S = S_1 \cup S_2 = \left\{\begin{pmatrix} 1 & 2 & 3 \\ 1 & 2 & 3 \end{pmatrix}, \begin{pmatrix} 1 & 2 & 3 \\ 1 & 3 & 2 \end{pmatrix} \text{ in } S_3\right\}$$
$\cup \{[0, a] / a \in \{0, 10, 20, 30\} \subseteq Z_{40}\} \subseteq G_1 \cup G_2$; S is only a quasi interval subbigroup of G and is not normal in G.

If $G = G_1 \cup G_2$ is a quasi interval bigroup which is non commutative but all its quasi interval bisubgroups are commutative then we define G to be a quasi commutative quasi interval bigroup.

The quasi interval bigroup G given in example 1.3.15 is only a quasi commutative quasi interval bigroup; but G is clearly non commutative.

Now having seen the substructures as in case of interval bigroup G, if G is a quasi interval bigroup of finite order say if $G = G_1 \cup G_2$ then $|G| = |G_1| \cdot |G_2|$ and all classical theorems true for usual groups hold good for quasi interval bigroups also. Now we just give some examples.

*Example 1.3.16:* Let $G = G_1 \cup G_2 = \{S_4\} \cup \{[0, a] / a \in Z_{150}\}$ be a quasi interval bigroup of finite order. $|G| = \underline{4} \times 150 = 3600$. It is easily verified every quasi interval bisubgroup of G divides the order of G hence the Lagrange theorem is true.

Further G is not a quasi commutative quasi interval bigroup. For take in G; $P = P_1 \cup P_2 = A_4 \cup \{[0, a] / a \in \{0, 10, 20, 30, \ldots, 140\} \subseteq G_1 \cup G_2$; we see P is quasi interval bisubgroup but clearly P is non commutative, hence G is not a quasi commutative quasi interval bigroup. G has also commutative interval subbigroups, for take $T = T_1 \cup T_2 = \left\{\begin{pmatrix} 1 & 2 & 3 & 4 \\ 1 & 2 & 3 & 4 \end{pmatrix}, \begin{pmatrix} 1 & 2 & 3 & 4 \\ 2 & 3 & 4 & 1 \end{pmatrix}, \begin{pmatrix} 1 & 2 & 3 & 4 \\ 3 & 4 & 1 & 2 \end{pmatrix}, \begin{pmatrix} 1 & 2 & 3 & 4 \\ 4 & 1 & 2 & 3 \end{pmatrix}\right\} \cup \{[0, a] / a \in \{0, 30, 60, 90, 120\} \subseteq Z_{150}\} \subseteq G_1 \cup G_2$; T is commutative quasi interval subbigroups. Consider an element
$$x = \begin{pmatrix} 1 & 2 & 3 & 4 \\ 2 & 3 & 4 & 1 \end{pmatrix} \cup \{[0, 20] / 20 \in Z_{150}\}$$
in G. We see $\begin{pmatrix} 1 & 2 & 3 & 4 \\ 1 & 2 & 3 & 4 \end{pmatrix} \cup \{[0, 0]\}$ is the biidentity element of G. We see $x^{-1} = \begin{pmatrix} 1 & 2 & 3 & 4 \\ 4 & 1 & 2 & 3 \end{pmatrix} \cup \{[0, 130]\} \in G$ is the biinverse of x. Consider $y = \begin{pmatrix} 1 & 2 & 3 & 4 \\ 2 & 3 & 4 & 1 \end{pmatrix} \cup \{[0, 50]\} \in G$.

We see $y^{12} = \begin{pmatrix} 1 & 2 & 3 & 4 \\ 1 & 2 & 3 & 4 \end{pmatrix} \cup \{[0, 0]\}$. It is easily verified Cauchy theorem is true. Also one can easily check Sylow theorems are true for G.

Now consider the following examples.

***Example 1.3.17:*** Let $G = G_1 \cup G_2 = <g / g^{17} = 1> \cup \{[0, a] / a \in Z_{19}\}$ be a quasi interval bigroup. We see G is a strongly simple quasi interval bigroup for both $G_1$ and $G_2$ does not contain subgroups. If only one of $G_1$ or $G_2$ has subgroups we call G to be a simple quasi interval bigroup.

***Example 1.3.18:*** Let $G = G_1 \cup G_2 = <g / g^{23} = 1 \cup \{[0, a] / a \in Z_{49}\}$ be a quasi interval bigroup. Clearly G is only a simple quasi interval bigroup and is not a strongly simple quasi interval bigroup for $G_2$ has subgroup with respect to addition modulo 49.

***Example 1.3.19:*** Let $G = G_1 \cup G_2 = \{D_{29}\} \cup \{[0, a] / a \in Z_{43}\}$ be a quasi interval bigroup. G is only a simple quasi interval bigroup for $G_2$ has no proper interval subgroups. Thus G is not a doubly simple quasi interval bigroup.

We have the following theorem, the proof of which is direct and is left as an exercise to the reader.

**THEOREM 1.3.1:** *A simple quasi interval bigroup is not strongly simple quasi interval bigroup.*

Now we can define biinterval group - semigroup (semigroup - group).

**DEFINITION 1.3.3:** *Let $G = G_1 \cup G_2$ where $(G_1, *)$ is an interval group and $(S_2, o)$ is an interval semigroup. G with operation '.' such that for every $x = [0, a] \cup [0, b]$ and $y = [0, c] \cup [0, d]$ in G we define*

$$x.y = ([0, a] \cup [0, b]) \cdot ([0, c] \cup [0, d])$$
$$= [0, a] * [0, c] \cup [0, b] o [0, d]$$
$$= [0, a * c] \cup [0, b \, o \, d] \in G.$$

*($G_1$, .) is defined as the biinterval group - semigroup. We call G a biinterval group - semigroup as elements in G are biintervals.*

We illustrate this situation by some examples.

***Example 1.3.20:*** Let $G = G_1 \cup G_2 = \{[0, a] |\, a \in Z_{40}, +\} \cup \{[0, b] \,|\, b \in Z_{20}$ under multiplication modulo 20$\}$ be the biinterval group-semigroup of finite order. Clearly $|G| = |G_1|\, |G_2| = 40 \times 20 = 800$. Suppose $x = [0, 9] \cup [0, 10]$ and $y = [0, 1] \cup [0, 7]$ are in G.

$$x.y = ([0, 9] \cup [0, 10]) \cdot ([0, 1] \cup [0, 7])$$
$$= ([0, 9] + [0, 1]) \cup ([0, 10] \cdot [0, 7])$$
$$= [0, 10] \cup [0, 10] \in G.$$

$[0, 0] \cup [0, 1]$ is the biidentity of G.

In general every biinterval group - semigroup need not contain the biidentity. This is evident from the following example.

***Example 1.3.21:*** Let $G = G_1 \cup G_2 = \{[0, a] \,/\, a \in Z_{45}\} \cup \{[0, a] \,/\, a \in 2Z^+ \cup \{0\}$ under multiplication$\}$ be a biinterval group - semigroup. Clearly order of G is infinite and G has no biidentity for the interval semigroup is not a monoid.

***Example 1.3.22:*** Let $G = G_1 \cup G_2 = \{[0, a] \,/\, a \in Z_7 \setminus \{0\}\} \cup \{[0, b] \,/\, b \in Z_{19}\}$ be an interval group - semigroup. We see $\{[0, 1] \cup [0, 1]\}$ is the biidentity of G.

When a biinterval group - semigroup (semigroup - group) has no biinterval subgroup- subsemigroup then we call G to be a simple biinterval group - semigroup.

We will first illustrate this by some examples.

***Example 1.3.23:*** Let $G = G_1 \cup S_1 = \{[0, a] / a \in Z_{20}\} \cup \{[0, b] / b \in Z_{40}\}$ be a biinterval group - semigroup. Consider $H = H_1 \cup H_2 = \{[0, a] / a \in \{0, 2, 4, 6, 8, \ldots, 18\} \cup \{[0, b] / b \in \{0, 10, 20, 30\}$ (multiplication modulo 40)$\} \subseteq G = G_1 \cup S_1$ is a biinterval subgroup - subsemigroup of S.

***Example 1.3.24:*** Let $G = G_1 \cup S_1 = \{[0, a] / a \in Z_{11}\} \cup \{[0, b] / b \in Z_7\}$ be a biinterval group-semigroup. We see $G_1$ has no interval subgroup where as $S = \{[0, a] / 0, 1, 6\} \subseteq S$ is an interval subsemigroup.

We under these conditions define the following.

**DEFINITION 1.3.4:** *Let $G = G_1 \cup S_1$ be a biinterval group - semigroup. If only one of $G_1$ or $S_1$ has interval subgroup or interval subsemigroup, then we define G to be a quasi subbiinterval group - semigroup.*

***Example 1.3.25:*** Let $G = G_1 \cup G_2 = \{S(\langle X \rangle);$ interval symmetric semigroup$\} \cup \{[0, a] / a \in Z_{42}\}$ be the biinterval semigroup-group of finite order which is non commutative.

***Example 1.3.26:*** Let $G = S_1 \cup G_1 = \{[0, a] / a \in Z_{17}\} \cup \{[0, a] / a \in Z_{41}\}$ be a biinterval semigroup-group. G is of finite order and order of G is 17.41.

Now we define quasi biinterval semigroup - group (group - semigroup) if only one of them is group or the interval semigroup or an interval group or interval semigroup.

We will illustrate this situation by some examples.

***Example 1.3.27:*** Let $G = G_1 \cup S_1 = \{Z_{15}, +\} \cup \{[0, a] / a \in Z_7, \times\}$ be a quasi biinterval group-semigroup of finite order. $|G| = 15.7 = 105$.

***Example 1.3.28:*** Let $G = G_1 \cup S_1 = \{S_9\} \cup \{[0, a] / a \in Z_{12}, *\}$ be a quasi biinterval group-semigroup.

*Example 1.3.29:* Let $G = S_1 \cup G_1 = \{S(\langle X \rangle) / \langle X \rangle = \langle ([0, x_1], \ldots, [0, x_7]) \rangle\} \cup \{Z_{15}, +\}$ be a quasi biinterval semigroup - group.

We will now give some substructures of them.

*Example 1.3.30:* Let $G = S_1 \cup G_1 = \{[0, a] / a \in Z_{45}, \times\} \cup \{Z_{40}, +\}$ be a quasi biinterval semigroup - group. Consider $H = H_1 \cup H_2 = \{[0, a] / a \in \{0, 5, 10, 15, 20, 25, 30, 35, 40\} \subseteq Z_{45}, \times\} \cup \{0, 10, 20, 30\} \subseteq Z_{40}, +\} \subseteq S_1 \cup G_1 = G$ is a quasi biinterval subsemigroup-subgroup (quasi subbiinterval semigroup- group).

*Example 1.3.31:* Let $G = S_1 \cup G_1 = \{Z_{120}, \times\} \cup \{[0, a] / a \in Z_{19} \setminus \{0\}, \times\}$ be a quasi biinterval semigroup-group. Choose $H = H_1 \cup H_2 = \{\{0, 10, 20, \ldots, 110\} \subseteq Z_{120}, \times\} \cup \{[0, a] / a \in \{1, 18\} \subseteq Z_{19}, \times\} \subseteq S_1 \cup G_1$ is a quasi subbiinterval semigroup-group of G.

*Example 1.3.32:* Let $G = G_1 \cup S_1 = \{Z_7, +\} \cup \{[0, a] / a \in Z_3 \setminus \{0\}\}$ be a quasi biinterval group - semigroup. G has no quasi subbiinterval group - semigroup.

In view of this we have the following results the proof of which is direct.

**THEOREM 1.3.2:** *Let $G = G_1 \cup S_1 = \{Z_p, +\} \cup \{any\ interval\ semigroup\}$ be a quasi biinterval group-semigroup. G is a quasi simple biinterval group-semigroup.*

**THEOREM 1.3.3:** *Let $G = G_1 \cup S_1 = \{[0, a] / a \in Z_p, +, p\ a\ prime\} \cup \{any\ semigroup\}$ be a quasi biinterval group-semigroup. G is a quasi simple biinterval group - semigroup.*

Now we cannot have all classical theorems to be true in case of quasi biinterval group - semigroup.

This will be illustrated by some examples.

***Example 1.3.33:*** Let $G = S_1 \cup G_1 = \{Z_{16}, \times\} \cup \{[0, a] / a \in Z_7 \setminus \{0\}, \times\}$ be a quasi biinterval semigroup - group of order $16 \times 6 = 96$.

Consider $H = H_1 \cup H_2 = \{\{0, 4, 8, 12, 1\} \subseteq Z_{16}, \times\} \cup \{[0, a] / a \in \{1, 6\} \subseteq Z_7 \setminus \{0\}\} \subseteq S_1 \cup G_1 = G$; H is a quasi subbiinteval semigroup - group of order $5 \times 2 = 10$ we see $10 \nmid 96$.

Thus the classical Lagrange's theorem for finite groups is not true in general for biinterval semigroup - group or quasi biinterval semigroup - group. This is evident from the above example.

***Example 1.3.34:*** Let $G = G_1 \cup S_1 = \{Z_{11} \setminus \{0\}, \times\} \cup \{[0, a] / a \in Z_9, \times\}$ be a quasi biinterval group - semigroup of order $10 \times 9 = 90$. $H = \{[0, a] / a \in \{1, 10\} \subseteq Z_{11} \setminus \{0\}, \times\} \cup \{[0, a] / a \in \{0, 3\} \subseteq Z_9, \times\} = H_1 \cup H_2 \subseteq G_1 \cup S_1 = G$; H is a quasi - subbiinterval group-semigroup of order 4. Clearly $4 \times 90$. Further Cauchy theorem is also not true for $x = \{10\} \cup \{[0, 8]\} \in G$ is such that $x^4 = \{1\} \cup \{[0, 1]\}$ identity element of G and $4 \nmid 90$. Thus in general Cauchy theorem for finite groups is not true in case of quasi interval group - semigroup.

***Example 1.3.35:*** Let $G = G_1 \cup S_1 = \{S_3\} \cup \{[0, a] / a \in Z_{11}, \times\}$ be a quasi - biinterval group semigroup. Consider $H = H_1 \cup H_2 = \{A_3\} \cup \{[0, 0], [0, 1], [0, 10] \in S_1\} \subseteq G_1 \cup S_1$ is a quasi sub biinterval group - semigroup of G. Now $o(G) = 6 \times 11 = 66$ and $o(H) = 3 \times 3 = 9$ and $9 \nmid 66$.

Thus all classical theorems except Cayleys theorem (when modified for semigroups) does not hold good for quasi biinterval group - semigroup.

Having defined quasi biinterval group - semigroups we can define also quasi biinterval group - groupoid (groupoid - group). We wish to state here that the term biinterval is used in an appropriate way but only to signify the structure under consideration is a bistructure. We also define the notion of biinterval group - groupoid.

Let $G = G_1 \cup G_2$ where $\{G_1, \text{'o'}\}$ is an interval group and $\{G_2, *\}$ is an interval groupoid we define '.' on G which operations is described in the following.

$G = G_1 \cup G_2 = \{[0, a] \cup [0, b] \,/\, [0,a] \in G_1 \text{ and } [0,b] \in G_2\}$.
Let $x = [0, a] \cup [0, b]$ and $y = [0, c] \cup [0, d]$ be in G. Now

$$\begin{aligned} x.y &= ([0, a] \cup [0, b]) \cdot ([0, c] \cup [0, d]) \\ &= [0, a] \text{ o } [0, c] \cup [0, b] * [0, d] \\ &= [0, a \text{ o } c] \cup [0, b * d] \in G \end{aligned}$$

as $[0,a], [0,c] \in G_1$ and $[0,c], [0,d] \in G_2$. We define $(G_1, .)$ to be a biinterval group groupoid.

If both $G_1$ and $G_2$ are of finite order we say G is of finite order and $|G| = |G_1| \times |G_2|$, even if one of $G_1$ or $G_2$ is of infinite order then we define G to be of infinite order. If both $G_1$ and $G_2$ are commutative we say G is commutative otherwise non commutative.

We will illustrate this situation by some examples.

***Example 1.3.36:*** Let $G = G_1 \cup G_2 = \{[0, a] \,/\, a \in Z_{42}, +\} \cup \{[0, b] \,/\, b \in Z_9, * (2, 4)\}$ be a biinterval group - groupoid. Clearly G is of finite order. For $|G| = |G_1| \, |G_2| = 42.9 = 378$. G is non commutative as $G_2$ is non commutative.
Let $x = [0, 24] \cup [0, 4]$ and $y = [0, 7] \cup [0, 3] \in G = G_1 \cup G_2$.

$$\begin{aligned} x.y &= ([0, 24] \cup [0, 4]) \cdot ([0, 7] \cup [0, 3]) \\ &= ([0, 24] + [0, 7]) \cup ([0, 4] * [0, 3]) \\ &= [0, 31] \cup [0, 8+12 \text{ (mod 9)}] \\ &= [0, 31] \cup [0, 2] \in G. \end{aligned}$$

We in general cannot talk about biidentity or biinverse for $G_2$ the groupoid may or may not have identity hence inverse may or may not exist.

*Example 1.3.37*: Let $G = G_1 \cup G_2 = \{[0, a] / a \in Z_{23}, *, (3, 2)\}$ $\cup \{[0, a] / a \in Z_{25}, +\}$ be an biinterval groupoid - group. Clearly G is of finite order and non commutative. G has no biidentity.

Further in general all the classical theorems for finite groups may not be true in case of finite biinterval group-groupoids (groupoids - group).

*Example 1.3.38:* Let $G = G_1 \cup G_2 = \{[0, a] / a \in Z_{19} \setminus \{0\}\} \cup \{[0, a] / a \in Z_{40}, *, (7, 9)\}$ be a biinterval group - groupoid of order $18 \times 40$. Take $H = H_1 \cup H_2 = \{[0, a] / a \in \{1, 18\} \subseteq Z_{19} \setminus \{0\}\} \cup \{[0, a] / a \in \{0, 10, 20, 30\} \subseteq Z_{40}, *, (7, 9)\} \subseteq G_1 \cup G_2 = G$ is a subbiinterval group-groupoid of order $2 \times 4 = 8$. We see $o(H) / o(G)$.

Thus we may have some sub biinterval group - groupoids whose biorder divides the biorder of the biinterval group - groupoid.

*Example 1.3.39*: Let $G = G_1 \cup G_2 = \{[0, a] / a \in Z_{25}, +\} \cup \{[0, a] / a \in Z_{120}, *, (3, 7)\}$ be a biinterval group – groupoid. Clearly G is non commutative. Biorder of G is $25 \times 120$. Consider $H = H_1 \cup H_2 = \{[0, a] / a \in \{0, 5, 10, 15, 20\} \subseteq Z_{25}, +\} \cup \{[0, a] / a \in \{0, 10, 20, 30, 40, 50, …, 110\} \subseteq Z_{120}, *, (3, 7)\} \subseteq G_1 \cup G_2$, H is a subbiinterval group - groupoid of G. $o(H) = 5 \times 12 = 60$. We see $o(H) / o(G)$.

Now we have also examples of finite quasi biinterval group-groupoids G such that it contains subbiinterval group - groupoid H such that $o(H) / o(G)$.

We will give some examples of them.

*Example 1.3.40*: Let $G = G_1 \cup G_2 = \{Z_{20}, +\} \cup \{[0, a] / a \in Z_9, *, (5, 3)\}$ be a biinterval group - groupoid of finite order $o(G) = 20 \times 9$. Consider $= \{0, 5, 10, 15\} \cup \{1, 2, 4, 5, 7, 8\} = H_1 \cup H_2 \subseteq G_1 \cup G_2$ be a subbiinterval group-groupoid. Clearly $o(H) = 4.6 = 24$. $o(G) = 20 \times 9$. $24 \nmid 180$. Hence the order does not divide.

We can have such examples. This distinguishes usual interval bigroups from biinterval group-groupoids.

Next we proceed onto describe formally quasi biinterval group-groupoids.

Let $G = G_1 \cup G_2$ only one of $G_1$ or $G_2$ is an interval group or interval groupoid we define '.' on G so that (G, .) is defined as the quasi biinterval group-groupoid or quasi interval groupoid-group.

Let (G, .) be a group and $(G_2, *)$ be an interval groupoid then $G = G_1 \cup G_2 = \{g \cup [0, a] \,/\, g \in G_1, [0, a] \in G_2\}$, define '.' on G as follows. $x = g \cup [0, a]$ and $y = h \cup [0, b]$ in G.

$$x.y = (g \cup [0, a]) \cup (h \cup [0, b])$$
$$= g.h \cup [0, a] * [0, b]$$
$$= g.h \cup [0, a * b] \in G.$$

Thus (G, .) is the quasi interval group - groupoid.

*Example 1.3.41:* Let $G = G_1 \cup G_2 = \{g \,/\, g^{12} = 1\} \cup \{[0, a] \,/\, a \in Z_{42}, *, (9, 8)\}$ be a quasi interval group - groupoid of G order $12 \times 42$.

*Example 1.3.42:* Let $G = G_1 \cup G_2 = \{[0, a] \,/\, a \in Z_{57}, *, (17, 11)\} \cup \{S_{10}\}$ be a quasi interval groupoid-group of order $57 \times o(S_{10}) = 57 \times \underline{10}$.

*Example 1.3.43:* Let $G = G_1 \cup G_2 = \{[0, a] \,/\, a \in Z_{409}, +\} \cup \{Z_9, (2, 4)\}$ be a quasi interval group-groupoid. Clearly order of G is $409 \times 9$.

*Example 1.3.44:* Let $G = G_1 \cup G_2 = \{[0, a] \,/\, a \in Z_{43} \setminus \{0\}, \times\} \cup \{Z_{45}, (8, 11)\}$ be a quasi interval group-groupoid of order $42 \times 45$.

Now having seen examples of quasi interval group - groupoids we now proceed onto give examples of their substructures.

***Example 1.3.45:*** Let $G = G_1 \cup G_2 = S_3 \cup \{[0, a] / a \in Z_9, *, (5, 3)\}$ be a quasi interval group-groupoid. Clearly G is non commutative. $o(G) = 6 \times 9 = 54$. Let $H = H_1 \cup H_2 = A_3 \cup \{1, 2, 4, 5, 7, 8 \in Z_9, *, (5, 3)\} \subseteq G_1 \cup G_2$ is a quasi interval subgroup - groupoid or quasi subbiinterval group - groupoid. $o(H) = 3 \times 6 = 18$. Clearly $o(H) / o(G)$.

***Example 1.3.46:*** Let $G = G_1 \cup G_2 = \{g / g^{24} = 1\} \cup \{[0, a] / a \in Z_6, *, (2, 2)\}$ be a quasi interval group-groupoid. Clearly G is commutative. $o(G) = 24 \times 6$.

***Example 1.3.47:*** Let $G = G_1 \cup G_2 = \{g / g^{20} = 1\} \cup \{[0, a] / a \in Z_{20}, *, (4, 4)\}$ be a quasi interval group-groupoid. G is commutative and is of order $20 \times 20$.

Now we have the following result which is left as an exercise for the reader.

**THEOREM 1.3.4:** *Let $G = G_1 \cup G_2 = \{g / g^n = 1\} \cup \{[0, a] / a \in Z_n, *, (t, t)\}$ be a quasi interval group - groupoid G is commutative.*

**THEOREM 1.3.5:** *Let $G = G_1 \cup G_2 = \{S_n\} \cup \{[0, a] / a \in Z_n, *, (t, u); t \neq u, t, u \in Z_n \setminus \{0\}\}$ be a quasi interval group - groupoid. G is non commutative.*

***Example 1.3.48:*** Let $G = G_1 \cup G_2 = \{[0, a] / a \in Z_{19} \setminus \{0\}, \times\} \cup \{Z_{19}, (3, 3)\}$ be a quasi interval group-groupoid. G is non commutative of order $18 \times 19$.

***Example 1.3.49:*** Let $G = G_1 \cup G_2 = \{[0, a] / a \in Z_{11}, +\} \cup \{Z_{12}, (3, 7)\}$ be a quasi interval group - groupoid, G is non commutative and has no substructures.

In view of this we have the following theorem the proof of which is direct.

**THEOREM 1.3.6:** *Let $G = G_1 \cup G_2 = \{[0, a] / a \in Z_p, +\} \cup \{Z_{19}, (3, 2)\}$ is a quasi interval group - groupoid p, a prime has no substructures.*

**THEOREM 1.3.7:** *Let $G = G_1 \cup G_2 = \{[0, a] / a \in Z_n$, n a composite number, $+\} \cup \{Z_m, (t, u)\}$ be a quasi interval group - groupoid. G has substructures.*

*Example 1.3.50:* Let $G = G_1 \cup G_2 = \{[0, a] / a \in Z_{12}, +\} \cup \{Z_{10}, (5, 6)\}$ is a quasi interval group - groupoid, G has substructures.

It is pertinent to mention that almost all classical theorems for finite groups in general is not true for these quasi interval group - groupoids. Further certain concepts cannot be even extended to these structures. Thus we have limitations however these structures can find its applications in appropriate fields. Now having seen these semi associative and non associative bistructures with single binary operation we now proceed onto define interval biloops.

1.4 Interval Biloops and their Generalization

In this section we introduce the notion of interval biloops quasi intervals biloops, interval loop - group, interval group - loop, interval groupoid - loop interval loop-semigroup and so on and describe them.

We now define and describe these structures. These structures also do not in general satisfy the classical theorems for finite groups.

**DEFINITION 1.4.1:** *Let $L = L_1 \cup L_2$ where both $L_1$ and $L_2$ are two distinct interval loops (L, .) with '.' an operation inherited from both $L_1$ and $L_2$ is a loop called the biinterval loop or interval biloop.*

We will illustrate this situation by some examples.

*Example 1.4.1*: Let $L = L_1 \cup L_2 = \{[0, a] \mid a \in \{e, 1, 2, \ldots, 9\}, m = 8, *$, where $[0, a] * [0, b] = [0, 7b + 8a \pmod 9)\} \cup \{[0, b] / b \in \{e, 1, 2, 3, \ldots, 15\}, m = 8, *\}$ be a biinterval loop.

Suppose $x = [0, 6] \cup [0, 9]$ and $y = [0, 4] \cup [0, 10]$ in $L$.

$$\begin{aligned} x \cdot y &= ([0, 6] \cup [0, 9] \cdot ([0, 4] \cup [0, 10]) \\ &= [0, 6] * [0, 4] \cup [0, 9] * [0, 10] \\ &= [0, 7.6 \, ; 4.8 \pmod 9)] \cup [0, 7.9 + 8.10 \pmod{15})] \\ &= [0, 2] \cup [0, 8] \in L. \end{aligned}$$

$L$ is of finite order and is of even order given by $|L_1| \cdot |L_2| = 10.16 = 160$.

We can construct several such interval biloops.

*Example 1.4.2*: Let $L_1 \cup L_2 = \{[0, a] \mid a \in \{e, 1, 2, \ldots, 11\}, 6, *\} \cup \{[0, b] / b \in \{e, 1, 2, .., 5\}, 3, *\}$ be a biinterval loop of order $12 \times 6 = 72$.

It is easily verified $L$ is a commutative interval biloop.

*Example 1.4.3*: Let $L = L_1 \cup L_2$ where $L_1$ and $L_2$ are given by the following tables.

Table of $L_1$

| *      | [0, e] | [0, 1] | [0, 2] | [0, 3] | [0, 4] | [0, 5] |
|--------|--------|--------|--------|--------|--------|--------|
| [0, e] | [0, e] | [0, 1] | [0, 2] | [0, 3] | [0, 4] | [0, 5] |
| [0, 1] | [0, 1] | [0, e] | [0, 3] | [0, 5] | [0, 2] | [0, 4] |
| [0, 2] | [0, 2] | [0, 5] | [0, e] | [0, 4] | [0, 1] | [0, 3] |
| [0, 3] | [0, 3] | [0, 4] | [0, 1] | [0, e] | [0, 5] | [0, 2] |
| [0, 4] | [0, 4] | [0, 3] | [0, 5] | [0, 2] | [0, 3] | [0, 1] |
| [0, 5] | [0, 5] | [0, 2] | [0, 4] | [0, 1] | [0, 3] | [0, e] |

Clearly $L_1$ is of order 6 built using the loop $L_5(2)$.

Now the table for the interval loop $L_2$ is as follows

| * | [0, e] | [0, 1] | [0, 2] | [0, 3] | [0, 4] | [0, 5] |
|---|---|---|---|---|---|---|
| [0, e] | [0, e] | [0, 1] | [0, 2] | [0, 3] | [0, 4] | [0, 5] |
| [0, 1] | [0, 1] | [0, e] | [0, 4] | [0, 2] | [0, 5] | [0, 3] |
| [0, 2] | [0, 2] | [0, 4] | [0, e] | [0, 5] | [0, 3] | [0, 1] |
| [0, 3] | [0, 3] | [0, 2] | [0, 5] | [0, e] | [0, 1] | [0, 4] |
| [0, 4] | [0, 4] | [0, 5] | [0, 3] | [0, 1] | [0, e] | [0, 2] |
| [0, 5] | [0, 5] | [0, 3] | [0, 1] | [0, 4] | [0, 2] | [0, e] |

We see $L_1$ is non commutative interval loop of order 6 where as $L_2$ is a commutative interval loop of order 6. Thus $L = L_1 \cup L_2$ is a non commutative biinterval loop of order $6 \times 6 = 36$.

We have several properties associated with them which we will be discussing. For more about loops refer [5, 9, 11].

Let us show by examples the properties satisfied by these bi interval loops.

**Example 1.4.4**: Let $L = L_1 \cup L_2 = \{[0, a] / a \in \{e, 1, 2, ..., 7\}, *, '6'\} \cup \{[0, b] / b \in \{e, 1, 2, ..., 13\}, *, 9\}$ be a biinterval loop of order $8 \times 14 = 112$. We see L is a S-biinterval loop, for $P = P_1 \cup P_2 = \{[0, e], [0, 6], *\} \cup \{[0, e], [0, 11], *\} \subseteq L_1 \cup L_2$ is an interval bigroup of order $2 \times 2 = 4$. Thus L is a Smarandache interval biloop and is non commutative.

**Example 1.4.5**: Let $L = L_1 \cup L_2 = \{[0, a] / a \in \{e, 1, 2, ..., 19\}, *, 8\} \cup \{[0, b] | b \in \{e, 1, 2, ..., 23\} *, 10\}$ be an interval biloop. This has no proper interval bisubloop but has only interval bisubgroups.

These types of biinterval biloops which has no biinterval subloops but only biinterval subgroups will be defined as Smarandache biinterval subgroup loop or Smarandache interval bisubgroup biloops.

We have the following theorem.

**THEOREM 1.4.1**: *Let $L = L_1 \cup L_2 = \{[0, a] / a \in \{e, 1, 2, ..., p\}$, p a prime $m \neq e$ or 1 or p, *$\} \cup \{[0, b] / b \in \{e, 1, 2, ..., q\}$ where q is a prime $m' \neq e$ or 1 or q, *) be a biclass of interval loops. (This is a class for as m' can vary from $2 \leq m \leq p - 1$ and $2 \leq m' \leq q - 1$ respectively). This class of biinterval loop is a S-biinterval bisubgroup biloop.*

*Proof*: Given p and q are primes so L has $(p + 1) \cdot (q + 1)$ elements (L given in the theorem). By the very construction every bielement $[0, a] \cup [0, b]$ in L generates an interval bigroup of order two. Hence L has no biinterval subloops only has biinterval subgroups. Hence the claims.

We can also define for biinterval loops principal isotope, as in case of loops [5, 9, 11]. We shall illustrate them by examples.

*Example 1.4.6*: Let $L = L_1 \cup L_2 = \{[0, a] / a \in \{e, 1, 2, 3, 4, 5\}, 2, *\} \cup \{[0, b] / b \in \{e, 1, 2, 3, 4, 5\}, 3, *\}$ be biinterval biloop given by the following tables.

| * | [0, e] | [0, 1] | [0, 2] | [0, 3] | [0, 4] | [0, 5] |
|---|---|---|---|---|---|---|
| [0, e] | [0, e] | [0, 1] | [0, 2] | [0, 3] | [0, 4] | [0, 5] |
| [0, 1] | [0, 1] | [0, e] | [0, 3] | [0, 5] | [0, 2] | [0, 4] |
| [0, 2] | [0, 2] | [0, 5] | [0, e] | [0, 4] | [0, 1] | [0, 3] |
| [0, 3] | [0, 3] | [0, 4] | [0, 1] | [0, e] | [0, 5] | [0, 2] |
| [0, 4] | [0, 4] | [0, 3] | [0, 5] | [0, 2] | [0, e] | [0, 1] |
| [0, 5] | [0, 5] | [0, 2] | [0, 4] | [0, 1] | [0, 3] | [0, e] |

| * | [0, e] | [0, 1] | [0, 2] | [0, 3] | [0, 4] | [0, 5] |
|---|---|---|---|---|---|---|
| [0, e] | [0, e] | [0, 1] | [0, 2] | [0, 3] | [0, 4] | [0, 5] |
| [0, 1] | [0, 1] | [0, e] | [0, 4] | [0, 2] | [0, 5] | [0, 3] |
| [0, 2] | [0, 2] | [0, 4] | [0, e] | [0, 5] | [0, 3] | [0, 1] |
| [0, 3] | [0, 3] | [0, 2] | [0, 5] | [0, e] | [0, 1] | [0, 4] |
| [0, 4] | [0, 4] | [0, 5] | [0, 3] | [0, 1] | [0, e] | [0, 2] |
| [0, 5] | [0, 5] | [0, 3] | [0, 1] | [0, 4] | [0, 2] | [0, e] |

Now let $S = S_1 \cup S_2 = \{[0, a] / a \in \{e, 1, 2, 3, 4, 5\}, \otimes, 2\} \cup \{[0, b] / b \in \{e, 1, 2, 3, 4, 5\}, \otimes, 3\}$ be the principal bi isotopes of $L = L_1 \cup L_2$ given by the following tables.

| $\otimes$ | [0, e] | [0, 1] | [0, 2] | [0, 3] | [0, 4] | [0, 5] |
|---|---|---|---|---|---|---|
| [0, e] | [0, 3] | [0, 2] | [0, 5] | [0, e] | [0, 1] | [0, 4] |
| [0, 1] | [0, 5] | [0, 3] | [0, 4] | [0, 1] | [0, e] | [0, 2] |
| [0, 2] | [0, 4] | [0, e] | [0, 3] | [0, 2] | [0, 5] | [0, 1] |
| [0, 3] | [0, e] | [0, 1] | [0, 2] | [0, 3] | [0, 4] | [0, 5] |
| [0, 4] | [0, 2] | [0, 5] | [0, 1] | [0, 4] | [0, 3] | [0, e] |
| [0, 5] | [0, 1] | [0, 4] | [0, e] | [0, 5] | [0, 2] | [0, 3] |

and the table of $S_2$ is as follows :

| $\otimes$ | [0, e] | [0, 1] | [0, 2] | [0, 3] | [0, 4] | [0, 5] |
|---|---|---|---|---|---|---|
| [0, e] | [0, 4] | [0, 5] | [0, 3] | [0, 1] | [0, e] | [0, 2] |
| [0, 1] | [0, 3] | [0, 2] | [0, 5] | [0, 3] | [0, 1] | [0, 4] |
| [0, 2] | [0, 5] | [0, 3] | [0, 1] | [0, 4] | [0, 2] | [0, e] |
| [0, 3] | [0, 2] | [0, 4] | [0, e] | [0, 5] | [0, 3] | [0, 1] |
| [0, 4] | [0, e] | [0, 1] | [0, 2] | [0, 3] | [0, 4] | [0, 5] |
| [0, 5] | [0, 1] | [0, e] | [0, 4] | [0, 2] | [0, 5] | [0, 3] |

Now we can define as in case of bi interval loops define Smarandache simple interval biloops. This task is also left to the reader.

***Example 1.4.7***: Let $L = L_1 \cup L_2 = \{[0, a] / a \in \{e, 1, 2, …, n\}$, *, m such that $(m – 1, n) = 1 = (m, n), 1 < m < n\} \cup \{[0, b] / b \in \{e, 1, 2, …, t\}$, *, s such that $\{s – 1, t) = (t, s) = 1$ and $1 < s < t\}$ be an interval biloop. It is easily verified L is a Smarandache interval bisimple loop.

We have the following theorem which guarantees the existence is a class of such interval biloops.

**THEOREM 1.4.2**: *$L = L_1 \cup L_2 = \{[0, a] / a \in \{e, 1, 2, ..., m\}, *, t$ where $(t, m) = (t - 1, m) = 1; 1 < t < m\} \cup \{[0, b] / b \in \{e, 1, 2, ..., n\}, *, s$ such that $(s, n) = (s-1, n) = 1$ and $1 < s < n\}$ ($m \neq n$ or if $m = n$ then $s \neq t$) be an interval biloop. L is a Smarandache simple interval biloop.*

The proof is direct and is left as an exercise for the reader to prove.

In fact we have for every t in L such that $(t, m) = 1$ and $(t - 1, m) = 1$ we have an interval loop hence we have a class of interval loops associated with each t which we can denote by $\{L_m[0,a](t)\}$. Similarly for $L_2$ we have a class denoted by $\{L_n [0, b] (s)\}$. So this class $\{L_m [0,a] (t)\} \cup \{L_n [0, b] (s)\}$ is a Smarandache interval bisimple loop.

We will define an interval biloop to be A Smarandache weakly Lagrange biloop if there exist atleast one interval bisubgroup $H = H_1 \cup H_2$ in $L = L_1 \cup L_2$ such that $o(H)/o(L)$ that is $|H_1| |H_2| / |L_1| |L_2|$.

We will first illustrate this situation by an example.

***Example 1.4.8***: Let $L = L_1 \cup L_2 = \{[0, a] / a \in \{e, 1, 2, ..., 15\}, *, 2\} \cup \{[0, b] / b \in \{e, 1, 2, ..., 7\}, *, 6\}$ be an interval biloop.

Clearly L is only a Smarandache weakly Lagrange biinterval loop.

For take $H = H_1 \cup H_2 \{[0, e], [0, 8], *, 2\} \cup \{[0, e], [0, 3], *, 6\} \subseteq L_1 \cup L_2$. H is an biinterval subgroup of G. $o(H) = 2 \times 2 = 4$. Now $o(L) = |L_1| |L_2| = 16.8$. $4 / 16.8$.

Consider $T = T_1 \cup T_2 = \{[0, a] / a \in \{e, 1, 2, 5, 8, 11, 14\} \subseteq \{e, 1, 2, 3, ..., 15\} *, 2\} \cup \{[0, e], [0, 4], *, 6\} \subseteq L_1 \cup L_2$, T is such that $|T_1| \cdot |T_2| = 6.2 \nmid 16.8$. Hence L is only a Smarandache weakly Lagrange interval biloop.

Infact we have a non empty class of Smarandache weakly Lagrange interval biloop.

**THEOREM 1.4.3**: *$L = \{L_n [0, a] (t)\} \cup \{L_m [0, b] (s) \}$ ($m \neq n$) be a class of biinterval loop (interval biloops). Every interval* *biloop in L is a Smarandache weakly Lagrange interval bilooop.*

The proof is straight forward and is hence left as an exercise to the reader [9].

We can define Smarandache Lagrange interval biloop in a similar way [9] and it is also easy to verify that every Smarandache Lagrange interval biloop is a Smarandache weakly Lagrange interval biloop [9].

Now we can define as in case of usual loops the notion of Smarandache Cauchy interval biloop.

Let $L = L_1 \cup L_2$ be an interval biloop. Let $x = [0, a] \cup [0, b] \in L$ we say x is a Smarandache-Cauchy biinterval element of L if $x^r = [0, a]^{r_1} \cup [0, b]^{r_2} = [0, 1] \cup [0, 1]$ and $r_1 > 1$, $r_2 > 2$ with $r_1 r_2 / |L_1| |L_2|$, otherwise x is not a Smarandache Cauchy bi interval element of L.

*Example 1.4.9*: Let $L = L_1 \cup L_2 = \{[0, a] / a \in \{e, 1, 2, …, 11], *, 3\} \cup \{0, a] / a \in \{e, 1, 2, …, 13\}, *, 7\}$ be an interval biloop.

Clearly $o(L) = |L_1| \cdot |L_2| = 12.14$. Consider $x = [0,2] \cup [0, 4] \in L$. $x * x = [0, 2] \times [0, 2] \cup [0, 4] [0, 4] = [0, 2 * 2] \cup [0, 4 * 4] = [0, e] \cup [0, e]$. Thus $r_1 = r_2 = 2$. Now $2.2/|L_1| |L_2|$. Thus x is a S-Cauchy biinterval element of L.

Recall if every bielement in L is a S-Cauchy bi interval element of L then we define L to be a Smarandache Cauchy interval biloop (S-Cauchy interval biloop) we will show we have a class of interval biloops which are S-Cauchy biinterval loops.

**THEOREM 1.4.4**: *Let $L = \{L_1\} \cup \{L_2\} = \{L_n[0, a](t)\} \cup \{L_m [0, b](s)\}$ $n \neq m$, be a class of biinterval loops. Every interval biloop in L is a Smarandache Cauchy interval biloop.*

Proof is direct for we see in every interval biloop every bi interval element $x = [0, a] \cup [0, b] \in L$ is of biorder 2.2 and since every interval loop $L_n [0, a] (t)$ is of even order $2.2 / |L_1| \cdot |L_2|$. One can extend the notion of Smarandache pseudo Lagrange loops to interval biloops.

Let $L = L_1 \cup L_2$ be an interval biloop of finite order. If the biorder of every interval S-subbiloop (S-biinterval subloop or S-interval subbiloop) divides the biorder $|L_1| |L_2|$ then we say L is a Smarandache pseudo Lagrange biinterval loop or Smarandache pseudo Lagrange interval biloop.

If L has atleast one S-interval subbiloop $K = K_1 \cup K_2 \subseteq L_1 \cup L_2$ such that $|K_1| . |K_2| / |L_1| . |L_2|$ then we say L is a Smarandache weakly pseudo Lagrange interval biloop [5, 9].

It is easily verified that every S-pseudo biinterval Lagrange loop is a S-weakly pseudo Lagrange biinterval loop [5,9,11]. Interested reader is expected to construct examples to this effect.

We can define Smarandache p-Sylow interval subbiloops as in case of loops [5, 9, 11]. We also can as in case of loops define Smarandache p-Sylow subgroup [5, 9, 11].

Further for Smarandache interval bisubgroup biloop we can define strong Sylow substructures. Let $L = L_1 \cup L_2$ be a Smarandache biinterval subgroup loop of finite biorder; if every interval subbigroup is either of a prime power biorder and that bidivides $o(L) = |L_1| |L_2|$ then we call L to be a Smarandache strong interval p-Sylow biloop (S-strong interval p-Sylow biloop) [5, 9, 11].

We have a class of Smarandache strong 2-Sylow biloops.

***Example 1.4.10:*** $L = L_1 \cup L_2 = L_5 [0,a] (3) \cup L_7 [0, b] (3)$ be an interval biloop given by the following tables.

Table of $L_5 [0, a] (3)$

| *     | [0, e] | [0, 1] | [0, 2] | [0, 3] | [0, 4] | [0, 5] |
|-------|--------|--------|--------|--------|--------|--------|
| [0, e]| [0, e] | [0, 1] | [0, 2] | [0, 3] | [0, 4] | [0, 5] |
| [0, 1]| [0, 1] | [0, e] | [0, 4] | [0, 2] | [0, 5] | [0, 3] |
| [0, 2]| [0, 2] | [0, 4] | [0, e] | [0, 5] | [0, 3] | [0, 1] |
| [0, 3]| [0, 3] | [0, 2] | [0, 5] | [0, e] | [0, 1] | [0, 4] |
| [0, 4]| [0, 4] | [0, 5] | [0, 3] | [0, 1] | [0, e] | [0, 2] |
| [0, 5]| [0, 5] | [0, 3] | [0, 1] | [0, 4] | [0, 2] | [0, e] |

Table for $L_7$ [0,b] (3)

| * | [0, e] | [0, 1] | [0, 2] | [0, 3] | [0, 4] | [0, 5] | [0, 6] | [0, 7] |
|---|---|---|---|---|---|---|---|---|
| [0, e] | [0, e] | [0, 1] | [0, 2] | [0, 3] | [0, 4] | [0, 5] | [0, 6] | [0, 7] |
| [0, 1] | [0, 1] | [0, e] | [0, 4] | [0, 7] | [0, 3] | [0, 6] | [0, 2] | [0, 5] |
| [0, 2] | [0, 2] | [0, 6] | [0, e] | [0, 5] | [0, 1] | [0, 4] | [0, 7] | [0, 3] |
| [0, 3] | [0, 3] | [0, 4] | [0, 7] | [0, e] | [0, 6] | [0, 2] | [0, 5] | [0, 1] |
| [0, 4] | [0, 4] | [0, 2] | [0, 5] | [0, 1] | [0, e] | [0, 7] | [0, 3] | [0, 6] |
| [0, 5] | [0, 5] | [0, 7] | [0, 3] | [0, 6] | [0, 2] | [0, e] | [0, 1] | [0, 4] |
| [0, 6] | [0, 6] | [0, 5] | [0, 1] | [0, 4] | [0, 7] | [0, 3] | [0, e] | [0, 2] |
| [0, 7] | [0, 7] | [0, 3] | [0, 6] | [0, 2] | [0, 5] | [0, 1] | [0, 4] | [0, e] |

We see every interval bielement x = [0, a] $\cup$ [0, b] $\in$ L is such that $x^2$ = [0, 1] $\cup$ [0, 1]. Thus L is a Smarandache strong 2-Sylow interval biloop.

**THEOREM 1.4.5**: *Let = {$L_p$ [0, a] (t)} $\cup$ {$L_q$ [0,b] (s)} (p and q two distinct primes) be a Smarandache interval biloop of order |p+1|. |q+1|. Then every interval biloop in L is a Smarandache strong biinterval 2-Sylow loop.*

Proof follows from the fact the biorder of L = |$L_1$| . |$L_2$| = |p+1| |q+1| = even number × even number and as p and q are two distinct primes every interval bielement in L is of biorder 2.2. Hence each interval biloop in L is a Smarandache strong 2-Sylow biloop.

Now we can define as in case of usual S-loops the notion of Smarandache interval biloop homomorphism and S-interval biloop isomorphism [5, 9, 11]. Interested reader can substantiate this with examples. As in case of general loops we define Smarandache commutative interval biloop [5, 9, 11].

*Example 1.4.11*: Let L = $L_1$ $\cup$ $L_2$ = {[0, a] / a $\in$ {e, 1, 2, …, 7}, *, 3} $\cup$ {[0, b] / b $\in$ {e, 1, 2, …, 11}, *, 3} be a biinterval loop. Clearly L is a S-commutative interval biloop.

We say an interval biloop $L = L_1 \cup L_2$ to be a Smarandache strongly commutative (S-strongly commutative) interval billoop if every proper bisubset $S = S_1 \cup S_2 \subseteq L_1 \cup L_2$ which is a biinterval group (interval bigroup) is a commutative interval bigroup.

Also it is clear from the very definition if $L = L_1 \cup L_2$ is S-strongly commutative interval biloop then L is a Smarandache commutative interval biloop.

*Example 1.4.12*: Let $L = L_1 \cup L_2 = \{[0, a] / a \in \{e, 1, 2, …, 19], *, 8\} \cup \{[0, b] / b \in \{e, 1, 2, …, 23\}, *, 9\}$ be an interval biloop. It is easily verified that L is a S-substrongly commutative interval biloop.

*Example 1.4.13*: Let $L = L_1 \cup L_2 = \{[0, a] / a \in \{e, 1, 2, …, 21], 11, *\} \cup \{[0, b] / b \in \{e, 1, 2, …, 15\}, 8, *\}$ be a bi interval loop.

Clearly L is S-commutative biinterval biloop.

Now we have the following interesting theorem.

**THEOREM 1.4.6**: *Let $L = \{L_p [0, a] (t)\} \cup \{L_q [0, b] (s)\}$, p and q two distinct primes be a class of biinterval loops. Every interval biloop in L is a S-strongly commutative interval biloop.*

Follows from the fact every proper bisubset $A = A_1 \cup A_1 \subseteq L$ which is an interval subbigroup is of biorder 2.2, and has no other interval subbiloops. Hence the claim.

Now as in case of loops we can in case of interval bilooops also define the notion of Smarandache cyclic interval biloop or Smarandache bicyclic interval biloop or Smarandache cyclic biinterval loop [5, 9, 11]. Also the notion of Smarandache strong cyclic biinterval loop as in case of loops. It is easily verified that every S-strong biinterval cyclic loop is a S-cyclic interval biloop.

We will illustrate this by some examples.

*Example 1.4.14*: $L = L_1 \cup L_2 = \{[0, a] / a \in \{e, 1, 2, 3,4,5\}, 4, *\} \cup \{[0, b] / b \in \{e, 1, 2, …, 29\}, *, 7\}$ be a biinterval loop. It is easily verified that L is S-strongly cyclic interval biloop.

We have a very large classes of class of S-strongly cyclic interval biloop which is stated in the following theorem the proof of which is left as an exercise to the reader.

**THEOREM 1.4.7**: *Let $L = \{L_p [0, a] (t)\} \cup \{L_q [0, b] (s)\}$, p and q two distinct primes be a class of biinterval loop. Then every biinterval loop in L is a S-strongly cyclic interval biloop.*

In fact by varying the primes p and q over the set of primes we can get an infinite class of S-strongly cyclic interval biloops.

We have the following interesting theorem which gives the number of strictly non commutative interval biloops which are S-strongly commutative interval biloops and S-strongly cyclic interval biloops.

**THEOREM 1.4.8**: *Let $L = \{L_n [0, a] (t)\} \cup \{L_m [0, b] (s)\}$ ($m \neq n$, $m > 3$, $n > 3$) where $n = p_1^{\alpha_1} \cdots p_k^{\alpha_k}$ and $m = q_1^{t_1} \cdots q_s^{t_s}$ with $\alpha_i \geq 1$, $t_j \geq 1$, $1 \leq i \leq k$ and $1 \leq j \leq s$. Then L contains exactly $F_n \cdot F_m$ interval biloops which are strictly non commutative and they are*
1. *S-strongly commutative interval biloops and*
2. *S-strongly cyclic interval bilooops where*

$$F_n = \prod_{i=1}^{K} (p_i - 3) p_i^{\alpha_i - 1} \quad , \quad F_m = \prod_{j=1}^{s} (q_j - 3) q_j^{t_j - 1}.$$

We can as in case of usual biloops define Smarandache pseudo commutative interval biloop and Smarandache strongly pseudo commutative interval biloop. We can also define the notion of Smarandache commutator interval bisubloop denoted by $L^S = L_1^s \cup L_2^s$. We have the following interesting result viz. if $L = L_1 \cup L_2$ be an S-interval biloop which has no S-interval subbiloops then $L' = L^s = L_1' \cup L_2' = L_1^s \cup L_2^s$ where $L' = L_1' \cup L_2' = \langle\{[0,x] \in L_1 / [0, x] = ([0, y], [0,z])$ for some $[0,y], [0,z] \in L_1\}\rangle \cup \langle\{[0,y] \in L_2 / [0,y] = ([0,s], [0,t]$ for some $[0,s], [0,t]$ in $L_2\}\rangle$ and $L^s = L_1^s \cup L_2^s = \{ L_1^s$ is the interval subloop generated by all the interval commutators $A_1$, $A_1$ a S-interval subloop of the interval loop $L_1$} ∪ { $L_2^s$ is the interval subloop generated by all the interval commutators $A_2$, $A_2$ a S-interval subloop of the interval loop $L_2$}. We have the following theorem.

**THEOREM 1.4.9**: *Let $L = \{L_n [0,a] (t) \mid n$ is a prime; $0 < t < n\}$ ∪ $\{L_m [0, b] (s) \mid m$ is a prime, $0 < s < m\}$ be a class of non commutative interval biloops. Then we have $L' = \{L'_l [0,a] (t)\}$ ∪ $\{L'_m [0, b](s)\} = L^s = \{L_n^s [0,a] (t)\}$ ∪ $\{L_m^s [0, a] (s)\}$ for every pair taken in L.*

Here $L'_n [0, a] (t) = \langle\{[0, d] \in L_n [0,a] (t) \mid [0, d] = ([0, b], [0, c])$ for some $[0, b], [0, c]$ in $L_n [0, a] (t)\}\rangle$. Similarly $L'_m [0, b] (s)$ is defined.

The proof is direct for more information please refer [5].

We can as in case of loops define for interval biloops the notion of Smarandache associative interval biloops [5].

*Example 1.4.15*: Let $L = L_1 \cup L_2 = \{[0, a] / a \in \{e, 1, 2, \ldots, 13\}, *, 5\} \cup \{[0,b] \mid b \in \{e, 1, 2, \ldots, 17\} *, 7\}$ be an interval biloop. Clearly L is not a S-associative interval biloop.

We have a very large class of interval biloops which are not S-associative interval biloops. This is evident from the following theorem, the proof of which is straight forward.

**THEOREM 1.4.10**: *Let $L = \{L_n [0, a] (t)\} \cup \{L_m [0, b] (s)\}$ (n and m are two distinct primes $1 < t < n$ and $1 < s < m$) be class of interval biloops. None of the interval biloops in this class of $(n+1)(m+1)$ biloops is an S-associative interval biloop.*

The notion of Smarandache strongly pairwise associative interval biloops can be defined in an analogous way for interval biloops. We will illustrate this by an example.

*Example 1.4.16* : Let $L = L_1 \cup L_2 = \{[0, a] / a \in \{e, 1, 2, \ldots, 19\}, t = 7, *, 1 < t < 19\} \cup \{[0, b] / b \in \{e, 1, 2, \ldots, 11\}, *, s = 5, 1 < s < 11\}$ m and n are distinct be a S-strongly pairwise associative interval biloop.

**THEOREM 1.4.11**: *Let $L = L_1 \cup L_2 = \{L_n[0, a](t) / 0 < t < n, *\} \cup \{L_m [0, b](s) \mid 0 < s < m, *\}$ ($m \neq n$, $m > 3$, $n > 3$; $m$ and $n$ positive integers) be a class of interval biloops. $L$ is a class of S-strongly pairwise associative interval biloops.*

We give only hint of the proof.

Let $a = [0, x] \cup [0, y] \in L_1^1 \cup L_2^2$ where $L_1^1$ and $L_2^2$ are interval loops from the class of loops $L_1$ and $L_2$ respectively. Similarly $b = [0, p] \cup [0, s] \in L_1^1 \cup L_2^2$.

Now clearly using some simple number theoretic techniques we have

$$\begin{align} (a\,b)a &= (([0, x] \cup [0, y])([0, p] \cup [0, s]) \times ([0, x] \cup [0, y]) \\ &= [0, \{xp) x] \cup [0\ (ys)\ y] \\ &= [0, x(px)] \cup [0, y(sy)] \\ &= a(ba) \ [\ ]. \end{align}$$

We can also define the notion of Smarandache associator interval subbiloop denoted by $L^A = L_1^{A_1} \cup L_2^{A_2} = \{<$interval subloop generated by all the interval associators in $A_1$, where $A_1$ is a S-interval subloop of $L_1$; that is $A_1 \subseteq L_1>\} \cup \{<$interval subloop generated by all the interval associators in $A_2$, where $A_2$ is a S-interval subloop of $L_2$ that is $A_2 \subseteq L_2>\}$.

If $L$ is a S-interval biloop which has no S-interval bisubloops then we have $A(L) = L^A = A_1 (L_1) \cup A_2 (L_2) = L_1^{A_1} \cup L_2^{A_2}$ that is the associator interval subbiloop of $L$ coincides with the S-associator interval subbiloop [ ].

We have the following theorem which guarantees such class of interval biloops.

**THEOREM 1.4.12**: *Let $L = L_1 \cup L_2 = \{L_n [0, a] (t) ; *, 1 < t < n\} \cup \{L_m [0,a] (s), *, 1 < s < m\}$ be a class of interval biloops where $L_i$'s are S-interval loops and has no S-interval subloops, $i = 1, 2$. Then $A\{L\} = A(L_1) \cup A(L_2) = L_1^A \cup L_2^A$.*

(Here we mean every pair of interval biloops in the class satisfies the condition which is represented by the class $L = L_1 \cup L_2$).

In general we cannot say this for S-interval biloops which has S-interval subbiloops. Interested reader can construct examples to this effect [5, 9, 11]. Consequent of this we have the following theorem.

**THEOREM 1.4.13**: *Let $L = L_1 \cup L_2$ be an S-interval biloop having a S-interval subbiloop $A = A_1 \cup A_2$, then $A(L) \neq L^A$ that is $A(L) = A_1(L_1) \cup A_2(L_2) \neq L_1^{A_1} \cup L_2^{A_2}$.*

This proof is by constructing counter examples to this effect [5, 9, 11]. We can as in case of loops define for interval biloops the notion of S-first normalizer and S-second normalizer and obtain the condition for the S-first normalize to be equal to S-second normalize for interval biloops.

***Example 1.4.17***: Let $L = L_1 \cup L_2 = \{[0, a] / a \in \{e, 1, 2, \ldots, 21\}, *, 8\} \cup \{[0, b] / b \in \{e, 1, 2, \ldots, 33\}, *, 5\}$ be an interval biloop. Clearly $H_1^1(7)$ is the interval subloop of L is given by the following table $H_1^1(7) = \{[0, e], [0, 1], [0, 8], [0, 15]\}$.

| * | [0, e] | [0, 1] | [0, 8] | [0, 15] |
|---|---|---|---|---|
| [0, e] | [0, e] | [0, 1] | [0, 8] | [0, 15] |
| [0, 1] | [0, 1] | [0, e] | [0, 15] | [0, 8] |
| [0, 8] | [0, 8] | [0, 15] | [0, e] | [0, 1] |
| [0, 15] | [0, 15] | [0, 8] | [0, 1] | [0, e] |

$H_1^2(11) = \{[0, e], [0, 11], [0, 12], [0, 23]\}$ is an interval subloop of $L_2$ given by the following table:

| * | [0, e] | [0, 1] | [0, 12] | [0, 23] |
|---|---|---|---|---|
| [0, e] | [0, e] | [0, 1] | [0, 12] | [0, 23] |
| [0, 1] | [0, 1] | [0, e] | [0, 23] | [0, 12] |
| [0, 12] | [0, 12] | [0, 23] | [0, e] | [0, 1] |
| [0, 23] | [0, 23] | [0, 12] | [0, 1] | [0, e] |

Consider the interval bisubloop $H_1^1(7) \cup H_1^2(11) \subseteq L_1 \cup L_2$. It is easily verified that $SN_1(H_1^1(7)) \cup SN_1(H_1^2(11)) = SN_2(H_1^1(7)) \cup SN_2(H_1^2(11))$.

In view of this we have the following important theorem.

**THEOREM 1.4.14**: *Let $L = L_1 \cup L_2 = \{L_n[0, a](t)\} \cup \{L_m[0, b](s)\}$ be the class of interval biloops (this forms a class of interval biloops). For any pair of interval biloops from L say $L_n^1 \cup L_m^2 \in L_1 \cup L_2$, let $H_{in}^1([0, a](p)) \cup H_{im}^2([0,b](q)) \subseteq L_n' \cup L_m^2$ be its interval S-subbiloop.*
*Then*
$$SN_1(H_{in}^1([0, a](p)) \cup SN_2(H_{im}^2([0,b](q))$$
$$= SN_2(H_{in}^1([0, a](p)) \cup SN_2(H_{im}^2([0,b](q))$$
*if and only if*
$(t^2 - t + 1, p) = (2t - 1, p)$ and $(s^2 - s + 1, q) = (2s - 1, q)$.

*Proof*: Let $L_n^1 \cup L_m^2$ be as in theorem and $H_{in}^1([0, a](p)) \cup H_{im}^2([0,b](q)) \subseteq L_n' \cup L_m^2$ be an S-interval bisubloop of $L_n^1 \cup L_m^2$. First we show that first interval S-binormalizer

$SN_1(H_{in}^1([0, a](p))) \cup SN_1(H_{im}^2([0,b](q))) = H_{in}^1(k_1) \cup H_{im}^2(k_2)$ where $k = k_1 \cup k_2 = p/d_1 \cup q/d_2$ and $d_1 \cup d_2 = (2t-1, p) \cup (2s-1, q)$, we use only simple number theoretic arguments and the definition $SN_1(H_{in}^1([0, a](p)) \cup SN_2(H_{im}^2([0,b](q)) = \{[0,j_1] \in L_n' / [0, j_1] H_{in}^1([0, a](p)) = H_{in}^1([0, a](p))[0, j_1]\} \cup \{[0,j_2] \in L_n^2 / [0, j_2] H_{im}^2([0,b](q)] = H_{im}^2([0,b](q))[0, j_2]\}$ is the S-first interval binormalizer of $H_{in}^1([0, a](p)) \cup H_{im}^2([0,b](q))$.

It is left for the reader to verify $[0, j_1] H_{in}^1([0, a](p)) \cup H_{in}^1([0,a](p)).[0, j_1]$ if and only if $(2t-1)(i - j_1) \equiv 0 \mod (p)$ and $[0, j_1] H_{im}^2([0, a](q)) = H_{im}^2([0, b](q))[0, j_2]$ if and only if $(2s - 1)(i - j_2) \equiv 0 \pmod{q}$.

For $[0, j_1] \notin H^1_{in}([0, a](p))$ and $[0, j_2] \notin H^2_{im}([0,b](q))$ further if $[0, j_1] \in H^1_{in}([0, a](p))$ we have $[0, j_1] H^1_{in}([0,a](p))$ = $H^1_{in}([0,a](p))[0, j_1]$.

Similar argument holds for $H^2_{im}([0,b](q))$ and $[0, j_2]$.

Now for the other part reasoning is done as in case of usual loops. Please refer [ ].

We have an interesting result relating $L_n([0, a)(t)) \cup L_m([0, b](s))$ where n and m are two distinct primes and $SN(L_n([0, a](t)) \cup SN(L_m([0, b](s))$.

**THEOREM 1.4.15**: *Let $L = L_1^n \cup L_2^m = \{L_n([0, a](t)\} \cup \{L_m([0, b])(s)\}$ be a class of interval biloops where n and m are two distinct primes. Then for every pair of interval biloops $L_n^1([0, a])(t) \cup L_m^2([0, b])(s) \in L$ we have $SN(L_n^1([0, a])(t)) \cup SN(L_m^2([0, b])(s)) = \{e\} \cup \{e\}$.*

For proof refer [5, 9, 11].

Analogously one can derive for interval biloops with appropriate changes.

We can for interval biloops define S-Moufang bicenter. This is easily done by suitably extending to interval biloops.

We have the following theorems the proofs can be obtained as in [ 5, 9, 11] with appropriate modifications.

**THEOREM 1.4.16**: *Let $L = L_1 \cup L_2 = \{L_n[0, a](t)\} \cup \{L_m[0, b](s)\}$, n and m are two distinct primes be a class of interval bilooops, then S-Moufang bicenter of every interval biloop from L say $L_n^1([0, a])(t) \cup L_m^2([0, b])(s)$ for fixed t and s is either $\{e\} \cup \{e\}$ or $L_n^1([0, a])(t) \cup L_m^2([0, b])(s)$.*

**THEOREM 1.4.17**: *Let $L = L_1 \cup L_2 = \{L_n[0, a](t)\} \cup \{L_m[0, b](s)\}$ be a class of interval biloops where n and m are two distinct primes.*

*Then $NZ(L_n^1([0, a])(t)) \cup NZ(L_m^2([0, b])(s)) = Z(L_n^1([0, a])(t)) \cup Z(L_m^2([0, b])(s)) = \{e\} \cup \{e\}$ for every pair of*

*interval biloops $L_n^1$ ([0, a]) (t)) $\cup$ $L_m^2$ ([0, b]) (s) in $L_1 \cup L_2$ for a fixed t and s.*

For proof refer [5, 9] and obtain the proof with proper modifications for interval biloops. We can as in case of loops define direct product in case of interval biloops and derive their related properties. Now using interval loops we can define interval group - loop, interval semigroup - loop, quasi interval biloops and interval groupoid - loop.

Now we proceed onto define these structures and study some of their related properties.

**DEFINITION 1.4.2**: *Let $L = L_1 \cup L_2$ where $L_1$ be an interval loop and $G_2$ is just a loop. We call L a quasi interval biloop. The operations from $L_1$ and $G_2$ are carried over to L.*

We will illustrate this with examples.

***Example 1.4.18***: Let $L = L_1 \cup L_2$ where $L_1$ is $L_5(3)$ and $L_2 = \{[0,a] / a \in \{e, 1, 2, …, 7\}, *, 6\}$ be a quasi interval biloop. We will just show how on L operations are carried out. Let $x = 2 \cup [0,4]$ and $y = 4 \cup [0, 6]$ be in $L = L_1 \cup L_2$,

x . y   =   $(2 \cup [0, 4]) . (4 \cup [0,6])$
        =   $2 * 4 \cup [0, 4] * [0, 6]$
        =   $(4 \times 3 - 2 \times 2) \pmod 5 \cup \{[0, 4 \times 6]\}$
        =   $\{12 + 1 \pmod 5\} \cup \{(0, \{36 - 4 \times 5\} \pmod 7)\}$
        =   $3 \cup [0, 2] \in L_1 \cup L_2 = L.$

$e \cup [0, e] \in L$ acts as the identity element. This quasi interval biloop is of order $6 \times 8 = 48$.

***Example 1.4.19***: Let $L = L_1 \cup L_2 = \{[0,a] \mid a \in \{e, 1, …, 11\}, *, 9\} \cup L_7(3)$ be a quasi interval biloop of order $20 \times 8 = 160$. $L = \{[0, a] \cup b / a \in \{e, 1, 2, …, 11\}$ and $b \in \{e, 1, 2, …, 7\}\}$. Operations on them can be carried out as shown in example 1.4.18. Now we can define quasi interval subbiloop, quasi interval S-biloops and so on. Several of the theorems proved for interval biloops can be derived also for quasi interval biloops with some simple changes.

We will illustrate these situations by some examples.

***Example 1.4.20***: L = $L_1 \cup L_2$ = {[0,a] / a ∈ {e, 1, 2, …, 21}, *, 5} $\cup$ $L_{33}$(8) be a quasi interval biloop. L is a S-quasi interval biloop. For take A = $A_1 \cup A_2$ = {[0, e], [0, 10]} $\cup$ {e, 6} $\subseteq L_1 \cup L_2$ is a quasi interval bigroup, hence L is a S-quasi interval biloop. Consider $H_1(7)$ = {[0, e], [0,1], [0,8], [0,15]} and $H_1(11)$ = {e, 1, 12, 23} subloops of $L_1$ and $L_2$ respectively. H = $H_1(7) \cup H_1(11) \subseteq L_1 \cup L_2$ is not only a quasi interval subbiloop but also H is a S-quasi interval bisubloop of order $4 \times 4 = 16$.

***Example 1.4.21***: Let . L = $L_1 \cup L_2$ = {[0,a] / a ∈ {e, 1, 2, …, 19}, *, 8} $\cup$ $L_{17}$ (3) be a quasi interval biloop. L is a S-quasi interval biloop.

For take A = $A_1 \cup A_2$ = {[0, e], [0, 3]} $\cup$ {e, 9} $\subseteq L_1 \cup L_2$ ; A is a quasi interval bigroup, hence L is a S-quasi interval biloop. Infact we have a class of S-quasi interval biloop.

**THEOREM 1.4.18**: *Let . L = $L_1 \cup L_2$ = $L_n(m) \cup$ {[0,a] / a ∈ {e, 1, 2, ..., t}, *, s; 1 < s < t (t, s) = 1 = (s – 1, t)} be a quasi interval biloop. L is a Smarandache quasi interval biloop.*

The proof is straight forward and hence is left as an exercise for the reader. In fact a class of quasi interval biloops exists. For in the theorem m and s can vary and we have a class. If we vary n and t we get classes of quasi interval biloops of finite order.

We have a class of quasi interval biloops which are bisimple. We will first illustrate by an example.

***Example 1.4.22***: Let L = $L_9$ (8) $\cup$ {[0,a] / a ∈ {e, 1, 2, …, 19}, *, 11} be a quasi interval biloop which is clearly bisimple.

In fact we have a class of bisimple quasi interval biloops.

**THEOREM 1.4.19**: *L = $L_1 \cup L_2$ = $L_n \cup \{L_t[0,a](s)\}$ be a class of quasi interval biloops. Clearly every pair of quasi interval biloops are simple.*

Proof is straight forward as L has no non-trivial quasi interval normal bisubloops. Hence the claim.

*Example 1.4.23*: Let $L = L_1 \cup L_2 = \{L_{19}(3)\} \cup \{[0,a] / a \in \{e, 1, 2, 3, \ldots, 23\}, *, 5\}$ be quasi interval biloop. Clearly L is a Smarandache quasi interval bisubgroup biloops.

We have a class of S-quasi interval bisubgroup biloop.

**THEOREM 1.4.20**: *Let $L = L_1 \cup L_2 = L_n \cup \{\{L_t[0,a] (s)\}$ where n and t are distinct primes, be a class of quasi interval biloops. L is a S-quasi interval subgroup biloop.*

Now we have a class of quasi interval S-Cauchy biloops.

*Example 1.4.24*: Let $L = L_{15}(8) \cup \{[0,a] / a \in \{e, 1, 2, \ldots, 21\}, 11, *\}$ be a quasi interval biloop. Clearly L is a quasi interval S-Cauchy biloop or S-Cauchy quasi interval biloop.

**THEOREM 1.4.21**: *Let $L = L_1 \cup L_2 = \{L_m\} \cup \{L_t [0,a] (s) / 1 < s < t, *\}$ be a class of quasi interval biloops. Every quasi interval biloop in this class is a S-Cauchy quasi interval biloop.*

The proof is left as an exercise to the reader.

*Example 1.4.25*: Let $L = L_1 \cup L_2 = L_{15}(2) \cup \{[0,a] / a \in \{e, 1, 2, \ldots, 15\}, *, 2\}$ be a quasi interval biloop of order $16 \times 16$.

We see L has quasi interval S-subbiloops. $H = H_1 \cup H_2 = \{e, 2, 5, 8, 11, 14\} \cup \{[0, a] / a \in e, 2, 5, 8, 11, 14\} \subseteq L_1 \cup L_2$ is a quasi interval S-subbiloop of order $6 \times 6$. Clearly o(H) $\nmid$ {L} that is 6.6 $\nmid$ 16.16.

**THEOREM 1.4.22**: *Let $L_n \cup \{L_m [0,a] (t)\} = L$ be a class of quasi interval biloops. Every quasi interval biloop in L is a S-weakly Lagrange quasi interval biloop.*

Proof is straight forward and is left as an exercise for the reader.

We have a class of quasi interval biloops which are Smarandache strong quasi interval 2-Sylow biloops.

*Example 1.4.26*: Let $L = L_1 \cup L_2 = L_7(3) \cup \{[0,a] / a \in \{e, 1, 2, \ldots, 29\}, *, 12\}$ be a quasi interval biloop which is a S-strong quasi interval 2-Sylow biloop.

We have the following theorem which guarantees the existence of a class of quasi interval Smarandache strong 2-Sylow biloops.

**THEOREM 1.4.23**: *Let $L = L_n \cup \{[0,a] / a \in \{e, 1, 2, ..., q\}, *, t, 1 < t < q\}$; n and q primes be a class of quasi interval biloops. Every quasi interval biloop in L is a Smarandache strong quasi interval 2-Sylow biloop.*

The proof is straight forward and hence is left as an exercise to the reader.

*Example 1.4.27*: Let $L = L_1 \cup L_2 = L_5(4) \cup \{[0,a] / a \in L_7(2)\}$ is a quasi interval biloop.
    Clearly L is a S-strongly cyclic quasi interval biloop.

**THEOREM 1.4.24**: *Let $L = L_1 \cup L_2 = L_n \cup \{[0,a] / a \in L_m(t); *; 1 < t < m\}$ ($n > 3$, $m > 3$) n and m are primes be a quasi interval biloop. Then every quasi interval biloop in L is bicyclic so L is a S-strongly cyclic quasi interval biloop.*

The proof is direct and hence left for the reader to prove. Now we will give the theorem which gives even the number of quasi interval bicyclic biloops.

**THEOREM 1.4.25**: *Let $L = L_n \cup \{L_m ([0,a] (s)\}$ ($n > 3$, $m > 3$) be a class of quasi interval biloops. If $n = p_1^{\alpha_1} p_2^{\alpha_2} \cdots p_t^{\alpha_t}$ and $m = q_1^{t_1} q_2^{t_2} \cdots q_p^{t_p}$ ($\alpha_i \geq 1$, $t_j \geq 1$; $1 \leq i \leq t$; $1 \leq j \leq p$) then it contains exactly $F_n \cdot F_m$ quasi interval biloops which are strictly noncommutative and they are*
  1) *S-strongly commutative quasi interval biloops and*
  2) *S-strongly quasi interval cyclic biloops where*
  $$F_n = \prod_{i=1}^{t} (p_i - 3) \, p_i^{\alpha_i - 1} \text{ and } F_m = \prod_{i=1}^{p} (q_i - 3) \, q_i^{t_i - 1}$$

Proof is direct for information related to it refer [5, 9].

Also we can say when the S-commutative quasi interval biloop coincides with the quasi interval biloop.

To this effect we first give an example and then give the related theorem.

**THEOREM 1.4.26**: *Let $L = L_n \cup \{L_m ([0,a] (t)\} \mid 1 < t < m\}$ be a class of quasi interval biloops n and m distinct primes and L is noncommutative. Every quasi interval biloop $P = P_1 \cup P_2 \subseteq L_n \cup \{L_n ([0, a]), (t)\}$ is such that $P' = P'_1 \cup P'_2 = P^s = P^s_1 \cup P^s_2 = P = P_1 \cup P_2$.*

For proof in an analogous way please refer [5, 9, 11].

We have a class of quasi interval biloops which are non associative. We will first illustrate this by an example.

***Example 1.4.28***: Let $L = L_1 \cup L_2 = L_n \cup \{[0,a] \mid a \in \{e, 1, 2, …, m\}, *, t; 1 < t < m\}$ where n and m are two distinct primes. For varying elements s and t between 1 to n and 1 to m respectively we get a class of quasi interval biloops. Every biloop in this class is not a S-associative quasi interval biloops.

We say quasi interval biloops $P = P_1 \cup P_2$ is Smarandache strongly pair wise associative if for all $[0, x] \cup [0, y]$, $[0, a] \cup [0, b]$ in $P_1 \cup P_2$ we must have $([0, x] [0, a] ([0, x]) \cup ([0,y] [0,b]) [0, y] = [0, x] ([0, a] [0,x]) \cup [0, y] ([0, b] [0, y])$. We show we have a class of S-strongly associative quasi interval biloops.

Let $L = L_1 \cup L_2 = L_9 (8) \cup \{ ([0,a] \mid a \in \{e, 1, 2, …, 15\}, *, 14\}$ be a quasi interval biloop which is clearly S-strongly associative quasi interval biloop.

**THEOREM 1.4.27**: *Let $L = L_1 \cup L_2 = L_n \cup \{L_m([0,a] (s)\} \mid 1 < s < m\}$ be a class of quasi interval biloops. Every quasi interval biloop in L is S-strongly associative quasi interval biloop.*

Proof given in [5] can be analogously used for $L = L_1 \cup L_2$.

***Example 1.4.29***: Let $L = L_1 \cup L_2 = L_{11}$ (5) $\cup$ { ([0,a] } a $\in$ {e, 1, 2, …, 25}, * 24} be a quasi interval biloop. Clearly L is S-biloop and has no S-subloops. Further $A(L) = A(L_1 \cup L_2) = A(L_1) \cup A(L_2) = L^A = L_1^A \cup L_2^A = L_1 \cup L_2$.

We have a class of quasi interval biloops which satisfy the following condition.

**THEOREM 1.4.28**: *Let $L = L_1 \cup L_2 = \{L_n\} \cup \{[0, a] \mid a \in \{e, 1, 2, .., m\}, *, 1 < t < m\}$ (m and n distinct odd numbers greater than three) be a class of quasi interval biloops. Every biloop $P = P_1 \cup P_2$ in L is such that*
1) *P is S-quasi interval biloop.*
2) *P has no S-quasi interval subbiloops*
3) *$A(P) = A(P_1) \cup A(P_2)$*
   *$P^A = P_1^A \cup P_2^A = P_1 \cup P_2$*
   *for every P in L.*

The proof is direct using definitions and simple number theoretic techniques.

We leave it as an exercise to the reader to prove that we have a class of quasi interval biloops for which the first normalize is equal to second normalilzer under the condition $(m^2 - m + 1, t) = (2m - 1, t)$ where $L_n(m) \in L_n$ and t/n and $(2p - p + 1, s) = (2p - 1, t)$ an interval loop $L_q$ in which s/q and $L_q$ is built using p where $(p, q) = (p, q - 1) = 1$.
Further when both the interval loop and loop are built using $Z_q$ and $Z_p$, where p and q are primes we see the quasi interval biloop $L = L_p$ (t) $\cup$ {$L_q$ ([0, a]) s} is such that SN (L) = SN ($L_p$ (t)) $\cup$ SN ($L_q$ ([0, a]) (s)) = {e} $\cup$ {e} where $1 < t < p$ and $1 < s < q$.
Further for all the quasi interval biloops given above we see the S-Moufang bicentre is {e} $\cup$ {e} or total of L.

Also in this case when p and q are two distinct primes we see the quasi interval biloops we have NZ (L) = NZ ($L_p$ (t)) $\cup$ NZ ($L_q$ [0, a] (s))) = Z (L) = Z ($L_p$ (t)) $\cup$ Z ($L_q$ [0, a] (s)) = {e} $\cup$ {e}.

We can also define the notion of direct product to obtain more classes of quasi interval biloops.

**DEFINITION 1.4.3**: *Let $L = L_1 \cup G$ where $L_1$ is an interval loop and $G$ is an interval group, then $L$ is a interval loop - group. The operations on $L_1$ and $G$ are carried on $L$ component wise.*

We will illustrate this situation by some examples.

*Example 1.4.30*: Let $L = L_1 \cup G_1 = \{[0, a] / a \in \{e, 1, 2, \ldots, 7\}, 3, *\} \cup \{[0, b] / b \in Z_{14}, +\}$ be a interval loop group of finite order.
Order of L is $8.14 = 112$.

*Example 1.4.31*: Let $L = G_1 \cup L_1 = \{Z_{19} \setminus \{0\}, \times\} \cup \{[0, a] / a \in \{e, 1, 2, \ldots, 25\}, *, 24\}$ be an interval group-loop of finite order.

*Example 1.4.32*: Let $L = G_1 \cup L_1 = \{Z_{19} \setminus \{0\}\} \cup \{[0, a] / a \in \{e, 1, 2, \ldots, 23\}, *, 12\}$ be an interval group-loop of order $18 \times 24$. If one of them is alone an interval structure we define then as quasi interval loop-group. Clearly L is a commutative quasi interval group-loop.

We can define quasi interval subgroup - subloop as in case of other algebraic structures.

*Example 1.4.33*: Let $G = G_1 \cup L_1 = \{Z_{25}, +\} \cup \{[0, a] / a \in \{e, 1, 2, \ldots, 15\}, *, 8\}$ be an quasi interval group - loop.

*Example 1.4.34*: Let $V = V_1 \cup V_2 = \{Z_{16}, +\} \cup \{[0, a] / a \in \{e, 1, 2, \ldots, 15\}, 8, *\}$ be an quasi interval group - loop. Consider $P = P_1 \cup P_2 = \{\{0, 4, 8, 12\}, +\} \cup \{[0, a] / a \in \{e, 1, 6, 11\}, *, 8\} \subseteq V_1 \cup V_2$ given by the following tables.

| + | 0 | 4 | 8 | 12 |
|---|---|---|---|----|
| 0 | 0 | 4 | 8 | 12 |
| 4 | 4 | 8 | 12 | 0 |
| 8 | 8 | 12 | 0 | 4 |
| 12 | 12 | 0 | 4 | 8 |

and

| * | [0, e] | [0, 1] | [0, 6] | [0, 11] |
|---|---|---|---|---|
| [0, e] | [0, e] | [0, 1] | [0, 6] | [0, 11] |
| [0, 1] | [0, 1] | [0, e] | [0, 11] | [0, 6] |
| [0, 6] | [0, 6] | [0, 11] | [0, e] | [0, 1] |
| [0, 11] | [0, 11] | [0, 6] | [0, 1] | [0, e] |

P is a quasi interval subgroup subloop of finite order. Clearly o (P) / o (L).

***Example 1.4.35***: Let $L = L_1 \cup L_2 = \{Z_{23}, +\} \cup \{[0, a] / a \in \{e, 1, 2, …, 47\}, *, 9\}$ be a quasi interval group-loop. L has no proper interval group-loop.

***Example 1.4.36***: Let $L = L_1 \cup G_2 = L_9 (8) \cup \{[0, a] / a \in \{Z_{12}\}, +\}$ be a quasi interval loop-group.

Clearly L has substructures.

***Example 1.4.37***: Let $L = L_1 \cup G_1 = \{[0, a] / a \in \{e, 1, 2, …, 17\}, *, 10\} \cup \{[0, b] / b \in Z_{23} \setminus \{0\}, \times\}$ be an interval loop-group of order $18 \times 22$. L has interval bigroup as and has no interval loop-group. $A = A_1 \cup A_2 = \{[0, e], [0, 7], *, 10\} \cup \{[0, 1], [0, 22]\} \subseteq L_1 \cup G_1$ is an interval bigroup of L.

***Example 1.4.38***: Let $G = L_1 \cup G_2 = \{[0, a] / a \in \{Z_{43}\}, +\} \cup \{[0, a] / a \in \{e, 1, 2, …, 19\}, *, 10\}$ be an interval group - loop. Clearly G has no interval substructures.

Infact we have a class of interval group - loop which has no substructures.

**THEOREM 1.4.29**: *Let $L = L_1 \cup L_2 = \{[0, a] / a \in Z_p, +; p$ a prime$\} \cup \{[0, b] / b \in Z_q, *, t ; 1 < t < q, q$ a prime$\}$ be an interval group-loop. L has no interval substructures.*

Proof is direct and hence is left as an exercise to reader.

***Example 1.4.39***: Let $L = L_1 \cup G_1 = L_{15}(8) \cup \{[0, a] / a \in Z_{40}, +\}$ be a quasi interval loop - group of finite order. Clearly $P = P_1 \cup \{[0, a] / a \in \{0, 4, 8, 12, 16, 20, 24, 28, 32, 36\} \subseteq Z_{40}, +\} \subseteq L_1 \cup G_1$ where $P_1 = \{\{e, 1, 6, 11\}$ is a subloop of $L_1$ so $P$ is a quasi interval subloop-subgroup of finite order.

***Example 1.4.40***: Let $L = L_1 \cup L_2 = \{L_{19}(8)\} \cup \{[0, a] / a \in Z_{49}, +\}$ be a quasi interval loop- group of order $20 \times 49$. Clearly $L$ has proper quasi interval bigroup.

The notion of Smarandache cannot arise as we are involving groups in these quasi structures.

Now we can define interval loop - semigroup and quasi interval loop - semigroup.

**DEFINITION 1.4.4**: *Let $L = S_1 \cup L_1$ where $S_1$ is an interval semigroup and $L_1$ is an interval loop; operations carried out on $L$ using operations of $S_1$ and $L_1$. We define $L$ to be an interval semigroup - loop.*

We will illustrate this by some examples so that one can easily follow how the operations on $L$ are carried out.

***Example 1.4.41***: Let $L = S_1 \cup L_1 = \{[0, a] / a \in Z_{24}, \times\} \cup \{[0, b] / b \in \{e, 1, 2, \ldots, 11\}, 8\}$ be an interval semigroup-loop of order $24 \times 12$.

Let $x = [0, 5] \cup [0, 7]$ and $y = [0, 12] \cup [0, 5] \in L$.
$$\begin{aligned}x.y &= ([0, 5] \cup [0, 7]) \,.\, )[0, 12] \cup [0, 5]) \\ &= [0, 5] \times [0, 12] \cup [0, 7] * [0, 5] \\ &= [0, 5 \times 12 \,(\text{mod } 24)] \cup [0, 40\text{-}50 \,(\text{mod } 11)] \\ &= [0, 12] \cup [0, 7+10 \,(\text{mod } 11)] \\ &= [0, 12] \cup [0, 6] \in L.\end{aligned}$$

Now it is easily verified that $L$ is an interval semigroup-loop.

***Example 1.4.42***: Let $L = L_1 \cup S_1 = \{[0, a] / a \in \{e, 1, 2, …, 15, *, 8\} \cup \{[0, b] / b \in Z_{28}, \times\}$ be an interval semigroup-loop of finite order. Clearly o (L) = $16 \times 28$. L is commutative.

***Example 1.4.43***: Let $L = S_1 \cup L_2 = \{S(X) / X = ([0, a_1], [0, a_2], [0, a_3], [0, a_4])\} \cup \{[0, b] \mid b \in \{e, 1, 2, …, 19\}, *, 2\}$ be an interval semigroup-loop. L is of finite order and is non commutative.

***Example 1.4.44***: Let $L = L_1 \cup S_1 = \{[0, a] / a \in \{e, 1, 2, …, 23\}, 4, *\} \cup \{S (\langle X \rangle) / \langle X \rangle = \langle \{[0, a_1], [0, a_2], [0, a_3], [0, a_4], [0, a_5]\} \rangle\}$ is an interval loop semigroup which is non commutative. We say $L = S_1 \cup L_1$ is a quasi interval semigroup-loop if only one of $S_1$ or $L_1$ is an interval structure and the other is a usual structure.

***Example 1.4.45***: Let $L = L_1 \cup S_1 = \{[0, a] / a \in Z_{40}, \times\} \cup L_{11}(3)$ be a quasi interval semigroup-loop of finite order. L is a non commutative structure.

***Example 1.4.46***: Let $L = L_1 \cup S_1 = \{[0, a] \mid a \in Z_{25}, \times\} \cup L_{19}(10)$ be a quasi interval semigroup-loop of finite order. Clearly L is a commutative structure.

We have a class of commutative quasi interval semigroup-group which is evident from the following theorems, the proof of which are direct.

**THEOREM 1.4.30**: *Let $L = S_1 \cup L_1 = \{[0, a] / a \in Z_m, \times\} \cup$*

$$\left\{L_n\left(\frac{n+1}{2}\right), *, \frac{n+1}{2}\right\}$$

*be a quasi interval semigroup-loop. L is commutative.*

**THEOREM 1.4.31**: *Let $L = L_1 \cup S_1 = \{[0, a] / a \in \{e, 1, 2, …, p\}, *, \frac{p+1}{2}\} \cup \{Z_n, \times\}$ be a quasi interval loop - semigroup of finite order. L is a commutative structure.*

Now we can also have a class of non commutative quasi interval loop semigroups which is evident from the following theorems the proof of which is direct.

**THEOREM 1.4.32**: *Let $S = S(n) \cup \{[0, a] / a \in \{e, 1, 2, ..., n\}, *, t \neq \frac{n+1}{2}\}$ be a quasi interval semigroup - loop. S is non commutative.*

**THEOREM 1.4.33**: *Let $S = S_1 \cup L_1 = \{S(X) / X = ([0, a_1], ..., [0, a_n])\} \cup \{L_n(m)\}$ ($m \neq \frac{n+1}{2}$) be a quasi interval semigroup - loop. S is non commutative.*

Now we can define substructures which is a matter of routine. We will however give examples of them.

*Example 1.4.47*: Let $G = S_1 \cup L_1$ where $S_1 = \{S(\langle X \rangle) / X = ([0, a_1], [0, a_2], ..., [0, a_8])\}$ be the special interval symmetric semigroup and $L_1 = \{[0, a] / a \in \{e, 1, 2, ..., 35\}, *, 9\}$ be the interval loop. G is a interval semigroup-loop of finite order. Take $P = P_1 \cup P_2 = \{S(X) \subseteq S(\langle X \rangle)\} \cup \{[0, a] / a \in \{3e, 1, 8, 15, 22, 29\}, *, 9\} \subseteq S_1 \cup L_1$, P is an interval subsemigroup - subloop of L.

*Example 1.4.48*: Let $G = \{Z_{50}, \times\} \cup \{[0, a] / a \in \{e, 1, 2, ..., 15\}, 8, *\}$ be a quasi interval semigroup - loop. Choose $H = H_1 \cup H_2 = \{\{0, 10, 20, 30, 40\} \subseteq Z_{50}, \times\} \cup \{[0, a] / a \in \{e, 1, 6, 11\} \subseteq \{e, 1, 2, ..., 15\}, 8, *\} \subseteq G$, H is a quasi interval subsemigroup-subloop of G.

*Example 1.4.49*: Let $G = G_1 \cup S_1 = L_{21}(5) \cup \{[0, a] / a \in Z_{30}, \times\}$ be a quasi interval semigroup-loop. $H = H_1 \cup H_2 = \{a \in \{e, 1, 8, 15\}, *, 5\} \cup \{[0, a] / a \in \{0, 2, 4, 6, 8, ..., 26, 28\}, \subseteq Z_{30}, \times\} \subseteq G_1 \cup S_1$, is a quasi interval subsemigroup - subloop of L.

We can define Smarandache interval smeigroup-loop and Smarandache quasi interval semigroup - loop.

We will leave the task of defining the Smarandache structure to the reader as it is a matter of routine, however we give examples of them.

***Example 1.4.50***: Let $L = S_1 \cup L_1 = \{Z_{18}, \times\} \cup \{[0, a] / a \in \{e, 1, 2, \ldots, 21\}, *, 11\}$ be a quasi interval semigroup - loop. Clearly L is a S-quasi interval semigroup-loop for L contains $A = A_1 \cup A_2 = \{\{1, 17\} \subseteq Z_{18}, \times\} \cup \{[0, e], [0, 19], *, 10\} \subseteq L = S_1 \cup L_1$ which is a quasi interval bigroup.

***Example 1.4.51***: Let $G = S_1 \cup L_1 = \{S (\langle X \rangle) / X = ([0, a_1], [0, a_2], [0, a_3])\} \cup \{L_{19} (7)\}$ be a quasi interval semigroup loop. Clearly G is a S-quasi interval semigroup loop for $H = H_1 \cup H_2 = S_X \cup \{\{e, 11\} \subseteq L_{19} (7)\}$ is a quasi interval bigroup of G.

***Example 1.4.52***: Let $G = S_1 \cup L_1 = \{[0, a] / a \in Z_{48}, \times\} \cup L_{15} (8)$ is a quasi interval semigroup-loop. Consider $H = H_1 \cup H_2 = \{[0, a] / a \in \{0, 4, 8, 12, \ldots, 44\}, \times\} \cup \{\{e, 1, 6, 11\} / 8, *\} \subseteq S_1 \cup L_1$ is a quasi interval subsemigroup-subloop of G. Also $P = P_1 \cup P_2 = \{[0, 1], [0, 4], \times\} \cup \{\{e, 12\}, *, 8\} \subseteq S \cup L_1$ is a quasi interval bigroup of G.

***Example 1.4.53***: Let $G = S_1 \cup L_1 = \{[0, a] / a \in Z_{20}, \times\} \cup \{[0, b] / b \in \{e, 1, 2, \ldots, 7\}, *, 4\}$ be a interval semigroup-loop. Clearly G is a S-interval semigroup-loop but G has no S-interval subsemigroup - subloop.

We have a class of interval semigroup loop which has no S-interval subsemigroup - subloop.

**THEOREM 1.4.34**: *Let $S = G_1 \cup L_1 = \{L_n / n \text{ a prime}\} \cup \{Z_p, +\}$ be a class of interval loop - semigroup for varying primes p. Clearly S has no S-quasi interval subloop - subsemigroup.*

The proof is direct and hence it left as an exercise for the reader.

**THEOREM 1.4.35**: *Let $L = S \cup L_1 = \{[0, a] / a \in Z_p$, p a prime p varying over the set of all primes, under $\times\} \cup \{[0, a] / a \in \{e, 1, 2, …, n\}$, n a prime, t, $1 < t < n$, \*, for varying t between $(1, n)\}$ be a class of interval semigroup - loop L has no S-interval subsemigroup - subloop.*

The proof is an easy consequence of the definition and hence is left as an exercise for the reader.

This is a quasi associative interval algebraic structure.

Now we proceed onto define interval loop-groupoids.

**DEFINITION 1.4.5**: *Let $G = G_1 \cup L_1$ where $G_1$ is an interval groupoid and $L_1$ is an interval loop and G inherits the operation from $G_1$ and $L_1$ denote the operation by '.', (G, .) is defined as the interval groupoid-loop.*

We will illustrate this by some examples.

*Example 1.4.54*: Let $G = G_1 \cup L_1 = \{[0, a] / a \in Z_8, \{3, 7\}, \otimes\} \cup \{[0, a] / a \in \{e, 1, 2, …, 7\}, *, 3\}$ be an interval groupoid-loop.

Suppose $x = [0, a] \cup [0, b]$ and $y = [0, c] \cup [0, d]$ be in G.
$$\begin{aligned}x.y &= ([0, a] \cup [0, b]) . ([0, c] \cup [0, d]) \\ &= [0, a] \otimes [0, c] \cup [0, b] * [0, d] \\ &= [0, a \otimes c] \cup [0, b * d] \\ &= [0, 3a + 7c \pmod 8] \cup [0, 3d - 2b \pmod 7] \text{ is in G.}\end{aligned}$$

*Example 1.4.55*: Let $G = G_1 \cup L_1 = \{[0, a] / a \in Z_{15}, *, (8, 4)\} \cup \{[0, d] / d \in \{e, 1, 2, …, 19\}, *, 4\}$ be a interval groupoid-loop.

*Example 1.4.56*: Let $G = G_1 \cup L = \{[0, a] / a \in Z^+, *, (3, 2)\} \cup \{[0, b] / b \in \{e, 1, 2, …, 43\}, 24, *\}$ be an interval groupoid - loop. Clearly G is of infinite order.

All other examples of interval groupoid - loop given are only of finite order.

We can define substructures in them, this task is left to the reader. We give only examples.

**Example 1.4.57**: Let $G = G_1 \cup L_1 = \{[0, a] / a \in Z_6, *, (4, 5)\} \cup \{[0, b] / b \in \{e, 1, 2, …, 35\}, *, 9\}$ be interval groupoid-loop of order $6 \times 36$.

Take $H = H_1 \cup H_2 = \{[0, 2], [0, 4], 0 / 0, 4, 2 \in Z_6, *, (4, 5)\} \cup \{[0, b] / b \in \{e, 1, 8, 15, 22, 29\} \subseteq \{e, 1, 2, …, 35\}, *, 9\} \subseteq G_1 \cup L_1$, H is an interval subgroupoid-subloop of G.

The order of H is 3.6 clearly o(H) / o(G).

**Example 1.4.58**: Let $G = H_1 \cup G_1 = \{[0, a] / a \in Z_8, *, (2, 6)\} \cup \{[0, b] / b \in \{e, 1, 2, …, 15\}, *, 8\}$ be an interval groupoid-loop of order 8.16. Consider $A = A_1 \cup A_2 = \{[0, a] / a \in \{0, 2, 4, 6\} \subseteq Z_8, *, (2, 6)\} \cup \{[0, b] / b \in \{e, 1, 6, 11\} \subseteq \{e, 1, 2, 3, …, 14, 15\}, *, 8\} \subseteq H_1 \cup G_1$. A is an interval subgroupoid - subloop of G and o (A) = 4.4 and o (A) / o (G), that is 4.4 / 8.16

Now we can define S-interval groupoid - loop. Further this algebraic structure is a non associative structure.

**Example 1.4.59**: Let $G = H_1 \cup H_2 = \{[0, a] / a \in Z_8, *, (2, 6)\} \cup \{[0, b] / b \in \{e, 1, 2, …, 21\}, *, 11\}$ be an interval groupoid - loop of order 8.22 = 176. Let $A = A_1 \cup A_2 = \{[0, a] / a \in \{0, 2, 3, 4, 6\} \subseteq Z_8, *, (2, 6)\} \cup \{[0, b] / b \in \{e, 1, 8, 15\}, *, 11\} \subseteq H_1 \cup H_2$ ; A is an interval subgroupoid - subloop of G. Clearly o (H) $\not|$ o (G) for o (H) = 5.4 and 5.4 $\not|$ 8.22. We can define Smarandache structures in these algebraic structure.

An interval groupoid - loop G is said to be a Smarandache interval groupoid - loop if G has a proper subset $H = H_1 \cup H_2$ where $H_1$ is an interval semigroup and $H_2$ is an interval subgroup.

We will illustrate this structure by an example.

***Example 1.4.60***: Let $G = G_1 \cup G_2 = \{[0, a] / a \in Z_5, *, (3, 3)\} \cup \{[0, a] / a \in \{e, 1, 2, …, 19\}, *, 8\}$ be an interval groupoid - loop. Consider $H = H_1 \cup H_2 = \{[0, 1] / 1 \in Z_5, *, (3, 3)\} \cup \{[0, e] [0, 9] / e, 9 \in \{e, 1, 2, …, 19\}, *, 8\} \subseteq G_1 \cup G_2$, is an interval semigroup-group, contained in G. Hence G is a S-interval groupoid - loop.

We have a class of such S-interval groupoid - loops.

***Example 1.4.61***: Let $H = H_1 \cup H_2 = \{[0, a] / a \in Z_6, *, (4, 5)\} \cup \{[0, b] / b \in \{e, 1, 2, …, 19\}, *, 3\}$ be a S-interval groupoid loop.

**THEOREM 1.4.36**: *Let $S = G_1 \cup G_2 = \{[0, a] / a \in Z_{2p}, p$ a prime, $*, (1, 2)\} \cup \{[0, b] / b \in \{e, 1, 2, ..., n\}; n > 3, 1 < m < n, (m, n) = (1-m, n) = 1, *\}$ be a class of interval groupoid - loop for varying p. Clearly every interval groupoid-loop in S is a Smarandache groupoid-loop.*

Proof is direct and is left as an exercise to the reader.

We call an interval groupoid - loop to be a S-Bol interval-groupoid loop if both the interval-groupoid and interval loop are S-Bol. Similarly S-Moufang, S-alternative and S-idempotent interval gropoid-loop.

Interested reader is expected to supply examples of these structures. However since both the interval structures are non - associative we can define interval S-quasi Bol groupoid - loop, or S-quasi Bol loop - groupoid, interval S-quasi Moufang loop - groupoid or interval S-quasi P-groupoid and so on. We will also give the related theorems. We will further give examples of them.

***Example 1.4.62***: Let $G = L_1 \cup G_1$ where $L_1 = \{[0, a] / a \in \{e, 1, 2, …, 19\}, *, 12\}$ be a S-strong interval cyclic loop and $G_1 = \{[0, b] / b \in Z_{28}, *, (7, 3)\}$ be an interval groupoid. Clearly $G_1$ is a not a S-strongly interval cyclic groupoid. Hence G the interval loop - groupoid is only a interval quasi S-strongly cyclic loop-groupoid.

***Example 1.4.63***: Let $G = G_1 \cup L_1 = \{[0, a] / a \in Z_{20}, *, (7, 8)\}$ $\cup \{[0, b] / b \in \{e, 1, 2, \ldots, 23\}, *, 8\}$ be an interval groupoid - loop. G is only a interval quasi S-strongly commutative groupoid loop as only $L_1$ the interval loop is S-strongly commutative interval loop where as $G_1$ is not a S-strongly interval commutative groupoid.

Several related results and properties can be derived with appropriate modifications. This can be treated as a matter of routine and carried out by the interested reader.

**Chapter Two**

# n-Interval Algebraic Structures with Single Binary Operation

In this chapter we introduce several new types of n-interval algebraic structures (n>2) with single binary operation. These structures are so mixed and they are utilized in places of appropriate applications. This chapter has five sections. Sections one deals with n-interval semigroups and analysis their properties. Section two introduces the notions of n-interval groupoids (n>3) and generalizes them. The notion of n-interval group and quasi n-interval structures using groups, semigroups groupoids are introduced for the first time in section three. The four section defines the notion of n-interval loops and describes a few properties related with them. The final section introduces the notion of mixed n-interval algebraic structures.

## 2.1 n-Interval Semigroups

In this section we introduce n-interval semigroups (n > 2) and describe a few properties related with them.

**DEFINITION 2.1.1**: *Let $S = S_1 \cup S_2 \cup \ldots \cup S_n$, $(n > 2)$ where each $S_i$ is an interval semigroup, $S_i \not\subseteq S_j$ for any $i \neq j$; $1 \leq i, j \leq n$. The operation on S is the component wise operation on each $S_i$ carried out in a systematic way and denoted by '.'; $1 \leq i \leq n$.*

*Thus any element $s \in S$ is represented as $s = s_1 \cup s_2 \cup \ldots \cup s_n$ where $s_i \in S_i$; $1 \leq i \leq n$ and $(S, .)$ is defined as the n-interval semigroup.*

*If $n = 2$ we call it as the biinterval semigroup or interval bisemigroup.*

If the order of every $S_i$ is finite S will be finite $1 \leq i \leq n$.
Even if one of the $S_i$'s is of infinite order, S will be of infinite order; $1 \leq i \leq n$.

We will first illustrate this situation by examples.

***Example 2.1.1***: Let $S = S_1 \cup S_2 \cup S_3 \cup S_4 \cup S_5 = \{[0, a] / a \in Z_{40}, \times\} \cup \{[0, b] / b \in Z^+ \cup \{0\}, +\} \cup \{[0, c] / c \in Z_{25}, +\} \cup \{[0, d] / d \in Z_{17}, \times\} \cup \{[0, g] / g \in Z_{12}, \times\}$ be a 5-interval semigroup.

Take
$$x = [0, 2] \cup [0, 4] \cup [0, 3] \cup [0, 7] \cup [0, 8]$$
and
$$y = [0, 1] \cup [0, 5] \cup [0, 20] \cup [0, 4] \cup [0, 4]$$
in S.

$$\begin{aligned}
x.y &= ([0, 2] \cup [0, 4] \cup [0, 3] \cup [0, 7] \cup [0, 8]) ([0, 1] \cup [0, 5] \cup [0, 20] \cup [0, 4] \cup [0, 4]) \\
&= ([0, 2].[0, 1] \cup [0, 4] [0, 5] \cup [0, 3] . [0, 20] \cup [0, 7] [0, 4] \cup [0, 8] [0, 4]) \\
&= [0, 2] \cup [0, 20] \cup [0, 10] \cup [0, 11] \cup [0, 8]
\end{aligned}$$
is in S.

Thus $(S, .)$ is a 5-interval semigroup of infinite order. Clearly S is a commutative S-interval semigroup as each $S_i$ is a commutative semigroup, $1 \leq i \leq 5$.

***Example 2.1.2***: Let $V = V_1 \cup V_2 \cup V_3 \cup V_4 =$

$$\left\{ \begin{bmatrix} [0,a] & [0,b] \\ [0,c] & [0,d] \end{bmatrix} \bigg/ a, b, c, d \in Z_{20}, \times \right\} \cup$$

$$\{([0, a], [0, b], [0, c], [0, d]) / a, b, c, d, \in Z_{15}\} \cup$$

$$\left\{ \begin{bmatrix} [0,a] \\ [0,b] \\ [0,c] \\ [0,d] \end{bmatrix} \bigg| a, b, c, d \in Z_{12} \right\} \cup \{[0, a] / a \in Z_{47}\}$$

be a 4-interval semigroup. Clearly V is of finite order but V is non commutative.

***Example 2.1.3***: Let $S = S_1 \cup S_2 \cup S_3 \cup S_4 \cup S_5 \cup S_6 =$

$\{([0, a], [0, b]) / a, b \in Z_{10}, \times\} \cup$

$$\left\{ \begin{bmatrix} [0,a] \\ [0,b] \\ [0,c] \\ [0,d] \end{bmatrix} \bigg| a, b, c, d \in Z_{15}, + \right\}$$

$\cup \{S(X) / X = ([0, a], [0, b], [0, c])\} \cup \{$All $3 \times 3$ interval matrices with intervals of the form $[0, a]$ where $a \in Z_{12}\} \cup \{[0, a] / a \in Z_{19}\} \cup \{$All $2 \times 4$ interval matrices with intervals of the form $[0, a]$ where $a \in Z_8, +\}$ be a 6-interval semigroup.
 Clearly S is of finite order and S is non commutative.

***Example 2.1.4***: Let $S = S_1 \cup S_2 \cup S_3 \cup S_4 = \{\sum_{i=0}^{5}[0,a]x^i \mid a \in Z_{12}, +\} \cup \{([0, b], [0, a], [0, c]) / a, b, c \in Z_{14}, +\} \cup$

$$\left\{ \begin{bmatrix} [0,a] \\ [0,b] \\ [0,c] \end{bmatrix} \bigg/ a, b, c \in Z_{15}, + \right\} \cup$$

$\{3 \times 5$ interval matrices with intervals of the form $[0, a]$ where a $\in Z_{20}, +\}$ be a 4-interval semigroup. S is of finite order and is commutative.

Now having seen examples of n-interval semigroups we give examples of n-interval subsemigroups and ideals in n-interval semigroups.

The task of giving definition is a matter of routine and hence is left as an exercise to the reader.

***Example 2.1.5***: Let $S = S_1 \cup S_2 \cup S_3 \cup S_4 = \{[0, a] / a \in Z^+ \cup \{0\}\} \cup \{([0, a], [0, b], [0, c]) / a, b, c \in Z_{20}\} \cup$

$$\left\{ \begin{bmatrix} [0,a] \\ [0,b] \\ [0,c] \end{bmatrix} \bigg/ a, b, c \in Z^+ \cup \{0\} \right\} \cup \left\{ \sum_{i=0}^{5} [0,a]x^i \bigg/ a \in Z_{40}, + \right\}$$

be a 4-interval semigroup.

Consider $A = A_1 \cup A_2 \cup A_3 \cup A_4 = \{[0, a] / a \in 3Z^+ \cup \{0\}\} \cup \{([0, a], 0, [0, b]) / a, b \in Z_{20}\} \cup$

$$\left\{ \begin{bmatrix} [0,a] \\ [0,b] \\ [0,c] \end{bmatrix} \bigg/ a, b, c \in 5Z^+ \cup \{0\} \right\} \cup \left\{ \sum_{i=0}^{8} [0,a]x^i \; \bigg| \right.$$

$a \in \{2, 0, 4, 8, \ldots, 36, 38\} \subseteq Z_{40}, +\} \subseteq S_1 \cup S_2 \cup S_3 \cup S_4$; A is a 4-interval subsemigroup of S. It is easily verified, A is not an ideal of S.

***Example 2.1.6***: Let $S = S_1 \cup S_2 \cup S_3 \cup S_4 \cup S_5 = \{[0, a] / a \in Z^+ \cup \{0\}\} \cup \{$All $3 \times 3$ interval matrices with intervals of the form $[0, a] / a \in Z^+ \cup \{0\}\} \cup \{$all $1 \times 5$ row interval matrices with intervals of the form $[0, a] / a \in Z^+ \cup \{0\}, \times\} \cup$

$$\left\{ \begin{bmatrix} [0,a] & [0,b] \\ [0,c] & [0,d] \end{bmatrix} \mid a, b, c, d \in Z^+ \cup \{0\} \right\} \cup \left\{ \sum_{i=0}^{8} [0,a]x^i \mid a \in Z^+ \cup \right.$$

$\{0\}, \times$ with $x^9 = 1\}$ be an 5-interval semigroup.

Consider $P = P_1 \cup P_2 \cup P_3 \cup P_4 \cup P_5 = \{[0, a] / a \in 5Z^+ \cup \{0\}\} \cup \{$All $3 \times 3$ interval matrices with intervals of the form $[0, a]$ where $a \in 7Z^+ \cup \{0\}, \times\} \cup \{$All $1 \times 5$ row interval matrices with intervals of the form $[0, a]$ where $a \in 19Z^+ \cup \{0\}, \times\} \cup$

$$\left\{ \begin{bmatrix} [0,a] & [0,b] \\ [0,c] & [0,d] \end{bmatrix} \mid a, b, c, d \in 3Z^+ \cup \{0\} \right\} \cup$$

$$\{ \sum_{i=0}^{8} [0,a]x^i / x^9 = 1, a \in 5Z^+ \cup \{0\}, \times\}$$

$\subseteq S_1 \cup S_2 \cup S_3 \cup S_4 \cup S_5 = S$ be a 5-interval subsemigroup of S. It is easily verified P is also a 5-interval ideal of S.

However it is interesting note the following result.

**THEOREM 2.1.1**: *Let $S = S_1 \cup S_2 \cup ... \cup S_n$ be a n-interval semigroup. Let $P = P_1 \cup P_2 \cup ... \cup P_n \subseteq S_1 \cup S_2 \cup ... \cup S_n$ be a n-interval subsemigroup. P in general is not a n-interval ideal of S. Further every n-interval ideal of a n-interval semigroup is a n-interval subsemigroup of S.*

The proof is direct and hence is left as an exercise for the reader.

Now having seen the notion of n-ideals in an n-interval semigroups and n-interval subsemigroups we now proceed onto define S-n-interval semigroups.

**DEFINITION 2.1.2**: *Let $S = S_1 \cup S_2 \cup ... \cup S_n$ be a n-interval semigroup. If each $S_i$ is a Smarandache semigroup then we define S to be a Smarandache n-interval semigroup (S-n-interval semigroup) $1 \leq i \leq n$.*

We will illustrate this situation by some examples.

***Example 2.1.7***: Let $S = S_1 \cup S_2 \cup S_3 \cup S_4$ be a 4-interval semigroup, where $S_1 = \{[0, a] / a \in Z_{12}, \times\}$, $S_2 = \{([0, a], [0, b] / a, b \in Z_9\}$, $S_3 = \{([0, a], [0, b], [0, c], [0, d]) / a \in Z_{11}\}$ and

$$S_4 = \left\{ \begin{bmatrix} [0,a] \\ [0,b] \\ [0,c] \\ [0,d] \end{bmatrix} / a, b, c, d \in Z_{18}, \text{ under '+'} \right\}.$$

Now choose $A = A_1 \cup A_2 \cup A_3 \cup A_4$

$= \{[0, 1], [0, 11] / 1, 11 \in Z_{12}\} \cup \{([0, 1] [0, 1]), ([0, 8], [0, 8]), ([0, 8], [0, 1]) ([0, 1], [0, 8])\} \cup \{([0, 1], [0, 1], [0, 1], [0, 1]), ([0, 10], [0, 10], [0, 10], [0, 10]) \cup$

$$\left\{ \begin{bmatrix} [0,a] \\ [0,b] \\ [0,c] \\ [0,d] \end{bmatrix} \mid a, b, c, d \in \{0, 3, 6, 9, 12, 15\} \subseteq Z_{18}, + \right\}$$

$\subseteq S_1 \cup S_2 \cup S_3 \cup S_4$. It is easily verified A is a 4-interval group in S.

Hence each $S_i$ is a S-interval semigroup, $1 \leq i \leq 4$. Thus S is a S-4-interval semigroup.

It is important and interesting to note that in general all n-interval semigroups need not be S-n-interval semigroups.

We will illustrate this situation by an example.

***Example 2.1.8***: Let $S = S_1 \cup S_2 \cup S_3 = \{[0, a] / a \in Z^+ \cup \{0\}, \times\} \cup \{\begin{bmatrix} [0,a] \\ [0,b] \end{bmatrix} / a, b \in Z^+ \cup \{0\}, +\} \cup \{([0, a] [0, b]) / a, b \in Q^+ \cup \{0\}\}$ be a 3-interval semigroup. It is easily verified S is not a S-3- interval semigroup.

We can define as a matter of routine the notion of S-interval subsemigroup. We will only give an example of it.

***Example 2.1.9***: Let $S = S_1 \cup S_2 \cup S_3 \cup S_4 = \{S(X) \mid X = ([0, a_1], [0, a_2], [0, a_3])\} \cup \{S(\langle X \rangle) / X = ([0, x_1], [0, x_2], [0, x_3])\} \cup \{[0, a] / a \in Z_{20}, +\} \cup \{([0, a], [0, b], [0, c], [0, d]) / a, b, c, d \in Z_{15}, +\}$ be a 4-interval semigroup. Consider $V = V_1 \cup V_2 \cup V_3 \cup V_4 = \{S_X\} \cup \{S_{\langle X \rangle}\} \cup \{[0, a] / a \in \{2, 0, 4, 6, 8, 10, \ldots, 18\} \subseteq Z_{20}, +\} \cup \{([0, a], [0, b], [0, c], [0, d]) / a, b, c, d \in \{0, 3, 6, 9, 12\} \subseteq Z_{15}, +\} \subseteq S_1 \cup S_2 \cup S_3 \cup S_4$.

It is easily verified V is a 4-interval subsemigroup of S. Consider $A = A_1 \cup A_2 \cup A_3 \cup A_4 = \{A_X\} \cup \{A_{\langle X \rangle}\} \cup \{[0, a] / a \in \{0, 4, 8, 12, 16\} \subseteq Z_{20}, +\}, \{([0, a] [0, a] [0, a]) / a \in \{0, 2, 4, 6, 8, 10, \ldots, 18\} \subseteq Z_{20}, +\} \subseteq V_1 \cup V_2 \cup V_3 \cup V_4 \subseteq S_1 \cup S_2 \cup S_3 \cup S_4$. Clearly each $A_i$ is a interval group in $S_i$; $1 \leq i \leq 4$. Thus A is a 4-interval group. Hence V is a S-4-interval subsemigroup.

***Example 2.1.10***: Let $S = S_1 \cup S_2 \cup S_3 \cup S_4 = \{[0, a] / a \in Z_{11}, \times\} \cup \{[0, b] / b \in Z_{13}, \times\} \cup \{[0, c] / c \in Z_{19}, \times\} \cup \{[0, d] / d \in Z_{23}, \times\}$ be a 4-interval semigroup. Clearly S is a S-4-interval semigroup as $A = A_1 \cup A_2 \cup A_3 \cup A_4 = \{[0, 1], [0, 10], \times\} \cup \{[0, 1], [0, 12], \times\} \cup \{[0, 1], [0, 18], \times\} \cup \{[0, 1], [0, 22], \times\} \subseteq S_1 \cup S_2 \cup S_3 \cup S_4$ is a 4-itnerval group.

***Example 2.1.11***: Let $S = S_1 \cup S_2 \cup S_3 \cup S_4 \cup S_5 = \{[0, a] / a \in Z_6, \times\} \cup \{[0, a] / a \in Z_q, \times\} \cup \{[0, a] / a \in Z_{16}, \times\} \cup \{[0, a] / a \in Z_{25}, \times\} \cup \{[0, a] / a \in Z_{36}, \times\}$ be a 5-interval semigroup.

Consider $P = P_1 \cup P_2 \cup P_3 \cup P_4 \cup P_5 = \{[0, a] / a \in \{0, 2, 4\} \subseteq Z_6, \times\} \cup \{[0, a] / a \in \{0, 3, 6\} \subseteq Z_9, \times\} \cup \{[0, a] / a \in \{0, 4, 8, 12\} \subseteq Z_{16}, \times\} \cup \{[0, a] / a \in \{0, 5, 10, 15, 20\} \subseteq Z_{25}, \times\} \cup \{[0, a] / a \in \{0, 6, 12, 18, 24\} \subseteq Z_{36}, \times\} \subseteq S_1 \cup S_2 \cup S_3 \cup S_4 \cup S_5$. P is a 5-interval subsemigroup of S. But P is not a S-5-interval subsemigroup of S. However S is a S-5-interval semigroup as $A = A_1 \cup A_2 \cup A_3 \cup A_4 \cup A_5 = \{[0, 1], [0, 5], \times\} \cup \{[0, 1], [0, 8], \times\} \cup \{[0, 1], [0, 24], \times\} \cup \{[0, 1], [0, 35], \times\} \subseteq S_1 \cup S_2 \cup S_3 \cup S_4 \cup S5$ is a 5-interval group. Thus S is a S-5-interval semigroup but every 5-interval subsemigroup of S need not be a S-5-interval subsemigroup.

In view of this we have the following theorem the proof of which is direct.

**THEOREM 2.1.2**: *Let $S = S_1 \cup S_2 \cup ... \cup S_n$ be a n-interval semigroup. If S has a S-n-interval subsemigroup then S is a S-n-interval semigroup. Suppose S is a S-n-interval semigroup then every n-interval subsemigroup of S in general is not a S-n-interval subsemigroup.*

Now having seen examples of these situations we leave the task of defining S-n-interval ideal and illustrate them by examples.

Now as in case of usual n-semigroups we can in case of n-interval semigroups also define the notion of n-zero divisors, n-units and n-idempotents and quasi n-zero divisors, quasi n-units and quasi n-idempotents.

**DEFINITION 2.1.3**: *Let $S = S_1 \cup S_2 \cup ... \cup S_n$ be a n-interval semigroup. Suppose for $x = x_1 \cup x_2 \cup ... \cup x_n \in S$ there exists a $y = y_1 \cup y_2 \cup ... \cup y_n \in S$ such that $x.y = x_1 y_1 \cup x_2 y_2 \cup ... \cup x_n y_n = 0 \cup 0 \cup ... \cup 0$ in S then we call x to be a n-interval zero divisor of S.*

We will illustrate this situation by some examples.

*Example 2.1.12*: Let $S = S_1 \cup S_2 \cup S_3 \cup S_4 = \{[0, a] | a \in Z_{12}, \times\}$ $\cup \{[0, b] / b \in Z_{15}, \times\} \cup \{([0, a], [0, b])$ where $a, b \in Z_{24}, \times\} \cup \left\{ \begin{bmatrix} [0,a] \\ [0,b] \end{bmatrix} \mid a, b \in Z_{10}, + \right\}$

be a 4-interval semigroup. Let

$$x = \{[0, 4] \cup [0, 3] \cup ([0, 12], [0, 6])\} \cup \begin{bmatrix} [0,5] \\ [0,7] \end{bmatrix} \} \in S.$$

We have

$$y = [0, 3] \cup [0, 5] \cup ([0, 2]\, [0, 4]) \cup \begin{bmatrix} [0,5] \\ [0,3] \end{bmatrix} \in S$$

is such that

$$\begin{aligned}
x.y \quad &= \quad ([0, 4] \cup [0, 3] \cup ([0, 12], [0, 6]) \cup (\begin{bmatrix}[0,5]\\ [0,7]\end{bmatrix})\\
&\qquad ([0, 3] \cup [0, 5] \cup ([0, 2], [0, 4]) \cup \begin{bmatrix}[0,5]\\ [0,3]\end{bmatrix}\\
&= \quad [0, 4] \times [0, 3] \cup [0, 3] \times [0, 5] \cup ([0, 12], [0, 6]) \times\\
&\qquad ([0, 2], [0, 4]) \cup \begin{bmatrix}[0,5]\\ [0,7]\end{bmatrix} + \begin{bmatrix}[0,5]\\ [0,3]\end{bmatrix}\\
&= \quad [0, 12] \cup [0, 15] \cup ([0, 24], [0, 24]) \cup \begin{bmatrix}[0,10]\\ [0,10]\end{bmatrix}\\
&= \quad 0 \cup 0 \cup 0 \cup 0
\end{aligned}$$

is a 4-interval zero divisor in S.

***Example 2.1.13***: Let $S = S_1 \cup S_2 \cup S_3 = \{[0, a] / a \in Z^+ \cup \{0\}\} \cup \{[0, a] / a \in Z_{19} \setminus \{0\}\} \cup \{[0, a] / a \in Q^+ \cup \{0\}\}$ be a 3-interval semigroup. S has no interval zero divisors.

***Example 2.1.14***: Let $G = G_1 \cup G_2 \cup G_3 \cup G_4 = \{([0, a], [0, b], [0, c], [0, d]) / a, b, c, d \in Z^+ \cup \{0\}\} \cup \{[0, a] / a \in Z_{16}, \times\} \cup \{([0, a], [0, b], [0, c], [0, d], [0, e]) / a, b, c, d, e, \in Q^+ \cup \{0\}\} \cup \{[0, a] / a \in Z_{40}, \times\}$ be a 4-interval semigroup. G has non trivial 4 - interval zero divisors. For take $x = ([0, 0], [0, 8], [0, 0], [0, 12]) \cup \{[0, 8]\} \cup ([0, 2], [0, 1/3], [0, 0], [0, 8/9], [0, 0]) \cup [0, 10] \in G$, we have $y = ([0, 2], [0, 0], [0, 12], [0, 0]) \cup \{[0, 4]\} \cup \{([0, 0], [0, 0], [0, 9/7], [0, 0], [0, 11/3])\} \cup [0, 8]$ in G such that

$$\begin{aligned}
x.y \quad &= \quad \{([0, 0], [0, 8], [0, 0], [0, 12]) \cup [0, 8] \cup ([0, 2],\\
&\qquad [0, 1/3], [0, 0], [0, 8/9], [0, 0]) \cup [0, 10]\} .\{[0, 2],\\
&\qquad [0, 0], [0, 12], [0, 0]) \cup [0, 4] \cup ([0, 0], [0, 0],\\
&\qquad [0, 9/7], [0, 0], [0, 11/3]) \cup [0, 8]\}\\
&= \quad ([0, 0], [0, 8], [0, 0], [0, 12]) \times ([0, 2], [0, 0],\\
&\qquad [0, 12], [0, 0]) \cup [0, 8] [0, 4] \cup ([0, 2], [0, 1/3],\\
&\qquad [0, 0], [0, 8/9], [0, 0]) \times ([0, 0], [0, 0], [0, 9/7],\\
&\qquad [0, 0], [0, 11/3]) \cup [0, 10] \times [0, 8]
\end{aligned}$$

$$
\begin{aligned}
=\ &([0, 0], [0, 2], [0, 8], [0, 0], [0, 0], [0, 12], [0, 12], [0, 0]) \\
&\cup\ [0, 32]\ \cup\ ([0, 2], [0, 0], [0, 1/3], [0, 0], [0, 0], \\
&[0, 9/7], [0, 8/9], [0, 0], [0, 0], [0, 11/3])\ \cup\ [0, 80] \\
=\ &([0, 0], [0, 0], [0, 0], [0, 0])\ \cup\ [0, 0]\ \cup\ ([0, 0], [0, 0], \\
&[0, 0], [0, 0], [0, 0])\ \cup\ [0, 0].
\end{aligned}
$$

Thus x is a 4 interval zero divisor in S.

Now one can define S-n-interval zero divisor analogous to S-zero divisors and illustrate by examples.

Now we proceed onto give examples of n-interval idempotents in a n-interval semigroup S.

***Example 2.1.15***: Let $S = S_1 \cup S_2 \cup S_3 \cup S_4 \cup S_5 = \{[0, a] / a \in Z_6, \times\} \cup \{[0, a] / a \in Z_{20}, \times\} \cup \{[0, a] / a \in Z_{18}, \times\} \cup \{[0, a] / a \in Z_{24}, \times\} \cup \{[0, a] / a \in Z_{10}\}$ be a 5-interval semigroup. Consider $x = [0, 3] \cup [0, 5] \cup [0, 9] \cup [0, 9] \cup [0, 5] \in S$.
Clearly
$$
\begin{aligned}
x^2\ =\ &([0, 3] \cup [0, 5] \cup [0, 9] \cup [0, 9] \cup [0, 5])\,([0, 3] \cup \\
&[0, 5] \cup [0, 9] \cup [0, 9] \cup [0, 5]) \\
=\ &[0, 3].[0, 3] \cup [0, 5]\,[0, 5] \cup [0, 9].[0, 9] \cup [0, 9] \\
&[0, 9] \cup [0, 5]\,[0, 5] \\
=\ &[0, 9\ (\mathrm{mod}\ 6)] \cup [0, 25\ (\mathrm{mod}\ 20)] \cup [0, 81\ (\mathrm{mod}\ 18)] \cup \\
&[0, 81\ (\mathrm{mod}\ 24)] \cup [0, 25\ (\mathrm{mod}\ 10)] \\
=\ &[0, 3] \cup [0, 5] \cup [0, 9] \cup [0, 9] \cup [0, 5] \\
=\ &x.
\end{aligned}
$$
Thus x is a 5-interval idempotent in S.

***Example 2.1.16***: Let $S = S_1 \cup S_2 \cup S_3 \cup S_4 = \{[0, a] / a \in Z^+ \cup \{0\}\} \cup \{[0, a] / a \in Q^+ \cup \{0\}\} \cup \{[0, a] / a \in R^+ \cup \{0\}\} \cup \{[0, a] / a \in Z_{43}\}$ be a 4 - interval semigroup. It is easily verified S has only trivial 4 - interval idempotents like $([0, 1] \cup [0, 1] \cup [0, 1] \cup [0, 1])$ or $[0, 0] \cup [0, 0] \cup [0, 0] \cup [0, 0]$ or elements of the form $[0, 1] \cup [0, 0] \cup [0, 1] \cup [0, 1]$ and so on. We call all those idempotents constructed using $[0, 0]$ and $[0, 1]$ as trivial interval idempotents.

*Example 2.1.17*: Let $S = S_1 \cup S_2 \cup S_3 \cup S_4 = \{[0, a] / a \in Z^+\} \cup \{[0, a] / a \in Q^+\} \cup \{[0, a] / a \in R^+\} \cup \{[0, a] / a \in Z_{43} \setminus \{0\}\}$ be a 4-interval semigroup. S has only non trivial 4-interval idempotents.

Let $S = S_1 \cup S_2 \cup \ldots \cup S_n$ be a n-interval semigroup we can define the notion of Smarandache Lagrange semigroup, Smarandache p-Sylow semigroup and Smarandache weakly Lagrange semigroup.

We see if $S = S_1 \cup S_2 \cup \ldots \cup S_n$ be n-interval semigroup say each $S_i$ is of order $m_i$ then $|S| = m_1 m_2 \ldots m_n$.

We say a n-interval subsemigroup P of S divides the order of P if $o(P) / o(S)$.

We define S-Lagrange interval semigroup and S-weakly Lagrange interval semigroup in an analogous way [10-3].

We leave this routine task to the reader but give some examples of them.

*Example 2.1.18*: Let $S = S_1 \cup S_2 \cup S_3 \cup S_4 = \{S(X)$ where $X = \{([0, a_1], [0, a_2], [0, a_3], [0, a_4])\}\} \cup \{S(Y) / Y = \{([0, a_1], [0, a_2], [0, a_3])\}\} \cup \{S(A) / A = \{([0, a_1], [0, a_2])\} \cup \{S(B) / B = \{([0, a_1], [0, a_2], \ldots, [0, a_7])\}\}$ be a 4-interval semigroup.

$$A = A_1 \cup A_2 \cup A_3 \cup A_4 =$$

$$\left\{ \begin{pmatrix} [0,a_1] & [0,a_2] & [0,a_3] & [0,a_4] \\ [0,a_1] & [0,a_2] & [0,a_3] & [0,a_4] \end{pmatrix}, \right.$$
$$\begin{pmatrix} [0,a_1] & [0,a_2] & [0,a_3] & [0,a_4] \\ [0,a_2] & [0,a_3] & [0,a_4] & [0,a_1] \end{pmatrix},$$
$$\begin{pmatrix} [0,a_1] & [0,a_2] & [0,a_3] & [0,a_4] \\ [0,a_3] & [0,a_4] & [0,a_1] & [0,a_2] \end{pmatrix},$$
$$\left. \begin{pmatrix} [0,a_1] & [0,a_2] & [0,a_3] & [0,a_4] \\ [0,a_4] & [0,a_1] & [0,a_2] & [0,a_3] \end{pmatrix} \right\} \cup$$

$$\left\{ \begin{pmatrix} [0,a_1] & [0,a_2] & [0,a_3] \\ [0,a_1] & [0,a_2] & [0,a_3] \end{pmatrix}, \begin{pmatrix} [0,a_1] & [0,a_2] & [0,a_3] \\ [0,a_2] & [0,a_1] & [0,a_3] \end{pmatrix}, \right.$$

$$\left. \begin{pmatrix} [0,a_1] & [0,a_2] & [0,a_3] \\ [0,a_3] & [0,a_1] & [0,a_2] \end{pmatrix} \right\} \cup$$

$$\left\{ \begin{pmatrix} [0,a_1] & [0,a_2] \\ [0,a_1] & [0,a_2] \end{pmatrix}, \begin{pmatrix} [0,a_1] & [0,a_2] \\ [0,a_2] & [0,a_1] \end{pmatrix} \right\} \cup$$

$$\left\{ \begin{pmatrix} [0,a_1] & [0,a_2] & \ldots & [0,a_7] \\ [0,a_1] & [0,a_2] & \ldots & [0,a_7] \end{pmatrix}, \right.$$

$$\begin{pmatrix} [0,a_1] & [0,a_2] & \ldots & [0,a_7] \\ [0,a_2] & [0,a_3] & \ldots & [0,a_1] \end{pmatrix}, \ldots,$$

$$\left. \begin{pmatrix} [0,a_1] & [0,a_2] & \ldots & [0,a_7] \\ [0,a_7] & [0,a_1] & \ldots & [0,a_6] \end{pmatrix} \right\}$$

$\subseteq S_1 \cup S_2 \cup S_3 \cup S_4$, A is a 4-intrval group of S. o(A) / o(S). We have o(S) = $4^4 \cdot 3^3 \cdot 2^2 \cdot 7^7$ and o(A) = 4.3.2.7 so o(A) / o(S).

We see S is only a S-weakly Lagrange 4-interval semigroup. Consider $H = H_1 \cup H_2 \cup H_3 \cup H_4 = S_X \cup S_Y \cup S_A \cup S_B \subseteq S_1 \cup S_2 \cup S_3 \cup S_4$ is a 4-interval subgroup.

Clearly o(H) $\not|$ o(S); where o(H) = 4! × 3! × 2! × 7! and o(H) $\not|$ o(S). Thus S cannot be a S-Lagrange 4-interval semigroup only a S-weakly Lagrange 4-itnerval semigroup.

Having seen an example of a S-weakly Lagrange 4-interval semigroup.

*Example 2.1.19*: Let $S = S_1 \cup S_2 \cup S_3 = \{[0, a] / a \in Z_{12}, \times\} \cup \{[0, a] / a \in Z_{10}, \times\} \cup \{[0, a] / a \in Z_6, \times\}$ be a 3-interval semigroup. Consider $H = H_1 \cup H_2 \cup H_3$ any 3-interval group in $S = S_1 \cup S_2 \cup S_3$. We see S is a S-weakly Lagrange 3-interval semigroup. Take $A = A_1 \cup A_2 \cup A_3 = \{[0, a] / a \in \{1, 11\} \subseteq Z_{12}\} \cup \{[0, a] / a \in \{1, 9\} \subseteq Z_{10}\} \cup \{[0, a] / a \in \{1, 5\} \subseteq Z_6\} \subseteq S_1 \cup S_2 \cup S_3$, A is a 3-interval group and o (A) = 2.2.2 /

12.10.6. Consider $S_1$, the subgroups in $S_1$ are $A_1 = \{[0, 1], [0, 11]\}$, $B_1 = \{[0, 1], [0, 5]\}$ and $C_1 = \{[0, 1], [0, 7]\}$. The interval groups in $S_3$ are as follows: $A_3 = \{[0, 1], [0, 5]\}$ is the only interval subgroup of $S_3$.

Now the subgroups of $S_2$ are $A_2 = \{[0, 1], [0, 9]\}$, $B_2 = \{[0, 1], [0, 3], [0, 9], [0, 7]\}$ is given by the following table.

| × | [0, 1] | [0, 3] | [0, 7] | [0, 9] |
|---|---|---|---|---|
| [0, 1] | [0, 1] | [0, 3] | [0, 7] | [0, 9] |
| [0, 3] | [0, 3] | [0, 9] | [0, 1] | [0, 7] |
| [0, 7] | [0, 7] | [0, 1] | [0, 9] | [0, 3] |
| [0, 9] | [0, 9] | [0, 7] | [0, 3] | [0, 1] |

Clearly $B_2$ is a interval subgroup of $S_2$ but $o(B_2) \nmid o(S_2)$.

$S_2$ has only two subgroups. Only one of them divide the order of $S_2$. Thus S is only a S-weakly Lagrange interval semigroup and is not a S-Lagrange interval semigroup.

*Example 2.1.20*: Let $S = S_1 \cup S_2 \cup S_3 \cup S_4 = \{[0, a] / a \in Z_{61}\} \cup \{[0, a] / a \in Z_7\} \cup \{[0, a] / a \in Z_{11}\} \cup \{[0, a] / a \in Z_{13}\}$ be a 4-interval semigroup.

It is easily verified S is a S-4-interval semigroup as $A = A_1 \cup A_2 \cup A_3 \cup A_4 = \{[0, a] / a = 1 \text{ and } 60, ×\} \cup \{[0, a] / a = 1 \text{ and } 6, ×\} \cup \{[0, a] / a = 1 \text{ and } 10, ×\} \cup \{[0, a] / a = 1 \text{ and } 12, ×\} \subseteq S_1 \cup S_2 \cup S_3 \cup S_4$ is a 4-interval subgroup of S. Thus S is a Smarandache 4-interval semigroup. Clearly $o(A) \nmid o(S)$ as $o(A) = 2^4$ and $o(S) = 61.7.11.13$. Further S is only a S-weakly Lagrange semigroup.

Now in view of this we have the following theorem.

**THEOREM 2.1.3**: *Let $S = S_1 \cup S_2 \cup ... \cup S_n$ be a n-interval semigroup where each $S_i = \{[0, a] / a \in Z_p, p \text{ a prime}, ×\}; 1 \leq i \leq n$. Clearly S is a S-n-interval semigroup and S is only a S-weakly Lagrange n-interval semigroup.*

The proof is left as an exercise to the interested reader.

***Example 2.1.21***: Let $S = S_1 \cup S_2 \cup S_3 \cup S_4 = \{[0, a] / a \in Z_6, \times\} \cup \{[0, b] / b \in Z_{12}, \times\} \cup \{S(X) / X = ([0, a_1], [0, a_2], [0, a_3])\} \cup \{[0, a] / a \in Z_8, \times\}$ be a 4-interval semigroup. It is easily verified S is a S-weakly cyclic 4-interval semigroup.

***Example 2.1.22***: Let $S = S_1 \cup S_2 \cup S_3 \cup S_4 = \{[0, a] / a \in Z_{18}, \times\} \cup \{[0, b] / b \in Z_{40}, \times\} \cup \{[0, c] / c \in Z_{64}, \times\} \cup \{[0, d] / d \in Z_{72}, \times\}$ be a 4-interval semigroup. Clearly S is a S-weakly cyclic 4-interval semigroup.

In view of this we have the following theorem which guarantees a classes of S-n-interval semigroups which are S-weakly cyclic n-interval semigroups.

**THEOREM 2.1.4**: *Let $S = S_1 \cup S_2 \cup ... \cup S_n = \{S(X_1) / X_1 = ([0, a_1], ..., [0, a_{m_1}])\} \cup \{S(X_2) / X_2 = ([0, a_1], ..., [0, a_{m_2}])\} \cup ... \cup \{S(X_n) / X_n = (([0, a_1], ..., [0, a_{m_n}])\}$ be a n-interval semigroup. Clearly S is a S-weakly cyclic n-interval semigroup.*

This proof is also direct and hence is left as an exercise to the reader.

***Example 2.1.23***: Let $S = S_1 \cup S_2 \cup S_3 \cup S_4 = \{[0, a] / a \in Z_{10}, \times\} \cup \{[0,a] | a \in Z_{32}, \times\} \cup \{[0, a] / a \in Z_{42}, \times\} \cup \{[0, a] / a \in Z_{28}, \times\}$ be a 4-interval semigroup. Clearly S is a S-2-Sylow 4-interval semigroup.

Inview of this we have a theorem which gurantees the existence of a class of S-2-Sylow n-interval semigroups.

**THEOREM 2.1.5**: *Let $S = S_1 \cup S_2 \cup ... \cup S_n = \{[0, a] / a \in Z_{2m_1}, \times\} \cup \{[0, a] / a \in Z_{2m_2}, \times\} \cup ... \cup \{[0, a] / a \in Z_{2m_n} \times\}$ be a n interval semigroup. Clearly S is a S-2-Sylow n-interval semigroup.*

The proof is left as an exercise.

***Example 2.1.24***: Let $S = S_1 \cup S_2 \cup S_3 \cup S_4 = \{[0, a] / a \in Z_{12}\} \cup \{[0, b] / b \in Z_{30}, \times\} \cup \{[0, a] / a \in Z_{40}, \times\} \cup \{[0, a] / a \in Z_{16}, \times\}$ be a 4-interval semigroup. S has S-Cauchy elements.

However every n-interval semigroup need not have S-Cauchy elements.

***Example 2.1.25***: Let $S = S_1 \cup S_2 \cup S_3 \cup S_4 = \{[0, a] / a \in Z_{11}, \times\} \cup \{[0, b] / b \in Z_{11}, \times\} \cup \{[0, a] / a \in Z_{19}, \times\} \cup \{[0, a] / a \in Z_{13}, \times\}$ be a 4-interval semigroup. S has S-Cauchy elements.

Inview of this we have a theorem which guarantees the existence of n-interval semigroups which have no S-Cauchy elements.

**THEOREM 2.1.6**: *Let $S = S_1 \cup S_2 \cup \ldots \cup S_n = \{[0, a] / a \in Z_{p_1}$, $p_1$ a prime, $\times\} \cup \{[0, a] / a \in Z_{p_2}$, $p_2$ a prime, $\times\} \cup \ldots \cup \{[0, a] / a \in Z_{p_n}$, $p_n$ a prime, $\times\}$ be a n-interval semigroup where $p_1, p_2, \ldots, p_n$ are n distinct primes. S has no S-Cauchy elements.*

The proof is straight forward and hence is left as an exercise to the reader.

***Example 2.1.26***: Let $S = S_1 \cup S_2 \cup S_3 \cup S_4 = \{S(X) / X = ([0, a_1], [0, a_2], \ldots, [0, a_6])\} \cup \{S(Y) / Y = ([0, a_1], [0, a_2], \ldots, [0, a_{10}])\} \cup \{S(A) / A = ([0, a_1], [0, a_2], \ldots, [0, a_{15}])\} \cup \{S(B) / B = ([0, a_1], [0, a_2], \ldots, [0, a_{21}])\}$ be a 4-interval semigroup. S is a S-2-Sylow 4-interval semigroup.
Also S is a S(3, 5, 3, 7) - Sylow 4-interval semigroup. Further S is also a S - (3, 2, 3, 3) Sylow 4-interval semigroup. S is also a S (3, 5, 5, 7) - Sylow 4-interval semigroup.

Thus if $S = S_1 \cup S_2 \cup \ldots \cup S_n$ is such that each $S_i = S(X_i)$ with $X_i = ([0, a_1], \ldots, [0, a_{m_i}])$ where $m_i = p_1^i \ldots p_{t_i}^i$ ($p_j^i$ distinct primes, $1 \leq j \leq t_i$) for $i = 1, 2, \ldots, n$ be a n-interval semigroup.

Then S is a S-$\left(p_{t_1}^1, p_{t_2}^2, ..., p_{t_n}^n\right)$ Sylow n-interval semigroup where $p_{t_{k_i}}^i$ can vary in $p_1^i$ ... $p_{t_i}^i$ i=1,2, ..., n and $1 < k_i < i$. Thus we get several S-$\left(p_{t_1}^1, p_{t_2}^2, ..., p_{t_n}^n\right)$ Sylow n-interval semigroups from S.

It is easily verified these S - ($p_1^i$ ... $p_{t_i}^i$) - Sylow n-interval semigroups in general need not be conjugate.

When the n-interval semigroups are constructed using interval matrix semigroups or interval polynomial semigroups calculations become very difficult. At this juncture it is suggested that a nice program in general be made so that calculations become easy.

Now we can define n-interval homomorphisms of n-interval semigroups in four ways.

(1) We take two n-interval semigroups $S = S_1 \cup ... \cup S_n$ and $P = P_1 \cup P_2 \cup ... \cup P_n$ and define n-homomorphism from $\eta: S \to P$ by assigning to each $S_i$ a unique $P_j$, $1 < i, j < n$ where $\eta = \eta_1 \cup \eta_2 \cup ... \cup \eta_n$ and $\eta_i : S_i \to P_j$ such that each $\eta_i$ is an interval semigroup homomorphism.

(2) Another way of defining $\eta : S \to P$ where $\eta = \eta_1 \cup \eta_2 \cup ... \cup \eta_n : S \to P$ is such that $\eta_i : S_i \to P_j$ to a $S_i$ any $P_j$ is assigned that for more than one $S_i$ the same $P_j$ may be assigned.

(3) Suppose $S = S_1 \cup S_2 \cup ... \cup S_n$ and $P = P_1 \cup P_2 \cup ... \cup P_m$ be any n-interval semigroup and m-interval semigroup respectively. n < m.

Let $\eta = \eta_1 \cup \eta_2 \cup ... \cup \eta_n: S \to P$ is such that $\eta_i : S_i \to P_j$ each $P_j$ is distinct or we can assign to more than one $S_i$ same $P_j$'s.

If m < n then $\eta = \eta_1 \cup \eta_2 \cup \ldots \cup \eta_n : S_1 \cup S_2 \cup \ldots \cup S_n \to P_1 \cup \ldots \cup P_m$ where $\eta_i : S_i \to P_j$ maps more than one $S_i$ to same $P_j$'s.

Thus when we define the n-homomorphism of n-interval semigroups we can have several n-homomorphism for each $S_i$ can be mapped on to any one of the $P_j$'s.

Interested reader can analyse the properties of n-interval homomorphism of n-interval semigroups.

Now we proceed onto define quasi n-interval semigroups or quasi (s, r) - interval semigroup.

**DEFINITION 2.1.4**: *Let $S = S_1 \cup \ldots \cup S_n$ where s of the semigroups are distinct interval semigroups and n - s of them are just n-s distinct semigroups. Then we define S to be a quasi n-interval semigroup or quasi (s, n-s) - interval semigroup.*

We will first illustrate this situation by some examples.

*Example 2.1.27*: Let $S = S_1 \cup S_2 \cup S_3 \cup S_4 \cup S_5 = \{[0, a] / a \in Z_{12}, \times\} \cup \{Z_{19}, \times\} \cup \{[0, a] / a \in Z_{25}, \times\} \cup \{Z_{36}, \times\} \cup \{[0, a] / a \in Z_{30}, \times\}$ be a quasi 5-interval semigroup or quasi (3, 2) - interval semigroup.

*Example 2.1.28*: Let $S = S_1 \cup S_2 \cup S_3 \cup S_4 \cup S_5 \cup S_6 = \{[0, a] / a \in Z_8, \times\} \cup \{Z_{20}, \times\} \cup \{Z_{48}, \times\} \cup \{[0, a] / a \in Z_{11}, \times\} \cup \{Z_{18}, \times\} \cup \{[0, a] / a \in Z_{12}, \times\}$ be a quasi 6-interval semigroup or a quasi (3, 3)-interval semigroup.

*Example 2.1.29*: Let $S = S_1 \cup S_2 \cup S_3 = (Z_{48}, \times) \cup \{[0, a] / a \in Z_{12}, \times\} \cup \{Z_{20}, \times\}$ be a quasi 3-interval semigroup or quasi (2, 1) interval semigroup.

Now having seen examples of quasi (r, s)-interval semigroups we now proceed onto give examples of substructures. The task of defining these substructures is left as an exercise to the reader.

***Example 2.1.30***: Let $S = S_1 \cup S_2 \cup S_3 \cup S_4 = \{Z_{40}, \times\} \cup \{[0, a] / a \in Z_{18}, \times\} \cup \{Z_{48}, \times\} \cup \{[0, a] / a \in Z_{20}, \times\}$ be a quasi 4-interval semigroup.

Consider $V = V_1 \cup V_2 \cup V_3 \cup V_4 = \{\{0, 2, 4, 8, …, 38\}, \times\} \cup \{[0, a] / a \in \{0, 3, 6, 9, 12, 15\}, Z_{18}, \times\} \cup \{\{0, 6, 12, 18, …, 42\} \subseteq Z_{48}, \times\} \cup \{[0, a] / a \in \{0, 5, 10, 15\} \subseteq Z_{20}, \times\} \subseteq S_1 \cup S_2 \cup S_3 \cup S_4$; V is a quasi 4-interval subsemigroup of S. Clearly V is also a quasi 4-interval ideal of S.

***Example 2.1.31***: Let $S = S_1 \cup S_2 \cup S_3 = \{[0, a] / a \in R^+ \cup \{0\}, \times\} \cup \{[0, a] / a \in Q^+ \cup \{0\}, \times\}$ a quasi 3-interval semigroup.

Consider $X = X_1 \cup X_2 \cup X_3 = \{[0, a] / a \in Q^+ \cup \{0\}\} \cup \{[0, a] / a \in Z^+ \cup \{0\}\} \cup \{Q^+\} \subseteq S_1 \cup S_2 \cup S_3 = S$; X is only a quasi 3-interval subsemigroup but is not a quasi 3-interval ideal of S.

Thus every quasi n-interval subsemigroup of a quasi n-interval semigroup need not be a quasi n-interval ideal.

Inview of this we have the following theorem the proof of which is direct.

**THEOREM 2.1.7**: *Let $S = S_1 \cup S_2 \cup … \cup S_n$ be a quasi n-interval semigroup. Every quasi n-interval ideal of S is a quasi n-interval subsemigroup of S, but in general a quasi n-interval subsemigroup need not be a quasi n-interval ideal of S.*

We will illustrate by some examples S-quasi n-interval semigroups.

It is pertinent to mention here that every quasi n-interval semigroup need not in general be a S-quasi n-interval semigroup.

***Example 2.1.32***: Let $S = S_1 \cup S_2 \cup S_3 \cup S_4 = \{Z_{45}, \times\} \cup \{[0, a] / a \in Z_{20}\} \cup \{Z_{17}, \times\} \cup \{[0, a] / a \in Z_{80}\}$ be a quasi 4-interval semigroup. Consider $H = H_1 \cup H_2 \cup H_3 \cup H_4 = \{\{1, 44\} \subseteq Z_{45}\} \cup \{[0, 1], [0, 19] / 1, 19 \in Z_{20}\} \cup \{1, 16\} \cup$ $\{[0, 1] [0, 79] / 1, 79 \in Z_{80}\} \subseteq S_1 \cup S_2 \cup S_3 \cup S_4$. H is a quasi 4-interval group. So S is a S-quasi 4-interval semigroup.

*Example 2.1.33*: Let $S = S_1 \cup S_2 \cup S_3 = \{[0, a] / a \in Z^+ \cup \{0\}\} \cup \{[0, a] / a \in R^+, +\} \cup \{3Z^+ \cup \{0\}, +\}$ be a quasi 3-interval semigroup. Clearly S is not a S-quasi 3-interval semigroup.

Now we will give classes of quasi n-interval semigroups which are S-quasi n-interval semigroups.

**THEOREM 2.1.8**: *Let $S = S_1 \cup S_2 \cup ... \cup S_n$ where $S_i = \{[0, a] / a \in Z_{n_i}, \times\}$ and $S_j = \{Z_{m_j}, \times\}$, $1 < i, j < n$ be a quasi n-interval semigroup. S is a S-quasi n-interval semigroup.*

The proof is direct however we give a small hint which makes the proof obvious.
Let $A = A_1 \cup A_2 \cup ... \cup A_n$ where $A_i = \{[0, 1], [0, n_i-1]\} \subseteq S_i$ and $A_j = \{1, m_j-1\} \subseteq S_j$ are subgroups and A is a quasi interval n-group so S is a S-quasi n-interval semigroup.

**THEOREM 2.1.9**: *Let $S = S_1 \cup S_2 \cup ... \cup S_n = S(X_1) \cup S(m_2) S(X_2) \cup ... \cup S(m_2) \cup ... \cup S(X_n)$ where $S(X_i)$ is the interval symmetric semigroup group on $n_i$ intervals and $S(m_j)$ are just symmetric semigroups on $m_j$ elements, $1 \leq i, j \leq n$. Thus S is a quasi n-interval semigroup. S is a S-quasi n-interval semigroup.*

For this theorem also we only give an hint.
Every $S_i$ if $S_i$ is a interval symmetric semigroup then $S_i = S(X_i)$ has $A_i = S_{X_i}$ to be the symmetric interval group in $S_i$. Similarly for $S_j = S(m_j)$ we have $A_j = S_{m_j} \subseteq S(m_j)$ is the symmetric group, $1 \leq i, j \leq n$. So $A = A_1 \cup A_2 \cup ... \cup A_n \subseteq S_1 \cup S_2 \cup ... \cup S_n$ is the quasi n-interval group, hence S is a S-quasi n-interval semigroup.

Now we will give examples of S-quasi interval subsemigroups.

***Example 2.1.34***: Let $S = S_1 \cup S_2 \cup S_3 \cup S_4 \cup S_5 = \{Z_{20}, \times\} \cup \{[0, a] / a \in Z_{15}\} \cup \{Z_{19}, \times\} \cup \{[0, a] / a \in Z_{30}\} \cup \{Z_{12}, \times\}$ be a quasi 5-interval semigroup.

$P = P_1 \cup P_2 \cup P_3 \cup P_4 \cup P_5 = \{\{1, 19\} \subseteq Z_{20}, \times\} \cup \{[0, 1], [0, 14] / \{1, 14\} \in Z_{15}\} \cup \{1, 18 / 1, 18 \in Z_{19}\} \cup \{[0, 1], [0, 29] / 1, 29 \in Z_{30}\} \cup \{1, 11 / 1, 11 \in Z_{12}, \times\} \subseteq S_1 \cup S_2 \cup S_3 \cup S_4 \cup S_5$ is a quasi 5-interval group. So S is a S-quasi 5-interval semigroup. Now consider the quasi 5-interval subsemigroup $T = T_1 \cup T_2 \cup T_3 \cup T_4 \cup T_5 = \{\{0, 10\}, \times\} \cup \{[0, 5], [0, 10], [0, 0], \times\} \cup \{[0, 0], [0, 1]\} \cup \{[0, 0], [0, 10]\} \cup \{0, 6\} \subseteq S_1 \cup S_2 \cup S_3 \cup S_4 \cup S_5$ is a quasi 5-interval subsemigroup of S. Clearly T is not a S-quasi 5-interval subsemigroup. However S is a S-quasi 5-interval semigroup.

Inview of this we have the following theorem.

**THEOREM 2.1.10**: *Let $S = S_1 \cup S_2 \cup \ldots \cup S_n$ be a S-quasi n-interval semigroup. Then if $T = T_1 \cup T_2 \cup \ldots \cup T_n \subseteq S_1 \cup S_2 \cup \ldots \cup S_n = S$ be a quasi n-interval subsemigroup of S. T in general need not be a S-quasi n-interval subsemigroup of S.*

The proof is by counter example. Example 2.1.34 will prove the theorem.

**THEOREM 2.1.11**: *Let $S = S_1 \cup S_2 \cup \ldots \cup S_n$ be a quasi n-interval semigroup. If $P = P_1 \cup P_2 \cup \ldots \cup P_n \subseteq S_1 \cup S_2 \cup \ldots \cup S_n$ be a S-quasi n-interval subsemigroup of S then S is a S-quasi n-interval semigroup.*

The proof is direct and hence is left as an exercise to the reader.

Now one can as in case of n-interval semigroups define the notion of n-interval zero divisors, n-interval units and n-interval idempotents [10-3].

We will illustrate this situation by some examples.

***Example 2.1.35***: Let $S = S_1 \cup S_2 \cup S_3 \cup S_4 = \{[0, a] / a \in Z_{20}\} \cup \{Z_{12}, \times\} \cup \{[0, a] / a \in Z_{15}\} \cup \{Z_{24}, \times\}$ be a quasi 4-interval semigroup. Now take x = {[0, 10]} ∪ {6} ∪ {[0, 5]} ∪ {12} ∈ S. Choose y = {[0, 6]} ∪ {2} ∪ {[0, 6]} ∪ {6} ∈ S. We see xy ={0} ∪ {0} ∪ {0} ∪ {0}. Thus S has 4-interval zero divisor.

Consider x = [0, 19] ∪ {11} ∪ {[0, 14]} ∪ {23} ∈ S. We see $x^2$ = [0, 1] ∪ {1} ∪ {[0, 1]} ∪ {1} so x is a 4-interval unit in S.

x = {[0, 5]} ∪ {9} ∪ {[0, 10]} ∪ {9} in S is such that $x^2$ = [0, 25 (mod 20)] ∪ {81 (mod 12)} ∪ [0, 100 (mod 15)} ∪ {81 (mod 24)} = [0, 5] ∪ {9} ∪ {[0, 10]} ∪ {0} = x. Thus x is a 4-interval idempotent of S.

***Example 2.1.36***: Let S = $S_1$ ∪ $S_2$ ∪ $S_3$ = {[0, a] / a ∈ $Z^+$, ×} ∪ {[0, b] / b ∈ $Q^+$} ∪ {$R^+$, ×} be a quasi 3-interval semigroup. S has no 3-interval units, no 3-interval zero divisors and no three interval idempotents. Thus we have 3-interval semigroups which have none of the special elements. It is further important to note that S is not Smarandache quasi 3-interval semigroup.

Further S has no quasi 3-interval ideals. However S has quasi 3-interval subsemigroup say P = $P_1$ ∪ $P_2$ ∪ $P_3$ = {$3Z^+$, ×} ∪ {[0, b] / b ∈ $Z^+$} ∪ {[0, b] / b ∈ $Q^+$} ⊆ $S_1$ ∪ $S_2$ ∪ $S_3$. Infact S has infinite number of quasi 3-interval subsemigroups none of them are quasi 3-interval ideals.

Now as in case of n-interval semigroups or quasi interval bisemigroups discuss and study the concept of zero divisors, S-zero divisors, idempotents and S-idempotents. We give examples of quasi n-interval semigroups using interval matrices and interval polynomials.

***Example 2.1.37***: Let S = $S_1$ ∪ $S_2$ ∪ $S_3$ ∪ $S_4$ ∪ $S_5$ ∪ $S_6$ ∪ $S_7$ = $\left\{ \sum_{i=0}^{\infty} [0,a]x^i \;\middle|\; a \in Z^+ \cup \{0\}, \times \right\}$ ∪ {All 8 × 4 interval matrices with intervals of the form [0, a] / a ∈ $Z_{40}$, +} ∪ {[0, $a_1$], [0, $a_2$], …, [0, $a_9$] / $a_i$ ∈ $Z_{23}$, 1 ≤ i ≤ 9, ×} ∪

$$\left\{ \begin{bmatrix} [0,a_1] \\ [0,a_2] \\ [0,a_3] \\ [0,a_4] \end{bmatrix} \middle| a_i \in Z_{120}, +, 1 \le i \le 4 \right\}$$

∪ {All 7 × 7 upper triangular interval matrices with intervals of the form [0, a] where a ∈ $Q^+ \cup \{0\}$, +} ∪

$$\left\{ \sum_{i=0}^{\infty} [0,a]x^i \middle| a \in Z_{48}, \times \right\} \cup \left\{ \sum_{i=0}^{27} [0,a]x^i \middle| a \in Z_{144}, + \right\}$$

is a 7-interval semigroup of infinite order which is clearly non commutative.

***Example 2.1.38***: Let $S = S_1 \cup S_2 \cup S_3 \cup S_4 \cup S_5 = \{([0, a_1], [0, a_2], ..., [0, a_{12}]) / a_i \in Z_{50}, 1 \le i \le 12, \times\}$ ∪ {All 3 × 3 matrices with entries from $Z_{250}$, ×} ∪ {all 7 × 2 interval matrices with intervals of the form [0, a] / a ∈ $Z_{420}$, +} ∪

$$\left\{ \sum_{i=0}^{\infty} a_i x^i \middle| a_i \in Q^+ \cup \{0\}, \times \right\} \cup$$

{S(8)} be a quasi 5-interval semigroup.

We see in these two examples the interval semigroups and semigroups are of varying types.

Now having seen such hectrogeneous examples now we proceed onto describe n-interval groupoids and quasi n-interval groupoids.

## 2.2 n-Interval Groupoids

In this section we proceed onto describe the notion of n-interval groupoids, quasi n-interval groupoids and (m, n) interval semigroup-groupoids (interval m-semigroup-n-groupoids) and its quasi analogue.

**DEFINITION 2.2.1**: *Let $G = G_1 \cup G_2 \cup G_3 \cup \ldots \cup G_n$ be such that each $G_i$ is an interval groupoid and the $G_i$'s are distinct ($G_i \not\subseteq G_j$, if $i \neq j$ $1 \leq i, j \leq n$). We define '.' on G which takes operations on each of the $G_i$'s; $1 \leq i \leq n$. We define (G, .) to be the n-interval groupoid (n > 2). If n = 2 we call G as interval bigroupoid.*

We will illustrate this by examples.

*Example 2.2.1*: Let $G = G_1 \cup G_2 \cup G_3 \cup G_4 = \{[0, a] \mid a \in Z_{14}, *, (3, 5)\} \cup \{[0, b] \mid b \in Z_{20}, *, (7, 0)\} \cup \{[0, c] \mid c \in Z_{19}, *, (4, 4)\} \cup \{[0, d] / d \in Z_{17}, *, (1, 3)\}$ be a 4-interval groupoid. We define '.' on G as follows. Suppose $x = [0, 3] \cup [0, 8] \cup [0, 1] \cup [0, 5]$ and $y = [0, 2] \cup [0, 1] \cup [0, 8] \cup [0, 10] \in G$. Then

$$
\begin{aligned}
x.y &= ([0, 3] \cup [0, 8] \cup [0, 1] \cup [0, 5]) \cdot ([0, 2] \cup [0, 1] \cup [0, 8] \cup [0, 10]) \\
&= [0, 3] * [0, 2] \cup [0, 8] * [0, 1] \cup [0, 1] \cup [0, 8] \cup [0, 5] * [0, 10] \\
&= [0, 3.3 + 2.5 \pmod{14}] \cup [0, 7.8 + 0.1 \pmod{20}] \cup [0, (1.4 + 1.8) \bmod 19] \cup [0, (1.5 + 3.10) \pmod{17}] \\
&= [0, 5] \cup [0, 16] \cup [0, 12] \cup [0, 1] \in G.
\end{aligned}
$$

Thus (G, .) is a 4-interval groupoid of finite order.

*Example 2.2.2*: Let $S = S_1 \cup S_2 \cup S_3 \cup S_4 \cup S_5 = \{[0, a] / a \in Z^+ \cup \{0\}, *, (3, 2)\} \cup \{[0, a] / a \in Z_{45}, *, (2, 1)\} \cup \{[0, a] / a \in Z_{20}, *, (13, 0)\} \cup \{[0, a] / a \in R^+ \cup \{0\}, *, (3, 0)\} \cup \{[0, a] / a \in Z_{48}, *, (0, 7)\}$ be a 5-interval groupoid.

*Example 2.2.3*: Let $S = S_1 \cup S_2 \cup S_3 \cup S_4 = \{$All $3 \times 3$ interval matrices with intervals of the form $[0, a]$; $a \in Z_{12}, *, (7, 0)\} \cup$

$$\left\{\left\{\sum_{i=0}^{\infty}[0,a]x^i \,\Big|\, a \in Z_{15}, *, (0,3)\right\} \cup \left\{\sum_{i=0}^{\infty}[0,a]x^i \,\Big|\, a \in Z_{19}, *, (4,0)\right\}\right\}$$

$\cup \{([0, a_1], [0, a_2], \ldots, [0, a_6]) \mid a_i \in Z_{120}, *, (8, 5)\}$ be a 4-interval groupoid where we describe the operations in each of the $S_i$, $1 \leq i \leq 4$.

Let
$$A = \begin{pmatrix} [0,3] & [0,1] & [0,7] \\ [0,1] & [0,8] & [0,0] \\ [0,2] & [0,0] & [0,5] \end{pmatrix}, B = \begin{pmatrix} [0,2] & [0,0] & [0,0] \\ [0,7] & [0,5] & [0,0] \\ [0,1] & [0,3] & [0,8] \end{pmatrix}$$
in $S_1$.

$$A * B = \begin{pmatrix} [0,3] & [0,1] & [0,7] \\ [0,1] & [0,8] & [0,0] \\ [0,2] & [0,0] & [0,5] \end{pmatrix} * \begin{pmatrix} [0,2] & [0,0] & [0,0] \\ [0,7] & [0,5] & [0,0] \\ [0,1] & [0,3] & [0,8] \end{pmatrix} =$$

$$\begin{pmatrix} [0,(3.7+0.2)\bmod 12] & [0,(0.7+0.1)\bmod 12] & [0,(0.7+0.7)\bmod 12] \\ [0,(7.7+0.1)\bmod 12] & [0,(5.7+8.0)\bmod 12] & [0,(7.0+0.0)\bmod 12 \\ [0,(7.1+0.2)\bmod 12] & [0,(8.3+0.0)\bmod 12] & [0,(8.7+0.5)\bmod 12 \end{pmatrix}$$

$$= \begin{pmatrix} [0,9] & [0,0] & [0,0] \\ [0,1] & [0,11] & [0,0] \\ [0,7] & [0,0] & [0,8] \end{pmatrix}.$$

We now describe the operation in $S_2$.

$$p(x) = [0, 5] x^7 + [0, 2] x^3 + [0, 3] x + [0, 1]$$
and
$$q(x) = [0, 1] x^6 + [0, 9] x + [0, 7] x^3 + [0, 4] ;$$

$$\begin{aligned} p(x) * q(x) &= [0, 0] x^7 + [0, 3] x^6 + [0, (2.0+7.3) \bmod 15] x^3 + [0, (3.0 + 9.3) \bmod 15] x + \\ & \quad [0, (0.1 + 4.3) \bmod 15] \\ &= [0, 3] x^6 + [0, 9] x^3 + [0, 9] x. \end{aligned}$$

Thus this is the way * on the interval polynomial groupoid $S_2$ is defined.

In a similar way the operation on $S_3$ is carried out. Now consider the operation * on $S_4$. Take x = ([0, 3], [0, 2], [0, 1], [0, 0], [0, 5], [0, 9]) and y = ([0, 9], [0, 12], [0, 20], [0, 40], [0, 8], [0, 1]) in $S_4$. Now x*y = ([0, (24 + 45) mod 120],

[0, (16+60) mod 120], [0, (8+120) mod 120], [0, (0+40×5) mod 120], [0, (40+40) mod 120], [0, (72+5) mod 120])
= ([0, 69], [0, 76], [0, 8], [0, 0], [0, 80], [0, 77]).

***Example 2.2.4***: Let $S = S_1 \cup S_2 \cup S_3 \cup S_4 \cup S_5 \cup S_6 = \{$All $5 \times 5$ interval matrices with intervals of the form [0, a] where $a \in Z_7$, *, $(3, 4)\} \cup \{([0, a_1], [0, a_2], \ldots, [0, a_{10}]) / a_i \in Z_{11}, 1 \leq i \leq 10$, *, $(8, 3)] \cup \{[0, a] / a \in Z_6$, *, $(4, 0)\} \cup \{[0, a] / a \in Z_{14}$, *, $(0, 3)\} \cup \{\sum_{i=0}^{9}[0,a]x^i$ , *, $a \in Z_{10}, (3, 2)\} \cup \{[0, a] / a \in Z_4$, *, $(2, 2)\}$ be a 6-interval groupoid.

We have seen examples of n-interval groupoids.
Now we proceed onto illustrate the substructure of an n-interval groupoid. As the definition is easy the reader is left with the task of defining the substructures.

***Example 2.2.5***: Let $S = S_1 \cup S_2 \cup S_3 \cup S_4 = \{[0, a] / a \in Z^+ \cup \{0\}, (3, 0), *\} \cup \{([0, a_1], [0, a_2], [0, a_3]) / a_i \in Z_{10}, 1 \leq i \leq 3$, *, $(7, 2)\} \cup \{[0, b] / b \in Z_{40}$, *, $(10, 13)\} \cup \{[0, c] / c \in Z_{12}$, *, $(8, 0)\}$ be a 4-interval groupoid. Consider $V = V_1 \cup V_2 \cup V_3 \cup V_4 = \{[0, a] / a \in 7Z^+ \cup \{0\}$, *, $(3, 0)\} \cup \{([0, a_1], 0, [0, a_2]) / a_1, a_2 \in Z_{10}$, *, $(7, 2)\} \cup \{[0, b] / b \in \{0, 10, 20, 30\} \subseteq Z_{40}$, *, $(10, 13)\} \cup \{[0, c] / c \in \{0, 2, 4, 6, 8, 10\} \subseteq Z_{12}$, *, $(8, 0)\} \subseteq S_1 \cup S_2 \cup S_3 \cup S_4$; V is a 4-interval subgroupoid of S.

***Example 2.2.6***: Let $S = S_1 \cup S_2 \cup S_3 \cup S_4 \cup S_5 = \{([0, a_1], [0, a_2], [0, a_3], [0, a_4], [0, a_5]) / a_i \in Z_{12}$, *, $1 \leq i \leq 5, (3, 2)\} \cup \{[0, a] / a \in Z_{35}$, *, $(7, 0)\} \cup \{$All $2 \times 2$ interval matrices with intervals of the form [0, a] where $a \in Z_{40}$, *, $(10, 2)\} \cup$

$$\left\{\begin{bmatrix}[0,a]\\[0,b]\\[0,c]\end{bmatrix} \middle| a,b,c \in Z_{15}, *, (3,2)\right\} \cup \left\{\sum_{i=0}^{7}[0,a]x^i \middle| a \in Z_{24}, *, (3,11)\right\}$$

be a 5-interval groupoid. $P = P_1 \cup P_2 \cup P_3 \cup P_4 \cup P_5 = \{([0, a_1], [0, a_2], \ldots, [0, a_5]) / a_i \in \{0, 2, 4, 6, 8, 10\} \subseteq Z_{12}$, *, $(3, 2)\}$

$\cup \{[0, a] / a \in \{0, 5, 10, \ldots, 30\} \subseteq Z_{35}, *, (7, 0)\} \cup \{$All $2 \times 2$ interval matrices with intervals of the form $[0, a]$ where $a \in \{0, 2, 4, \ldots, 38\} \subseteq Z_{40}, (10, 2)\} \cup$

$$\left\{ \begin{bmatrix} [0,a] \\ [0,b] \\ [0,c] \end{bmatrix} \middle| a,b,c \in \{0,5,10\} \subseteq Z_{15}, *, (3,2) \right\} \cup$$

$$\left\{ \sum_{i=0}^{3} [0,a]x^i \middle| a \in Z_{24}, *, (3,11) \right\}$$

$\subseteq S_1 \cup S_2 \cup S_3 \cup S_4 \cup S_5$ is a 5-interval subgroupoid of S.

***Example 2.2.7***: Let $S = S_1 \cup S_2 \cup S_3 \cup S_4 = \{$All $2 \times 5$ interval matrices with intervals of the form $[0, a]$, $a \in Z_8, *, (3,1)\} \cup$

$$\left\{ \sum_{i=0}^{9} [0,a]x^i \middle| a \in Z_{10}, (2,0), * \right\}$$

$\cup \{[0, a] / a \in Z_{12}, *, (4, 3)\} \cup \left\{ \sum_{i=0}^{4} [0,a]x^i \middle| a \in Z_{40}, *, (0,7) \right\}$

be a 4-interval groupoid.

Consider $G = G_1 \cup G_2 \cup G_3 \cup G_4 = \{[0, a] / a \in \{0, 2, 4, 6\} \subseteq Z_8, *, (3, 1)\} \cup$

$$\left\{ \sum_{i=0}^{9} [0,a]x^i \middle| a \in \{0,2,4,6,8\} \subseteq Z_{10}, *, (2,0) \right\} \cup$$

$\{[0, a] / a \in \{0, 3, 6, 9\} \subseteq Z_{12}, (4, 3)\} \cup$

$$\left\{ \sum_{i=0}^{4} [0,a]x^i \middle| a \in \{0,4,8,12,\ldots,36\} \subseteq Z_{40}, *, (0,7) \right\}$$

$\subseteq S_1 \cup S_2 \cup S_3 \cup S_4 = S$ is a 4-interval subgroupoid of S.

Now having seen n-interval subgroupoids now we proceed onto describe other properties like S-n-interval groupoids and special identities satisfied by the n-interval groupoids. We will call a n-interval groupoid $S = S_1 \cup S_2 \cup \ldots \cup S_n$ to be a S-n-interval groupoid if S has a proper subset $A = A_1 \cup A_2 \cup \ldots \cup A_n$ where each $A_i$ is an interval semigroup under the operations of S, $1 \leq i \leq n$, that is if S contains a n-interval semigroup then we call S to be a Smarandache n-interval groupoid.

We will illustrate this situation by some examples.

***Example 2.2.8:*** Let $S = S_1 \cup S_2 \cup S_3 \cup S_4 \cup S_5 = \{[0, a] / a \in Z_{10}, *, (5, 6)\} \cup \{[0, a] / a \in Z_{12}, *, (3, 9)\} \cup \{[0, a] / a \in Z_{12}, *, (3, 4)\} \cup \{[0, a] / a \in Z_4, *, (2, 3)\} \cup \{[0, a] / a \in Z_6, *, (4, 3)\}$ be a 5-interval groupoid it is easily verified S is a Smarandache 5-interval groupoid of finite order.

***Example 2.2.9***: Let $S = S_1 \cup S_2 \cup S_3 \cup S_4 = \{[0, a] / a \in Z_6, *, (3, 5)\} \cup \{[0, a] / a \in Z_{14}, *, (7, 8)\} \cup \{[0, a] / a \in Z_{12}, *, (5, 10)\} \cup \{[0, a] / a \in Z_{12}, *, (2, 10)\}$ be a 4-interval groupoid and S is a Smarandache 4-interval grouped of finite order.

***Example 2.2.10***: Let $G = G_1 \cup G_2 \cup G_3 \cup G_4 \cup G_5 \cup G_6 \cup G_7 = \{[0, a] / a \in Z_9, *, (7, 3)\} \cup \{[0, a] / a \in Z_{16}, *, (7, 10)\} \cup \{[0, a] / a \in Z_{20}, *, (10, 11)\} \cup \{[0, a] / a \in Z_{19}, *, (17, 3)\} \cup \{[0, a] / a \in Z_{14}, *, (8, 7)\} \cup \{[0, a] / a \in Z_{22}, *, (12, 11)\} \cup \{[0, a] / a \in Z_{27}, *, (23, 5)\}$ is a 7-interval groupoid which is a S-7-interval groupoid.

In view of this we have the following theorem which guarantees the existence of a class of S-n-interval groupoids.

**THEOREM 2.2.1**: *Let $G = G_1 \cup G_2 \cup G_3 \cup \ldots \cup G_n$ where $G_i = \{[0, a] / a \in Z_{m_i}, *, (t_i, u_i)$ where $(t_i, u_i) = 1$ and $t_i + u_i \equiv 1 \pmod{m_i}\}$ true for $i = 1, 2, \ldots, n$ and each $m_i > 5$. Then G is a Smarandache n-interval groupoid of order $m_1, m_2, \ldots, m_n$.*

Proof is straight forward as every n-element. $x = [0, a_1] \cup [0, a_2] \cup \ldots \cup [0, a_n] \in G$ such that $x*x = x$ so $\{x\}$ is a n-interval semigroup in G.

Hence the claim.

We say an n-interval groupoid $G = G_1 \cup G_2 \cup \ldots \cup G_n$ is said to be an n-interval idempotent groupoid if for every $x = x_1 \cup x_2 \ldots \cup x_n$ in G we have $x * x = x$.

We will first illustrate this by some examples.

***Example 2.2.11***: Let $S = S_1 \cup S_2 \cup S_3 \cup S_4 = \{[0, a] / a \in Z_{10}, (7, 4), *\} \cup \{[0, a] / a \in Z_{13}, *, (9, 5)\} \cup \{[0, a] / a \in Z_{20}, *, (11, 10)\} \cup \{[0, a] / a \in Z_{15}, *, (9, 7)\}$ be a 4-interval groupoid. It is easily verified G is a 4-interval idempotent groupoid. For take $x = [0, 5] \cup [0, 7] \cup [0, 12] \cup [0, 9]$ in S.

$$\begin{aligned} x^2 &= [0, 5] * [0, 5] \cup [0, 7] * [0, 7] \cup [0, 12] * [0, 12] \cup [0, 9] * [0, 9] \\ &= [0, (35 + 20) \bmod 10] \cup [0, (63+35) \bmod 13] \cup [0, (132+120) \bmod 20] \cup [0, (81+63) \bmod 15] \\ &= [0, 5] \cup [0, 7] \cup [0, 12] \cup [0, 9] \\ &= x. \end{aligned}$$

It is easily verified $x^2 = x$ for every $x \in S$ is a 4-interval groupoid. It is easily verified G is a 4-interval idempotent groupoid.

In view of this we have the following theorem which guarantees the existence of n-interval idempotent groupoid.

**THEOREM 2.2.2**: *Let $G = G_1 \cup G_2 \cup G_3 \cup \ldots \cup G_n$ be a n-interval groupoid $G_i = \{[0, a] / a \in Z_{m_i}, (t_i, u_i) = 1\ (t_i + u_i) \equiv 1 \bmod m_i, *\}$ is an interval groupoid for each $i = 1, 2, \ldots, n$.*
*G is a n-interval idempotent groupoid.*

The proof is direct and is left as an exercise to the reader.
Now we will give examples of n-interval groupoids which has S-n-interval subgroupoids.

***Example 2.2.12***: Let $S = S_1 \cup S_2 \cup S_3 \cup S_4 = \{[0, a] / a \in Z_{14}, *, (7, 8)\} \cup \{[0, a] / a \in Z_{12}, *, (5, 10)\} \cup \{[0, a] / a \in Z_{10}, *, (5, 6)\} \cup \{[0, a] / a \in Z_{12}, *, (1, 3)\}$ be a 4-interval groupoid of finite order. Consider $A = A_1 \cup A_2 \cup A_3 \cup A_4 = \{[0, a] / a \in \{0, 4\}, *, \{7, 8\}\} \cup \{[0, a] / a \in \{0, 6\} \subseteq Z_{12}, *, (5, 10)\} \cup \{[0, a] / a \in \{0, 5\} \subseteq Z_{10}, *, (5, 6)\} \cup \{[0, a] / a \in \{0, 6, 3, 9\} \subseteq Z_{12}, *, (1, 3)\} \subseteq S_1 \cup S_2 \cup S_3 \cup S_4$. Clearly a is a Smarandache 4-interval subgroupoid. As $P = P_1 \cup P_2 \cup P_3 \cup P_4 = \{[0, 4], *, (7, 8)\} \cup \{[0, 6] / 6 \in Z_{12}, *, (5, 16)\} \cup \{[0, 5] / 5 \in Z_{10}, *, (5, 6)\} \cup \{[0, 0], [0, 6] / 0, 6 \in Z_{12}, *, (1, 3)\} \subseteq A_1 \cup A_2 \cup A_3 \cup A_4 = A \subseteq S_1 \cup S_2 \cup S_3 \cup S_4 = S$ is a 4-interval semigroup in A.

Thus A is a S-4-interval subgroupoid of S.

Now we have the following nice theorem the proof of which is direct.

**THEOREM 2.2.3**: *Let $G = G_1 \cup G_2 \cup \ldots \cup G_n$ be a n-interval groupoid. If $A = A_1 \cup A_2 \cup \ldots \cup A_n$ is a n-interval subgroupoid of G which is a S-n-interval subgroupoid, then G is itself a S-n-interval groupoid. But if G is a S-n-interval groupoid then in general every n-interval subgroupoid of G need not be a S n-interval subgroupoid of G.*

To this effect the interested reader can give examples.

***Example 2.2.13***: Let $G = G_1 \cup G_2 \cup G_3 = \{[0, a] / a \in Z_4, *, (2,3)\} \cup \{[0, a] / a \in Z_{12}, (3,4), *\} \cup \{[0, a] / a \in Z_6, *, (2,3)\}$ be a Smarandache Bol 3-interval groupoid.

This can be easily verified by the reader. However G is not a S-strong Bol 3-groupoid.

***Example 2.2.14***: Let $G = G_1 \cup G_2 \cup G_3 = \{[0, a] / a \in Z_6, *, (4,3)\} \cup \{[0, a] / a \in Z_4, *, (2,3)\} \cup \{[0, a] / a \in Z_6, *, (3,5)\}$ be a Smarandache strong 3-interval P-groupoid. Hence G is a Smarandache 3-interval P-groupoid.

The following theorem guarantees the existence of Smarandache n-interval P-groupoid.

**THEOREM 2.2.4**: *Let $G = G_1 \cup G_2 \cup ... \cup G_n$ be a n-interval groupoid. If each $G_i = \{[0, a] / a \in Z_{m_i}, *, t_i + u_i \equiv 1 \pmod{m_i} (t_i, u_i)\}$ is an interval groupoid $1 \leq i \leq n$. G is a Smarandache alternative n-interval groupoid if and only if $t_i^2 = t_i \pmod{m_i}$ and $u_i^2 = u_i \pmod{m_i}$, for every i=1,2,...,n.*

The proof is direct for analogous methods refer [ ].

The next theorem guarantees the existence of a class of Smarandache n-interval P-groupoids.

**THEOREM 2.2.5**: *Let $G = G_1 \cup G_2 \cup ... \cup G_n$ where $G_i = \{[0, a] / a \in Z_{m_i}, *, (t_i, u_i) 1 + u_i = 1 \pmod{m_i}\}$, $1 \leq i \leq n$ be a n-interval groupoid. G is a Smarandache n-interval P-groupoid if and only if $t_i^2 = t_i \pmod{m_i}$ and $u_i^2 = u_i \pmod{m_i}$ for i = 1, 2, ..., n.*

Now we will give example of a Smarandache strong n-interval Moufang groupoid.

*Example 2.2.15*: Let $G = G_1 \cup G_2 \cup G_3 \cup G_4 \cup G_5 = \{[0, a] / a \in Z_{12}, *, (4, 9)\} \cup \{[0, a] / a \in Z_{10}, (5, 6), *\} \cup \{[0, a] / a \in Z_{15}, *, (6, 10)\} \cup \{[0, a] / a \in Z_{20}, (5, 16), *\} \cup \{[0, a] / a \in Z_{21}, (7, 15)\}$ be a 5-interval groupoid.
Clearly G is a Smarandache strong Moufang 5-interval groupoid.

We give a theorem which gurantees the existence of Smarandache strong Moufang n-interval groupoid.

**THEOREM 2.2.6**: *Let $G = G_1 \cup G_2 \cup ... \cup G_n = \{[0, a] / a \in Z_{m_1}, *, (t_1, u_1), t_1 + u_1 = 1 \pmod{m_1}\} \cup \{[0, a] / a \in Z_{m_2}, *, (t_2, u_2), t_2 + u_2 = 1 \pmod{m_2}\} \cup ... \cup \{[0, a] / a \in Z_{m_n}, *, (t_n, u_n), t_n + u_n = 1 \pmod{m_n}\}$ be a n-interval groupoid. G is a*

*Smarandache strong Moufang n-interval groupoid if and only if* $u_i^2 = u_i \pmod{m_i}$ *and* $t_i^2 = t_i \pmod{m_i}$, $i=1,2,\ldots,n$.

Now we give the existence of a class of Smarandache strong Bol n-interval groupoid.

**THEOREM 2.2.7**: *Let* $G = G_1 \cup G_2 \cup \ldots \cup G_n = \{[0, a] / a \in Z_{m_1}, *, (t_1, u_1); t_1 + u_1 \equiv 1 \pmod{m_1}\} \cup \{[0, a] / a \in Z_{m_2}, *, (t_2, u_2), t_2 + u_2 \equiv 1 \pmod{m_2}\} \cup \ldots \cup \{[0, a] / a \in Z_{m_n}, *, (t_n, u_n), *, t_n + u_n \equiv 1, \pmod{m_n}\}$ *be a n-interval groupoid. G is a Smarandache strong Bol n-interval groupoid if and only if* $t_i^2 = t_i \pmod{m_i}$ *and* $u_i^2 = u_i \pmod{m_i}$; $i=1,2,\ldots,n$.

The proof is straight forward for
$(([0, x_i] * [0, y_i]) * ([0, z_i]) * [0, x_i]$
$\quad = ([0, t_i x_i + u_i y_i] * [0, z_i]) * [0, x_i]$
$\quad = [0, t_i^2 + t_i u_i y_i + u_i z_i] * [0, x_i]$
$\quad = [0, t_i^3 x_i + t_i^2 u_i y_i + t_i u_i z_i + u_i x_i]$
$\quad = [0, t_i x_i + t_i u_i y_i + t_i u_i z_i + u_i x_i] \quad\quad - \text{I}$
$\quad\quad\quad\quad\quad\quad\quad\quad\quad ( t_i^3 = t_i \pmod{m_i})$

$[0, x_i] * (([0, y_i] * [0, z_i]) * [0, x_i])$
$\quad = [0, x_i] * ([0, t_i y_i + u_i z_i] * [0, x_i])$
$\quad = [0, x_i] * [0, t_i^2 u_i y_i + t_i u_i z_i + u_i^2 x_i]$
$\quad = [0, t_i x_i + t_i^2 u_i y_i + t_i u_i^2 x_i + u_i x_i] \quad\quad - \text{II}$
$\quad\quad\quad\quad\quad\quad\quad\quad\quad ( u_i^2 = u_i \pmod{m_i})$

I and II are equal for every $[0, x_i], [0, y_i], [0, z_i] \in G_i$; for $i = 1, 2, \ldots, n$. Hence G is a Smarandache strong Bol n-interval groupoid.

We give examples of Smarandache n-interval idempotent groupoids.

*Example 2.2.16*: Let $G = G_1 \cup G_2 \cup G_3 \cup G_4 \cup G_5 = \{[0, a] / a \in Z_{11}, *, (6,6)\} \cup \{[0, a] / a \in Z_{19}, *, (10,10)\} \cup \{[0, a] / a \in Z_{13}, *, (7,7)\} \cup \{[0, a] / a \in Z_{23}, *, (12,12)\} \cup \{[0, a] / a \in Z_{43},$

*, (22,22)} be a 5-interval groupoid. It is easily verified G is a Smarandache 5-interval idempotent groupoid.

We have a class of Smarandache n-interval groupoid, which is evident from the following theorems.

**THEOREM 2.2.8**: *Let $G = G_1 \cup G_2 \cup G_3 \cup \ldots \cup G_n = \{[0, a] / a \in Z_{p_1}, *, \left(\frac{p_1+1}{2}, \frac{p_1+1}{2}\right), *, p_1 \text{ a prime}\} \cup \{[0, a] / a \in Z_{p_2}, *, \left(\frac{p_2+1}{2}, \frac{p_2+1}{2}\right), *, p_2 \text{ a prime}\} \cup \ldots \cup \{[0, a] / a \in Z_{p_n}, *, \left(\frac{p_n+1}{2}, \frac{p_n+1}{2}\right), *, p_n \text{ a prime}\}$ (all the $p_i$'s are n-distinct primes, i=1,2,…,n). G is a Smarandache n-interval idempotent groupoid.*

The proof is straight forward and hence is left as an exercise for the reader to prove.

Now having seen the properties enjoyed by n-interval groupoids now we proceed onto define quasi n-interval groupoids and (t, s) interval semigroup - groupoid and (t, s) quasi interval semigroup- groupoid.

Let $G = G_1 \cup G_2 \cup \ldots \cup G_n$ where some t number $G_i$'s are distinct interval groupoids and the remaining n-t are just groupoids then we define G to be a quasi n-interval groupoid or quasi (t, (n-t)) - interval groupoid.

If in $G = G_1 \cup G_2 \cup \ldots \cup G_n$, t of the $G_i$'s are interval semigroups and n-t- of the $G_j$'s are interval groupoids and groupoids we define G to be a quasi n-interval semigroup-groupoid.

We will illustrate this situation by some examples.

***Example 2.2.17***: Let $G = G_1 \cup G_2 \cup G_3 \cup G_4 \cup G_5 \cup G_6 \cup G_7$ = $\{Z_9 (3, 2)\} \cup \{[0, a] / a \in Z_{40}, *, (19,11)\} \cup \{[0, a] / a \in Z_{24}, *, (11,13)\} \cup \{Z_{12} (7, 5)\} \cup \{[0, a] / a \in Z_{42}, *, (8,11)\} \cup \{Z_{27} (3, 1)\} \cup \{Z_{45} (11, 13)\}$ be a quasi 7-interval groupoid.

***Example 2.2.18***: Let $G = G_1 \cup G_2 \cup G_3 \cup G_4 = \{$All $3 \times 3$ interval matrices with intervals of the form $[0, a]$ where $a \in Z_7$, $(2, 3)\} \cup \{(a_1, a_2, a_3, a_4) / a_i \in Z_5, *, (3, 3)\} \cup \{[0, a] / a \in Z_9, *, (2, 7)\} \cup \{Z_{11}(3, 2)\}$ be a quasi 4-interval groupoid.

Let
$$x = \begin{pmatrix} [0,3] & [0,2] & [0,1] \\ 0 & [0,4] & 0 \\ [0,3] & 0 & [0,2] \end{pmatrix} \cup (2, 1, 0, 3) \cup \{[0, 7]\} \cup \{3\}$$

and
$$y = \begin{pmatrix} [0,3] & 0 & [0,1] \\ 0 & [0,5] & 0 \\ [0,5] & [0,4] & [0,2] \end{pmatrix} \cup (3, 1, 2, 4) \cup [0, 5] \cup \{9\} \in G.$$

x.y is calculated as follows;

$$x.y = \begin{pmatrix} [0,3] & [0,2] & [0,1] \\ 0 & [0,4] & 0 \\ [0,3] & 0 & [0,2] \end{pmatrix} * \begin{pmatrix} [0,3] & 0 & [0,1] \\ 0 & [0,5] & 0 \\ [0,5] & [0,4] & [0,2] \end{pmatrix}$$
$$\cup (2, 1, 0, 3) * (3, 1, 2, 4) \cup [0, 7] * [0, 5] \cup \{3 * 9\}$$

$$= \begin{pmatrix} [0,6+9(\bmod 7)] & [0,4+0(\bmod 7)] & [0,2+3(\bmod 7)] \\ [0,0+0(\bmod 7)] & [0,8+15(\bmod 7)] & [0,0] \\ [0,6+15(\bmod 7)] & [0,12(\bmod 7)] & [0,4+6(\bmod 7)] \end{pmatrix}$$
$$= [(6+9(\bmod 5), (3+3)(\bmod 5), (0+6 \bmod 5), (9+12) \bmod 5]$$
$$\cup [0, 14 + 35 (\bmod 9)] \cup [9+18 (\bmod 11)]$$

$$= \begin{bmatrix} [0,1] & [0,4] & [0,5] \\ [0,0] & [0,9] & [0,0] \\ [0,0] & [0,5] & [0,3] \end{bmatrix} \cup (0, 4, 1, 1) \cup [0, 4] \cup \{5\}$$

is in G. Thus G is a quasi 4-interval groupoid of finite order.

*Example 2.2.19*: Let $G = G_1 \cup G_2 \cup G_3 \cup G_4 \cup G_5 = Z_{10} (3, 2) \cup \{[0, a] / a \in Z_{11}, *, (8, 9)\} \cup \{Z_{45} (7, 2)\} \cup \{[0, a] / a \in Z_9, *, (2, 4)\} \cup \{\sum_{i=0}^{7}[0,a]x^i / a \in Z_{40}, *, (3, 2)\}$ be a quasi 5-interval groupoid of finite order.

*Example 2.2.20*: Let $G = G_1 \cup G_2 \cup G_3 \cup G_4 = \{[0, a] / a \in Z_{10}, \times\} \cup \{[0, a] / a \in Z_{12}, \times\} \cup \{[0, a] / a \in Z_{19}, \times\} \cup \{[0, a] / a \in Z_{20}, \times, (11, 7)\}$ be a 4-interval semigroup-groupoid of finite order.

*Example 2.2.21*: Let $S = S_1 \cup S_2 \cup S_3 \cup S_4 \cup S_5 = \{S (X) / X = ([0, a_1], [0, a_2], [0, a_3], [0, a_4])\} \cup \{[0, a] / a \in Z^+ \cup \{0\}, (9, 3)\} \cup \{[0, a] / a \in Z_{120}, *, (43, 29)\} \cup \{[0, a] / a \in Z_{40}, \times\} \cup \{[0, a] / a \in Z_{25}, *, (3, 7)\}$ be a 5-interval groupoid semigroup of infinite order.

Having seen examples of these structures it is a matter of routine to define substructures, however we give examples of them.

*Example 2.2.22*: Let $S = S_1 \cup S_2 \cup S_3 \cup S_4 = \{[0, a] / a \in Z_4, (2, 3), *\} \cup \{Z_{12} (3, 4)\} \cup \{[0, a] / a \in Z_{10}, *, (5, 6)\} \cup Z_{12} (1, 3)$ be a 4-interval groupoid. $H = H_1 \cup H_2 \cup H_3 \cup H_4 = \{\{[0, 0], [0, 2] / 0, 2 \in Z_4, *, (2, 3)\} \cup \{0, 4 \in Z_{12}, *, (3, 4)\} \cup \{[0, 0], [0, 2] / 0, 2 \in Z_{10}, *, (5, 6)\} \cup \{0, 3, 6, 9, \in Z_{12}, *, (1, 3)\} \subseteq S_1 \cup S_2 \cup S_3 \cup S_4$; is a quasi 4-interval subgroupoid of S.

*Example 2.2.23*: Let $G = G_1 \cup G_2 \cup G_3 \cup G_4 = \{[0, a] / a \in Z_{12}, *, (1, 3)\} \cup \{Z_6 (4, 5)\} \cup \{[0, a] / a \in Z_8, *, (2, 6)\} \cup Z_{12} (10, 8)$ be a quasi 4-interval groupoid of finite order. Consider $A = A_1 \cup A_2 \cup A_3 \cup A_4 = \{[0, a] / a \in \{0, 3, 6, 9\} \subseteq Z_{12}, *, (1, 3)\} \cup \{\{1, 3, 5\} \subseteq Z_6, *, (4, 5)\} \cup \{[0, a] / a \in \{0, 2, 4, 6\} \subseteq Z_8, *, (2, 6)\} \cup \{\{6, 2, 10\} \subseteq Z_{12}, *, (10, 8)\} \subseteq S_1 \cup S_2 \cup S_3 \cup S_4$, A is a quasi 4-interval subgroupoid of G.

Now we can as in case of n-interval groupoids mention special identities satisfied by quasi n-interval groupoids.

In case of quasi n-interval groupoids also all theorems given for n-interval groupoids can be proved with appropriate modifications. We now proceed onto describe substructure of n-interval semigroup groupoid.

*Example 2.2.24:* Let $G = G_1 \cup G_2 \cup G_3 \cup G_4 \cup G_5 = \{[0, a] / a \in Z_{24}, \times\} \cup \{[0, a] / a \in Z_9, *, (5, 3)\} \cup \{[0, a] / a \in Z_{40}, \times\} \cup \{[0, a] / a \in Z_{12}, *, (3, 9)\} \cup \{[0, a] / a \in Z_{12}, \times\}$ be a 5-interval semigroup-groupoid. Let $A = A_1 \cup A_2 \cup A_3 \cup A_4 \cup A_5 = \{[0, a] / a \in \{0, 2, 4, \ldots, 22\} \subseteq Z_{24}, \times\} \cup \{[0, a] / a \in \{1, 2, 4, 5, 7, 8\} \subseteq Z_9, *, (5, 3)\} \cup \{[0, a] / a \in \{0, 10, 20, 30\} \subseteq Z_{40}, \times\} \cup \{[0, a] / a \in \{0, 3, 6, 9\} \subseteq Z_{12}, *, (3, 9)\} \cup \{[0, a] / a \in \{0, 3, 6, 9\} \subseteq Z_{12}, \times\} \subseteq G_1 \cup G_2 \cup G_3 \cup G_4 \cup G_5$, A is a 5-interval subsemigroup - subgroupoid of G.

*Example 2.2.25*: Let $G = G_1 \cup G_2 \cup G_3 \cup G_4 = \{[0, a] / a \in Z_{24}, \times\} \cup \{[0, a] / a \in Z_{12}, *, (3, 9)\} \cup \{[0, a] / a \in Z_{10}, *, (5, 6)\} \cup \{S(X) / X = ([0, a_1], [0, a_2], [0, a_3], [0, a_4])\}$ be a 4-interval semigroup-groupoid. It is easily verified G is a Smarandache 4-interval semigroup-groupoid.

Now having seen examples of these new structures now we proceed onto define the notion of quasi n-interval semigroup-groupoid. $G = G_1 \cup G_2 \cup \ldots \cup G_n$ is a quasi n-interval semigroup-groupoid if some $G_i$'s are interval groupoids some $G_j$'s are groupoids some $G_k$'s are interval semigroups and the rest are semigroups.

We will describe them by some examples.

*Example 2.2.26*: Let $V = V_1 \cup V_2 \cup V_3 \cup V_4 \cup V_5 \cup V_6 = \{Z_{25}, \times\} \cup \{Z_7 (3, 2)\} \cup \{[0, a] / a \in Z_{40}, \times\} \cup \{[0, a] / a \in Z_{23}, *, (3, 2)\} \cup \{\begin{bmatrix} a & b \\ c & d \end{bmatrix} / a, b, c, d \in Z_{12}, \times\} \cup \{([0, a], [0, b], [0, c]) / a, b, c \in Z_{43}, *, (7, 0)\}$ be a quasi 6-interval semigroup-groupoids.

If both of interval semigroups and semigroups are Smarandache and both groupoids and interval groupoids are Smarandache then we define the quasi n-interval semigroup-groupoid to be a Smarandache quasi n-interval semigroup-groupoid.

Interested reader can construct examples of them. The notion of zero divisors, idempotents, S-zero divisors and S-idempotent are defined in case of these structures also. However it is only semiassociative so one cannot deal with special identities. All other results can be derived and illustrated with examples by any interested reader. In the next section we proceed on to define and describe n-interval groups and their generalizations.

## 2.3 n-Interval Groups and their Properties

In this section we proceed onto describe n-interval groups, quasi n-interval groups, n-interval group-semigroups and n-interval groupoid - groups and enumerate a few of the properties related with them.

**DEFINITION 2.3.1**: *Let $G = G_1 \cup G_2 \cup ... \cup G_n$ be such that each $G_i$ is an interval group and $G_i \neq G_j$ or $G_i \not\subseteq G_j$; if $i \neq j$, $1 \leq i, j \leq n$. Then G obtains the operation '.', componentwise inherited from $G_i$'s so G with this operation is defined as the n-interval group. If n = 2 we get the interval bigroup of biinterval group.*

We will illustrate this situation by some examples.

*Example 2.3.1*: Let $G = G_1 \cup G_2 \cup G_3 \cup G_4 = \{[0, a] / a \in Z_{12}, +\} \cup \{[0, a] / a \in Z_{19} \setminus \{0\}, \times\} \cup \{[0, a] / a \in Z_{17}, +\} \cup \{[0, a] / a \in Z_{43} \setminus \{0\}, \times\}$ be a 4-interval group of finite order. Clearly G is a commutative 4-interval group. Consider $x = \{[0, 3] \cup [0, 2] \cup [0, 7] \cup [0, 40]\}$ and $y = [0, 7] \cup [0, 10] \cup [0, 12] \cup [0, 10]$ in G.

x.y    =    ([0, 3] + [0, 7]) $\cup$ ([0, 2] × [0, 10]) $\cup$ ([0, 7] + [0, 12]) $\cup$ ([0, 40] × [0, 10])

$$= [0, 10 \,(\text{mod } 12)] \cup [0, 20 \,(\text{mod } 19)] \cup [0, 19 \,(\text{mod } 17)] \cup [0, 400 \,(\text{mod } 43)]$$
$$= [0, 10] \cup [0, 1] \cup [0, 2] \cup [0, 13] \in G.$$

It is easily verified G has a group structure.

*Example 2.3.2*: Let $G = G_1 \cup G_2 \cup G_3 \cup G_4 \cup G_5 = \{S_x \,/\, x = ([0, a_1], [0, a_2], [0, a_3], [0, a_4])\} \cup \{\sum_{i=0}^{\infty}[0,a]x^i \,/\, a \in Z_{29} \setminus \{0\}, \times\} \cup \{[0, a] \,/\, a \in Z_{45}, +\} \cup \{[0, a] \,/\, a \in Z_{420}, +\} \cup \{[0, a] \,/\, a \in Z_{11} \setminus \{0\}, \times\}$ be a 5-interval group. G is non commutative infinite 5-interval group.

*Example 2.3.3*: Let $G = G_1 \cup G_2 \cup G_3 \cup G_4 \cup G_5 \cup G_6 = \{[0, a] \,/\, a \in Z_{15}, +\} \cup \{[0, a] \,/\, a \in Z_{23} \setminus \{0\}, \times\} \cup \{\sum_{i=0}^{9}[0,a]x^i \,/\, a \in Z_{45}, +\} \cup \{S_x \,/\, x = ([0, a_1] \,[0, a_2] \,[0, a_3])\}\} \cup \{A = \begin{bmatrix} [0,a] & [0,b] & [0,e] \\ [0,c] & [0,d] & [0,f] \end{bmatrix} \,/\, a, b, c, d, e, f, \in Z_{420}, +\} \cup \{\Sigma [0, a] \,/\, a \in Z_{19} \setminus \{0\}, \times\}$ be a 6-interval group. $e = [0, 0] \cup [0, 1] \cup [0, 0] \cup \begin{pmatrix} [0,a_1] & [0,a_2] & [0,a_3] \\ [0,a_1] & [0,a_2] & [0,a_3] \end{pmatrix} \cup \left\{ \begin{bmatrix} 0 & 0 & 0 \\ 0 & 0 & 0 \end{bmatrix} \right\} \cup [0, 1]$ acts as the 6-interval identity of G.

Now having seen examples of n-interval groups we now proceed onto describe substructures in them and illustrate it with examples.

**DEFINITION 2.3.2**: *Let $G = G_1 \cup G_2 \cup G_3 \cup \ldots \cup G_n$ be an n-interval group under the operation '.'. Let $H = H_1 \cup H_2 \cup H_3 \cup \ldots \cup H_n \subseteq G_1 \cup G_2 \cup G_3 \cup \ldots \cup G_n$; if each $H_i$ is an interval subgroup of $G_i$, i=1,2,…, n then (H, .) is defined as the n-interval subgroup of G. If each $H_i$ is normal in $G_i$ for i = 1, 2, 3, …, n then we call H to be a n-interval normal subgroup of G.*

*We will say G is n-interval simple group if G has no n-interval normal subgroup.*

Now these definitions are described by the following examples.

***Example 2.3.4:*** Let $G = G_1 \cup G_2 \cup G_3 \cup G_4 \cup G_5 = \{[0, a] / a \in Z_{40}, +\} \cup \{[0, a] / a \in Z_{29} \setminus \{0\}, \times\} \cup \{S_x$ where $X = ([0, a_1]$
$[0, a_2], [0, a_3])\} \cup \{\begin{bmatrix}[0,a]\\[0,b]\end{bmatrix} / b, a \in Z_{45}, +\} \cup \{\begin{bmatrix}[0,a] & [0,b]\\[0,c] & [0,d]\end{bmatrix} /$
a, b, c, d $\in Z_{200}, +\}$ be a 5-interval group. Choose $H = H_1 \cup H_2 \cup H_3 \cup H_4 \cup H_5 = \{[0, a] / a \in \{0, 2, 4, \ldots, 38\}, +\} \cup \{[0, 1], [0, 28] / 1, 28 \in Z_{29} \setminus \{0\}, \times\} \cup$

$$\left\{\begin{pmatrix}[0,a_1] & [0,a_2] & [0,a_3]\\[0,a_1] & [0,a_2] & [0,a_3]\end{pmatrix}, \begin{pmatrix}[0,a_1] & [0,a_2] & [0,a_3]\\[0,a_2] & [0,a_1] & [0,a_3]\end{pmatrix}\right\} \cup$$

$$\left\{\begin{bmatrix}[0,a]\\[0,b]\end{bmatrix} \mid a, b \in \{0, 3, 6, 9, 12, \ldots, 42, +\}\right\} \cup \left\{\begin{bmatrix}[0,a] & 0\\ 0 & [0,b]\end{bmatrix} /\right.$$

a, b $\in Z_{200}, +\} \subseteq G_1 \cup G_2 \cup G_3 \cup G_4 \cup G_5$, H is a 5-interval subgroup of G. Every 5-interval subgroup in G is not a 5-interval normal subgroup of G. Infact G is not a 5-interval simple group.

***Example 2.3.5***: Let $G = G_1 \cup G_2 \cup G_3 \cup G_4 \cup G_5 = \{A_x / x = ([0, a_1], [0, a_2], [0, a_3], [0, a_8])\} \cup \{[0, a] / a \in Z_{11}, +\} \cup \{[0, 1], [0, 26]/ 1, 26 \in Z_{27}, \times\} \cup \{A_x / x = \{[0, a_1], [0, a_2], \ldots, [0, a_{10}]\} \cup \{[0, a] / a \in Z_5, +\}$ be a 5-interval group. G is a 5-interval simple group. Infact G has no 5-interval subgroup.

***Example 2.3.6***: Let $G = G_1 \cup G_2 \cup G_3 \cup G_4 \cup G_5 = \{[0, a] / a \in Z_{45}, +\} \cup \{([0, a_1], [0, a_2], \ldots, [0, a_9]) / a_i \in Z_{20}, 1 \le i \le 9, +\}$
$\cup \{\begin{bmatrix}[0,a_1]\\[0,a_2]\\[0,a_3]\end{bmatrix} / a_i \in Z_{42}, +. 1 \le i \le 3\} \cup \{\begin{bmatrix}[0,a] & [0,a_1]\\[0,a_1] & [0,a]\end{bmatrix} / a, a1$
$\in Z_{24}, +\} \cup \{[0, a] / a \in Z_{19} \setminus \{0\}, \times\}$ be a 5-interval group. This has 5-interval subgroups. G is commutative and is of finite order.

Now having seen examples of substructures in n-interval groups we can extend all the classical theorems for groups to n-interval groups with out any difficulty.

Further if $G = G_1 \cup G_2 \cup \ldots \cup G_n$ is a n-interval group where each $G_i$ is of finite order then $o(G) = o(G_1) o(G_2) \ldots o(G_n) = |G_1| |G_2| \ldots |G_n|$.

We have several classical theorems like Lagrange, Sylow, Cauchy, Cayley which are true for n-interval groups.

We now proceed onto describe quasi n-interval groups. Let $G = G_1 \cup G_2 \cup \ldots \cup G_n$ be such that some of the groups $G_i$ are groups and the rest of the $G_j$ are interval groups. We define G to be quasi n-interval group.

We give examples of this structure.

*Example 2.3.7*: Let $G = G_1 \cup G_2 \cup G_3 \cup G_4 \cup G_5 = S_3 \cup \{[0, a] / a \in Z_{12}, +\} \cup \{[0, a] / a \in Z_{17} \setminus \{0\}, \times\} \cup D_{26} \cup \{[0, a] / a \in Z_{19} \setminus \{0\}, \times\}$ be a quasi 5-interval group. Clearly G is of finite order and is non commutative. Order of $G = 6 \times 12 \times 16 \times 12 \times 18$.

*Example 2.3.8*: Let $G = G_1 \cup G_2 \cup G_3 \cup G_4 = \{[0, a] / a \in Z_{45}, +\} \cup \{[0, a] / a \in Z_{19} \setminus \{0\}\} \cup <g / g^5 = 1> \cup (Z_{10}, +)$ be a quasi 4-interval group $o(G) = o(G_1) o(G_2) o(G_3) o(G_4) = 45 \times 18 \times 5 \times 10$. Clearly G is commutative.

*Example 2.3.9*: Let $G = G_1 \cup G_2 \cup G_3 \cup G_4 \cup G_5 = \{A_9\} \cup \{[0, a] / a \in Z_{19}, +\} \cup \{[0, a] / a \in Z_{11}, +\} \cup \{A_7\} \cup \{g / g^7 = 1\}$ be a quasi 5-interval group.

Clearly G is non commutative and $o(G) = \frac{\lfloor 9}{2} \times 19 \times 11 \frac{\lfloor 7}{2} \times 7$. Further G has no quasi 5-interval normal subgroup hence G is simple.

*Example 2.3.10*: Let $G = G_1 \cup G_2 \cup G_3 \cup G_4 = \{S_4\} \cup \{[0, a] / a \in Z_{12}, +\} \cup \{[0, a] / a \in Z_{11} \setminus \{0\}, \times\} \cup \{A_5\}$ be a quasi 4-interval group of order $\lfloor 4 \times 12 \times 10 \times \lfloor 5/2$. Clearly G is non commutative. Consider $H = H_1 \cup H_2 \cup H_3 \cup H_4 = \{A_4\} \cup \{[0, a] / a \in \{0, 2, 4, 6, 8, 10\}, +\} \cup \{[0, 1], [0, 10]/ 1, 10 \in Z_{11} \setminus \{0\}, \times\} \cup \{<\begin{pmatrix} 1 & 2 & 3 & 4 & 5 \\ 2 & 3 & 4 & 5 & 1 \end{pmatrix}>\}$. $G_1 \cup G_2 \cup G_3 \cup G_4$. H is a quasi 4-interval subgroup of G. Clearly G has no quasi 4-interval normal subgroups as $A_5$ has no normal subgroup.

o (H) = $\lfloor 4/2 \cdot 6 \times 2 \times 5 = 720$.

Interested reader can construct any number of examples. All classical theorem for groups are true in case of quasi n-interval groups. One can easily verify this claim.

Now we proceed onto define n-interval semigroup group. Let $G = G_1 \cup G_2 \cup \ldots \cup G_n$ where some of the $G_i$'s are interval semigroups and the rest are interval groups. We define G with inherited operations from G to be a n-interval semigroup - group.

We will illustrate this situation by some examples.

*Example 2.3.11*: Let $G = G_1 \cup G_2 \cup G_3 \cup G_4 \cup G_5 = \{[0, a] / a \in Z_{10}, \times\} \cup \{[0, a] / a \in Z_{11} \setminus \{0\}, \times\} \cup \{[0, b] / b \in Z_{40}, \times\} \cup \{[0, a] / a \in Z_{15}, +\} \cup \{([0, a] [0, b] [0, c])/ a, b, c \in Z_{24}, \times\}$ be a 5-interval group-semigroup of finite order. Clearly G is commutative 5-interval group semigroup.

*Example 2.3.12*: Let $G = G_1 \cup G_2 \cup G_3 \cup G_4 = \{[0, a] / a \in Z_{16}, +\} \cup \{([0, a], [0, b], [0, c])/ a, b, c \in Z_{15}, \times\} \cup \{$All $4 \times 4$ interval matries with intervals of the form [0, a] where $a \in Z_{30}$, $\times\} \cup \{[0, a] / a \in Z_3 \setminus \{0\}, \times\}$ is a 4-interval semigroup-group of finite order.

We can have substructures in them.

*Example 2.3.13*: Let $G = G_1 \cup G_2 \cup G_3 \cup G_4 \cup G_5 = \{[0, a] / a \in Z_{40}, +\} \cup \{S_{(X)} / X = ([0, a_1] \ldots [0, a_5])$ where $a_i \in I$, $1 \le i \le 5\} \cup \{S_Y / Y = ([0, a_1], \ldots, [0, a_4]) / a_i \in I_m$ $1 \le i \le 4\} \cup \{[0, a] / a \in Z_{13} \setminus \{0\}, \times\} \cup \{[0, a] / a \in Z_{16}, +\}$ be a 5-interval semigroup – group of finite order. Clearly G is non commutative.

***Example 2.3.14***: Let $G = G_1 \cup G_2 \cup G_3 \cup G_4 = \{S (X) / X = ([0, 1], [0, 2], [0, 3])\} \cup \{[0, a] / a \in Z_{24}, \times\} \cup \{[0, a] / a \in Z_7 \setminus \{0\}, \times\} \cup \{[0, a] / a \in Z_{14}, +\}$ be a 4-interval group-semigroup of finite order.

Clearly G is non commutative and G is not simple. G has 4-interval subgroup - subsemigroup.

$H = H_1 \cup H_2 \cup H_3 \cup H_4 = \{S_x\} \cup \{[0, a] / a \in 2Z_{24}, \times\} \cup \{[0, 1], [0, 6] \times\} \cup \{[0, a] / a \in \{0, 7\} +\} \subseteq G_1 \cup G_2 \cup G_3 \cup G_4$ is a 4-interval subsemigroup subgroup of G.

Now we will proceed onto define the notion of quasi interval semigroup - group.

Let $G = G_1 \cup G_2 \cup \ldots \cup G_n$ where $G_1, G_2, \ldots, G_n$ are set of semigroups, interval semigroups, groups and interval groups. G with the inherited operations from $G_1, G_2, \ldots, G_n$ forms the quasi n-interval semigroup-group.

We will illustrate this situation by some examples.

***Example 2.3.15***: Let $S = S_1 \cup S_2 \cup S_3 \cup S_4 \cup S_5 = \{S (3)\} \cup \{[0, a] / a \in Z_{24}, \times\} \cup \{S_5\} \cup \{[0, a] / a \in Z_{19} \setminus \{0\}, \times\} \cup \left\{\begin{bmatrix}[0,a]\\[0,b]\end{bmatrix} / a, b \in Z_{25}, +\right\}$ be a 5-interval semigroup - group. Clearly S is of finite order but non commutative.

Now we give an example of a quasi n-interval semigroup - group of infinite order and non commutative.

***Example 2.3.16***: Let $G = G_1 \cup G_2 \cup G_3 \cup G_4 \cup G_5 = \{Z, +\} \cup \{[0, a] / a \in Z^+ \cup \{0\}, \times\} \cup \{S_{26}\} \cup S (46) \cup \{[0, a] / a \in Z_{45}, \times\}$ be the quasi 5-interval semigroup - group. Clearly G is non commutative but is of infinite order.

***Example 2.3.17***: Let $G = G_1 \cup G_2 \cup G_3 \cup G_4 = \{[0, a] / a \in Z_7 \setminus \{0\}, \times\} \cup \{([0, a], [0, b], [0, c]) / a, b, c \in Z_{25}, \times\} \cup \{G_3 = <g / g^{27} = 1>\} \cup \{Z_{28}, \times\}$ be a quasi 4-interval semigroup-group of finite order.

Defining substructures in them is a matter of routine we will give one or two examples of them.

***Example 2.3.18***: Let $G = G_1 \cup G_2 \cup G_3 \cup G_4 \cup G_5 = \{S_7\} \cup \{S(3)\} \cup \{[0, a] / a \in Z_{240}, \times\} \cup \{[0, a] / a \in Z_{43} \setminus \{0\}, \times\} \cup \left\{\begin{bmatrix}[0,a]\\ [0,b]\\ [0,c]\end{bmatrix}\right| a, b, c \in Z_{12}, +\}$ be a 5-quasi interval group-semigroup of finite order which is non commutative. Consider $W = W_1 \cup W_2 \cup W_3 \cup W_4 \cup W_5 = A_7 \cup \{S_3\} \cup \{[0, a] / a \in 2Z_{240}, \times\} \cup \{[0, a] / a = 1, 42 \in Z_{43} \setminus \{0\}, \times\} \cup \left\{\begin{bmatrix}[0,a]\\ [0,b]\\ [0,c]\end{bmatrix}\right/ a, b, c \in 2Z_{12}, +\} \subseteq G_1 \cup G_2 \cup G_3 \cup G_4 \cup G_5$; clearly W is a quasi 5-interval subgroup -subsemigroup of G.

***Example 2.3.19***: Let $G = G_1 \cup G_2 \cup G_3 \cup G_4 = \{S(X)$ where $X = ([0, a_1], [0, a_2] [0, a_3])\} \cup \{S_7\} \cup <g / g^{120} = 1> \cup \{[0, a] / a \in Z_{47} \setminus \{0\}, \times\}$ be quasi 4-interval group - semigroup. Take $H = H_1 \cup H_2 \cup H_3 \cup H_4 = \{A(X) / X = ([0, a_1], [0, a_2], [0, a_3])\} \cup \{A_7\} \cup \{<g^{2n} / n = 0, 2, \ldots, 118$ when $g^{120} = 1>\} \cup \{[0, a] / a = 1$ and $46 \in Z_{47} \setminus \{0\}, \times\} \subseteq G_1 \cup G_2 \cup G_3 \cup G_4$. H is a quasi 4-interval subgroup-subsemigroup.

It is pertinent to mention here that all classical theorems in group theory are not in general true in case of quasi 4-interval group-semigroup.

This sort of contradictions can be easily derived by constructing counter examples.

***Example 2.3.20***: Let $G = G_1 \cup G_2 \cup G_3 \cup G_4 \cup G_5$ be a quasi 5 interval semigroup - group. If even only one of the $G_i$ is S (n) or S(X) symmetric semigroup or symmetric interval semigroup we see Lagrange theorem for finite groups, Sylow theorems for finite groups and Cauchy theorem for finite groups cannot hold good.

We will leave the task of providing counter examples to the reader. Now in the next section we proceed onto define n-interval loops, quasi n-interval loops and their generalizations.

## 2.4 n-Interval Loops

In this section we proceed onto describe and define n-interval loops, quasi n-interval loops and n-interval loop-group and so on.

**DEFINITION 2.4.1:** *Let $G = G_1 \cup G_2 \cup G_3 \cup ... \cup G_n$ be such that each $G_i$ be an interval loop. G with the component wise operations from each $G_i$ is defined to be the n-interval loop.*

We will illustrate this by some simple examples.

*Example 2.4.1*: Let $L = L_1 \cup L_2 \cup L_3 \cup L_4 \cup L_5 = \{[0, a] / a \in \{e, 1, 2, ..., 19\}$ (3), *\} $\cup \{[0, a] / a \in \{e, 1, 2, ..., 7\}, *, 5\} \cup \{[0, a] / a \in \{e, 1, 2, 3, ..., 23\}, 10, *\} \cup \{[0, a] / a \in \{e, 1, 2, ..., 15\}, 8, *\} \cup \{[0,a] \mid a \in \{e, 1, 2, ..., 17\}, *, 8\}$ be a 5-interval loop. Clearly L is of finite order and L is non commutative.

*Example 2.4.2*: Let $L = L_1 \cup L_2 \cup L_3 \cup L_4 = \{[0, a] / a \in \{e, 1, 2, ..., 29\}, *, 20\} \cup \{[0, a] / a \in \{e, 1, 2, ..., 11\}, *, 6\} \cup \{[0, a] / a \in \{e, 1, 2, ..., 9\}, *, 8\} \cup \{[0, a] / a \in \{e, 1, 2, ..., 31\}, *, 9\}$ be a 4-interval loop.

*Example 2.4.3*: Let $L = L_1 \cup L_2 \cup L_3 \cup L_4 \cup L_5 = \{[0, a] / a \in \{e, 1, 2, ..., 43\}, *, 8\} \cup \{[0, a] / a \in \{e, 1, 2, ..., 43\}, *, 40\} \cup \{[0, a] / a \in \{e, 1, 2, ..., 43\},. *, 32\} \cup \{[0, a] / a \in \{e, 1, 2, ..., 43\}, *, 25\} \cup \{[0,a] \mid a \in \{e, 1, 2, ..., 43\}, *, 15\}$ be a 5-interval loop of order $(43)^5$. Clearly L is non commutative.

*Example 2.4.4*: Let $L = L_1 \cup L_2 \cup L_3 \cup L_4 = \{[0, a] / a \in \{e, 1, 2, ..., 17\}, 9, *\} \cup \{[0, a] / a \in \{e, 1, 2, ..., 27\}, 14, *\} \cup \{[0, a] / a \in \{e, 1, 2, ..., 19\}, *, 10\} \cup \{[0, a] / a \in \{e, 1, 2, ..., 21\},$

11, *} be a 4-interval loop. Clearly L is a commutative 4-interval loop.

***Example 2.4.5***: Let L = $L_1 \cup L_2 \cup L_3$ = {[0, a] / a ∈ {e, 1, 2, …, 9}, 8, *} ∪ {[0, a] / a ∈ {e, 1, 2, …, 15}, *, 8} ∪ {[0, a] / a ∈ {e, 1, 2, …, 25}, *, 12} be a 3-interval loop. Consider H = $H_1 \cup H_2 \cup H_3$ = {[0, a] / a ∈ {e, 1, 4, 7} ⊆ {e, 1, 2, 3, …, 9}, *, 8} ∪ {[0, a] / a ∈ {e, 1, 6, 11} ⊆ {e, 1, 2, …, 15}, *, 8} ∪ {[0, a] / a ∈ {e, 1, 6, 11, 16, 21}} ⊆ $L_1 \cup L_2 \cup L_3$ is a 3-interval subloop of L.

Now having given one example of n-interval subloop interested reader can construct n-interval subloops and study their properties.

We will give some of the properties of n-interval loops.

***Example 2.4.6***: Let L = $L_1 \cup L_2 \cup L_3 \cup L_4 \cup L_5$ = {[0, a] / a ∈ {e, 1, 2, …, 19}, *, 10} ∪ {[0, a] / a ∈ {e, 1, 2, …, 19}, *, 10} ∪ {[0, a] / a ∈ {e, 1, 2, …, 27} *, 17} ∪ {[0, a] / a ∈ {e, 1, 2, …, 43}, *, 25} ∪ {[0, a] / a ∈ {e, 1, 2, …, 47}, *, 32} be a 5-interval loop of finite order.

***Example 2.4.7***: let L = $L_1 \cup L_2 \cup L_3 \cup L_4$ = {[0, a] / a ∈ {e, 1, 2, …, 45}, *, 8} ∪ {[0, a] / a ∈ {e, 1, 2, …, 9}, *, 8} ∪ {[0, a] / a ∈ {e, 1, 2, …, 25}, *, 12} ∪ {[0, a] / a ∈ {e, 1, 2, …, 15}, *, 2} be a 4-interval loop.

H = $H_1 \cup H_2 \cup H_3 \cup H_4$ = {[0, a] / a ∈ {e, 1, 16, 31} ⊆ {e, 1, 2, …, 45} *, 8} ∪ {[0, a] / a ∈ {e, 1, 4, 7} ⊆ {e, 1, 2, …, 9}, *, 8} ∪ {[0, a] / a ∈ {e, 1, 6, 11, 16, 21} ⊆ {e, 1, 2, 3, …, 25}, *, 12} ∪ {[0, a] / a ∈ {e, 1, 6, 11} ⊆ {e, 1, 2, …, 15} *, 2} ⊆ $L_1 \cup L_2 \cup L_3 \cup L_4$ is 4-interval subloop of L.

The reader is expected to define n-interval subloops of a n-interval loop.

We now give a theorem which guarantees the existence of n-interval subloops of an n-interval loop.

**THEOREM 2.4.1**: *Let $L = L_1 \cup L_2 \cup ... \cup L_n = \{[0, a] / a \in \{e, 1, 2, ..., m_1\}, *, s_1\} \cup \{[0, a] / a \in \{e, 1, 2, ..., m_2\}, *, s_2\} \cup ... \cup \{[0, a] / a \in \{e, 1, 2, ..., m_n\}, *, s_n\}$ be a n-interval loop where each $H = H_1(t_1) \cup H_2(t_2) \cup ... \cup H_n(t_n) = \{e, i, i + t_1, ..., (k_1 - t_1)t_1\} \cup \{e, j, j + t_2, ..., (k_2 - t_2)t_2\} \cup ... \cup \{e, p, p + t_n, ..., (k_n - t_n)t_n\}$ where $k_i = n_i / t_i$ $i = 1, 2, ..., n$. H is a n-interval subloop infact a S-n-interval subloop of L $(1 < t_1, j < t_2, ..., p < t_n)$.*

The proof is left to the reader as an exercise. This proves the existence of n-interval subloops and S-n-interval subloops of an n-interval loop.

Now we can define Smarandache n-interval loop. Let $L = L_1 \cup L_2 \cup L_3 \cup ... \cup L_n$ be a n-interval loop. If each $L_i$ is a Smarandache interval loop then we call L to be Smarandache n-interval loop (S-n-interval loop).

We will illustrate this by some examples.

*Example 2.4.8*: Let $L = L_1 \cup L_2 \cup L_3 \cup L_4 = \{[0, a] / a \in \{e, 1, 2, ..., 13\}, *, 8\} \cup \{[0, a] / a \in \{e, 1, 2, ..., 15\}, *, 8\} \cup \{[0, a] / a \in \{e, 1, 2, ..., 17\}, *, 7\} \cup \{[0, a] / a \in \{e, 1, 2, ..., 19\}, *, 6\}$ be a 4-interval loop. L is a S-4-interval loop for $A = A_1 \cup A_2 \cup A_3 \cup A_4 = \{[0, e], [0, 7], *, 8\} \cup \{[0, e] [0, 5], *, 8\} \cup \{[0, e], [0, 15], *, 7\} \cup \{[0, e] [0, 11], *, 6\} \subseteq L_1 \cup L_2 \cup L_3 \cup L_4$ is a 4-interval group in L hence L is a S-4-interval loop.

*Example 2.4.9*: Let $L = L_1 \cup L_2 \cup L_3 \cup L_4 \cup L_5 = \{[0, a] / a \in \{e, 1, 2, ..., 19\}, *, 7\} \cup \{[0, a] / a \in \{e, 1, 2, ..., 19\}, *, 6\} \cup \{[0, a] / a \in \{e, 1, 2, ..., 19\}, *, 3\} \cup \{[0, a] / a \in \{e, 1, 2, ..., 19\}, *, 2\} \cup \{[0, a] / a \in \{e, 1, 2, ..., 19\}, *, 9\}$ be a 5-interval loop. Clearly L is a S-5-interval loop.

Now using these examples we can have the following theorem which guarantees the existence of S-n-interval loops.

**THEOREM 2.4.2**: *Let $L = L_1 \cup L_2 \cup ... \cup L_n = \{[0, a] / a \in \{e, 1, 2, ..., m_1\}, *, p_1\} \cup \{[0,a] / a \in \{e, 1, 2, ..., m_2\}, *, p_2\} \cup ...$*

$\cup \{[0,a] \mid a \in \{e, 1, 2, ..., m_n\}, *, p_n\}$ be a n-interval loop. L is a S-n-interval loop.

The proof is direct for in each $L_i$, $A_i = \{[0,e], [0,a]\}$ is an interval subgroup.
Hence the theorem follows.

**COROLLARY 2.4.1**: *If $L = L_1 \cup L_2 \cup ... \cup L_n$ be a S-n-interval loop then every n-interval subloop of L is a S-n-interval subloop of L.*

Further we give examples of S-n-interval simple loops.

***Example 2.4.10***: Let $L = L_1 \cup L_2 \cup L_3 \cup L_4 = \{[0, a] / a \in \{e, 1, 2, ..., 23\}, *, 14\} \cup \{[0, a] / a \in \{e, 1, 2, ..., 29\}, *, 28\} \cup \{[0, a] / a \in \{e, 1, 2, ..., 31\}, *, 17\} \cup \{[0, a] / a \in \{e, 1, 2, ..., 37\}, *, 19\}$ be a 4-interval loop. It can be easily verified that L has no S-4-interval normal subloops. Hence L is a 4-interval S-simple loop.

In view of this we have a theorem which establishes the existence of a class of n-interval S-simple loops.

**THEOREM 2.4.3**: *Let $L = L_1 \cup L_2 \cup ... \cup L_n = \{[0, a] / a \in \{e, 1, 2, ..., m_1\}, *, p_1\} \cup \{[0, a] / a \in \{e, 1, ..., m_2\}, *, p_2\} \cup ... \cup \{[0, a] / a \in \{e, 1, ..., m_n\}, *, p_n\}$ be an n-interval loop. L is a n-interval S-simple loop.*

The proof is left as an exercise to the reader. Now we have seen a class of n-interval S-simple loops. Now we will give an example of a Smarandache n-interval subgroup-loop.

***Example 2.4.11***: Let $L = L_1 \cup L_2 \cup L_3 \cup L_4 \cup L_5 \cup L_6 = \{[0, a] / a \in \{e, 1, 2, ..., 23\}, *, 4\} \cup \{[0, a] / a \in \{e, 1, 2, ..., 29\}, *, 5\} \cup \{[0, a] / a \in \{e, 1, 2, ..., 17\}, *, 8\} \cup \{[0, a] / a \in \{e, 1, 2, ..., 47\}, *, 9\} \cup \{[0, a] / a \in \{e, 1, 2, ..., 53\}, *, 12\} \cup \{[0, a] / a \in \{e, 1, 2, ..., 41\}, *, 14\}$ be a 6-interval loop. Clearly L is a 6-interval S-subgroup loop.

Inview of this example we will now proceed onto give a theorem which gurantees the existence of n-interval S-subgroup loop.

**THEOREM 2.4.4**: *Let $L = L_1 \cup L_2 \cup ... \cup L_n = \{[0, a] / a \in \{e, 1, 2, ..., p_1\}, *, m_1\} \cup \{[0, a] / a \in \{e, 1, 2, ..., p_2\}, *, m_2\} \cup ... \cup \{[0, a] / a \in \{e, 1, 2, ..., p_n\}, *, m_n\}$ be a n-interval loop where $p_1, p_2, ..., p_n$ are n primes then L is a n-interval S-subgroup loop.*

The proof is straight forward and is left as an exercise for the reader to prove.
We will give an example of Smarandache Cauchy n-interval loop.

***Example 2.4.12***: Let $L = L_1 \cup L_2 \cup L_3 \cup ... \cup L_6 = \{[0, a] / a \in \{e, 1, 2, ..., 25\}, *, 7\} \cup \{[0, a] / a \in \{e, 1, 2, ..., 37\}, *, 11\} \cup \{[0, a] / a \in \{e, 1, 2, ..., 85\}, *, 7\} \cup \{[0, a] / a \in \{e, 1, 2, ..., 93\}, *, 17\} \cup \{[0, a] / a \in \{e, 1, 2, ..., 23\}, *, 10\} \cup \{[0, a] / a \in \{e, 1, 2, ..., 43\}, *, 27\}$ be a 6-interval loop. Clearly L is a 6-interval S-Cauchy loop.

Infact we have a class of n-interval S-Cauchy loops which is evident from the following theorem.

**THEOREM 2.4.5**: *Let $L = L_1 \cup L_2 \cup ... \cup L_n = \{[0, a] / a \in \{e, 1, 2, ..., m_1\}, *, t_1\} \cup \{[0, a] / a \in \{e, 1, 2, ..., m_2\}, *, t_2\} \cup ... \cup \{[0, a] / a \in \{e, 1, 2, ..., m_n\}, *, t_n\}$ be n-interval loop. L is a n-interval S-Cauchy loop.*

*Proof:* Every loop $L_i$ in L is of even order. Further $|L| = |L_1| ... |L_n| = 2^n M$, where

$$M = \left|\frac{L_1}{2}\right|\left|\frac{L_2}{2}\right|...\left|\frac{L_n}{2}\right| ; 1 \leq i \leq n.$$

Since every interval loop $L_i$ is a S interval loop and $2 / |L_i|$; $1 \leq i \leq n$, we see every interval loop $L_i$ is a S-Cauchy interval loop, hence L is a n-interval S-Cauchy loop.

Now we proceed onto describe by an example of a 4-n interval Lagrange loop.

***Example 2.4.13***: Let $L = L_1 \cup L_2 \cup L_3 \cup L_4 = \{[0, a] / a \in \{e, 1, 2, …, 23\}, *, 14\} \cup \{[0, a] / a \in \{e, 1, 2, …, 29\}, *, 14\} \cup \{[0, a] / a \in \{e, 1, 2, …, 37\}, *, 14\} \cup \{[0, a] / a \in \{e, 1, 2, …, 43\}, *, 14\}$ be a 4-interval loop. Clearly L is a 4-interval Lagrange loop.

Since each $L_i$ is an interval loop of order a prime number plus one. The only interval subgroups are of order 2, $1 \leq i \leq 4$. Thus L is a 4-interval S-Lagrange loop.

Inview of this we have the following theorem which gives a class of S-n-interval Lagrange loop.

**THEOREM 2.4.6**: *Let $G = G_1 \cup G_2 \cup … \cup G_n = \{[0, a] / a \in \{e, 1, 2, …, p_1\}, *, m_1\} \cup \{[0, a] / a \in \{e, 1, 2, …, p_2\}, *, m_2\} \cup … \cup \{[0, a] / a \in \{e, 1, 2, …, p_n\}, *, m_n\}$ be n-interval loop where $p_1, …, p_n$ are primes. Then G is a S-n-interval Lagrange loop.*

The proof is left as an exercise to the reader.

**THEOREM 2.4.7**: *Let $L = L_1 \cup L_2 \cup … \cup L_n = \{[0, a] / a \in \{e, 1, 2, …, m_1\}, t_1, *\} \cup \{[0, a] / a \in \{e, 1, 2, …, m_2\}, t_2, *\} \cup … \cup \{[0, a] / a \in \{e, 1, 2, …, m_n\}, t_n, *\}$ be a n-interval loop; $m_1, …, m_n$ are non primes. L is not a S-n-interval Lagrange loop but only a S-n-interval weakly Lagrange loop.*

This proof is also direct and hence left as an exercise to the reader. However we give examples of them.

***Example 2.4.14***: Let $L = L_1 \cup L_2 \cup L_3 \cup L_4 = \{[0, a] / a \in \{e, 1, 2, …, 15\}, *, 2\} \cup \{[0, a] / a \in \{e, 1, 2, …, 15\}, *, 8\} \cup \{[0, a] / a \in \{e, 1, 2, …, 15\}, *, 14\} \cup \{[0, a] / a \in \{e, 1, 2, …, 63\}, *, 14\}$ is a 4-interval loop, which can be easily checked to be a S-n-interval weakly Lagrange loop.

We next proceed onto show Smarandache strong n-interval p-Sylow loop by examples and give a theorem which gurantees the existence of such n-interval loops.

*Example 2.4.15*: Let $L = L_1 \cup L_2 \cup L_3 \cup L_4 = \{[0, a] / a \in \{e, 1, 2, ..., 13\}, *, 10\} \cup \{[0, a] / a \in \{e, 1, 2, ..., 43\}, *, 20\} \cup \{[0, a] / a \in \{e, 1, 2, ..., 23\}, *, 20\} \cup \{[0, a] / a \in \{e, 1, 2, ..., 19\}, *, 12\}$ be a 4-interval loop. It is easily verified L is a Smarandache strong 4-interval 2-Sylow loop.

In this we observe each of the interval loops $L_i$ given in example 2.4.14 is of order a prime + 1, ($p_i+1$). Thus we have the following theorem which proves the existence of a class of n-interval loops which are Smarandache strong 2-Sylow n-interval loops.

**THEOREM 2.4.8**: *Let $L = L_1 \cup L_2 \cup ... \cup L_n = \{[0, a] / a \in \{e, 1, 2, ..., p_1\}, *, t_1\} \cup \{[0, a] / a \in \{e, 1, 2, ..., p_2\}, *, t_2\} \cup \{[0, a] / a \in \{e, 1, 2, ..., p_3\}, *, t_3\} \cup ... \cup \{[0, a] / a \in \{e, 1, 2, ..., p_n\}, *, t_n\}$ be a n-interval loop where $p_1, p_2, ..., p_n$ are primes. Then L is a Smarandache strong n-interval 2-Sylow loop.*

Proof is direct and hence is left for the reader to prove [ ]. It is pertinent to mention here that the interval loop $L_i$ in L can be many for varying $t_i$ and for the fixed $p_i$. Thus we get a class $L_{p_i}$ ([0, a] ($t_i$)) for every i, i = 1, 2, ..., n. Hence we have a class of n-interval loops which are Smarandache strong n-interval 2-Sylow loops. Now we will give examples of Smarandache commutative n-interval loops.

*Example 2.4.16*: Let $L = L_1 \cup L_2 \cup L_3 \cup L_4 \cup L_5 = \{[0, a] / a \in \{e, 1, 2, ..., 21\}, *, 11\} \cup \{[0, a] / a \in \{e, 1, 2, ..., 15\}, *, 8\} \cup \{[0, a] / a \in \{e, 1, 2, ..., 25\}, *, 12\} \cup \{[0, a] / a \in \{e, 1, 2, ..., 27\}, *, 26\} \cup \{[0, a] / a \in \{e, 1, 2, ..., 9\}, *, 2\}$ be a 5-interval loop. Clearly L is a Smarandache commutative 5-interval loop.

We have a class of S-commutative n-interval loops.

**THEOREM 2.4.9**: *Let $L = L_1 \cup L_2 \cup L_3 \cup \ldots \cup L_n = \{[0, a] / a \in \{e, 1, 2, \ldots, m_1\}, *, t_1\} \cup \{[0, a] / a \in \{e, 1, 2, \ldots, m_2\}, *, t_2\} \cup \ldots \cup \{[0, a] / a \in \{e, 1, 2, \ldots, m_n\}, *, t_n\}$ be a n-interval loop. L is a Smarandache commutative n-interval loop.*

The proof is direct [5].

*Example 2.4.17*: Let $L = L_1 \cup L_2 \cup L_3 \cup L_4 \cup L_5 = \{[0, a] / a \in \{e, 1, 2, \ldots, 19\}, *, 8\} \cup \{[0, a] / a \in \{e, 1, 2, \ldots, 11\}, *, 7\} \cup \{[0, a] / a \in \{e, 1, 2, \ldots, 23\}, *, 12\} \cup \{[0, a] / a \in \{e, 1, 2, \ldots, 31\}, *, 13\} \cup \{[0, a] / a \in \{e, 1, 2, \ldots, 37\}, *, 20\}$ be a 5-interval loop. It is easily verified that L is a Smarandache strongly commutative 5-interval loop.

We now give a class of n-interval loops which are Smarandache strongly commutative n-interval loop. We first observe the entries of the interval loops given in example 2.4.17 is built using e and $Z_{p_i}$, $p_i$ are primes.

**THEOREM 2.4.10**: *Let $L = L_1 \cup L_2 \cup L_3 \cup \ldots \cup L_n = \{[0, a] / a \in \{e, 1, 2, \ldots, p_1\}, *, m_1\} \cup \{[0, a] / a \in \{e, 1, 2, \ldots, p_2\}, *, m_2\} \cup \ldots \cup \{[0, a] / a \in \{e, 1, 2, \ldots, p_n\}, *, m_n\}$ be a n-interval loop, where $p_1, p_2, \ldots, p_n$ are primes. L is a Smarandache strongly commutative n-interval loop.*

The proof easily follows from the result [5, 11-2].

Now we give examples of S-strongly n-cyclic loop.

*Example 2.4.18*: Let $L = L_1 \cup L_2 \cup L_3 \cup L_4 \cup L_5 = \{[0, a] / a \in \{e, 1, 2, \ldots, 19\}, *, 3\} \cup \{[0, a] / a \in \{e, 1, 2, \ldots, 23\}, *, 4\} \cup \{[0, a] / a \in \{e, 1, 2, \ldots, 29\}, *, 3\} \cup \{[0, a] / a \in \{e, 1, 2, \ldots, 31\}, *, 4\} \cup \{[0, a] / a \in \{e, 1, 2, \ldots, 37\}, *, 3\}$ be a 5-interval loop. Clearly L is a Smarandache strong 5-interval cyclic loop. All the loops are built using $Z_{p_i} \cup \{e\}$, $p_i$'s are prime, $1 \leq i \leq 4$.

We will illustrate this situation by a theorem. This theorem will prove the existence of a class of n-interval loops which are Smarandache strongly cyclic.

**THEOREM 2.4.11**: *Let $L = L_1 \cup L_2 \cup L_3 \cup ... \cup L_n = \{[0, a] / a \in \{Z_{p_1} \cup \{e\}\}, *, t_1\} \cup \{[0, a] / a \in \{Z_{p_2} \cup \{e\}\}, *, t_2\} \cup ... \cup \{[0, a] / a \in \{Z_{p_n} \cup \{e\}\}, *, t_n\}$ be a n-interval loop. L is a Smarandache strongly cyclic n-interval loop.*

For proof refer [5].

Now we can derive most of the results which are true for interval biloops in case of n-interval loops like S-right nucleus, S-left nucleus, S-center, S-Moufang center and so on.

We will give examples of n-interval loops in which Smarandache first normalizer is equal to Smarandache second normalizer in S-n-interval subloop.

***Example 2.4.19***: Let $L = L_1 \cup L_2 \cup L_3 \cup L_4 = \{[0, a] / a \in \{e, 1, 2, ..., 9\}, *, 8\} \cup \{[0, a] / a \in \{e, 1, 2, ..., 25\}, *, 7\} \cup \{[0, a] / a \in \{e, 1, 2, ..., 49\}, *, 9\} \cup \{[0, a] / a \in \{e, 1, 2, ..., 121\}, *, 3\}$ be a 4-interval loop. Clearly if we take $H = H_1 \cup H_2 \cup H_3 \cup H_4 = \{[0, a] / a \in \{e, 1, 4, 7\}, *, 8\} \cup \{[0, a] / a \in \{e, 1, 6, 11, 16, 21\}, *, 7\} \cup \{[0, a] / a \in \{e, 1, 8, 19, 40\}, *, 9\} \cup \{[0, a] / a \in \{e, 1, 12, 23, 34, 45, 56, 67, 78, 89, 100, 111\}, *, 3\} \subseteq L_1 \cup L_2 \cup L_3 \cup L_4$ is a 4-interval S-subloop of L.

We see $SN_1(H) = SN_2(H)$ that is $SN_1(H_1(3)) \cup SN_1(H_1(5)) \cup SN_1(H_1(7)) \cup SN_1(H_1(11)) = SN_2(H_1(3)) \cup SN_2(H_1(5)) \cup SN_2(H_1(7)) \cup SN_2(H_1(11))$.

In view of this we have the following theorem.

**THEOREM 2.4.12**: *Let $L = L_1 \cup L_2 \cup L_3 \cup ... \cup L_n = \{[0, a] / a \in \{e, 1, 2, ..., m_1\}, *, p_1\} \cup \{[0, a] / a \in \{e, 1, 2, ..., m_2\}, *, p_2\} \cup ... \cup \{[0, a] / a \in \{e, 1, 2, ..., m_n\}, *, p_n\}$ be a n-interval loop. Suppose $H = H_i^1(t_1) \cup H_i^2(t_2) \cup ... \cup H_i^n(t_n) \subseteq L_1 \cup L_2 \cup ... \cup L_n$ be a S-subloop of L; $SN_1(H) = SN_2(H)$ if and only if $(m_i^2 - m_i + 1, t_i) = (2m_i - i, t_i)$ for every i, i = 1, 2, ..., n.*

For analogous proof refer [5, 11-2].

***Example 2.4.20***: Let $L = L_1 \cup L_2 \cup L_3 \cup L_4 \cup L_5 = \{[0, a] / a \in \{e, 1, 2, ..., 7\}, *, 4\} \cup \{[0, a] / a \in \{e, 1, 2, ..., 13\}, *, 4\} \cup \{[0, a] / a \in \{e, 1, 2, ..., 19\}, *, 14\} \cup \{[0, a] / a \in \{e, 1, 2, ..., 23\}, *, 20\}$ be a 4-interval loop. Clearly SN (L) = $\{[0, e] \cup [0, e] \cup [0, e] \cup [0, e]\}$.

In view of this we have the following theorem. We first make a small observation each $L_i$ is built using $Z_p$, p a prime viz in this; example 2.4.20, 7, 13, 19 and 23 are primes, so that S(L) is the identity.

**THEOREM 2.4.13**: *Let $L = L_1 \cup L_2 \cup ... \cup L_n = \{[0, a] / a \in \{e, 1, 2, ..., p_1\}, *, t_1\} \cup \{[0, a] / a \in \{e, 1, 2, ..., p_2\}, *, t_2\} \cup ... \cup \{[0, a] / a \in \{e, 1, 2, ..., p_n\}, *, t_n\}$ be a n-interval loop. $p_i$'s are primes, $1 \leq i \leq n$. We see SN (L) = $[0, e] \cup ... \cup [0, e]$ (Infact we get a class of n-interval loop $\{L\} = \{L_1\} \cup ... \cup \{L_n\}$ as $L_i$ can be constructed using any $t_i$; $1 < t_i < p_i$, i=1, 2, ..., n).*

*Proof:* Given $p_1,..., p_n$ are prime. Hence L has no S-n- interval subloop but L is a S-n-interval loop.

Hence N (L) = $\{[0, e] \cup ... \cup [0, e]\}$ further we know SN (L) = N (L). Hence we have used analogous results for interval loops. SN (L) = $[0, e] \cup [0, e] \cup ... \cup [0, e]$.

***Example 2.4.21***: Let $L = L_1 \cup L_2 \cup L_3 = \{[0, a] / a \in \{e, 1, 2, ..., 23\}, *, 14\} \cup \{[0, a] / a \in \{e, 1, 2, ..., 19\}, *, 11\} \cup \{[0, a] / a \in \{e, 1, 2, ..., 13\}, *, 10\}$ be a 3-interval loop. It is easily verified S-Moufang centre of L is either $[0, e] \cup [0, e] \cup [0, e]$ or $L = L_1 \cup L_2 \cup L_3$.

We have the following theorem.

**THEOREM 2.4.14**: *Let $L = L_1 \cup L_2 \cup ... \cup L_n = \{[0, a] / a \in \{e, 1, 2, ..., p_1\}, *, t_1\} \cup \{[0, a] / a \in \{e, 1, 2, ..., p_2\}, *, t_2\} \cup ... \cup \{[0, a] / a \in \{e, 1, 2, ..., p_n\}, *, t_n\}$ be a n-interval loop where $p_1, p_2, ..., p_n$ are primes. Then S-Moufang n-center of L is $\{[0, e] \cup [0, e] \cup ... \cup [0, e]\}$ or $L_1 \cup L_2 \cup ... \cup L_n$.*

For analogous proof refer [5].

**THEOREM 2.4.15**: *Let $L = L_1 \cup L_2 \cup L_3 \cup L_4 \cup L_5 = \{[0, a] / a \in \{e, 1, 2, ..., p_1\}, *, t_1\} \cup \{[0, a] / a \in \{e, 1, 2, ..., p_2\}, *, t_2\} \cup ... \cup \{[0, a] / a \in \{e, 1, 2, ..., p_n\}, *, t_n\}$ be a n-interval loop where $p_1, ..., p_n$ are n primes. $NZ(L) = Z(L) = \{[0, e], [0, e], ..., [0, e]\}$.*

The proof is direct and hence left as an exercise to the interested reader.

All other properties enlisted in interval biloops can be analogously obtained for n-interval loops. It is left as an exercise to the reader.

Now we proceed onto define and describe quasi n-interval loops.

**DEFINITION 2.4.2**: *Let $L = L_1 \cup L_2 \cup ... \cup L_n$ where some $L_i$'s are interval loops and the rest are just loops. L inherits the operations from each $L_i$; $1 \leq i \leq n$, denoted by '.', (L, .) is defined as the quasi n-interval loop.*

We will illustrate this situation by some examples.

***Example 2.4.22***: Let $L = L_1 \cup L_2 \cup L_3 \cup L_4 \cup L_5 = \{L_9(8)\} \cup \{[0, a] / a \in \{e, 1, 2, ..., 11\}, *, 7\} \cup \{L_{11}(3)\} \cup \{[0, a] / a \in \{e, 1, 2, ..., 13\}, *, 9\} \cup \{[0, a] / a \in \{e, 1, 2, ..., 15\}, *, 8\}$ be a quasi 5-interval loop. Suppose $x = \{6 \cup [0, 10] \cup 7 \cup [0, 8] \cup [0, 12]\}$, $y = \{2 \cup [0, 7] \cup 3 \cup [0, 5] \cup [0, 10]\}$ are in L.

x.y =   $(6 * 2) \cup [(0, 10] * [0, 7]) \cup (7 * 3) \cup ([0, 8] * [0, 5]) \cup ([0, 12] * [0, 10])$
   =   $(16-42 \pmod 9) \cup \{[0, 49-60 \pmod{11})]\} \cup \{9 - 7 \times 2 \pmod{11}\} \cup \{[0, 45-64 \pmod{13})\} \cup \{[0, 80-84 \pmod{15})]\}$.
   =   $\{1\} \cup \{[0, 0]\} \cup \{6\} \cup \{[0, 3]\} \cup \{[0, 11]\} \in L$.

Thus L is a quasi 5-interval loop.

***Example 2.4.23***: Let $L = L_1 \cup L_2 \cup L_3 \cup L_4 = \{L_{23}(7)\} \cup \{[0, a] / a \in \{e, 1, 2, ..., 23\}, *, 20\} \cup \{L_{19}(6)\} \cup \{[0, a] / a \in$

{e, 1, 2, …, 19}, *, 12} be a quasi 4-interval loop of order $24^2\ 20^2$.

Clearly L is non commutative. Now we can define substructure in quasi n-interval loop, this task is left as an exercise to the reader. However we will provide this by some simple examples.

***Example 2.4.24***: Let $L = L_1 \cup L_2 \cup L_3 \cup L_4 \cup L_5 = \{[0, a] / a \in \{e, 1, 2, …, 15\}, *, 8\} \cup \{[0, a] / a \in \{e, 1, 2, …, 15\}, *, 14\} \cup L_{15}(2) \cup \{[0, a] / a \in \{e, 1, 2, …, 21\}, *, 11\} \cup L_{21}(5)$ be a quasi 5-interval loop. Consider $H = H_1 \cup H_2 \cup H_3 \cup H_4 \cup H_5 = \{[0, a] / a \in \{e, 1, 6, 11\}, *, 8\} \cup \{[0, a] / a \in \{e, 2, 7, 12\}, *, 14\} \cup \{[0, a] / a \in \{e, 1, 4, 7, 10, 13\}, *, 2\} \cup \{[0, a] / a \in \{e, 1, 8, 15\}, *, 11\} \cup \{\{e, 1, 4, 7, 10, 13, 16, 19\}, *, 5\} \subseteq L_1 \cup L_2 \cup L_3 \cup L_4 \cup L_5$ is a quasi 5-interval subloop of L.

***Example 2.4.25***: Let $L = L_1 \cup L_2 \cup L_3 \cup L_4 = \{L_{49}(10)\} \cup \{[0, a] / a \in \{e, 1, 2, …, 49\}, *, 12\} \cup \{[0, a] / a \in \{e, 1, 2, …, 121\}, *, 4\} \cup \{L_{121}(8)\}$ be a quasi 4-interval loop. Consider $H = H_1 \cup H_2 \cup H_3 \cup H_4 = \{\{e, 1, 8, 15, 22, 29, 36, 43\}, *, 10\} \cup \{[0, a] / a \in \{e, 2, 9, 16, 23, 30, 37, 44\}, *, 12\} \cup \{[0, a] / a \in \{e, 1, 12, 23, 34, 45, 56, 67, 78, 89, 100, 111\}, *, 4\} \cup \{\{e, 3, 14, 25, 36, 47, 58, 69, 80, 91, 102, 113\}, *, 8\} \subseteq L_1 \cup L_2 \cup L_3 \cup L_4$ is a quasi 4-interval subloop of L.

We say a quasi n-interval loop is a Smarandache quasi n-interval loop (S-quasi n-interval loop) if L has a proper subset $A = A_1 \cup A_2 \cup … \cup A_n \subseteq L_1 \cup L_2 \cup … \cup L_n$ such that A is a quasi n-interval group with respect to the operations on L.

We will illustrate first this by some examples.

***Example 2.4.26***: Let $L = L_1 \cup L_2 \cup L_3 \cup L_4 \cup L_5 = \{[0, a] / a \in \{e, 1, 2, …, 29\}, *, 8\} \cup L_{19}(8) \cup \{[0, a] / a \in \{e, 1, 2, …, 17\}, *, 3\} \cup L_{23}(9) \cup \{[0, a] / a \in \{e, 1, 2, …, 143\}, *, 15\}$ be a quasi 5-interval loop. $H = H_1 \cup H_2 \cup H_3 \cup H_4 \cup H_5 = \{[0, e], [0, 7], *, 8\} \cup \{e, 12, *\} \cup \{[0, e], [0, 10], *, 3\} \cup \{e, 15, *, 9\} \cup \{[0, e] [0, 141], *, 15\} \subseteq L_1 \cup L_2 \cup L_3 \cup L_4 \cup L_5$ is a S-quasi 5-interval loop, as H is a S-quasi 5-interval group.

*Example 2.4.27*: Let L = $L_1 \cup L_2 \cup L_3 \cup L_4 \cup L_5 \cup L_6$ = {[0, a] / a ∈ {e, 1, 2, ..., 47}, *, 10} $\cup L_{53}$ (2) $\cup L_{61}$ (5) $\cup$ {[0, a] / a ∈ {e, 1, 2, ..., 79}, *, 14} $\cup L_{19}$ (3) $\cup$ {[0, a] / a ∈ {e, 1, 2, ..., 101}, *, 42} be a quasi 6-interval loop. H = $H_1 \cup H_2 \cup ... \cup H_6$ = {[0, 1], [0, 46], *, 10} $\cup$ {e, 38, *, 2} $\cup$ {e, 48, *, 5} $\cup$ {[0, e], [0, 77], *, 14} $\cup$ {[e, 8, * 3} $\cup$ {[0, e] [0, 98], *, 42} $\subseteq L_1 \cup L_2 \cup ... \cup L_6$ is a quasi 6-interval group hence L is a S-quasi 6-interval loop.

Now having seen examples of S-quasi n-interval loops we now proceed onto give a class of S-quasi n-interval loops.

**THEOREM 2.4.16**: *Let L = $L_1 \cup L_2 \cup ... \cup L_n$ = $L_{m_1}$ ($t_1$) $\cup$ {[0, a] / a ∈ {e, 1, 2, ..., $m_2$}, *, $t_2$} $\cup$ { $L_{m_3}$ ($t_3$)} $\cup$ { $L_{m_4}$ ($t_4$)} $\cup$ ... $\cup$ {[0, a] / a ∈ {e, 1, 2, ..., $m_n$}, *, $t_n$} be a quasi n-interval loop. L is a S-quasi n-interval loop.*

Infact we get of class of such loops using $m_1, m_2, ..., m_n$ and appropriately varying $t_1, ..., t_n$ in $m_1, m_2, ..., m_n$ respectively.

The proof is easy and hence is left as an exercise to the reader.

*Example 2.4.28*: Let L = $L_1 \cup L_2 \cup L_3 \cup L_4 \cup L_5$ = {$L_9$ (5)} $\cup$ {$L_{15}$ (8)} $\cup$ {$L_{17}$ (3)} $\cup$ {[0, a] / a ∈ {e, 1, 2, ..., 11}, *, 4} $\cup$ {[0, a] / a ∈ {e, 1, 2, ..., 19}, *, 8} be a quasi 5-interval loop. We see A (L) = L.

For more refer [5].

Inview of this we have the following theorem which gives a class of loops {L} such that for each L in {L} are have A (L) = L.

**THEOREM 2.4.17**: *{L} = { $L_{m_1}$ } $\cup$ { $L_{m_2}$ } $\cup$ {{[0, a] / a ∈ {e, 1, ..., $m_i$}, *, $t_i$; $1 < t_i < m_i$}} $\cup$ ... $\cup$ { $L_{m_{i+r}}$ } $\cup$ ... $\cup$ {{[0, a] / a ∈ {e, 1, 2, ..., $m_n$}, *, $t_n$; $1 < t_n < m_n$}} be a class of quasi n-* *interval loops. Every quasi n-interval loop L in {L} is such that A (L) = L.*

Proof follows from the fact that every loop $L_i$ in L is such that $A(L_i) = L_i$ where $L_i$ is an interval loop or other wise.

***Example 2.4.29***: Let $L = L_1 \cup L_2 \cup L_3 \cup L_4 \cup L_5 = \{L_9 (5)\} \cup \{[0, a] / a \in \{e, 1, 2, …, 11\}, *, 6\} \cup \{L_{13} (7)\} \cup \{[0, a] / a \in \{e, 1, 2, …, 21\}, *, 11\} \cup L_{43} (22)$ be a quasi 5-interval loop.

We have the following theorem which guarantees the existence of quasi n-interval loop.

**THEOREM 2.4.18**: *Let $L = L_1 \cup L_2 \cup … \cup L_n = \{L_{m_1}\left(\frac{m_1+1}{2}\right)\} \cup \{L_{m_2}\left(\frac{m_2+1}{2}\right)\} \cup … \cup \{\{[0, a] / a \in \{e, 1, …, m_i\}, *, \left(\frac{m_i+1}{2}\right)\} \cup … \cup \{L_{m_n}\left(\frac{m_n+1}{2}\right)\}$ be a quasi n-interval loop. L is a commutative quasi n-interval loop.*

***Example 2.4.30***: Let $L = L_1 \cup L_2 \cup L_3 \cup L_4 \cup L_5 = \{L_{11} (3)\} \cup \{[0, a] / a \in \{e, 1, 2, …, 19\}, *, 3\} \cup L_{17} (8) \cup \{[0, a] / a \in \{e, 1, 2, …, 29\}, *, 24\} \cup L_{43} (8)$ be a quasi 5-interval loop. L has quasi 5-interval subloops. L is a Smarandache 5-interval subgroup loop.

**THEOREM 2.4.19**: *Let $L = L_1 \cup L_2 \cup … \cup L_n = \{L_{p_1}(m_1)\} \cup \{[0, a] / a \in \{e, 1, 2, ..., p_2\}, *, m_2\} \cup \{L_{p_3}(m_3)\} \cup \{[0, a] / a \in \{e, 1, 2, ..., p_4\}, *, m_4\} \cup \{[0, a] / a \in \{e, 1, 2, ..., p_5\}, *, m_5\} \cup … \cup \{L_{p_n}(m_n)\}$ be a quasi n-interval loop. L is a S-quasi n-interval subgroup loop.*

***Example 2.4.31***: Let $L = L_1 \cup L_2 \cup L_3 = \{L_7 (3)\} \cup \{[0, a] / a \in \{e, 1, 2, …, 15\}, *, 8\} \cup L_{25} (8)$ be a quasi 3-interval loop. It is easily verified L is a S-quasi 3-interval simple loop. (quasi 3-interval S-simple loop).

In view of this example we have the following theorem.

**THEOREM 2.4.20**: *Let $L = L_1 \cup L_2 \cup ... \cup L_n = \{L_{m_1}(t_1)\} \cup \{[0, a] / a \in \{e, 1, 2, ..., m_2\}, t_2, *\} \cup ... \cup \{L_{m_n}(t_n)\}$ be a quasi n-interval loop. L is a S-quasi n-interval simple loop (quasi n-interval S-simple loop).*

Proof is obvious.

*Example 2.4.32*: Let $L = L_1 \cup L_2 \cup L_3 \cup L_4 = \{L_9(5)\} \cup \{[0, a] / a \in \{e, 1, 2, ..., 17\}, *, 4\} \cup \{L_{15}(14)\} \cup L_{29}(3)$ be a quasi 4 interval loop. L is a S-quasi 4-interval Cauchy loop (quasi 4-interval S-Cauchy loop).

In view of this we have the following theorem.

**THEOREM 2.4.21**: *Let $L = L_1 \cup L_2 \cup ... \cup L_n = \{[0, a] / a \in \{e, 1, ..., m_1\}, *, t_1\} \cup \{L_{m_2}(t_2)\} \cup \{L_{m_3}(t_3)\} \cup ... \cup \{[0, a] / a \in \{e, 1, 2, ..., m_{n-1}\}, *, t_{n-1}\} \cup L_{m_n}(t_n)$ be a quasi n-interval loop. L is a S-Cauchy quasi n-interval loop.*

Proof follows from the fact order of each $L_i$ is of even order and $2^n/|L|$. Infact we get a class of such loops for appropriate $t_1$, ..., $t_n$.

*Example 2.4.33*: Let $L = L_1 \cup L_2 \cup L_3 \cup L_4 \cup L_5 = \{[0, a] / a \in \{e, 1, ..., 19\}, *, 8\} \cup L_{11}(9) \cup L_{12}(7) \cup \{[0, a] / a \in \{e, 1, ..., 23\}, *, 12\} \cup L_{43}(2)$ be a quasi 5-interval loop. L is a Smarandache Lagrange quasi 5-interval loop.

*Example 2.4.34*: Let $L = L_1 \cup L_2 \cup L_3 \cup L_4 = \{L_{15}(8)\} \cup \{[0, a] / a \in \{e, 1, ..., 25\}, *, 7\} \cup L_{49}(9) \cup \{[0, a] / a \in \{e, 1, ..., 21\}, *, 5\}$ be a quasi 4-interval loop. Clearly L is a Smarandache weakly Lagrange quasi 4-interval loop.

In view of these examples we have a class of S-weakly Lagrange quasi n-interval loops and S-Lagrange n-interval loop.

**THEOREM 2.4.22**: *Let $L = L_1 \cup L_2 \cup ... \cup L_n = \{[0, a] / a \in \{e, 1, ..., p_1\}, *, t_1\} \cup \{L_{p_2}(t_2)\} \cup \{L_{p_3}(t_3)\} \cup ... \cup \{[0, a] / a \in$*

$\{e, 1, 2, ..., p_{n-1}\}, *, t_{n-1}\} \cup L_{p_n}(t_n)$ *be a quasi n-interval loop, here $p_1, p_2, ..., p_n$ are n primes. L is a S-Lagrange quasi interval loop.*

Proof is direct hence left as an exercise to the reader.

**THEOREM 2.4.23**: *Let $L = L_1 \cup L_2 \cup ... \cup L_n = L_{m_1}(t_1) \cup \{[0, a] / a \in \{e, 1, ..., m_2\}, *, t_2\} \cup \{L_{m_3}(t_3)\} \cup ... \cup \{[0, a] / a \in \{e, 1, 2, ..., m_n\}, *, t_n\}$ be a quasi n-interval loop where $m_1, m_2, ..., m_n$ are not primes. Then L is only a S-weakly Lagrange quasi n-interval loop.*

***Example 2.4.35***: Let $L = L_1 \cup L_2 \cup L_3 \cup L_4 = \{L_{19}(8)\} \cup \{[0, a] / a \in \{e, 1, ..., 11\}, *, 3\} \cup L_{13}(9) \cup \{[0, a] / a \in \{e, 1, ..., 43\}, *, 29\}$ be a quasi 4-interval loop. Clearly L is a Smarandache strong quasi 4-interval 2-Sylow loop.

In view of this we have the following theorem the proof of which is easy. Only one thing to observe from example 2.4.35 is that all the $Z_i$'s used are prime fields of characteristic 19, 11, 13 and 29. Keeping this in mind we state the theorem.

**THEOREM 2.4.24**: *Let $L = L_1 \cup L_2 \cup ... \cup L_n = \{L_{p_1}(m_1)\} \cup \{[0, a] / a \in \{e, 1, ..., p_2\}, *, m_2\} \cup \{L_{p_3}(m_3)\} \cup ... \cup \{[0, a] / a \in \{e, 1, 2, ..., p_n\}, *, m_n\}$ be a quasi n-interval loop where $p_1, p_2, ..., p_n$ are primes. L is a Smarandache strongly quasi n-interval commutative loop.*

It is verified L mentioned in the above theorem is also a Smarandache strongly quasi n-interval cyclic loop.

Several other results obtained in case of n-interval loops can also be derived for quasi n-intervals. This task is left to the reader.

Now we proceed onto define n-interval loop - group.

**DEFINITION 2.4.3**: *Let $L = L_1 \cup L_2 \cup ... \cup L_n$ be such that some $L_i$'s are interval groups and rest are interval loops. L inherits* *the operations from $L_i.(1 \leq i \leq n)$. L is a n-interval group - loop.*

We will first illustrate this situation by some examples.

***Example 2.4.36***: Let $L = L_1 \cup L_2 \cup L_3 \cup L_4 = \{[0, a] / a \in Z_{28}, +\} \cup \{[0, a] / a \in \{e, 1, \ldots, 47\}, *, 9\} \cup \{[0, a] / a \in Z_{19} \setminus \{0\}, \times\} \cup \{[0, a] / a \in \{e, 1, \ldots, 43\}, *, 7\}$ be a 4-interval group-loop.

***Example 2.4.37***: Let $L = L_1 \cup L_2 \cup L_3 \cup L_4 \cup L_5 = \{[0, a] / a \in Z_{10}, +\} \cup \{[0, a] / a \in Z_{11} \setminus \{0\}, \times\} \cup \{[0, a] / a \in Z_{27}, +\} \cup \{[0, a] / a \in \{e, 1, 2, \ldots, 13\}, *, 7\} \cup \{[0, a] / a \in \{e, 1, 2, \ldots, 27\}, *, 11\}$ be a 5-interval group-loop.

We cannot talk of any of the usual properties studied for non associative structures. The only thing we can analyse is about substructures and Lagrange theorem, Cauchy element and Sylow theorems.

We will just give examples of substructures.

***Example 2.4.38***: Let $L = L_1 \cup L_2 \cup L_3 \cup L_4 \cup L_5 = \{[0, a] / a \in Z_{40}, +\} \cup \{[0, a] / a \in Z_{11} \setminus \{0\}, \times\} \cup \{[0, a] / a \in \{e, 1, 2, \ldots, 15\}, *, 8\} \cup \{[0, a] / a \in \{e, 1, 2, \ldots, 25\}, *, 12\} \cup \{[0, a] / a \in Z_{24}, +\}$ be a 5-interval group-loop.

Consider $H = H_1 \cup H_2 \cup H_3 \cup H_4 \cup H_5 = \{[0, a] / a \in \{0, 4, 8, 12, 16, \ldots, 36\}, +\} \cup \{[0, 1] [0, 11], \times\} \cup \{[0, a] / a \in \{e, 1, 6, 11\}, *, 8\} \cup \{[0, a] / a \in \{e, 1, 6, 11, 16, 21\}, *, 12\} \cup \{[0, a] / a \in \{0, 3, 6, 9, 12, 15, 18, 21\}, +\} \subseteq L_1 \cup L_2 \cup L_3 \cup L_4 \cup L_5$ is a 5-interval subgroup - subloop of L.

If L has no n-interval normal subgroup - normal subloop then we call L to be a n-interval S-simple group-loop.

***Example 2.4.39***: Let $L = L_1 \cup L_2 \cup L_3 \cup L_4 \cup L_5 = \{[0, a] / a \in Z_{11}, +\} \cup \{[0, a] / a \in \{e, 1, 2, \ldots, 11\}, *, 7\} \cup \{[0, a] / a \in \{e, 1, 2, \ldots, 13\}, *, 8\} \cup \{[0, a] / a \in Z_{17}, +\} \cup \{[0, a] / a \in \{e, 1, 2, \ldots, 23\}, *, 18\}$ be a 5-interval loop-group.

Clearly L is a 5-interval S-simple loop-group or S-simple 5 interval loop - group.

For we know very well the notion of Smarandache simple in case of group has no meaning.

***Example 2.4.40***: Let $L = L_1 \cup L_2 \cup L_3 \cup L_4 \cup L_5 = \{[0, a] / a \in Z_{23} \setminus \{0\}, \times\} \cup \{[0, a] / a \in Z_{47} \setminus \{0\}, \times\} \cup \{[0, a] / a \in \{e, 1, 2, …, 21\}, *, 11\} \cup \{[0, a] / a \in \{e, 1, 2, …, 15\}, *, 8\} \cup \{[0, a] / a \in Z_{37} \setminus \{0\}, \times\}$ be a 5-interval group - loop. Every element in L is a Cauchy element.

For o (L) $= |L_1| |L_2| |L_3| |L_4| |L_5|$
$= 22 \times 46 \times 22 \times 16 \times 36$.

In general we have Cauchy theorem to be true in case of n-interval loop-group for the special class of interval loops.

**THEOREM 2.4.25**: *Let $L = L_1 \cup L_2 \cup … \cup L_n$ be a n-interval loop - group where some $L_i$ are interval groups and the rest of the $L_j$'s are interval loops of the form $\{[0, a] / a \in \{e, 1, …, m_i\}, t_i, *\}$. Then every element in L is a Cauchy element.*

Proof is straight forward as every element in these type of interval loops are of order two and every interval loop of this type is of even order.

***Example 2.4.41***: Let $L = L_1 \cup L_2 \cup L_3 \cup L_4 \cup L_5 = \{[0, a] / a \in Z_{46}, +\} \cup \{[0, a] / a \in \{e, 1, 2, …, 25\}, *, 12\} \cup \{[0, a] / a \in Z_{13} \setminus \{0\}, \times\} \cup \{[0, a] / a \in \{e, 1, 2, …, 15\}, 8, \times\} \cup \{[0, a] / a \in Z_{20}, +\}$ be a 5-interval group - loop. Clearly L is a only S-weakly 5-interval group-loop.

***Example 2.4.42***: Let $L = L_1 \cup L_2 \cup L_3$ be a 3-interval group-loop where $L_1 = \{[0, a] / a \in \{e, 1, 2, …, 19, *, 14\}$, $L_2 = \{[0, a] / a \in Z_{24}, +\}$ and $L_3 = \{[0, a] / a \in Z_{23} \setminus \{0\}, \times\}$. It is easily verified that L is a S-Lagrange 3-interval group-loop.

When we say S-Lagrange the Smarandache qualifies only the loop and not the group as for groups Smarandache has no relevance.

We will now illustrate by some theorems which gurantees such n-interval group-loops.

**THEOREM 2.4.26**: *Let $L = L_1 \cup L_2 \cup \ldots \cup L_n$ where some $L_i$'s are interval groups of finite order and the rest of the $L_j$'s are interval loops of the form $\{[0, a] / a \in \{e, 1, 2, \ldots, m_j\}, *, t_j\}$ and ($m_j$'s are non prime odd numbers greater than three) be a n-interval group-loop. L is a Smarandache weakly Lagrange n-interval loop-group.*

The proof is direct and easy and hence is left as an exercise for the reader to prove.

**THEOREM 2.4.27**: *Let $L = L_1 \cup L_2 \cup \ldots \cup L_n$, where some $L_i$'s are interval groups of finite order and the rest of the $L_j$'s are interval loops of the form $\{[0, a] / a \in \{e, 1, 2, \ldots, p_j\}, t_j, *, 1 < t_j < p_j\}$, $p_j$'s are prime. L is a n-interval group-loop and L is a Smarandache Lagrange n-interval loop - group.*

Proof is direct hence left as an exercise.

We cannot prove p-Sylow theorems for n-interval group-loops. However only special p-Sylow theorems can be derived for these n-interval group-loops.

*Example 2.4.43*: Let $L = L_1 \cup L_2 \cup L_3 \cup L_4 = \{[0, a] / a \in Z_{24}, +\} \cup \{[0, a] / a \in \{e, 1, 2, \ldots, 19\}, *, 12\} \cup \{[0, a] / a \in \{e, 1, 2, \ldots, 13\}, *, 10\} \cup \{[0, a] / a \in Z_{19} \setminus \{0\}, \times\}$ be a 4-interval group-loop. Clearly L is a 4-interval Smarandache 2-Sylow loop-group.

In view of this we have the following theorem.

**THEOREM 2.4.28**: *Let $L = L_1 \cup L_2 \cup \ldots \cup L_n$ be a n-interval group - loop where some $L_i$'s are interval groups of finite order and the rest are interval loops of the form $\{[0, a] / a \in \{e, 1, 2, \ldots, p\} / *, t_i, 1 < t_i < p\}$, p a prime. L is a Smarandache strong 2-Sylow interval loop-group.*

Automatically interval groups are of finite order hence will be Sylow but the interval loops of only this type is a Smarandache strong 2-Sylow loops.

Similar results proved for n-interval loop can also be easily extended in case of n-interval loop-group. We see however all results cannot be extended for general interval loops.

Now having seen this type of n-interval loop-group we will now proceed onto define n-interval loop - semigroup.

**DEFINITION 2.4.4**: *Let $L = L_1 \cup L_2 \cup \ldots \cup L_n$ be such that some $L_i$'s are interval semigroups and the rest interval loops. The operation on L is inherited from each $L_i$, $1 \leq i \leq n$, L will be defined as n-interval loop-semigroup or n-interval semigroup-loop.*

We will illustrate this situation by some example.

*Example 2.4.44*: Let $L = L_1 \cup L_2 \cup L_3 \cup L_4 \cup L_5 = \{[0, a] / a \in Z_{24}, \times\} \cup \{[0, a] / a \in \{e, 1, 2, \ldots, 25\}, *, 8\} \cup \{[0, a] / a \in Z_{12}, \times\} \cup \{[0, a] / a \in \{e, 1, 2, \ldots, 27\}, *, 11\} \cup \{[0, a] / a \in Z_{40}, \times\}$ be a 5-interval semigroup-loop of finite order o (L) = $|L_1|$ $|L_2|$ $|L_3|$ $|L_4|$ $|L_5|$ = 24.25.12.27.40. Clearly L is non commutative.

*Example 2.4.45*: Let $L = L_1 \cup L_2 \cup L_3 = \{$ S (X) / X = $\{([0, a_1] [0, a_2] [0, a_3] [0, a_4])\} \cup \{[0, a] / a \in Z_{14}, \times\} \cup \{[0, a] / a \in Z_{42}, \times\}$ be a 3-interval loop semigroup of finite order and L is non commutative. Now o (L) = $4^4.14.42$.

Now as only one of the structure is associative and other is non associative we cannot proceed to arrive results for non associative structure.

We can define substructures, this task is left as an exercise.

*Example 2.4.46*: Let $L = L_1 \cup L_2 \cup L_3 \cup L_4 = \{[0, a] / a \in \{e, 1, 2, \ldots, 15\}, *, 8\} \cup \{[0, a] / a \in Z_{24}, \times\} \cup \{[0, a] / a \in \{e, 1, 2, \ldots, 25\}, *, 12\} \cup \{[0, a] / a \in Z_{40}, \times\}$ be a 4-interval loop semigroup. Consider $H = H_1 \cup H_2 \cup H_3 \cup H_4 = \{[0, a] / a \in \{e, 1, 4, 7, 10, 13\}, *, 8\} \cup \{[0, a] / a \in \{0, 2, 4, \ldots, 22\} \subseteq Z_{24}, \times\} \cup \{[0, a] / a \in \{e, 1, 6, 11, 16, 21\}, *, 12\} \cup \{[0, a] / a \in \{0, 10, 20, 30\}, \times\}$ be the 4-interval subloop - subsemigroup of L.

*Example 2.4.47*: Let $S = S_1 \cup S_2 \cup S_3 = \{[0, a] / a \in Z_{27}, \times\} \cup \{[0, a] / a \in Z^+ \cup \{0\}, \times\} \cup \{[0, a] / a \in \{e, 1, 2, \ldots, 11\}, *, 4\}$ be 3-interval loop-semigroup. Clearly S is of infinite order, but S-commutative and has 3-interval subloop-subsemigroup of only infinite order. Except when the subsemigroup of $S_2$ is taken as $\{0, 1\}$ under product. Thus S has only one finite 3-interval subloop - subsemigroup, but infinite number of infinite 3 interval subloop - subsemigroup. One more observation is infact the interval subloop is an interval group only.

Other properties like zero divisors, etc associated with interval semigroup cannot be studied as interval loops do not have such properties associated with it.

Next we proceed onto study describe and define n-interval loop-groupoid.

Let $G = G_1 \cup G_2 \cup \ldots \cup G_n$ be a n interval in which some $G_i$'s are interval loops and the rest are interval groupoids. G inherits operations of every $G_i$ done component wise denoted by '.' (G, .) is defined as the n-interval loop-groupoid.

It is nice to observe that both the algebraic structures are non associative hence they have several common properties enjoyed by them which will be discussed. We will illustrate this situation by some examples.

*Example 2.4.48*: Let $G = G_1 \cup G_2 \cup G_3 \cup G_4 \cup G_5 \cup G_6$ be a 6-interval loop groupoid where $G_1 = \{[0, a] / a \in Z_7, *, (3, 2)\}$, $G_2 = \{[0, a] / a \in \{e, 1, 2, \ldots, 7\}, *, 4\}$ $G_3 = \{[0, a] / a \in Z_{12}, *, (3, 7)\}$, $G_4 = \{[0, a] / a \in \{e, 1, 2, \ldots, 19\}, *, 10\}$, $G_5 = \{[0, a] / a \in Z_{18}, *, (0, 7)\}$ and $G_6 = \{[0, a] / a \in \{e, 1, 2, \ldots, 23\}, *, 7\}$. Clearly G is of finite order. $|G| = |G_1| |G_2| \ldots |G_6|$
$= 7 \times 8 \times 12 \times 20 \times 18 \times 24$.

*Example 2.4.49*: Let $L = L_1 \cup L_2 \cup L_3 \cup L_4 = \{[0, a] / a \in Z_{25}, *, (7, 0)\} \cup \{[0, a] / a \in \{e, 1, 2, \ldots, 101\}, *, 6\} \cup \{[0, a] / a \in Z_{32}, *\} \cup \{[0, a] / a \in \{e, 1, 2, \ldots, 43\}, *, 7\}$ be a 4-interval loop-groupoid of finite order which is clearly non commutative.

Now having seen examples of n interval loop-groupoids we now leave the task of defining substructures to the reader, but we give some examples of substructures in them.

***Example 2.4.50***: Let $L = L_1 \cup L_2 \cup L_3 \cup L_4 \cup L_5 = \{[0, a] / a \in \{e, 1, 2, \ldots, 15\}, *, 8\} \cup \{[0, a] / a \in \{e, 1, 2, \ldots, 25\}, *, 9\} \cup \{[0, a] / a \in \{e, 1, 2, \ldots, 21\}, *, 11\} \cup \{[0, a] / a \in \{0, 2, 4, 6, 8, 10\}, *, (3, 0)\} \cup \{[0, a] / a \in \{e, 3, 6, 9, 12\}, *, (0, 2)\} \subseteq L_1 \cup L_2 \cup L_3 \cup L_4 \cup L_5$. H is a 5-interval subgroupoid subloop. Both L and H are non commutative and is of finite order.

***Example 2.4.51***: Let $L = L_1 \cup L_2 \cup L_3 \cup L_4 = \{[0, a] / a \in Z_7, +\} \cup \{[0, a] / a \in \{e, 1, 2, \ldots, 11\}, *, 10\} \cup \{[0, a] / a \in Z_{13}, +\} \cup \{[0, a] / a \in \{e, 1, 2, \ldots, 17\}, *, 12\}$ be a 4-interval groupoid-loop.

L has no 4-interval subgroupoid-subloop. Thus L is a simple 4-interval groupoid-loop.

Now we will study some of the special identities satisfied n-interval groupoid-loops for this we need to know about Smarandache n-interval groupoid-loops.

We will call a n-interval groupoid-loop $G = G_1 \cup G_2 \cup \ldots \cup G_n$ to be Smarandache n-interval groupoid - loop if G contains a n-interval group-semigroup. It is interesting to note that in general all n-interval loop-groupoids are not Smarandache n-interval loop-groupoids.

We will illustrate this situation by some examples.

***Example 2.4.52***: Let $G = G_1 \cup G_2 \cup G_3 \cup G_4 = \{[0, a] / a \in \{e, 1, 2, \ldots, 11\}, *, 4\} \cup \{[0, a] / a \in \{e, 1, 2, \ldots, 13\}, *, 5\} \cup \{[0, a] / a \in Z_{12}, *, (2, 0)\} \cup \{[0, a] / a \in \{e, 1, 2, \ldots, 17\}, *, 10\}$ be a 4-interval loop-groupoid. G is a Smarandache 4-interval loop-groupoid.

***Example 2.4.53***: Let $L = L_1 \cup L_2 \cup L_3 \cup L_4 \cup L_5 = \{[0, a] / a \in Z_{15}, *, (3, 6)\} \cup \{[0, a] / a \in \{e, 1, 2, \ldots, 15\}, *, 8\} \cup \{[0, a] / a \in \{e, 1, 2, \ldots, 17\}, *, 3\} \cup \{[0, a] / a \in \{e, 1, 2, \ldots, 25\}, *, 9\} \cup \{[0, a] / a \in Z_{12}, *, (2, 1)\}$ be a 5-interval loop-groupoid. Clearly L is a Smarandache 5-interval loop-groupoid.

Now having seen examples of S-n-interval loop-groupoid interested reader can construct examples of n-interval loop-groupoids which satisfies special identities.

***Example 2.4.54***: Let $L = L_1 \cup L_2 \cup L_3 \cup L_4 \cup L_5 = \{[0, a] / a \in \{e, 1, 2, \ldots, 41\}, *, 20\} \cup \{[0, a] / a \in \{e, 1, 2, \ldots, 43\}, *, 38\} \cup \{[0, a] / a \in \{e, 1, 2, \ldots, 29\}, *, 14\} \cup \{[0, a] / a \in Z_{25}, *, (3, 0)\} \cup \{[0, a] / a \in Z_{30}, *, (4, 7)\}$ be a n-interval groupoid-loop. L is not S-5-interval Bol groupoid - loop. L is also not a Smarandache 5-interval Moufang loop-groupoid.

## 2.5 n-Interval Mixed Algebraic Structures

In this section we define the new notion of n-interval mixed algebraic structures and describe a few properties associated with them. These structures are so unique that they are mixture of associative and non associative algebraic structures.

**DEFINITION 2.5.1**: *Let $M = M_1 \cup M_2 \cup \ldots \cup M_n$ be a n-interval set, where some $M_i$'s are interval loops, some $M_j$'s are interval groupoids, some $M_k$'s are interval groups and the rest are interval semigroups $(1 \leq i, j, k < n)$; M obtains the operation '.' which is the componentwise operation of every $M_i$; $i=1, 2, \ldots n$. (M, .) is defined as the n-interval mixed algebraic structure or mixed n-interval algebraic structure.*

We will first accept even if three algebraic structures are involved then also M is a mixed n-interval algebraic structure. Thus what we demand is more than two algebraic structure must be present in the mixed n-interval algebraic structure. We see the n-interval algebraic structures discussed in earlier sections are not mixed n-interval algebraic structures as they do not contains more than two structures.

We will illustrate this situation by some examples.

***Example 2.5.1***: Let $G = G_1 \cup G_2 \cup G_3 \cup G_4 \cup G_5 = \{[0, a] / a \in Z_{20}, \times\} \cup \{[0, a] / a \in Z_{19} \setminus \{0\}, \times\} \cup \{[0, a] / a \in Z_{11}, +\} \cup$

$\{[0, a] / a \in Z_{40}, *(3, 1)\} \cup \{[0, a] / a \in \{e, 1, 2, \ldots, 27\} *, 8\}$ be a mixed 5-interval algebraic structure.

***Example 2.5.2***: Let $G = G_1 \cup G_2 \cup G_3 \cup G_4 = \{[0, a] / a \in Z_{40}, \times\} \cup \{[0, a] / a \in Z_{15}, +\} \cup \{[0, a] / a \in \{e, 1, \ldots, 29\}, *, 8\} \cup \{[0, a] / a \in Z_{12}, \times\}$ be mixed 4-interval algebraic structure.

***Example 2.5.3***: Let $G = L_1 \cup L_2 \cup L_3 = \{[0, a] / a \in Z_7, +\} \cup \{[0, a] / a \in Z_{16}, \times\} \cup \{[0, a] / a \in Z_7, *, (3, 2)\}$ be a mixed 3-interval algebraic structure.

***Example 2.5.4***: Let $G = L_1 \cup L_2 \cup L_3 \cup L_4 \cup L_5 \cup L_6 = \{[0, a] / a \in Z_{64}, \times\} \cup \{$All $3 \times 3$ interval matrices with intervals of the form $[0, a]$ where $a \in Z_{11}$ under matrix addition$\} \cup \{$All $1 \times 7$ row interval matrices with intervals of the form $[0, a]$ where $a \in Z_{11}, *, (3, 7)\} \cup \{[0, a] / a \in \{e, 1, 2, \ldots, 45\}, *, 17\} \cup \{\sum_{i=0}^{10}[0,a]x^i / a \in Z^+ \cup \{0\}$ under polylnomial addition$\} \cup \{\sum_{i=0}^{\infty}[0,a]x^i / a \in Z_{12}$, under polynomial multiplication$\}$ be a 6-interval mixed algebraic structure.

Clearly L is of infinite order and is non commutative. We can only define substructure of one type as this is a mixed n-interval algebraic structure. We leave the task of defining the substructure to the reader, but give examples of the same.

***Example 2.5.5***: Let $G = G_1 \cup G_2 \cup G_3 \cup G_4 = \{[0, a] / a \in \{e, 1, 2, \ldots, 27\}, *, 8\} \cup \{[0, a] / a \in Z_{12}, \times\} \cup \{[0, a] / a \in Z_{45}, +\} \cup \{[0, a] / a \in Z_{14}, *, (2, 0)\}$ be a 4-interval mixed algebraic structure. Consider $H = H_1 \cup H_2 \cup H_3 \cup H_4 = \{[0, a] / a \in \{e, 1, 4, 7, 10, 13, 16, 19, 22, 25\}, *, 8\} \cup \{[0, a] / a \in \{0, 2, 4, 6, 8, 10\}, \times\} \cup \{[0, a] / a \in \{0, 5, 10, 15, 20, 25, 30, 35, 40\}, +\} \cup \{[0, a] / a \in \{0, 2, 4, 6, 8, 10, 12\}, *, (2, 0)\} \subseteq G_1 \cup G_2 \cup G_3 \cup G_4$, H is a mixed 4-interval algebraic substructure of G.

We cannot define identities or other algebraic properties as they are mixed.

Now the following examples will show the way these mixed structures can be used in n-models.

*Example 2.5.6*: Let $L = L_1 \cup L_2 \cup L_3 \cup L_4 = \{$All $5 \times 3$ interval matrices with intervals of the form [0, a] where $a \in Z_{12}$, *, (5, 7)$\} \cup \{$All $4 \times 7$ interval matrices with intervals of the form [0, a] where $a \in Z_{20}$ under usual matrix addition$\} \cup \{$All $4 \times 4$ interval matrices with intervals of the form [0, a] where $a \in \{e, 1, 2, \ldots, 19\}$, *, 12$\} \cup \{$All $7 \times 5$ interval matrices with intervals of the form [0, a] where $a \in Z_{40}$ under matrix addition$\}$ be a mixed 4-interval algebraic structures which is non commutative and is of finite order.

*Example 2.5.7*: Let $L = L_1 \cup L_2 \cup L_3 = \{$All $5 \times 5$ interval matrices with intervals of the form [0, a] with $a \in Z_{42}$ under matrix addition$\} \cup \{$All $3 \times 7$ interval matrices with intervals of the form [0, a] with $a \in Z_{18}$ under matrix addition$\} \cup \{$All $6 \times 2$ interval matrices with intervals of the form [0, a] with $a \in Z_{27}$ with operation *, (3, 8)$\}$ be mixed 3-interval algebraic structures. This sort of mixed 3-interval matrices can be used in mathematical models.

We will just show how operations on this interval matrices are carried out on each of interval groupoids, interval semigroups, interval groups and interval loops.

Suppose B, A is a m × n interval matrix from an m × n interval matrix groupoid $G = \{([0, a_{ij}]) \,/\, a_{ij} \in Z_m, *, (t, u), t, u \in Z_m\}$

$$A = ([0, a_{ij}]) \quad 1 \leq i \leq m, \ 1 \leq j \leq n$$

and

$$B = ([0, b_{ij}]) \quad 1 \leq i \leq m, \ 1 \leq j \leq n$$

$$\begin{aligned} A*B &= ([0, a_{ij}] * [0, b_{ij}]) \\ &= ([0, ta_{ij} + u\, b_{ij} \,(\text{mod } m)]); \end{aligned}$$

\* in general is non associative. We follow the same operation even if m = n.

Suppose S is an interval matrix semigroup.

We can have two operations. If A, B ∈ S are two interval m × n matrices m ≠ n then A + B, + the usual matrix addition is the only operation on S. If A, B ∈ S are such that the interval matrices are square matrices then we can have usual interval matrix addition 'or' usual interval matrix multiplication.

In both cases this S has a semigroup structure. 'or' used only in the mutually exclusive sense.

Suppose L denotes the set of all interval m × n matrices with intervals of the form [0, a] with entries from {{e, 1, 2, …, m} m > 3, m odd, t, \*, where 1 < t < m with (t, m) = (t – 1, m) = 1 and for [0, a], [0, b] in this set.

$$[0, a] * [0, b] = [0, tb - (t-1) a \pmod m)].$$

Now
$$A = ([0, a_{ij}]) \quad 1 \leq i \leq m \quad 1 \leq j \leq n$$
and
$$B = ([0, b_{ij}]) \text{ in } L,$$
then
$$A * B = ([0, (tb_{ij} - (t-1) a_{ij}) \bmod m]).$$

This is the only operation L which is compatible and even if m = n we have only this operation.

If G is the set of all m × n interval matrices then usual interval matrix addition is the operation. If m ≠ n multiplication cannot be adopted for interval matrices. If m = n then the interval matrix can have either addition or multiplication.

Thus we have these mixed n-interval structures to function as dynamical systems in the mathematical modeling.

**Chapter Three**

# APPLICATIONS OF INTERVAL STRUCTURES AND N-INTERVAL STRUCTURES

These n-interval structures and biinterval structures are newly introduced in this book. A few of the probable applications are mentioned in this chapter.

1. Some of these interval structures can be adopted in finite interval analysis.

2. If the experts wants more than one approximate solution or if the experts are interested in appropriate solution with some possible flexibility these interval structures can be used. For the experts can have the liberty to choose the appropriate solution from the interval than forcefully accepting the approximate solution.

3. The introduction of n-interval matrix loops, (groupoids or groups or semigroups) gives the expert liberty to choose the suitable n-interval algebraic structure depending on the experiment / study. So non associative n-interval matrix algebraic structure would be a boon to them.

4. The n-interval m × s matrices with the operations defined on them can be used in interval stiffness matrix, which will yield an interval solution.

5. With the advent of computers several such models can be studied simultaneously using these n-interval matrix algebraic structures or mixed n-interval matrix algebraic structure.

**Chapter Four**

# SUGGESTED PROBLEMS

In this chapter we suggest around 295 problems for the reader some of which, are simple and some of them are difficult and some at research level.

1.  Obtain some interesting properties about interval bisemigroups.

2.  Let $S = S_1 \cup S_2$ be a finite interval bisemigroup.
    a.  Can every interval bisubsemigroup divide the order of S?
    b.  Does a interval bisemigroup have in general interval bisubsemigroups which are not interval biideals?

3.  Give an example of a non commutative interval bisemigroup of finite order.

4.  Give an example of a commutative interval bisemigroup.

5.  Does there exist a cyclic interval bisemigroup?

6.  Is $S = S_1 \cup S_2 = \{[0, a] / a \in Z^+ \cup \{0\}\} \cup \{S <X> / X = ([0, a_1], [0, a_2], [0, a_3])\}$ a finite interval bisemigroup?
    a. Find interval subbisemigroups in S.
    b. Prove S has interval biideal!
    c. Does S contain interval bisubsemigroups which are not ideals?

7.  Does there exist an interval bisemigroup which has no interval biideals?

8.  Does there exist an interval bisemigroup of prime order? Justify your answer.

9.  Does there exist an interval bisemigroup in which every interval bisubsemigroup is an interval biideal?

10. Find only left biideals of $S = S_1 \cup S_2 = \{S(\langle X \rangle) / X = ([0, a_1], [0, a_2], [0, a_3])\} \cup \{S(\langle Y \rangle)$ where $Y = \{([0, b_1], \ldots, [0, b_8])\}$ the interval bisemigroup (symmetric bisemigroup).
    a. Show every left biideal in general is not a right
    b. biideal.
    c. Can S have two sided interval biideals?
    d. Can S have interval bisubsemigroups which are not interval biideals?

11. Give an example of an interval bisemigroup which has no interval bizero divisors.

12. Does there exist interval bisemigroups which are interval idempotent bisemigroups?

13. Give an example of an interval bisemigroup which has no non trivial idempotents.

14. Give an example of a bisemigroup which has non trivial zero divisors.

15. Give an example of an interval bisemigroup which has no zero divisors and no idempotents.

16. Let $S = S_1 \cup S_2 = \{[0, a] / a \in Z_{450}, \times\} \cup \{[0, b] / b \in Z_{41}, \times\}$ be an interval bisemigroup.
    a. Can S have S-zero divisors?
    b. Is S a S-interval bisemigroup?
    c. Can S have S-interval biideals?
    d. Can S have S-interval bisubsemigroups?
    e. Give at least an interval biideal which is not an S-interval biideal.
    f. Can S have S-idempotent?

17. Give an example of a interval bisemigroup which is not an S-interval bisemigroup.

18. Prove the class of special symmetric biinterval semigroups are always S-interval semigroups.

19. Can $S = S_1 \cup S_2 = \{[0, a] / a \in Z_{19}, +\} \cup \{[0, b] / b \in Z_{43}, +\}$ the interval bisemigroup be a S-interval bisemigroup?

20. Show if in problem (19) + is replaced by × S is a S-interval bisemigroup.

21. Let $S = S_1 \cup S_2 = \{[0, a] / a \in Z_{420}, \times\} \cup \{[0, b] / b \in Z_{240}, \times\}$ be an interval bisemigroup.
    a. Find all bizero divisors of S.
    b. Does S have S-bizero divisors?
    c. Find all biidempotents of S.
    d. Can S have S-biidempotents?
    e. Find atleast 5 biideals of S.
    f. Can S have S-biideals?
    g. Can S have S-interval bisemigroups which are not S-ideals?
    h. Can S have S-interval bisubsemigroups?

      i. Can S have binilpotents?

22. Let $S = S_1 \cup S_2$ be a quasi interval bisemigroup. Obtain some properties enjoyed by these algebraic structures which are not enjoyed by the interval bisemigroups.

23. Let $S = S_1 \cup S_2 = \{Z_{30}, \times\} \cup \{[0, a] \,/\, a \in Z^+ \cup \{0\}\}$ be a quasi interval bisemigroup.
    a. Is S a S-quasi interval bisemigroup?
    b. Can S have S-interval biideals?
    c. Can S have S-interval bisubsemigroups?
    d. Find quasi interval bisubsemigroup in S.
    e. Find quasi interval biideals in S which are not S-ideals?
    f. Can S have S-zero divisors?
    g. Can S have S-idempotents?

24. Let S and G be any two interval bisemigroups. Define a interval bisemigroup homomorphism from S to G. Is it always possible to define bikernel of an interval bisemigroup homomorphism? Justify!

25. Let $S = S_1 \cup S_2 = \{[0, a] \,/\, a \in Z^+ \cup \{0\}, \times\} \cup \{[0, b] \,/\, b \in Z_{244}, \times\}$ be an interval bisemigroup.
    a. Define $\eta : S \to S$ so that $\eta$ is an interval bisemigroup homomorphism such that $\eta$ has a nontrivial bikernel.
    b. Define $\eta : S \to S$ so that $\eta$ is one to one but different from identity bihomomorphism.
    c. Is biker $\eta = \ker \eta_1 \cup \ker \eta_2$ where $\eta_i : S_i \to S_i$, $1 \le i \le 2$, an interval biideal of S. Illustrate this situation by an example.

26. Let $S = S_1 \cup S_2 = \{[0, a] \,/\, a \in 3Z^+ \cup \{0\}, \times\} \cup \{[0, b] \,/\, b \in 5Z^+ \cup \{0\}, \times\}$ be an interval bisemigroup under multiplication.
    a. Is S a S-interval bisemigroup?
    b. Does S have S-interval biideals?
    c. Can S have S-bizero divisors?

27. Let $G = G_1 \cup G_2 = S(X) \cup \{[0, a] / a \in Z_{43}, \times\}$ where $X = \{([0, a_1], [0, a_2], \ldots, [0, a_{11}])\}$ be an interval bisemigroup.
   a. Find the order of G.
   b. Does the biorder of every interval bisubsemigroup divide the biorder of the interval bisemigroup?
   c. Does G have S-interval biideals?
   d. Is G a S-interval bisemigroup?
   e. Can S have biideals which are not S-biideals?

28. Obtain some interesting properties related with the interval bisemigroup.
   (Properties like S-Lagrange, S-p-Sylow, S-Cauchy,…)

29. Does there exist an interval bisemigroup in which every element is S-Cauchy?

30. What is the marked difference between a S-interval bisemigroup and interval bisemigroup?

31. Let $G = G_1 \cup G_2$ be any S-interval matrix bisemigroup. Can S have S biideals?

32. Let $G = G_1 \cup G_2 = \left\{ \begin{bmatrix} [0,a] & [0,b] \\ [0,c] & [0,d] \end{bmatrix} \mid a, b, c, d \in Z^+ \cup \{0\}, \times \right\}$

   $\cup \left\{ \begin{bmatrix} [0,a] \\ [0,b] \\ [0,c] \\ [0,d] \end{bmatrix} \mid a, b, c, d \in Z^+ \cup \{0\}, + \right\}$ be an interval

   bisemigroup.
   a. Does S have S-interval biideals?
   b. Can S have S-zero divisors?
   c. Can S have interval bisemigroups which are not S-
   d. interval bisemigroups?

33. Let $S = S_1 \cup S_2 = \{([0, a_1], [0, a_2], \ldots, [0, a_{11}]) / a_i \in Z_{27}, \times\} \cup$
$$\left\{ \begin{bmatrix} [0,a_1] & [0,a_2] & [0,a_3] & [0,a_4] \\ [0,a_5] & [0,a_6] & [0,a_7] & [0,a_8] \end{bmatrix} \mid a_i \in Z_{45}, +, 1 \le i \le 8 \right\}$$ be an interval bisemigroup.
   a. Find the biorder of S.
   b. Can S satisfy the modified form of S-Lagrange theorem?
   c. Can S have S-biideals?
   d. Can S have S-zero divisors?

34. Let $G = G_1 \cup G_2 = \left\{ \begin{bmatrix} [0,a_1] & [0,a_2] \\ [0,a_3] & [0,a_4] \\ [0,a_5] & [0,a_6] \\ [0,a_7] & [0,a_8] \end{bmatrix} \mid a_i \in Z^+ \cup \{0\}, + \right\} \cup$
$$\left\{ \begin{bmatrix} [0,a_1] & [0,a_2] & \ldots & [0,a_9] \\ [0,a_{10}] & [0,a_{11}] & \ldots & [0,a_{18}] \end{bmatrix} \mid a_i \in Z^+ \cup \{0\}, + \right\}$$ be an interval matrix bisemigroup.
   a. Show G is of infinite order.
   b. Can G have S-biideals?
   c. Can G have S-zero divisors or just zero divisors?

35. Let $G = G_1 \cup G_2 = $ {Set of all $5 \times 5$ interval matrices with intervals of the form $[0, a_i]$ where $a_i \in Z_5$} $\cup$ {set of all $6 \times 6$ interval matrices with intervals of the form $[0, a_i]$ where $a_i \in Z_4$} be an interval matrix bisemigroup.
   a. Find the biorder of G.
   b. Is $P = P_1 \cup P_2 = $ {set of all $5 \times 5$ interval diagonal matrices with intervals of the form $[0, a]$ where $a \in Z_5$} $\cup$ {set of all diagonal $4 \times 4$ interval matrices with intervals of the form $[0, a_i]$ where $a_i \in Z_4$} $\subseteq G_1 \cup G_2$ an S biideal of G?

c. If diagonal interval matrices in P are replaced by upper triangular interval matrices will that structure be an interval biideal? Justify.
    d. Can G have S-bizerodivisors?
    e. Can G have biidempotents?
    f. Can G have S-binilpotents?
    g. Can G have S-Cauchy bielements?

36. Give some interesting properties enjoyed by interval matrix bisemigroups which are not in general true for interval bisemigroups.

37. Let $S = S_1 \cup S_2 = \left\{ \sum_{i=0}^{8}[0,a]x^i \big/ a \in Z_9, + \right\} \cup \left\{ \sum_{i=0}^{12}[0,b]x^i \big/ b \in Z_{12}, + \right\}$ be an interval polynomial bisemigroup.
    a. Find biideals if any in S.
    b. What is the biorder of S?
    c. Does S contain any Cauchy bielement?
    d. Is S a S-interval polynomial bisemigroup?
    e. Can S have S-biidempotents?

38. Let $G = \left\{ \sum_{i=0}^{\infty}[0,a]x^i \big/ a \in Z_{40}, '\times' \right\} \cup \left\{ \sum_{i=0}^{\infty}[0,b]x^i \big/ b \in Z_{25}, '\times' \right\}$ be an interval bisemigroup.
    a. Prove G is of infinite biorder.
    b. Find some interval biideals in G.
    c. Can G have S-biideals?
    d. Can G have S-bizero divisors?

39. Let $S = S_1 \cup S_2 = \left\{ \sum_{i=0}^{5} [0,a]x^i \,\middle/\, a \in Z^+ \cup \{0\}, + \right\} \cup \left\{ \sum_{i=0}^{\infty} [0,b]x^i \,\middle/\, b \in Z_{12} \right\}$ be an interval polynomial bisemigroup.
   a. Find the biorder of S.
   b. Find biideals and S-biideals in S.
   c. Find interval subbisemigroups which are not biideals of S.

40. Let $S = S_1 \cup S_2 = \{[0, a] \mid a \in Z_{40}, \times\} \cup \left\{ \sum_{i=0}^{3} [0,a]x^i \,\middle/\, a \in Z_{12}, '+' \right\}$ be an interval bisemigroup.
   a. Find the biorder of S.
   b. Does S satisfy the S-Lagrange theorem for
   c. bisemigroups?
   d. Can S satisfy the S-Cauchy theorem?
   e. Can S have S-bizero divisors?

41. Let $S = S_1 \cup S_2 = \left\{ \begin{bmatrix} [0,a] \\ [0,b] \\ [0,c] \\ [0,d] \end{bmatrix} \,\middle|\, a, b, c, d \in Z_{14}, + \right\} \cup \{[0, a] \,/\, a \in Z_{12}, \times\}$ be an interval bisemigroup.
   a. What is the order of S?
   b. Can S have S interval bisubsemigroups?
   c. Find S-bizero divisors if any in S.
   d. Can S have S-biidempotents?

42. Let $G = G_1 \cup G_2 = \{\text{All } 8 \times 8 \text{ interval matrices with intervals of the form } [0, a] \,/\, a \in Z_2, \times\} \cup \left\{ \sum_{i=0}^{9} [0,a_i]x^i \,/\, a_i \in Z_2, + \right\}$ be an interval bisemigroup.
   a. Prove G is of finite biorder.
   b. Find the biorder of G.

c. Can G have S-zero divisors?
   d. Can G have biideals which are not S-biideals?
   e. Can G have bizero divisors which are not S-bizero divisors?

43. Let $G = G_1 \cup G_2 = \{[0, a] / a \in Z_{12}, \times\} \cup \{$All $7 \times 1$ interval column matrices with entries from $Z_3\}$ be an interval bisemigroup.
   a. What is the biorder of G?
   b. Find atleast 4-interval sub bisemigroups.
   c. Is G a S-interval bisemigroup?
   d. Can G have S-interval biideals? Justify.

44. Let $S = S_1 \cup S_2 = \left\{ \sum_{i=0}^{\infty} [0,a]x^i \,\middle/\, a \in Z_2 \right\} \cup \{[0, a] / a \in Z_7\}$ be an interval bisemigroup. Enumerate the properties enjoyed by S.

45. Give an example of a bisimple interval bisemigroup.

46. Give an example of a S-bisimple interval bisemigroup.

47. Is a S-bisimple interval bisemigroup bisimple? Justify.

48. Let $S = S_1 \cup S_2 = \{([0, a], [0, b]) / a, b \in Z_{15}, \times\} \cup \{([0, a], [0, b], [0, c], [0, d]) / a, b, c, d \in Z_{12}\}$ be an interval semigroup.
   a. What is the order of S?
   b. Find bizero divisor and S-bizero divisor in S?
   c. Can S have S-biideals?
   d. Is every biideal a S-biideal? Justify!

49. Does there exist an interval bisemigroup of order 47? Justify!

50. Give an example of an interval bisemigroup of order 1048.

51. What are the special properties enjoyed by interval bigroupoids?

52. Give conditions under which an interval bigroupoid G is a S-interval bigroupoid.

53. Give an example of an interval bigroupoid which is not a S-interval bigroupoid.

54. Let $G = G_1 \cup G_2 = \{[0, a] / a \in Z_{12}, (3, 2), *\} \cup \{[0, b] / b \in Z_9, *, (8, 0)\}$ be an interval bigroupoid.
    a. What is the biorder of G?
    b. Find S-interval bisubgroupoid in G.
    c. Is every interval bisubgroupoid in G a S-interval
    d. subbigroupoid?

55. Let $G = \{[0, a] / a \in Z_7, *, (1, 4)\} \cup \{[0, b] / b \in Z_{11}, *, (2, 3)\}$ be an interval bigroupoid.
    a. Find the biorder of G.
    b. Does G have interval subbigroupoid?
    c. Can G have bizero divisors?
    d. Is G commutative?
    e. Is G a S-interval bigroupoid?
    f. Does the biorder of interval subbigroupoid divide
    g. biorder of G?

56. Let $G = G_1 \cup G_2 = \{[0, a] / a \in Z^+ \cup \{0\}, *, (3, 2)\} \cup \{[0, a] / a \in 3Z^+ \cup \{0\}, *, (0, 30)\}$ be an interval bigroupoid.
    a. Prove G is of infinite biorder.
    b. Find interval bisubgroupoids if any in G.
    c. Does G have S-interval subbigroupoids?
    d. Is G a S-interval bigroupoid?
    e. Can G satisfy any of the well known identities?

57. Let $G = G_1 \cup G_2 = \{[0, a] / a \in Z_{41}, *, (3, 9)\} \cup \{[0, b] / b \in Z_{43}, *, (3, 9)\}$ be an interval bigroupoid.

a. Find biorder of G.
b. Can G have interval bisubgroupoids?
c. Is G a S-interval bigroupoid?
d. Is G a S-interval P-bigroupoid?
e. Is G commutative?
f. Can G have S-interval subbigroupoids?

58. Let G = {[0, a] / a ∈ R⁺ ∪ {0}, *, (8, 3)} ∪ {[0, a] / a ∈ Q⁺ ∪ {0}, *, (7, 11)} = G₁ ∪ G₂ be an interval bigroupoid.
    a. Prove G is infinite.
    b. Enumerate all properties enjoyed by G.
    c. Is G a S-interval bigroupoid?
    d. Can G have S-interval bisubgroupoids?

59. Let G = G₁ ∪ G₂ = {[0, a] / a ∈ Z₁₉, (11, 9), *} ∪ {[0, b] / b ∈ Z₁₁, (7, 5), *} be an interval bigroupoid. Find all the special properties enjoyed by G.

60. Does there exist an interval bigroupoid which is not an S-interval bigroupoid?

61. Does there exists an interval bigroupoid in which all interval bisubgroupoids are S-interval subbigroupoids?

62. Give an example of an interval Bol bigroupoid.

63. Give an example of an S-strong interval Bol bigroupoid.

64. Give an example of a S-Bol interval bigroupoid.

65. Give an example of an interval bigroupoid of order 218.

66. Does there exist an interval bigroupoid of order 23? Justify your answer.

67. Let G = G₁ ∪ G₂ = {[0, a] / a ∈ Z₁₂, *, (2, 6)} ∪ {[0, b] / b ∈ Z₂₀, *, (10, 2)} be an interval bigroupoid. Can G have an interval subbigroupoid H so that o (H) / o (G) ?

68. Let $G = G_1 \cup G_2 = \{[0, a] / a \in Z_{19}, *, (3, 2)\} \cup \{[0, b] / b \in Z_{43}, *, (11, 3)\}$ be an interval bigroupoid. Does G contain an interval subbigroupoid H such that o(H) / o(G)?

69. Give an example of an interval bigroupoid which has no S-interval subbigroupoids.

70. Give an example of an interval bigroupoid which has no interval normal subbigroupoids.

71. Give an example of an interval bigroupoid which is an interval normal bigroupoid.

72. Give an example of an interval bigroupoid which is a Moufang interval bigroupoid.

73. Give an example of an interval bigroupoid G in which every interval subbigroupoid is a S-Moufang interval bigroupoid but G is not a Moufang interval bigroupoid.

74. Suppose $G = G_1 \cup G_2$ be an interval bigroupoid of order p. q, where p and q are two distinct primes, Can G be an interval alternative bigroupoid?

75. Determine some important and interesting properties enjoyed by interval bigroupoids of order pq where p and q are primes.

76. Let $G = G_1 \cup G_2 = \{[0, a] / a \in Z_{18}, *, (3, 0)\} \cup \{[0, b] / b \in Z_{36}, *, (0, 5)\}$ be an interval bigroupoid. What are the special properties enjoyed by G?

77. Let $G = G_1 \cup G_2 = \{[0, a] / a \in Z_{43}, *, (7, 0)\} \cup \{[0, b] / b \in Z_{47}, *, (0, 7)\}$ be an interval bigroupoid.
    a. Is G left alternative?
    b. Can G be alternative?
    c. Is G a P-interval bigroupoid?
    d. Can G be a S-interval bigroupoid?
        Justify all your claims.

78. Give an example of an interval bigroupoid which has no zero divisors.

79. Give an example of an interval bigroupoid which has every element to be an idempotent.

80. Give an example of an interval bigroupoid which has biidentity.

81. Does there exist an interval bigroupoid in which every interval subbigroupoid is an interval bisemigroup?

82. Does there exist an interval bigroupoid in which no interval subbigroupoid is an interval bisemigroup? If that is the case what is the speciality of that interval bigroupoid?

83. Define the notion of finite biorder of elements in an interval bigroupoid and illustrate them with an example.

84. Let $G = \{[0, a] / a \in Z_{20}, *, (3, 3)\} \cup \{[0, b] / b \in Z_{13}, *, (10, 10)\}$ be an interval bigroupoid.
    a. Is G an interval bisemigroup?
    b. Can G has non associative triples?
    c. Can G have interval subbigroupoids?
    d. Is G a S-interval bigroupoid?
    e. Does G satisfy any of the special identities?

85. Give an example of an interval bigroupoid which is right alternative but not left alternative.

86. Give an example of a commutative interval bigroupoid.

87. Give an example of an S-interval inner commutative bigroupoid.

88. Give an example of an interval P-bigroupoid.

89. Prove there exist an infinite class of interval bigroupoids of finite order.

90. Prove there does not exist an interval bigroupoid of prime order.

91. Let $G = G_1 \cup G_2 = \{[0, a] / a \in Z_{11}, *, (3, 2)\} \cup \{[0, b] / b \in Z_{42}, *, (7, 4)\}$ and $H = H_1 \cup H_2 = \{[0, a] / a \in Z_{11}, *, (5, 7)\} \cup \{[0, b] / b \in Z_{42}, *, (3, 2)\}$ be two interval bigroupoids.
    (i) How many distinct homomorphisms from G to H exist?
    (ii) Can ever G be isomorphic with H?

92. Let $G = G_1 \cup G_2 = \{[0, a] / a \in 3Z^+, *, (3, 2)\} \cup \{[0, b] / b \in 5Z^+, *, (3, 2)\}$. Is G an interval bigroupoid? Prove your claim.

93. Let $G = G_1 \cup G_2 = \{[0, a] / a \in Z_{49}, *, (7, 0)\} \cup \{[0, b] / b \in Z_{25}, *, (5, 0)\}$ be an interval bigroupoid. Does G enjoy any special property?

94. Let $G = G_1 \cup G_2 = \{[0, a] / a \in Z_7, *, (3, 2)\} \cup \{[0, b] / b \in Z_7, *, (2, 3)\}$ be an interval bigroupoid. Does G have any special property associated with it?

95. Let $G = G_1 \cup G_2 = \{[0, a] / a \in Z_{24}, *, (11, 13)\} \cup \{a / a \in Z_{40}, *, (7, 3)\}$ be an interval quasi bigroupoid.
    (i) Is G a S-quasi interval bigroupoid?
    (ii) Can G have S-quasi interval subbigroupoids?

96. Let $G = G_1 \cup G_2 = \{[0, a] / a \in Z^+ \cup \{0\}, *, (3, 2)\} \cup \{x / x \in Q^+ \cup \{0\}, *, (3, 29)\}$ be a quasi interval bigroupoid.
    a. Can G be a S-quasi interval bigroupoid?
    b. Does G satisfy any of the special identities?
    c. Does G contain quasi interval subbigroupoids?
    d. Give any of the special features enjoyed by G.

97. Let $G = G_1 \cup G_2 = \{[0, a] / a \in Z_9, *, (4, 4)\} \cup \{a / a \in R, *, (\sqrt{3}, 2)\}$ be a quasi interval bigroupoid.
   a. Does G have S-quasi interval S subbigroupoids?
   b. Can G satisfy any one of the special identities?

98. Let $G = G_1 \cup G_2 = \{([0, a_1], [0, a_2], [0, a_3], [0, a_4]) / a_i \in Z_{14}, *, (3, 7), 1 \leq i \leq 4\} \cup \{[0, b] / b \in Z_{43}, *, (3, 7)\}$ be an interval bigroupoid. Mention atleast two special features satisfied by G.

99. Let $G = G_1 \cup G_2 = \{\sum_{i=0}^{\infty}[0,a]x^i / a \in Z_{11}, *, (3, 2)\} \cup \{[0, a] / a \in Z_3, *, (1, 2)\}$ be an interval bigroupoid.
   (In $G_1$; $([0, a] x^i) * ([0, b] x^j) = [0, a * b] x^{i+j} = [0, 3a + 2b \pmod{11}] x^{i+j}$ extended for any sum).
   a. Is G a S-interval bigroupoid?
   b. Can G have interval subbigroupoids?
   c. Does G satisfy any special identities?

100. Let $G = G_1 \cup G_2$ be an interval semigroup-groupoid. Analyse the properties specially enjoyed by G.

101. Let $G = G_1 \cup G_2 = \{[0, a] / a \in Z_{40}, \times\} \cup \{[0, b] / b \in Z_{40}, *, (7, 11)\}$ be an interval semigroup-groupoid.
   a. Is G a S-interval semigroup-groupoid?
   b. Can G have interval subsemigroup-subgroupoid H such that o (H) / o (G)?

102. Give an example of a S-interval semigroup-groupoid.

103. Give an example of an interval semigroup-groupoid which is not a S-interval semigroup-groupoid.

104. Does there exist an interval semigroup-groupoid in which every interval subsemigroup-subgroupoid is a S-interval subsemigroup-subgroupoid?

105. Let $G = G_1 \cup G_2 = \{[0, a] / a \in Z_{14}, *, (3, 2)\} \cup \{S(X)\}$ be an interval groupoid - semigroup where $X = ([0, a_1], [0, a_2], \ldots, [0, a_7])$;
   a. Find the biorder of G.
   b. Is G a S-interval groupoid-semigroup?
   c. Does G have S-interval subgroupoid-subsemigroup?
   d. Does G contain interval bisubsemigroup?

106. Let $G = G_1 \cup G_2 = \{ \sum_{i=0}^{8}[0,a_i]x^i \mid a_i \in Z_{15}, +\} \cup \{[0, a] / a \in Z_{15}, *, (2, 7)\}$ be an interval semigroup-groupoid.
   a. What is the biorder of G?
   b. Can G have S-interval subsemigroup-subgroupoid?

107. Let $G = G_1 \cup G_2 = \{[0, a] / a \in Z_{16}, *, (3, 7)\} \cup \{[0, a] / a \in Z_{16}, \times\}$ be an interval semigroup-groupoid. Enumerate some of the special properties enjoyed by G.

108. Let $G = G_1 \cup G_2 = Z_{12} (3, 7) \cup \{[0, a] / a \in Z_{12}, \times\}$ be a quasi interval groupoid-semigroup.
   a. Is G a S-quasi interval groupoid-semigroup?
   b. What is the biorder of G?
   c. Does G have any proper bistructures?

109. Let $G = G_1 \cup G_2 = \{[0, a] / a \in \{e, 1, 2, \ldots, 23\}, *, 3\} \cup \{[0, a] / a \in \{e, 1, 2, \ldots, 11\}, *, 3\}$ be a biinterval loop or interval biloop.
   a. What is the biorder of G?
   b. Is G a S-interval biloop?
   c. Can G have S-interval subbiloop?
   d. Is G S-bisimple?
   e. Does G satisfy any one of the special identities?
   f. Does G satisfy in particular right alternative condition?

110. Let $G = G_1 \cup G_2 = \{[0, a] / a \in \{e, 1, 2, \ldots, 43\} *, 7\} \cup \{[0, b] / b \in \{e, 1, 2, 3, 4, \ldots, 13\}, *, 7\}$ be an interval biloop.
   a. What is the biorder of G?
   b. Is G commutative?

      c. Is G S-Lagrange?
      d. Is G a S-2-Sylow?
      e. Can G have substructures other than biorder 4? Justify all your claims.

111. Give an example of an interval biloop which is not a S-interval biloop.

112. Let G = $G_1 \cup G_2$ be an interval biloop where $G_1$ = {[0, a] / a ∈ {e, 1, 2, …, 45}, (23), *} and $G_2$ = {[0, b] / b ∈ {e, 1, 2, …, 55}, *, 13};
      a. Find all interval bisubloops of G.
      b. Prove all interval bisubloops of G are also S-interval bisubloops.
      c. Prove G is also an S-interval biloop.
      d. Does G satisfy any one of the special identities?

113. Let G = $G_1 \cup G_2$ be an interval biloop where $G_1 \cup G_2$ = {[0, a] / a ∈ {e, 1, 2, …, 19}, *, 18} ∪ {[0, b] / b ∈ {e, 1, 2, *, 43}, *, 12}.
      a. Does G satisfy any one of the special identities?
      b. Prove G has no S-interval bisubloops or interval bisubloops.
      c. Can G be interval biloop which is bisimple?
    Justify your answers.

114. Let G = $G_1 \cup G_2$ = {[0, a] / a ∈ {e, 1, 2, …, 47}, *, 12} ∪ $L_{25}(7)$ be a quasi interval biloop.
      a. Is G a S-quasi interval biloop?
      b. Is G a S-strong Moufang quasi interval biloop?

115. Give an example of a Moufang interval biloop.

116. Give an example of a right alternative interval biloop.

117. Give an example of a Jordan interval biloop.

118. Give an example of an interval Burck biloop.

119. Give an example of an interval P-biloop.

120. Give an example of an S-interval Moufang loop.

121. Give an example of a S-interval strongly cyclic biloop.

122. Is $L = L_1 \cup L_2 = \{[0, a] / a \in \{e, 1, 2, \ldots, 13\}, *, 10\} \cup \{[0, b] / b \in \{e, 1, 2, \ldots, 17\}, *, 9\}$ an interval biloop? Is L a S-strongly interval biloop?

123. Give an example of a S-strongly commutative interval biloop.

124. Give an example of a S-pseudo commutative interval biloop.

125. For the interval biloop $L = L_1 \cup L_2 = \{[0, a] / a \in \{e, 1, 2, \ldots, 29\}, *, 12\} \cup \{[0, b] / b \in \{e, 1, 2, \ldots, 31\}, *, 15\}$ find its principal biisotope.
Does it preserve the properties enjoyed by $L = L_1 \cup L_2$?

126. Let $L = L_1 \cup L_2$ be an interval biloop. Find $Z(L) = Z(L_1) \cup Z(L_2)$ where $L = \{[0, a] / a \in \{e, 1, 2, \ldots, 23\}, *, 22\} \cup \{[0, a] / a \in \{e, 1, 2, \ldots, 29\}, *, 28\}$.

127. Let $G = G_1 \cup G_2$ be an interval biloop, where $G_1 = \{[0, a] / a \in \{e, 1, 2, \ldots, 47\}, *, 43\}$ and $G_2 = \{[0, a] / a \in \{e, 1, 2, \ldots, 43\}, *, 42\}$. Is G a S-Lagrange interval biloop?

128. Let $G = G_1 \cup G_2 = L_7(4) \cup \{[0, a] / a \in \{e, 1, 2, \ldots, 7\}, *, 5\}$ be a quasi interval biloop. Find the right regular birepresentation of L.

129. Obtain some interesting properties enjoyed by interval biloops.

130. Obtain the special properties related with S-interval biloops.

131. Obtain some interesting properties related with quasi interval biloops.

132. Let $G = G_1 \cup G_2 = L_{23}(9) \cup \{[0, a] / a \in Z_{40}, \times\}$ be a quasi interval loop - semigroup.
   a. Find the biorder of G.
   b. Find substructure of G.

133. Let $G = G_1 \cup G_2 = \{[0, a] / a \in \{e, 1, 2, \ldots, 27\}, 11\} \cup \{[0, b] / b \in Z_{42}, \times\}$ be an interval loop - semigroup.
   a. Find interval subloop-subsemigroups of G.
   b. Is G a S-interval loop-semigroup?
   c. Does G have S-interval subloop-subsemigroup?

134. Let $G = G_1 \cup G_2 = L_{17}(3) \cup \{[0, a] / a \in \{e, 1, 2, \ldots, 43\}, *, 2\}$ be a quasi interval biloop.
   a. Find $S(N(G)) = S(N(G_1)) \cup S(N(G_2))$.
   b. Is $SZ(G) = SZ(G_1) \cup SZ(G_2)$?

135. Does there exist an interval biloop whose Smarandache binucleus is empty?

136. Let $L = L_1 \cup L_2 = \{[0, a] / a \in \{e, 1, 2, \ldots, 15\}, *, 2\} \cup \{[0, b] / b \in \{e, 1, 2, \ldots, 45\}, *, 8\}$ be an interval biloop.
   a. Find $SZ(L) = SZ(L_1) \cup SZ(L_2)$.
   b. Is $SN(L) = SN(L_1) \cup SN(L_2)$?
   c. Determine an interval S-bisubloop $A = A_1 \cup A_2$ of L and find $SN_1(A)$ and $SN_2(A)$.

137. Is $G = G_1 \cup G_2 = \{[0, a] / a \in Z_6, *, (4, 3)\} \cup \{[0, b] / b \in Z_6, *, (3, 5)\}$ the interval bigroupoid a S-P-interval bigroupoid?

138. Give an example of a S-strong P-interval bigroupoid.

139. Give an example of a S-Bol interval bigroupoid.

140. Give an example of a S-strong right alternative interval bigroupoid.

141. Give an example of a S-strong Moufang interval loop-groupoid.

142. Give an example of a S-Bol interval loop-groupoid.

143. Characterize those interval bigroupoids which are S-strong Bol bigroupoids.

144. Give some special properties enjoyed by S-strong Moufang biloops.

145. Enumerate the properties enjoyed by S-strong idempotent interval bigroupoid.

146. Does there exist a S-strong Bol interval matrix bigroupoid? Illustrate your claim if it exists.

147. Does there exists a S-strong Moufang interval polynomial bigroupoid?

148. Give an example of a S-right alternative interval matrix bigroupoid.

149. Give an example of a S-strong quasi interval P-bigroupoid.

150. Prove $G = G_1 \cup G_2 = \{[0, a] / a \in Z_n, *, (t, u)\} \cup \{[0, b] / b \in Z_m, *, (r, s)\}$ is a S-alternative interval bigroupoid if and only if $t^2 = t \pmod{n}$, $u^2 = u \pmod{n}$, $s^2 = s \pmod{m}$ and $r^2 = r \pmod{m}$ with $t + u = 1 \pmod{n}$ and $s + r = 1 \pmod{m}$. Illustrate this situation by an example.

151. Is the interval bigroupoid given in problem, (150) a S-strong Bol interval bigroupoid?

152. Is the interval bigroupoid $G = G_1 \cup G_2 = \{[0, a] / a \in Z_{2p}, *, (1, 2)\} \cup \{[0, b] / b \in Z_{2q}, *, (1, 2)\}$, p and q two distinct primes a S-interval bigroupoid ?

153. Give an example of a S-interval bigroupoid which is not a S-interval P-bigroupoid.

154. Give an example of an interval bigroupoid which is not a S-strong Bol interval bigroupoid.

155. Give an example of a S-interval bigroupoid which is not an S-interval idempotent bigroupoid.

156. Let $G = G_1 \cup G_2 = \{[0, a] / a \in Z_{20}, +\} \cup \{[0, b] / b \in Z_{42}, +\}$ be an interval bigroup.
    a. Verify whether Lagrange theorem for finite groups is satisfied by G.
    b. Find all interval subbigroups of G.
    c. Verify p-Sylow theorems for G.

157. Let $G = G_1 \cup G_2 = \{S_X, X = ([0, a_1], \ldots, [0, a_7])\} \cup \{[0, a] / a \in Z_{19} \setminus \{0\}, \times\}$ be an interval bigroup.
    a. Prove G is non commutative.
    b. Find normal interval bisubgroups of G.
    c. Is Cauchy theorem for groups satisfied by G?

158. Obtain some interesting properties enjoyed by interval groups.

159. Can Lagrange theorem in general be true for interval bigroups of finite biorder?

160. Give an example of an interval bigroup which is bicyclic.

161. Can all the Sylow theorem be true in case of finite interval bigroups? Justify your answer.

162. Obtain the classical homomorphism theorems in case of interval bigroups.

163. Determine the conditions under which we can extend Cayleys theorem in case of interval bigroups.

164. Enumerate all classical properties which are not true in case of interval bigroups.

165. Define interval bigroup automorphisms. Is the collection of all interval bigroup automorphisms a bigroup? Justify your claim.

166. Can we derive all properties associated with cosets in case of groups in case of interval bigroup?

167. Define permutation birepresentation of an interval bigroup.

168. Suppose $G = \{[0, a] / a \in Z_p, +\} \cup \{[0, b] / b \in Z_q, +\}$, p and q primes be an interval bigroup. What is the special property enjoyed by it?

169. Define interval matrix bigroup. Illustrate it by an example.

170. Let $G = G_1 \cup G_2 = S_{20} \cup \{[0, a] / a \in Z_{20}, \times\}$ be a quasi interval bigroup.
   a. What is the order of G?
   b. Find atleast 3 quasi bisubgroups of G.
   c. Can G have more than one quasi interval normal bisubgroup?
   d. Is Lagrange theorem for finite groups true in case of the quasi interval bigroup G?

171. Let $G = G_1 \cup G_2 = <g / g^{26} = 1> \cup \{[0, a] / a \in Z_{13} \setminus \{0\}, \times\}$ be a quasi interval bigroup.
   a. Find the order of G.
   b. Prove G has quasi interval subbigroups and all of them are binormal in G.

172. Let $G = S_5 \cup S_x$ where $X = \{[0, a_1], [0, a_2], [0, a_3]\}\}$ be the quasi interval bigroup.
   a. Find the biorder of G.
   b. Find all quasi interval normal subbigroups of G.
   c. Find the (p, q)-Sylow quasi interval bisubgroups of G ( p = 5 and q = 3).

173. Let $G = G_1 \cup G_2 = \{S_{10}\} \cup \left\{ \begin{bmatrix} [0,a] & [0,b] \\ [0,c] & [0,d] \end{bmatrix} \mid a, b, c, d \in Z_{20} \right.$

with $[0, ad] - [0, bc] \neq [0, 0]$ that is $[0, (ad-bc) \pmod{20}] \neq [0, 0]$.

    a. Is G a commutative quasi interval bigroup?
    b. Find the biorder of G.
    c. Can G have quasi interval normal bisubgroups?
    d. Find atleast 5 distinct quasi interval bisubgroups.
    e. What is the biorder of them? Will their biorder divide the biorder of G?

174. Let $G = G_1 \cup G_2 = \{A = \begin{bmatrix} a & b & c \\ d & e & f \\ g & h & i \end{bmatrix} \mid \det A \neq 0;$ a, b, c, d, e,

f, g, h, i $\in Z_{29}\} \cup \{A = \begin{bmatrix} [0,a_1] & [0,a_2] & [0,a_3] \\ [0,a_4] & [0,a_5] & [0,a_6] \\ [0,a_7] & [0,a_8] & [0,a_9] \end{bmatrix} \mid \det A \neq 0,$

$a_i \in Z_{19}, 1 \leq i \leq 9\}$. Is G a quasi interval bigroup?

175. Does their exist a quasi interval bigroup G of finite biorder which has a quasi interval subbigroup whose biorder does not divide the biorder of G? Justify your answer.

176. Let $G = G_1 \cup G_2 = \{[0, a] / a \in Z_{26}, +\} \cup \{[0, b] / b \in Z_{26}, \times\}$ be an interval group-semigroup.
    a. What is the biorder of G?
    b. Can Lagrange theorem for finite groups be true in case of this G?
    c. Find atleast 2 interval subgroup-subsemigroup of G.

177. Let $G = G_1 \cup G_2 = \{$All $3 \times 3$ interval matrices with intervals of the form $[0, a]$ where $a \in Z^+ \cup \{0\}\} \cup \left\{ \begin{bmatrix} [0,a_1] \\ [0,a_2] \\ [0,a_3] \\ [0,a_4] \end{bmatrix} \mid a_i \in Z_{20}, +, 1 \leq i \leq 4 \right\}$ be an interval semigroup-group. Prove $G$ has infinite number of interval subsemigroup-subgroup.

178. Let $G = G_1 \cup G_2 = \{ \sum_{i=0}^{7} [0,a]x^i / a \in Z_{15}, +\} \cup \{ \sum_{i=0}^{\infty} [0,a]x^i / a \in Z^+ \cup \{0\}\}$ be an interval group-semigroup. Prove $G$ in problem (177) is different from $G$ given in problem (178).

179. Let $G = G_1 \cup G_2 = \{ \sum_{i=0}^{\infty} [0,a]x^i \mid a_i \in Z^+ \cup \{0\}\} \cup \{ \sum_{i=0}^{\infty} [0,a]x^i \mid a_i \in Z_{48}\}$ be an interval bisemigroup. Find interval biideals of $G$. Does $G$ contain an interval subbisemigroup which is not an interval biideal of $G$?

180. Let $G = \{ \sum_{i=0}^{6} [0,a_i]x^i \mid a_i \in Z_{40}, x^7 = 1, \times\} \cup \{$All $8 \times 8$ matrices $A$ with entries from $Z_{20}$ where $|A| \neq 0\}$ be a quasi interval semigroup-group.
   a. What is the biorder of $G$?
   b. Find quasi interval subsemigroup-subsemigroups in $G$.

181. Let $G = G_1 \cup G_2 = \{<g / g^{20} = 1> \cup \{S_X$ where $X = ([0, a_1], [0, a_2], \ldots, [0, a_{10}])\}$ be a quasi interval bigroup.
   a. What is the biorder of $G$?
   b. Find quasi interval normal bisubgroups in $G$.
   c. Using the quasi interval binormal bisubgroup. $H = \{1, g^4, g^8, g^{12}, g^{16}\} \cup A(X) = H_1 \cup H_2 \subseteq G_1 \cup G_2$, define the quasi interval quotient bigroup $G/H = G_1/H_1 \cup G_2/H_2$. What is the biorder of $G/H$?

182. Describe any stricking property enjoyed by a quasi interval bigroupoid.

183. Find some dissimilarities between interval biloops and interval bigroupoids in general.

184. Determine some similarities between quasi interval bigroups and interval bigroups.

185. Find atleast one marked difference between a quasi interval semigroup and a quasi interval group semigroup.

186. Find some applications of interval bigroups.

187. How can interval biloops be used in the biedge colouring problems $K_{2n} \cup K_{2m}$?

188. Can interval bigroupoids be used in biautomaton?

189. Will interval bisemigroups be used in bisemiautomaton so as to yield a better result?

190. Let $G = G_1 \cup G_2 \cup G_3$ be a 3-interval groupoid or interval trigroupoid where $G_1 = \{[0, a] | a \in Z_{11}, *, (3, 2)\}$, $G_2 = \{[0, b] / b \in Z_{20}, *, (1, 4)\}$ and $G_3 = \{[0, c] / c \in Z_{40}, *, (0, 11)\}$. Find the triorder of the triinterval groupoid. Find substructures of G.

191. Obtain any interesting property associated with n-interval groupoids.

192. Give an example of a S-strong 5-interval Bol groupoid.

193. Give an example of a S-5-interval Bol groupoid.

194. Give an example of a 7-interval P-groupoid.

195. Does their exist a 8-interval groupoid which is not a S-8-interval groupoid?

196. Give an example of a 10-interval groupoid which is not a S-strong right alternative 10-interval groupoid.

197. Give an example of a 5-interval groupoid which is a S-left alternative 5-interval groupoid.

198. Give an example of a S-strong 4-interval Moufang groupoid.

199. Find a necessary and sufficient condition of a 7-interval groupoid to be a S-strong P-groupoid.

200. Give an example of a 3-interval idempotent groupoid.

201. Give an example of a 5-interval semigroup which has no 5-interval ideals.

202. Give an example of a 12-interval semigroup in which every 12-interval subsemigroup is a 12-interval ideal.

203. What makes the study of n-interval semigroups interesting?

204. Obtain some interesting results about n-interval groupoids (n > 2).

205. Characterize those n-interval groupoids which are not Smarandache n-interval groupoids.

206. Give an example of a 7-interval groupoid which is not a S-7-interval groupoid.

207. Give an example of a 6-interval groupoid in which every 6-interval subgroupoid is a S-6-interval subgroupoid.

208. Give an example of a 17-interval groupoid G which has no S - 17 interval subgroupoids but G is a S-17-interval groupoid.

209. Can these new structures be used in cryptography?

210. Give an example of a 3-interval semigroup which is a S-interval semigroup.

211. Give an example of a 5-interval semigroup which is not a S-5-interval semigroup.

212. Give an example of a n-interval semigroup (n > 2) in which every n-interval subsemigroup is a S-n interval subsemigroup.

213. Give an example of a n-interval semigroup (n>2) in which every n-interval subsemigroup is a S-n-interval ideal.

214. Let $G = G_1 \cup G_2 \cup G_3 \cup G_4$ be a 4-interval groupoid where $G_1 = \{[0, a] / a \in Z_8, *, (3, 0)\}$; $G_2 = \{[0, b] / b \in Z_{10}, *, (7, 0)\}$; $G_3 = \{[0, c] / c \in Z_8, *, (0, 3)\}$ and $G_4 = \{[0, d] / d \in Z_{10}, *, (0, 7)\}$.
    a. What is the 4-order of G?
    b. Find some 4-interval subgroupoids.
    c. Is G a S-4-interval groupoid?
    d. Does G satisfy any one of the special identities?
    e. Can G be have S-4-interval subgroupoids?

215. Let $G = G_1 \cup G_2 \cup G_3 = \{[0, a] / a \in Z_{12}, (3, 2), *\} \cup \{[0, b] / b \in Z_{12}, (7, 5), *\} \cup \{[0, c] / c \in Z_{12}, (8, 11), *\}$ be a 3-interval groupoid?
    a. Is G a S-3- interval groupoid?
    b. Can G have S-3- interval subgroupoids?
    c. Does G satisfy any one of the special identities?
    d. Is G a normal 3-interval groupoid?

216. Obtain some interesting properties enjoyed by n interval semigroups.

217. Does there exist a n-interval semigroup which is S-Lagrange n-interval semigroup?

218. Give an example of a n-interval semigroup which is not S-weakly Lagrange n-interval semigroup.

219. Let $S = S_1 \cup S_2 \cup S_3 \cup S_4 = \{[0, a] / a \in Z_{11}, \times\} \cup \{[0, a] / a \in Z_{13}, \times\} \cup \{[0, a] / a \in Z_{19}, \times\} \cup \{[0, a] / a \in Z_5, \times\}$ be a 4-interval semigroup.
   a. What is the order of S?
   b. Prove S is a S-4-interval semigroup.
   c. Is S a S-weakly Lagrange 4-interval semigroup?
   d. Is S a S-Lagrange 4 interval semigroup?
   e. Find all the 4-interval subgroups in S.

220. Let $S = S_1 \cup S_2 \cup S_3 \cup S_4 \cup S_5 = \{[0, a] / a \in Z_{22}, \times\} \cup \{[0, a] / a \in Z_{26}, \times\} \cup \{[0, a] / a \in Z_{42}, \times\} \cup \{[0, a] / a \in Z_{240}, \times\} \cup \{[0, a] / a \in Z_{36}, \times\}$ be a 5-interval semigroup.
   a. Find the order of S.
   b. Is S a S-5-interval semigroup?
   c. Is S a S-Lagrange 4-interval semigroup?
   d. Can S be only a S-weakly Lagrange 4-interval semigroup?
   e. Find 4-interval zero divisors in S.
   f. Does S have 4-interval idempotents?
   g. Does S have 4-interval units?

221. Illustrate by example a quasi n-interval semigroup which is S-quasi n-interval semigroup but has no S-quasi n-interval subsemigroup.

222. Find all the zero divisors and idempotents and units in $S = S_1 \cup S_2 \cup S_3 = \{[0, a] / a \in Z_{24}, \times\} \cup \{Z_{28}, \times\} \cup \{[0, b] | b \in Z_{42}, \times\}$, the quasi 3-interval semigroup.

223. Give an example of a quasi n-interval semigroup in which every quasi n-interval subsemigroup is a S-quasi n-interval subsemigroup.

224. Give an example of a quasi n-interval semigroup which has no S - quasi n-interval ideals.

225. Enumerate all the special properties enjoyed by a quasi n-interval semigroup.

226. Let $S = S_1 \cup S_2 \cup S_3 \cup S_4 \cup S_5 = \{Z_8, \times\} \cup \{[0, a] / a \in Z_9, \times\} \cup \{Z_{12}, \times\} \cup \{[0, a] / a \in Z_{15}, \times\} \cup \{Z_{18}, \times\}$ be a quasi 5-interval semigroup.
   a. What is the order of S?
   b. Is S a S-Lagrange quasi 5-interval semigroup?
   c. Is S only a S-weakly Lagrange quasi 5-interval semigroup?
   d. Is S a S-quasi 5-interval semigroup?
   e. Is S a S-quasi 5-interval cyclic semigroup?
   f. Can S be a S-p-Sylow 5-quasi interval semigroup?
   g. Find zero divisors, units and idempotents if any in S.
   h. Does S have S-zero divisors, S-units and S-idempotents?

227. Let $S = S_1 \cup S_2 \cup S_3 \cup S_4$ be a quasi 4-interval semigroup where $S_1 = S(3)$, $S_2 = S(7)$, $S_3 = S(5)$ and $S_4 = \{S(X) / X = ([0, a_1], [0, a_2], \ldots, [0, a_6])\}$.
   a. Is S a S-Lagrange quasi 4-interval semigroup?
   b. Is S a S-weakly Lagrange quasi 4-interval
   c. semigroup?
   d. Find S-quasi 4-interval subsemigroups of S (if any).
   e. Find S-quasi 4-interval ideals in S.
   f. Does S have zero divisors?
   g. Find the set of all units in S.
   h. Can S have S-units?
   i. What is the order of S?

228. Let $S = S_1 \cup S_2 \cup S_3 \cup S_4 \cup S_5 = \{[0, a] / a \in 3Z^+ \cup \{0\}\} \cup \{8Z^+ \cup \{0\}, \times\} \cup \{[0, a] / a \in 19Z^+\} \cup \{13Z^+, \times\} \cup \{11Z^+, \times\}$ be a quasi 5-interval semigroup.
   a. Can S have S-zero divisors?
   b. Can S have idempotents?
   c. Is S a S-quasi 5-interval semigroup?
   d. Can S have S-quasi 5-interval ideals?
   e. Find at least 2 quasi 5-interval subsemigroups in S.

229. Let $S = S_1 \cup S_2 \cup S_3 \cup S_4 \cup S_5 = \{[0, a] / a \in Z_{24}, \times\} \cup \{Z_{30}, \times\} \cup \{[0, a] / a \in Z_{45}, \times\} \cup \{Z_{120}, \times\} \cup \{Z_{240}, \times\}$ be a quasi 5-interval semigroup.
   a. Find the order of S.

b. Can S have S-zero divisors?
c. Find quasi 5-interval ideals of S.
d. Is S a S-quasi 5-interval semigroup?
e. Can S have S-quasi 5-interval ideals?
f. Can S have quasi 5-interval subsemigroups?

230. Let $S = S_1 \cup S_2 \cup S_3 \cup S_4 \cup S_5 = \{Z_{12}, \times\} \cup \{[0, a] / a \in Z_{36}, \times\} \cup S(5) \cup \{3 \times 3$ interval matrices with intervals of the form $[0, a]$ where $a \in Z_{24}\} \cup \{([0, a] [0, b], [0, c]) / a, b, c \in Z_{10}, \times\}$ be a quasi 5-interval semigroup.
    a. What is the order of S?
    b. Is S a S-quasi 5-interval semigroup?

231. Give an example of a Smarandache strong n-interval Bol loop.

232. Give an example of a S-quasi n-interval Bol loop.

233. Does there exist any n-interval Bruck loop of finite order?

234. Construct a n-interval m × s matrix loop which is a S-n-interval m × s matrix loop.

235. Show by an example all n-interval loop which is not a Smarandache n-interval loop.

236. Obtain some interesting applications of S-n-interval loop.

237. Derive some interesting properties about quasi n-interval loops of finite order.

238. Prove in general all n-interval loops are not S-strongly Lagrange n-interval loop.

239. Give an example of a S- Lagrange n-interval loop.

240. Is every n-interval loop a S-Cauchy loop?

241. Give an example of a 7-interval loop which is S-Cauchy loop.

242. Give an example of a quasi 8-interval loop which is a S-2-Sylow loop.

243. Give an example of a mixed 5-interval algebraic structure of finite order.

244. Let $L = L_1 \cup L_2 \cup \ldots \cup L_7$ be a 7-interval loop built using $L_7$ (5), $L_9$ (5), $L_{13}$ (7), $L_{21}$(11), $L_{43}$ (20), $L_{25}$ (9) and $L_{49}$ (9).
    a. What is the order of L?
    b. Is L a S-strongly Lagrange?
    c. Is L commutative?
    d. Can L have 7-inteval subloops?
    e. Is L a S-2-Sylow loop?
    f. Can L have 7-interval normal subloop?
    g. Is L S-simple.
    h. Is $SN_1(L) = SN_2(L)$?
    i. Does L satisfy any the well known identities?
    j. Find the S-Moufang center of L.
    k. Is L a power associative loop?
    l. Prove L is a S-7-interval loop.
    m. Can L have proper S-7-interval subloops?

245. Let $G = G_1 \cup G_2 \cup G_3 = \{3 \times 2$ interval matrices with intervals of the form [0, a], $a \in Z_{22}$, *, (10, 0)$\} \cup \{2 \times 4$ interval matrices with interval of the form [0, a] where $a \in Z_{12}$ under $+\} \cup \{3 \times 3$ interval matrices with intervals of the form [0, a], $a \in Z_8$ under $\times\}$ be a mixed 3 interval matrix algebraic structure.
    a. Is G commutative?
    b. What is the order of G?
    c. Find at least 3 substructures in G.
    d. Does the substructure satisfy any well known classical theorems for finite groups?

246. Give an example of a 5-interval group-loop of order 256.

247. Give an example of a 3-interval groupoid-loop G of infinite order and such that G is non commutative and has no substructures.

248. Give an example of a 8-interval groupoid-loop which is not a S-8-interval groupoid-loop.

249. Obtain some interesting applications of quasi n-interval groupoid-loops?

250. Find zero divisors and S-zero divisors of a 3-interval groupoid-semigroup G by constructing G of order 9.12.24.

251. Give an example of a n-interval semigroup which is not a S-n-interval semigroup.

252. Give an example of a quasi 18-interval semigroup with out zero divisors.

253. Can n-interval semigroups constructed using special symmetric semigroups have zero divisors? Justify your claim.

254. Give an example of a quasi 4-interval semigroup S using symmetric interval semigroups of order $3^3$, $4^4$, $2^2$ and $5^5$.

255. Does S in problem (254) S-Lagrange?

256. Does S in problem (254) S-quasi 4-interval symmetric semigroup?

257. Obtain atleast one quasi 5-interval semigroup which is not a S-quasi 5-interval semigroup.

258. State and prove any interesting result on mixed quasi n-interval algebraic structure.

259. Can any mixed quasi n-interval algebraic structure be S-quasi n-interval algebraic structure? Justify your claim.

260. Give an example of a mixed quasi n-interval algebraic structure which has no mixed quasi n-interval algebraic substructure.

261. Does there exists a 5-interval loop which is not a S-Cauchy 5-interval loop.

262. Give an example of S-Cauchy 6-interval loop.

263. Give an example of a S-strong Lagrange 7-interval loop.

264. Give an example of a 5-interval loop which is not a S-Lagrange 5-interval loop.

265. Let $L = L_1 \cup L_2 \cup L_3 \cup L_4 \cup L_5 = \{[0, a] / a \in \{e, 1, 2, ...., 17, *, 10\} \cup L_{13}(7) \cup Z_8(3, 4) \cup \{[0, a] / a \in Z_{18}, *, 90, 5\} \cup \{[0, a] / a \in Z_{12}, *, (0, 5)\}$ be a 5 quasi interval loop-groupoid.
   a. Find the order of L.
   b. Is L a S-5-quasi interval loop-groupoid?
   c. Does L satisfy any of the special identities?
   d. Is L a S-strongly Moufang?
   e. Obtain any other interesting property about L.

266. Obtain some interesting properties about n-quasi interval semigroup-groupoid.

267. Give an example of a quasi n-interval semigroup groupoid which is not a S-quasi interval semigroup-groupoid.

268. Can a n-interval loop of odd order exist?

269. Give an example of a 5-interval loop of order 251.

270. Give an example of a n-interval groupoid - loop of order $2^t$. (t is at the readers choice, t > n)

271. Let $L = L_1 \cup L_2 \cup L_3$ be a 3-loop of order 16.18.20 built using $Z_{15}$, $Z_{17}$ and $Z_{19}$.

a. Is L S-simple?
b. Is L S-Lagrange?
c. Prove L is S-Cauchy.
d. Is L S-strongly cyclic?
e. Does L have proper S-subloops?

272. Give an example of a 3-simple interval group.

273. Does there exists a 6-interval group which is non commutative of order $\underline{3}$. 12 $\underline{5}$.25.29.43?

274. Let $L = L_1 \cup L_2 \cup L_3 \cup L_4 = \{$All $2 \times 2$ interval matrices with intervals of the form [0, a] where $a \in Z_{12}$ under multiplication$\} \cup \{([0, a_1] [0, a_2] \ldots [0, a_5]) \mid a_i \in \{e, 1, 2, \ldots, 27\}, *, 11\} \cup \{ \begin{bmatrix} [0, a_1] \\ [0, a_2] \\ [0, a_3] \end{bmatrix} / a_i \in Z_{15}, *, (2, 3); 1 \le i \le 3\} \cup (Z_{20}, \times)$ be a quasi 4-interval mixed algebraic structure.
a. What is the order of L?
b. Find substructures of L.

275. Let $G = G_1 \cup G_2 \cup G_3 \cup G_4 \cup G_5 = \{3 \times 3$ interval matrices with intervals of the form [0, a] where $a \in Z_{12}$, under multiplication$\} \cup \{3 \times 5$ interval matrices with intervals of the form [0, a] where $a \in Z_{40}, +\} \cup \{1 \times 6$ interval matrices with intervals of the form [0, a] where $a \in Z_{25}, \times\} \cup \{4 \times 3$ interval matrices with intervals of the form [0, a] where $a \in Z_{120}, +\} \cup \{2 \times 2$ upper triangular interval matrices with intervals of the form [0, a] where $a \in Z_{21}, \times\}$ be a 5-interval semigroup.
a. What is the order of G?
b. Find substructures in G.
c. Find zero divisors in G.
d. Find S-zero divisors in G.

- e. Find S-ideals if any in G.
- f. Find S-subsemigroups which are not S-ideals if any
- g. in G.
- h. Find idempotents in G.

276. Give an example of a quasi mixed n-interval algebraic structure which does not have any substructure.

277. Give an example of a quasi mixed 5-interval algebraic structure which has atleast 3-substructures.

278. What is the order of $S = \{Z_{27}, \times\} \cup \{Z_{40}, +\} \cup \{[0, a] / a \in Z_9, *, (2, 3)\} \cup \{[0, a] / a \in \{e, 1, 2, …, 27\}, *, 5\} \cup L_{29}(7)$? Does S have substructures?

279. Give an example of a 5-interval S-simple loop.

280. Give an example of a 6-interval simple group.

281. What can one say about homomorphism of 3-interval group into a 4-interval group? Illustrate this situation by some examples.

282. Give any nice application of n-interval groupoids.

283. Can n-interval semigroups be used in automaton construction? (Here interval solution for machines to be used)?

284. What is the application of n-interval loops in colouring problem of $K_{2n}$?

285. Give an example of a n-interval loop which is has no S-subloops.

286. Give any interesting applications of mixed n-interval matrix algebraic structure.

287. What can one say about applications of quasi mixed n-interval matrix algebraic structure?

288. Give any possible applications of n-interval matrix groupoids.

289. Can the notion of n-interval matrix groups help in any applications in physics?

290. What can one say about the applications of quasi n-matrix interval semigroups?

291. Determine those n-matrix interval semigroups which has no S-n-ideal.

292. Find some interesting applications of n-matrix interval groupoids built using $Z_n$'s.

293. What is the benefit of using n-interval structures?

294. Find an example 7-interval groupoid-loop which is S-Moufang.

295. Give an example quasi 6-interval loop-groupoid which is a S-Bol.

296. Give an example of a 5-interval loop-groupoid which is Smarandache strong right alternative.

297. Obtain some special properties about S-strong Bol n-interval groupoid-loops which is not in general true for other n-interval groupoid-loops.

298. Give an example of a S-quasi 4-interval groupoid of order 420. How many such S-quasi 4-interval groupoids can be constructed using $S_i = ([0, a] / a \in Z_i, *, (p, q))$, $0 < i < \infty$?

# FURTHER READING

# INDEX

# P

# Q

**R**

**S**

# ABOUT THE AUTHORS

**Dr.W.B.Vasantha Kandasamy** is an Associate Professor in the Department of Mathematics, Indian Institute of Technology Madras, Chennai. In the past decade she has guided 13 Ph.D. scholars in the different fields of non-associative algebras, algebraic coding theory, transportation theory, fuzzy groups, and applications of fuzzy theory of the problems faced in chemical industries and cement industries. She has to her credit 646 research papers. She has guided over 68 M.Sc. and M.Tech. projects. She has worked in collaboration projects with the Indian Space Research Organization and with the Tamil Nadu State AIDS Control Society. She is presently working on a research project funded by the Board of Research in Nuclear Sciences, Government of India. This is her 55$^{th}$ book.

On India's 60th Independence Day, Dr.Vasantha was conferred the Kalpana Chawla Award for Courage and Daring Enterprise by the State Government of Tamil Nadu in recognition of her sustained fight for social justice in the Indian Institute of Technology (IIT) Madras and for her contribution to mathematics. The award, instituted in the memory of Indian-American astronaut Kalpana Chawla who died aboard Space Shuttle Columbia, carried a cash prize of five lakh rupees (the highest prize-money for any Indian award) and a gold medal.
She can be contacted at vasanthakandasamy@gmail.com
Web Site: http://mat.iitm.ac.in/home/wbv/public_html/
or http://www.vasantha.in

**Dr. Florentin Smarandache** is a Professor of Mathematics at the University of New Mexico in USA. He published over 75 books and 150 articles and notes in mathematics, physics, philosophy, psychology, rebus, literature.

In mathematics his research is in number theory, non-Euclidean geometry, synthetic geometry, algebraic structures, statistics, neutrosophic logic and set (generalizations of fuzzy logic and set respectively), neutrosophic probability (generalization of classical and imprecise probability). Also, small contributions to nuclear and particle physics, information fusion, neutrosophy (a generalization of dialectics), law of sensations and stimuli, etc. He can be contacted at smarand@unm.edu